\documentclass[sn-mathphys,Numbered]{sn-jnl}

\usepackage{graphicx}%
\usepackage{multirow}%
\usepackage{amsmath,amssymb,amsfonts}%
\usepackage{amsthm}%
\usepackage{mathrsfs}%
\usepackage[title]{appendix}%
\usepackage{xcolor}%
\usepackage{textcomp}%
\usepackage{manyfoot}%
\usepackage{booktabs}%
\usepackage{algorithm}%
\usepackage{algorithmicx}%
\usepackage{algpseudocode}%
\usepackage{listings}%
\usepackage{placeins}
\usepackage{accents}
\usepackage[utf8]{inputenc}
\usepackage[T1]{fontenc}

\makeatletter \@addtoreset{figure}{section} \makeatother
\renewcommand\figurename{$\mathbf{Fig.}$}
\renewcommand{\thefigure}{\arabic{section}.\arabic{figure}}
\numberwithin{equation}{section}

\def\v{\varepsilon}
\def\e{\eta}

\def\t{\theta}
\def\k{\kappa}

\def\g{\gamma}

\def\r{\rho}

\def\z{\zeta}

\def\f{\frac}
\def\pa{\partial}

\def\dis{\displaystyle}

\def\var{\varepsilon}
\newcommand{\ta}{\varpi}

\newcommand{\R}{\mathbb{R}}

\makeatletter
\newcommand{\rmnum}[1]{\romannumeral #1}
\newcommand{\Rmnum}[1]{\expandafter\@slowromancap\romannumeral #1@}
\makeatother

\theoremstyle{thmstyleone}%
\newtheorem{theorem}{Theorem}[section]
\newtheorem{lemma}[theorem]{Lemma}
\newtheorem{proposition}[theorem]{Proposition}%
\newtheorem{corollary}[theorem]{Corollary}

\theoremstyle{thmstyletwo}%
\newtheorem{remark}{Remark}[section]%

\theoremstyle{thmstylethree}%
\newtheorem{definition}{Definition}[section]%

\begin{document}

\title[Global Finite-Energy Solutions of the Compressible Euler-Poisson Equations]{\qquad\,\, Global Finite-Energy Solutions of \\ the Compressible Euler-Poisson Equations
for General Pressure Laws with Large Initial Data of Spherical Symmetry}

\author*[1,2,4]{\fnm{Gui-Qiang G.} \sur{Chen}}\email{chengq@maths.ox.ac.uk}

\author[2,3]{\fnm{Feimin} \sur{Huang}}\email{fhuang@amt.ac.cn}

\author[2]{\fnm{Tianhong} \sur{Li}}\email{thli@math.ac.cn}

\author[2]{\fnm{Weiqiang} \sur{Wang}}\email{wangweiqiang@amss.ac.cn}

\author[2,3]{\fnm{Yong} \sur{Wang}}\email{yongwang@amss.ac.cn}

\affil*[1]{\orgdiv{Mathematical Institute}, \orgname{University of Oxford}, \orgaddress{\city{Oxford}, \postcode{OX2 6GG}, \country{UK}}}

\affil[2]{\orgdiv{Academy of Mathematics and Systems Science}, \orgname{Chinese Academy of Sciences}, \orgaddress{ \city{Beijing}, \postcode{100190}, \country{China}}}

\affil[3]{\orgdiv{School of Mathematical Sciences}, \orgname{University of Chinese Academy of Sciences}, \orgaddress{ \city{Beijing}, \postcode{100049}, \country{China}}}

\affil[4]{\orgdiv{School of Mathematical Sciences}, \orgname{Fudan University}, \orgaddress{\city{Shanghai}, \postcode{200433}, \country{China}}}

\abstract{We are concerned with global finite-energy solutions of the three-dimensional compressible
Euler-Poisson equations with {\it gravitational potential} and {\it general pressure law},
especially including the constitutive equation
of {\it white dwarf stars}. In this paper, we construct global finite-energy solutions
of the Cauchy problem for the Euler-Poisson equations with large initial data of spherical symmetry
as the inviscid limit of the solutions of the corresponding Cauchy problem for the compressible
Navier-Stokes-Poisson equations.
The strong convergence of the vanishing viscosity solutions is achieved through
entropy analysis, uniform estimates in $L^p$, and a more general compensated compactness
framework via several new ingredients.
A key estimate is first established for the integrability of the density over unbounded domains
independent of the vanishing viscosity coefficient.
Then a special entropy pair is carefully designed via solving a Goursat problem for the entropy
equation such that a higher integrability of the velocity is established, which is a crucial step.
Moreover, the weak entropy kernel for the general pressure law and its fractional derivatives of
the required order near vacuum ($\r=0$) and far-field ($\r=\infty$) are carefully analyzed.
Owing to the generality of the pressure law, only the $W^{-1,p}_{\mathrm{loc}}$-compactness of
weak entropy dissipation measures with $p\in [1,2)$ can be obtained;
this is rescued by the equi-integrability of weak entropy pairs which can be established by the
estimates obtained above, so that the div-curl lemma still applies.
Finally, based on the above analysis of weak entropy pairs, the $L^p$ compensated compactness
framework for the compressible Euler equations with general pressure law is established.
This new compensated compactness framework and the techniques developed in this paper
should be useful for solving further nonlinear problems with similar features.}


\maketitle

\setcounter{tocdepth}{1}
\tableofcontents

\section{\,Introduction}
We are concerned with global finite-energy solutions of the three-dimensional (3-D) compressible Euler-Poisson equations (CEPEs)
that take the form:
\begin{equation}\label{1.1}
	\left\{
	\begin{aligned}
		&\partial_{t}\rho+
		\nabla\cdot\mathcal{M}=0,\\
		&\partial_{t}\mathcal{M}+
		\nabla\cdot\Big(\frac{\mathcal{M}\otimes \mathcal{M}}{\rho}\Big)+\nabla P+\r\nabla\Phi=0,\\
		&\Delta\Phi={k_{g}\r},
	\end{aligned}
	\right.
\end{equation}
for $(t,\mathbf{x}):=(t,x_1,x_2,x_3)\in \R_+^4:=\R_{+}\times \R^{3}=(0,\infty)\times \R^3$.
System \eqref{1.1} is used to model the motion of compressible gaseous stars under a self-consistent
gravitational field ({\it cf.} \cite{Chandrasekhar-1938}),
where $\rho$ is the density, $P=P(\rho)$ is the pressure, $\mathcal{M}\in \R^{3}$ is the momentum, $\Phi$
represents the gravitational potential of gaseous stars as ${k_{g}}>0$,
$\nabla=(\partial_{x_1}, \partial_{x_2}, \partial_{x_3})$,
and $\Delta=\partial_{x_1x_1}+\partial_{x_2x_2}+\partial_{x_3x_3}$.
Without loss of generality, by scaling, we take ${k_{g}}=1$ throughout this paper.

The constitutive pressure-density relation $P(\rho)$ depends on the types
of gaseous stars.
The class of polytropic gases, {\it i.e.},
\begin{equation}\label{1.6-1}
	P(\rho)=\kappa \rho^{\gamma} \qquad \text{ for $\kappa>0$ and $\gamma\in (1,3)$},
\end{equation}
has been widely investigated in mathematics.
From the point view of astronomy, the constitutive pressure $P(\rho)$ for certain gaseous stars
is
not of the polytropic form.
For example, the pressure law of a white dwarf star takes the following
form ({\it cf.} \cite{Chandrasekhar-1938, Strauss-Wu-2020}):
\begin{equation}\label{1.6}	
	P(\rho)=\mathcal{C}_{1}\int_{0}^{\mathcal{C}_{2}\rho^{\frac{1}{3}}}\frac{s^4}{\sqrt{\mathcal{C}_{3}+s^2}}\,\mathrm{d}s
	\qquad\, \text{ for $\rho>0$},
\end{equation}
where $\mathcal{C}_{1}, \mathcal{C}_{2}$, and $\mathcal{C}_{3}$ are positive constants.
It can be checked that $P(\rho)\cong \kappa_1\rho^\frac53$ as $\rho\to 0$
and $P(\rho)\cong \kappa_2\rho^\frac43$ as $\rho\to \infty$
for some positive constants $\kappa_1$ and $\kappa_2$.

In this paper, we consider a general
pressure law in which any pressure function
$P(\rho)$ satisfies the following conditions:
\begin{itemize}
	\item [(i)] The pressure function $P(\rho)$ is in $C^1([0,\infty))\cap C^4(\R_+)$ and satisfies the hyperbolic and genuinely nonlinear conditions:
	\begin{equation}\label{1.3}
		P'(\rho)>0,\quad 2P'(\rho)+\rho P''(\rho)>0\qquad\,\, \text{for }\rho>0.
	\end{equation}
	
	\smallskip	
	\item [(ii)] There exists a constant $\rho_{*}>0$ such that
	\begin{equation}\label{1.4}
		P(\rho)=\kappa_{1} \rho^{\gamma_1}\big(1+\mathcal{P}_1(\rho)\big) \qquad \text { for } \rho \in[0, \rho_{*}),
	\end{equation}
	with some constants $\gamma_1\in (1,3)$ and $\kappa_1>0$, and a function $\mathcal{P}_1(\rho)\in C^4(\R_+)$
	satisfying that
	$|\mathcal{P}_1^{(j)}(\rho)|\leq C_{*}\rho^{\gamma_1-1-j}$ for $\rho\in (0,\rho_{*})$
	and $j=0,\cdots,4$, where $C_{*}>0$ is a constant depending only on $\rho_{*}$.
	
	\smallskip
	\item [(iii)] There exists a constant $\rho^{*}> \rho_{*}>0$ such that
	\begin{equation}\label{1.5}
		P(\rho)=\kappa_2\rho^{\gamma_2}\big(1+\mathcal{P}_2(\rho)\big) \qquad \text { for } \rho \in[\rho^{*},\infty),
	\end{equation}
	with some constants $\gamma_2\in (\frac{6}{5},\gamma_{1}]$ and $\kappa_{2}>0$,
	and a function $\mathcal{P}_2(\rho)\in C^{4}(\R_{+})$ satisfying that
	$|\mathcal{P}_{2}^{(j)}(\rho)|\leq C^{*}\rho^{-\epsilon-j}$ for $\rho\in [\rho^{*},\infty)$ and $j=0,\cdots,4$,
	where $\epsilon>0$, and $C^{*}>0$ is a constant depending only on $\rho^{*}$.
\end{itemize}

It is direct to see that the polytropic gases in \eqref{1.6-1} satisfy assumptions \eqref{1.3}--\eqref{1.5}.
Moreover, the white dwarf star \eqref{1.6} is also included with
\begin{equation}\label{1.7}
	\gamma_1=\frac{5}{3},\quad\k_{1}=\f{1}{5\sqrt{\mathcal{C}_3}}\mathcal{C}_1\mathcal{C}_2^5,
	\quad\gamma_2=\frac{4}{3},\quad \k_{2}=\frac{1}{4}\mathcal{C}_1\mathcal{C}_2^4,\quad \epsilon=\frac{2}{3}.
\end{equation}
The restriction: $\gamma_2>\frac{6}{5}$ is necessary to ensure
the global existence of finite-energy
solutions with finite total mass. Such a condition is also needed for the existence of the Lane-Emden solutions; see \cite{Chandrasekhar-1938, Lin-1997}.

We consider the Cauchy problem of \eqref{1.1} with the initial data:
\begin{equation}\label{1.2}
	(\rho,\mathcal{M})(0,\mathbf{x})=(\rho_0,\mathcal{M}_0)(\mathbf{x})\,\rightarrow\, (0,\mathbf{0})\qquad\,\, \text{ as $|\mathbf{x}|\to \infty$},
\end{equation}
subject to the far field condition:
\begin{equation}\label{1.2-1}
	\Phi(t,\mathbf{x})\,\rightarrow\, 0\qquad\,\, \text{ as $|\mathbf{x}|\to \infty$}.
\end{equation}

The global existence of solutions of the Cauchy problem \eqref{1.1} and \eqref{1.2}--\eqref{1.2-1}
is a longstanding open problem.
Many efforts have been made for the polytropic gas case \eqref{1.6-1}.
Considerable progress has been made on the smooth or special solutions
under some restrictions on the initial data.
Among the most famous solutions of CEPEs \eqref{1.1} are the Lane-Emden steady solutions ({\it cf.} \cite{Lin-1997}),
which describe spherically symmetric gaseous stars in equilibrium and minimize the energy
among all possible configurations ({\it cf.} \cite{Lieb-Yau-1987}).
There exist expanding solutions for the non-steady CEPEs \eqref{1.1}.
Hadzi\'{c}-Jang \cite{Hadzic-Jang-2018} proved the nonlinear stability of
the affine solution (which is linearly expanding) under small spherically symmetric perturbations
for $\g=\frac{4}{3}$, while the stability problem for $\g\neq \frac{4}{3}$ is still widely open.
A class of linearly expanding solutions for $\g=1+\frac{1}{k}$ with $k\in \mathbb{N}\backslash \{1\}$,
or $\g\in (1,\frac{14}{13})$,
was further constructed in \cite{Hadzic-Jang-2019}.
For $1< \g<\frac43$,
the concentration (collapse) phenomena may happen. Indeed,
as $\gamma=\frac{4}{3}$, there exists an homologous concentration
solution;
see \cite{Fu-Lin-1998,Goldreich-Weber-1980,Makino-1992}.
More recently, Guo-Hadzi\'{c}-Jang \cite{Guo-Hadzic-Jang-2020}
observed a continued
concentration solution for $1< \g<\frac43$; see also \cite{Huang-Yao-2022}.
A kind of smooth radially symmetric self-similar solutions exhibiting gravitational collapse
for $1\leq \gamma<\frac{4}{3}$ can be found
in \cite{Guo-Hadzic-Jang-2021, Guo-Hadzic-Jang-Schrecker-2021}.
We refer to \cite{Luo-Xin-Zeng-2014, Makino-1986}
for the local well-posedness of smooth solutions.

Owing to the strong nonlinearity and hyperbolicity, the smooth solutions of \eqref{1.1} with \eqref{1.6-1}
may break down in a finite time, especially when the initial data are large
({\it cf.} \cite{Chen-Wang-1998, Makino-1992}).
Therefore, weak solutions have to be considered for large initial data.
For gaseous stars surrounding a solid ball, Makino \cite{Makino-1997} obtained
the local existence of weak solutions for $\g\in (1,\frac{5}{3}]$ with spherical symmetry;
also see Xiao \cite{Xiao-2016} for global weak solutions with a class of initial data.
For this case, the possible singularity at the origin is prevented since the domain was
considered outside a ball. Luo-Smoller \cite{Luo-Smoller-2008} proved the conditional stability
of rotating and non-rotating white dwarfs and rotating supermassive stars;
see also Rein \cite{Rein-2003} for the conditional nonlinear stability of
the Lane-Emden steady solutions.

Another fundamental question is whether global solutions can be constructed
via the vanishing viscosity limit
of the solutions of the compressible Navier-Stokes-Poisson  equations (CNSPEs):
\begin{equation}\label{1.8}
	\left\{
	\begin{aligned}
		&\partial_{t}\rho+
		\nabla\cdot\mathcal{M}=0,\\
		&\partial_{t}\mathcal{M}+
		\nabla\cdot\Big(\frac{\mathcal{M}\otimes \mathcal{M}}{\rho}\Big)
		+\nabla P+\r\nabla\Phi
    =\v
		\nabla\cdot\Big(\mu(\rho)D\big(\frac{\mathcal{M}}{\rho}\big)\Big)
		+\v\nabla\Big(\lambda(\rho)
		\nabla\cdot\big(\frac{\mathcal{M}}{\rho}\big)\Big),\\
		&\Delta\Phi=\r,
	\end{aligned}
	\right.
\end{equation}
where $D(\frac{\mathcal{M}}{\rho})=\frac{1}{2}\big(\nabla(\frac{\mathcal{M}}{\rho})+(\nabla(\frac{\mathcal{M}}{\rho}))^{\bot}\big)$
is the stress tensor, the Lam\'{e} (shear and bulk) viscosity coefficients $\mu(\rho)$ and $\lambda(\rho)$ depend on the density (that may vanish on the vacuum) and satisfy
\begin{equation}\label{1.8-1}
	\mu(\rho)\geq 0,\quad \mu(\rho)+3\lambda(\rho)\geq 0\qquad\,\,\, \text{for $\rho\geq 0$},
\end{equation}
and parameter $\v>0$ is the inverse of the Reynolds number. Formally, as $\v\to 0$, the sequence of the solutions of
CNSPEs \eqref{1.8} converges to a corresponding solution of CEPEs \eqref{1.1}.
However, the rigorous proof has been one of the most challenging problems
in mathematical fluid dynamics; see Chen-Feldman \cite{Chen-Feldman-2018} and Dafermos \cite{Dafermos-2016}.

The limit problem with vanishing physical viscosity dates back to the pioneering paper by Stokes \cite{Stokes-1998}.
Most of the known results were around the inviscid limit from the compressible Navier-Stokes
to the Euler equations for the polytropic gas case \eqref{1.6-1}.
The first rigorous proof of the vanishing viscosity limit from the Navier-Stokes
to the Euler equations was provided by Gilbarg \cite{Gilbarg-1951},
in which he established the existence and inviscid limit of the Navier-Stokes shock layers.
For the case of large data, due to the lack of $L^{\infty}$ uniform estimate, the $L^{\infty}$ compensated compactness
framework
\cite{Ding-Chen-Luo, Ding-Chen-Luo-1989, DiPerna-1983, Huang-Wang-2002, Lions-Perthame-Souganidis-1996, Lions-Perthame-Tadmor-1994}
fails to work directly in the inviscid limit of the compressible Navier-Stokes equations.
An $L^{p}$ compensated compactness framework was first studied
by LeFloch-Westdickenberg \cite{LeFloch-Westdickenberg-2007}
for the isentropic Euler equations for the case $\gamma\in (1,\frac{5}{3})$ in \eqref{1.6-1},
and was further developed by Chen-Perepelitsa \cite{Chen-Perepelitsa-2010} to all $\gamma>1$ for \eqref{1.6-1}
with a simplified proof;
see also \cite{Chen-Wang-2020} for spherically symmetric solutions of the M-D isentropic Euler equations.
We also refer to \cite{Schrecker-Schulz-2019, Schrecker-Schulz-2020}
for the 1-D case of asymptotically isothermal gas, {\it i.e.}, $\g_2=1$ in \eqref{1.5}.
More recently, Chen-He-Wang-Yuan \cite{Chen-He-Wang-Yuan-2021} established both
the strong inviscid limit of CNSPEs \eqref{1.8}
and the global existence of spherically symmetric solutions of CEPEs \eqref{1.1}
with large data for polytropic gases \eqref{1.6-1}.

The main purpose of this paper is to  establish the global existence of spherically symmetric finite-energy
solutions of \eqref{1.1} with general pressure law \eqref{1.3}--\eqref{1.5}:
\begin{equation}\label{1.10}
	\rho(t,\mathbf{x})=\rho(t,r),\quad \mathcal{M}(t,\mathbf{x})=m(t,r)\frac{\mathbf{x}}{r},\quad \Phi(t,\mathbf{x})=\Phi(t,r)
	\qquad\,\, \text{for }r=|\mathbf{x}|,
\end{equation}
subject to the initial condition:
\begin{equation}\label{1.11}
	(\rho,\mathcal{M})(0,\mathbf{x})=(\rho_0,\mathcal{M}_{0})(\mathbf{x})=(\rho_{0}(r),m_0(r)\frac{\mathbf{x}}{r})\,\to\, (0,\mathbf{0})
	\qquad\,\, \text{as }r\to \infty,
\end{equation}
and the asymptotic boundary condition:
\begin{equation}\label{1.11-1}
	\Phi(t,\mathbf{x})=\Phi(t,r)\,\rightarrow\, 0\qquad\,\, \text{as }r\to \infty.
\end{equation}
Systems \eqref{1.1} and \eqref{1.8} for spherically symmetric solutions take the following respective forms:
\begin{equation}\label{1.12}
	\left\{
	\begin{aligned}
		&\rho_{t}+m_{r}+\frac{2}{r}m=0,\\
		&m_{t}+\Big(\frac{m^2}{\rho}+P(\rho)\Big)_{r}+\frac{2}{r}\frac{m^2}{\rho}+\rho\Phi_{r}=0,\\
		&\Phi_{rr}+\frac{2}{r}\Phi_{r}=\rho,
	\end{aligned}
	\right.
\end{equation}
and
\begin{equation}\label{1.13}
	\left\{
	\begin{aligned}
		&\rho_{t}+m_{r}+\frac{2}{r}m=0,\\
		&m_{t}+\Big(\frac{m^2}{\rho}+P(\rho)\Big)_{r}+\frac{2}{r}\frac{m^2}{\rho}+\r\Phi_{r}
		\\		
&\qquad
=\v\Big((\mu(\r)+\lambda(\r))\big(\big(\frac{m}{\rho}\big)_r+\frac{2}{r}\frac{m}{\rho}\big)\Big)_{r}
		-\frac{2\v}{r}\mu(\rho)_{r}\frac{m}{\rho},\\
		&\Phi_{rr}+\frac{2}{r}\Phi_{r}=\rho.
	\end{aligned}
	\right.
\end{equation}

\smallskip
The study of spherically symmetric solutions is motivated by many important physical problems
such as stellar dynamics including gaseous stars and supernovae
formation \cite{Chandrasekhar-1938, Rosseland-1964, Whitham-1974}.
An important question is how the waves behave as they move radially inward near the origin,
especially under the self-gravitational force for gaseous stars.
The spherically symmetric solutions of the compressible Euler equations may blow up
near the origin \cite{Courant-Friedrichs-1962, Li-Wang-2006, Merle-Raphael-Rodnianski-Szeftel-2020, Whitham-1974}
at certain time in some situations.
Considering the effect of gravitation, a fundamental problem for CEPEs \eqref{1.1} is
whether a concentration (delta-measure) is formed at the origin.
This problem was answered in \cite{Chen-He-Wang-Yuan-2021} for polytropic gases in \eqref{1.6-1}
when the initial total-energy is finite
that no delta-measure is formed for the density at the origin
for the two cases:
(i) $\gamma>\frac{6}{5}$; (ii) $\g\in (\frac65, \frac43]$
and the initial total-energy is finite and the total mass is less than a critical mass.

In this paper, we establish the global existence of finite-energy solutions of the Cauchy problem \eqref{1.1}
and \eqref{1.11}--\eqref{1.11-1} with spherical symmetry as the inviscid
limits
of global weak solutions of CNSPEs \eqref{1.8} with general pressure law
\eqref{1.3}--\eqref{1.5}, especially including the white dwarf star \eqref{1.6}.
The $L^p$ compensated compactness framework for the general pressure is also established.
Moreover, it is proved that no delta-measure is formed for the density at the origin in the limit,
{and the critical mass for the white dwarf star is the same as
	the Chandrasekhar limit for
	the polytropic gas \eqref{1.6-1} with $\g=\f43$}. The precise statements of the main results
are given in \S 2.

To achieve these, the main strategy is to develop entropy analysis, uniform estimates in $L^p$,
and a more general compensated compactness framework
to prove that there exists a strongly convergent subsequence of solutions of CNSPEs \eqref{1.8}
and show that the limit is the finite-energy weak solution of CEPEs \eqref{1.1} with general pressure law.
This consists of the following three steps:
\begin{itemize}
	\item Establish the uniform $L^p$ estimates of the solutions of CNSPEs \eqref{1.8}
	independent of $\varepsilon$ for some $p>1$;
	\item Show the
	compactness for weak entropy dissipation measures;	
	\item Prove that the associated Young measure $\nu_{(t,r)}$ is the delta measure
	almost everywhere
	which leads to a
	subsequence of solutions of CNSPEs \eqref{1.8} strongly converging
	to the global finite-energy solution of CEPEs \eqref{1.1}.
\end{itemize}

The generality of pressure $P(\rho)$ causes essential difficulties in the analysis for all of the above steps.
We now describe these difficulties and show how they can be overcome:

\smallskip
(i) The crucial step in the $L^p$ estimates
is to show that $\rho |u|^3$ ($u:=\frac{m}{\rho}$ is the velocity) is uniformly bounded in $L^1_{\rm loc}$.
This estimate might be obtained through constructing appropriate entropy $\hat{\eta}$, which is a solution of $(\rho,u)$ to the entropy equation:
\begin{equation}\label{entropy-equation}
	\eta_{\rho\rho}-\frac{P'(\rho)}{\rho^{2}}\eta_{uu}=0,
\end{equation}
with corresponding entropy flux $\hat{q}$.
If $(\rho,u)$ is the solution of \eqref{1.13}, any entropy-entropy flux
pair (entropy pair, for short) $(\hat{\eta},\hat{q})$ satisfies
\begin{align*}
		&(\hat{\eta} r^{2})_{t}+(\hat{q}r^{2})_{r}
		+2r\, (-\hat{q}+\rho u\hat{\eta}_{\r}+\rho u^2\hat{\eta}_m)\nonumber\\
		&
		=\v\,r^{2}\big((\rho u_{r})_{r}+2\rho\big(\frac{u}{r}\big)_{r}\big)\hat{\eta}_{m}
	-\r\int_{a}^{r}\rho\, z^{2}\mathrm{d}z\,\hat{\eta}_{m};
\end{align*}
see \eqref{4.29} below.
For the polytropic gas case \eqref{1.6-1}, there is an explicit formula of the entropy
kernel $\chi(\rho,u)$ so that
$\chi* \psi$
is the entropy,
where $*$ denotes the convolution and $\psi(s)$ is any smooth function.
By choosing $\psi(s)=\frac12s|s|$ as in \cite{Chen-He-Wang-Yuan-2021},
the corresponding entropy flux $\hat{q}$ satisfies that $\hat{q}\ge c_0\rho |u|^3$ and
$-\hat{q}+\rho u\hat{\eta}_{\r}+\rho u^2\hat{\eta}_m\le 0$.
Then the uniform bound of
$\rho |u|^3r^2$  in $L^1_{\rm loc}$ follows ({\it cf.} \cite{Chen-He-Wang-Yuan-2021}).

However, there is no explicit formula of the entropy kernel $\chi$
for the general pressure satisfying \eqref{1.3}--\eqref{1.5}.
Even for the special entropy pair generated by $\psi(s)=\frac12s|s|$,
it is difficult to prove that $\hat{q}\ge c_0\rho|u|^3$ and
$-\hat{q}+\rho u\hat{\eta}_{\r}+\rho u^2\hat{\eta}_m\leq 0$,
due to the lack of explicit formula of the entropy kernel $\chi$.
Hence, the above approach does not apply directly, so we have to seek a
new method to establish the uniform local integrability of $\rho|u|^3$.
One of the novelties of this paper is that
a special entropy $\hat{\eta}$ is constructed by solving a Goursat problem
of the entropy equation \eqref{entropy-equation} in the domain:
$|u|\le k(\rho):=\int_{0}^{\rho}\sqrt{P'(y)}/{y}\,\mathrm{d} y$,
so that
$\hat{\eta}$ is chosen as the mechanical energy $\eta^*$
(see \eqref{1.25}) when $u\ge k(\rho)$, $-\eta^*$ when $u\le-k(\rho)$,
and the boundary condition
for the Goursat problem is given on
the characteristics curves: $u\pm k(\rho)=0$.
One advantage of such a special entropy pair $(\hat{\eta},\hat{q})$ is
that  $\hat{q}\ge c_0\r|u|^3$ as $|u|\geq k(\rho)$,
and  $|\hat{q}|\le C\rho^{\g_2+1}$ for large $\rho$
as $|u|\le k(\rho)$ via careful analysis for the Goursat problem;
see Lemma \ref{lem4.3} for details.
Moreover,
$-\hat{q}+\rho u\hat{\eta}_{\r}+\rho u^2\hat{\eta}_m$ vanishes
as $|u|\geq k(\rho)$.
Similarly, $|-\hat{q}+\rho u\hat{\eta}_{\r}+\rho u^2\hat{\eta}_m|\le C\rho^{\g_2+1}$
for large $\rho$ as $|u|\leq k(\rho)$.

To show $\rho |u|^3$ is uniformly bounded in $L^1_{\rm loc}$,
it remains to prove that
\begin{equation}\label{inte}
	\int_{0}^{T}\int_{d}^{\infty}\rho^{\g_2+1}\,r\mathrm{d}r\mathrm{d}t
\end{equation}
is uniformly bounded for any $T>0$ and $d>0$.
It should be noted that the local integrability
$\int_{0}^{T}\int_{d}^{D}\rho^{\g_2+1}\,\mathrm{d}r\mathrm{d}t\le C$ was obtained in \cite{Chen-He-Wang-Yuan-2021},
but it is not enough yet to obtain the uniform $L^1_{\rm loc}$ estimate for $\rho |u|^3$.
Fortunately, we can obtain even stronger estimate than \eqref{inte}, {\it i.e.},
\begin{equation}\label{inte1}
	\int_{0}^{T}\int_{d}^{\infty}\rho^{\g_2+1}\,r^2\mathrm{d}r\mathrm{d}t\le C,
\end{equation}
by an elaborate analysis; see Lemma \ref{lem2.4} and Corollary \ref{cor2.6} for details.

\smallskip
(ii)
For the polytropic gas case in \eqref{1.6-1},
Chen-Perepelitsa \cite{Chen-Perepelitsa-2010, Chen-Perepelitsa-2015}
and Chen-He-Wang-Yuan \cite{Chen-He-Wang-Yuan-2021} proved
the $H_{\rm loc}^{-1}$--compactness for weak entropy dissipation measures
via the explicit formula of the weak entropy kernel $\chi$ by convolution with any test function of compact support,
which also implies that the entropy pair $(\eta, q)$ is in $L^r_{\rm loc}, r>2$.
However, it is not clear how the $H_{\rm loc}^{-1}$--compactness
for the general pressure satisfying \eqref{1.3}--\eqref{1.5}
can be shown by using the expansions of the weak entropy kernel established
in \cite{Chen-LeFloch-2000,Chen-LeFloch-2003}.
Motivated by \cite{Schrecker-Schulz-2020},
we instead show the $W_{\mathrm{loc}}^{-1,p}$--compactness for $1\le p<2$,
so that an improved div-curl lemma ({\it cf}. \cite{Conti-Dolzmann-Muller-2011}) applies, which leads to the commutation identity for the entropy pairs.
In fact, we can show that the entropy flux function $q$
is bounded by $\rho^{\frac{\gamma_2+1}{2}}$ (see \eqref{6.43}) as $\rho$ is large by careful analysis
on the expansion of the entropy pair so that $q\in L^2_{\rm loc}$.
Then the interpolation compactness yields the $W^{-1,p}$ compactness for $1\le p<2$;
see Lemma \ref{lem6.6} for details.

\smallskip
(iii)
The argument for the reduction of the associated Young measure $\nu_{(t,r)}(\rho, u)$, introduced
in \cite{Chen-He-Wang-Yuan-2021, Chen-Perepelitsa-2010, Chen-Perepelitsa-2015},
for the polytropic gas case in \eqref{1.6-1},
can be roughly stated as follows:
Show first that
every connected subset of the support of the Young measure is a bounded interval;
then use the $L^\infty$ reduction technique introduced
in \cite{Chen-1986, Ding-Chen-Luo, DiPerna-1983, Lions-Perthame-Souganidis-1996}
for a bounded supported Young measure to show that the Young measure is either a delta measure or supported on the vacuum line.
This method essentially relies on the explicit formula of
the weak entropy kernel $\chi$.
For the general pressure law satisfying \eqref{1.3}--\eqref{1.5},
the above method does not apply directly, since it is difficult
to show that every connected subset of the support of the Young measure is a bounded interval.
Motivated by \cite{Chen-LeFloch-2000, Chen-LeFloch-2003, Lions-Perthame-Souganidis-1996, Schrecker-Schulz-2019,
	Schrecker-Schulz-2020},
we carefully analyze the singularities of $\partial^{\lambda_1+1} \chi$
with $\lambda_1=\frac{3-\g_1}{2(\g_1-1)}$ for large $\rho$ and
fully exploit the property: $(\rho^{\g_2+1},\rho|u|^3)\in L^1({\rm d}\nu_{(t,r)})$
so that the $\partial^{\lambda_1+1}-$derivatives can be operated
in the commutation relation; see Lemmas \ref{lem7.1}--\ref{lem7.3} for details.
Then we prove that the Young measure is either a delta measure or
supported on the vacuum line by similar arguments
as in \cite{Chen-LeFloch-2000, Chen-LeFloch-2003, Lions-Perthame-Souganidis-1996, Schrecker-Schulz-2020}.
This new compensated compactness framework and the techniques
developed in this paper
should be useful for solving further nonlinear problems with similar features.

\smallskip
Finally, we remark that there are some related results on CNSPEs \eqref{1.8}
and the compressible Euler equations.
For weak solutions of CNSPEs \eqref{1.8},
we refer to \cite{Ducomet-Feireisl-Petzelotva-Straskraba-2004, Jang-2010, Kobayashi-Suzuki-2008,
	Kong-Li-2017}
with constant viscosity,
and \cite{Duan-Li-2015,Ducomet-Necasova-Vasseur,Zhang-Fang-2009} with density-dependent viscosity.
Recently, Luo-Xin-Zeng \cite{Luo-Xin-Zeng-2014, Luo-Xin-Zeng-2016-1, Luo-Xin-Zeng-2016-2}
proved the large-time stability of the Lane-Emden solution for $\g\in (\frac{4}{3},2)$.
We also refer to the BD entropy developed in  \cite{Bresch-Desjardins-2002, Bresch-Desjardins-2004,
	Bresch-Desjardins-2007, Bresch-Desjardins-Lin-2003}, which provides
a new estimate for the gradient of the density.
For the compressible Euler equations, we refer
to \cite{Chen-1986, Chen-Schrecker-2018,
	Jang-Masmoudi-2015,
	Li-Wang-2006,
	Schrecker-2020}
and the references cited therein.

\smallskip
The rest of this paper is organized as follows:
In \S2, the finite-energy solutions of
the Cauchy problem \eqref{1.1} and \eqref{1.2}--\eqref{1.2-1} for CEPEs are introduced, and
the main theorems of this paper are given.
In \S 3, some elementary quantities and basic properties about the pressure and
related internal energy are provided, and then some remarks
on $M_{\rm c}$ are also given.
The entropy analysis for weak entropy pairs
for the general pressure satisfying \eqref{1.3}--\eqref{1.5}
is presented in \S 4,
especially a special entropy pair is constructed
by solving a Goursat problem for the entropy equation \eqref{1.25-1}.
In \S 5,  a free boundary problem \eqref{2.1}--\eqref{2.6} for \eqref{1.13} is analyzed,
and some uniform estimates of solutions are derived, including the basic energy estimate,
the BD-type entropy estimate, and the higher integrabilities of the density and the velocity.
In \S 6, the global existence of weak solutions of CNSPEs \eqref{1.8} is established,
and some uniform $L^p$ estimates in Theorem \ref{thm1.2} are also obtained.
In \S 7,
we prove the $W_{\mathrm{loc}}^{-1,p}$--compactness
of the entropy dissipation measures for the weak solutions of \eqref{1.13} and
complete
the proof of Theorem \ref{thm1.2}.
In \S 8, the $L^p$--compensated compactness framework for the general pressure
law \eqref{1.3}--\eqref{1.5}
(Theorem \ref{thm1.3}) is established, which leads to the proof of Theorem \ref{thm1.1}
by taking the inviscid
limit of weak solutions of CNSPEs \eqref{1.8} in \S 9.
Appendix \ref{AppendixB} is devoted to the presentation of both the sharp Sobolev inequality
that is used in \S 5 and some variants of Gr\"{o}nwall's inequality which are used in
the proof of several estimates in \S4.

\smallskip
{\bf Notations:} Throughout this paper, we denote $C^{\alpha} (\Omega), L^{p}(\Omega), W^{k, p}(\Omega)$,
and $H^{k}(\Omega)$ as the standard H\"{o}lder space, and the corresponding Sobolev spaces, respectively,
on domain $\Omega$ for $\alpha\in (0,1)$ and $p\in[1, \infty]$.
$C_{0}^{k}(\Omega)$ represents the space of
continuously differentiable functions up to the $k$th order
with compact support over $\Omega$,
and $\mathcal{D}(\Omega):=C_{0}^{\infty}(\Omega)$.
We also use $L^{p}(I ; r^{2}\mathrm{d}r)$
or $L^{p}([0, T) \times I; r^{2}\mathrm{d} r\mathrm{d} t)$
for an open interval $I \subset \mathbb{R}_{+}$ with measure $r^{2}\mathrm{d}r$
or $r^{2}\mathrm{d}r \mathrm{d}t$ correspondingly,
and $L_{\mathrm{loc}}^{p}([0, \infty) ; r^{2}\mathrm{d}r)$
to represent $L^{p}([0, R]; r^{2}\mathrm{d} r)$ for any fixed $R>0$.

\section{\,Mathematical Problem and Main Theorems}
The spherically symmetric initial data function $(\rho_{0},\mathcal{M}_{0})(\mathbf{x})$
given in \eqref{1.11} is
assumed to be of both finite initial total-energy:
\begin{align}\label{1.14}
	E_{0}&:=\int_{\mathbb{R}^{3}}\Big(\frac{1}{2}\Big|\frac{\mathcal{M}_{0}}{\sqrt{\rho_{0}}}\Big|^{2}
	+\rho_{0} e(\rho_{0})\Big)(\mathbf{x}) \,\mathrm{d}\mathbf{x}=\omega_3\int_0^\infty\Big(\frac{1}{2}\frac{m^2_{0}}{\rho_{0}}
	+\rho_{0} e(\rho_{0})\Big)(r)\, r^2\mathrm{d}r<\infty,
\end{align}
and initial total-mass:
\begin{equation}\label{1.14-1}
	M:=\int_{\R^3}\rho_0(\mathbf{x})\,\mathrm{d} \mathbf{x}=\omega_3\int_{0}^{\infty}\rho_0(r)\,r^{2}\mathrm{d} r<\infty,
\end{equation}
where the internal energy $e(\rho)$ is related to the pressure by
\begin{equation}\label{1.15}
	e'(\rho)=
	\frac{P(\rho)}{\rho^2}, \qquad e(0)=0,
\end{equation}
and $\omega_n:=\frac{2\pi^{\frac{n}{2}}}{\Gamma(\frac{n}{2})}$ denotes the surface area of the unit sphere in $\R^n$. The initial potential $\Phi_0(\mathbf{x})$ is determined by
\begin{equation}\label{1.15-1}
	\Delta \Phi_0(\mathbf{x})=\rho_0(\mathbf{x}), \qquad \lim_{|\mathbf{x}|\to \infty}\Phi_0(\mathbf{x})=0.
\end{equation}

For $\gamma_{2}\in (\frac{6}{5},\frac{4}{3}]$, we define the critical mass $M_{\rm c}$ as follows:

\smallskip
{(i) When $\gamma_2=\frac{4}{3}$,
	\begin{align}\label{1.18-0}
		M_{\rm c}:=M_{\rm ch},
	\end{align}
where $M_{\rm ch}$ is the Chandrasekhar limit that
is the total mass of the Lane--Emden steady solution $(\rho_{s}(|{\bf x}|),0)$
for $P(\rho)=\kappa_2\rho^{\frac{4}{3}}$:\,
$\rho_{s}(|{\bf x}|)$ has compact support and is determined by the equations:
\begin{align*}
\nabla_{{\bf x}}P(\rho_{s}(|{\bf x}|))+\rho_{s}(|{\bf x}|)\nabla_{{\bf x}}\Phi({\bf x})=0,\quad \Delta_{{\bf x}}\Phi({\bf x})=\rho_{s}(|{\bf x}|),
\quad P(\rho_{s}|{\bf x}|)=\kappa_2(\rho_{s}(|{\bf x}|))^{\frac{4}{3}},
\end{align*}
with the center density $\rho_{s}(0)=\varrho$.
It is well-known that $M_{\rm ch}$ is a uniform constant with respect to the center density $\varrho$ ({\it cf.} \cite{Chandrasekhar-1938}}).

\smallskip
(ii) When $\g_2\in (\frac{6}{5},\frac{4}{3})$,
	\begin{equation}\label{1.18-3}
		M_{\rm c}:=\sup_{\beta>0}M_{\rm c}(\beta)
	\end{equation}
	with
	\begin{align}
		(4-3\g_2)\Big(\frac{B_{\beta}}{3(\g_2-1)}\Big)^{-\frac{3(\g_2-1)}{4-3\g_2}}M_{\rm c}(\beta)^{-\frac{5\g_2-6}{4-3\g_2}}
         -\omega_{3}^{-1}\beta M_{\rm c}(\beta)=\frac{E_{0}}{\omega_{3}},\label{1.18-3b}
	\end{align}
and
\begin{equation}\label{1.18-4}
\begin{aligned}
		&
B_{\beta}:= \frac{2}{3} \omega_{4}^{-\frac{2}{3}} \omega_{3}^{\frac{4-3 \gamma_2}{3(\gamma_2-1)}}
		(C_{\max}(\beta))^{\frac{5\g_2-6}{3(\g_2-1)}},\\
&
C_{\max}(\beta):=\sup_{\rho\geq 0}\big(\rho^{\g_2-1}(\beta+e(\rho))^{-1}\big)^{\frac{1}{5\g_2-6}}>0.
\end{aligned}
\end{equation}

\smallskip
It is clear in \eqref{1.18-3}--\eqref{1.18-4} that $M_{\rm c}(\beta)$ is well determined  for $\beta>0$
and $\g_2\in (\frac{6}{5},\frac{4}{3})$.
Some useful properties of $M_{\rm c}:=\sup_{\beta>0}M_{\rm c}(\beta)$ will be presented
in Proposition \ref{prop3.1} below.  We also point out that $M_{\rm c}$ in \eqref{1.18-0} is strictly larger than
the one obtained in \cite[(2.8)]{Chen-He-Wang-Yuan-2021} for $\gamma_2=\frac{4}{3}$ ({\it cf}. \cite{Cheng-Cheng-Lin}).

\smallskip
For the spherically symmetric initial data $(\rho_{0}, m_{0},\Phi_{0})(r)$ imposed
in \eqref{1.10}--\eqref{1.11-1} satisfying \eqref{1.14}--\eqref{1.14-1},
using similar arguments as in \cite[Appendix A]{Chen-He-Wang-Yuan-2021},
we can construct a sequence of approximate initial data functions $(\rho_{0}^{\v},m_{0}^{\v}, \Phi_{0}^{\v})(r)$ satisfying
\begin{equation}\label{1.27}
	\begin{aligned}
		&\int_{0}^{\infty} \rho_{0}^{\varepsilon}(r)\, r^{2}\mathrm{d}  r=\frac{M}{\omega_3},
		\qquad \Phi_{0r}^{\v}=\frac{1}{r^2}\int_{0}^{r}\rho_{0}^{\v}(z)\,z^2\mathrm{d}z,\\
		&\,E_{0}^{\v}:=\omega_{3} \int_{0}^{\infty}
		\Big(\frac{1}{2}\Big|\frac{m_{0}^{\varepsilon}}{\sqrt{\rho_{0}^{\varepsilon}}}\Big|^{2}
		+\rho_{0}^{\varepsilon} e(\rho_{0}^{\varepsilon})\Big)\,r^{2}\mathrm{d}r\leq C(E_{0}+1)<\infty,\\
		&\,E_{1}^{\v}:=\varepsilon^{2}\omega_{3} \int_{0}^{\infty}\big|\partial_{r} \sqrt{\rho_{0}^{\varepsilon}(r)}\big|^{2}
		\,r^{2}\mathrm{d}r \leq C \varepsilon(M+1)<\infty.
	\end{aligned}
\end{equation}
Moreover, as $\v\to 0$,  $(E_{0}^{\varepsilon}, E_{1}^{\varepsilon}) \rightarrow (E_{0}, 0)$ and
\begin{align*}
	&(\rho_{0}^{\varepsilon}, \rho_{0}^{\varepsilon}u_{0}^{\varepsilon})(r)
	\rightarrow (\rho_{0}, \rho_{0}u_{0})(r) \qquad
	\text {in } L^{\tilde{q}}([0, \infty); r^{2}\mathrm{d} r) \times L^{1}([0, \infty); r^{2}\mathrm{d} r),\\
	&\Phi_{0r}^{\v}\to \Phi_{0r}\qquad \text{in }L^2([0,\infty);r^{2}\mathrm{d}r),
\end{align*}
where $\tilde{q}\in \{1,\gamma_{2}\}$.
Furthermore, there exists $\v_0\in (0,1]$ such that, for any $\v\in (0,\v_0]$,
\begin{equation}\label{1.27-2}
	M<M_{\rm c}^{\v}\qquad \text{for }\g_2\in (\frac{6}{5},\frac{4}{3}],
\end{equation}
where $M_{\rm c}^{\v}$ is defined in \eqref{1.18-0}--\eqref{1.18-4} by replacing $E_0$ with $E_0^{\v}$.

\smallskip
Now we introduce the weak entropy pairs
of the 1-D isentropic Euler system ({\it cf.} \cite{Chen-LeFloch-2000, Lax-1971}):
\begin{equation}\label{2.10a}
	\left\{
	\begin{aligned}
		&\rho_{t}+m_{r}=0,\\
		&m_{t}+\big(\frac{m^2}{\rho}+P(\rho)\big)_{r}=0.
	\end{aligned}
	\right.
\end{equation}
A pair of functions $(\eta(\rho,m),q(\rho,m))$ is called an entropy pair
of \eqref{2.10a} if
\begin{equation}\label{1.23}
	\nabla q(\rho,m)=\nabla \eta(\rho,m)\nabla \Big(\begin{matrix}m\\
		\frac{m^2}{\rho}+P(\rho)\end{matrix}\Big).
\end{equation}
Moreover, $\eta(\rho,m)$ is called a weak entropy if $\eta(\rho,m)\vert_{\rho=0}=0$,
and a convex entropy if $\nabla^2\eta(\rho, m)\ge 0$.
The mechanical energy and energy flux pair is defined as
\begin{equation}\label{1.25}
	\eta^{*}(\rho,m)=\frac{1}{2}\frac{m^2}{\rho}+\rho e(\rho),\qquad q^{*}(\rho,m)=\frac{1}{2}\frac{m^3}{\rho}+m(\rho e(\rho))',
\end{equation}
which is a convex weak entropy pair.
From \eqref{1.23}, any entropy satisfies
\begin{equation}\label{1.25-1}
	\eta_{\rho\rho}-\frac{P'(\rho)}{\rho^{2}}\eta_{uu}=0
\end{equation}
with $u=\frac{m}{\rho}$. It is known in \cite{Chen-LeFloch-2000, Chen-LeFloch-2003, Lions-Perthame-Souganidis-1996, Lions-Perthame-Tadmor-1994} that any regular weak entropy can be generated by the convolution of a smooth function $\psi(x)$ with
the fundamental solution $\chi(\rho,u,s)$ of the entropy equation \eqref{1.25-1}, {\it i.e.},
\begin{equation}\label{2.21e}
	\eta^{\psi}(\rho,u)=\int_{\R}\chi(\rho,u,s)\psi(s)\,\mathrm{d}s.
\end{equation}
The corresponding entropy flux is generated from the flux kernel $\sigma(\rho,u,s)$ (see \eqref{6.2}), {\it i.e.},
\begin{align}\label{2.21q}
	q^{\psi}(\rho,u)=\int_{\R}\sigma(\rho,u,s)\psi(s)\,\mathrm{d}s.
\end{align}

We first consider the Cauchy problem of CNSPEs \eqref{1.8} with approximate initial data:
\begin{equation}\label{1.19}
	(\rho,\mathcal{M},\Phi)\vert_{t=0}=(\rho_0^{\v},\mathcal{M}_{0}^{\v},\Phi_{0}^{\v})(\mathbf{x}):=(\rho_{0}^{\v}(r), m_{0}^{\v}(r)\frac{\mathbf{x}}{r}, \Phi_{0}^{\v}(r)),
\end{equation}
subject to the far field condition:
\begin{equation}\label{1.19-1}
	\Phi^{\v}(t,\mathbf{x})\longrightarrow 0\qquad \text{as }|\mathbf{x}|\to \infty.
\end{equation}
For concreteness, we take $\v\in (0,1]$ and the viscosity coefficients $(\mu(\rho),\lambda(\rho))=(\rho,0)$ in \eqref{1.8}.

\smallskip
\begin{definition}\label{def1.2}
A triple $(\rho^{\varepsilon}, \mathcal{M}^{\varepsilon},\Phi^{\v})(t,\mathbf{x})$ is said to be a weak solution of the Cauchy problem \eqref{1.8} and \eqref{1.19}  if
\begin{itemize}
\item [\rm (\rmnum{1})] $\rho^{\varepsilon}(t, \mathbf{x}) \geq 0$,
and $(\mathcal{M}^{\varepsilon}, \frac{\mathcal{M}^{\varepsilon}}{\sqrt{\rho^{\v}}})(t, \mathbf{x})=\mathbf{0}\,$
{\it a.e.} on  $\{(t, \mathbf{x})\,:\,\rho^{\varepsilon}(t, \mathbf{x})=0\}\,${\rm (}vacuum{\rm )},
\begin{equation*}
\begin{aligned}
&\rho^{\varepsilon} \in L^{\infty}(0, T ; L^{\gamma_2}(\mathbb{R}^{3})), \quad \nabla \sqrt{\rho^{\varepsilon}} \in L^{\infty}(0, T ; L^{2}(\mathbb{R}^{3})), \\
&\frac{\mathcal{M}^{\varepsilon}}{\sqrt{\rho^{\varepsilon}}} \in L^{\infty}(0, T ; L^{2}(\mathbb{R}^{3})),\quad \Phi^{\v}\in L^{\infty}(0,T;L^{6}(\R^3)),\quad \nabla\Phi^{\v}\in L^{\infty}(0,T;L^2(\R^3)).
\end{aligned}
\end{equation*}

\item [\rm (\rmnum{2})] For any $t_{2} \geq t_{1} \geq 0$ and any
		$\zeta(t, \mathbf{x}) \in C_{0}^{1}([0, \infty) \times \mathbb{R}^{3})$,
		the mass equation $\eqref{1.8}_{1}$ holds in the sense{\rm :}
		\begin{equation*}
			\int_{\mathbb{R}^{3}}(\rho^{\varepsilon} \zeta)(t_{2}, \mathbf{x})\,\mathrm{d} \mathbf{x}
			-\int_{\mathbb{R}^{3}}(\rho^{\varepsilon} \zeta)(t_{1}, \mathbf{x})\,\mathrm{d} \mathbf{x}
			=\int_{t_{1}}^{t_{2}} \int_{\mathbb{R}^{3}}(\rho^{\varepsilon} \zeta_{t}+\mathcal{M}^{\varepsilon} \cdot \nabla \zeta)(t, \mathbf{x}) \,\mathrm{d} \mathbf{x} \mathrm{d} t.
		\end{equation*}
		\item [\rm (\rmnum{3})] For any $\Psi=(\Psi_{1}, \Psi_{2}, \Psi_{3})(t,\mathbf{x}) \in (C_{0}^{2}([0, \infty) \times \mathbb{R}^{3}))^3$,
           the momentum equations $\eqref{1.8}_{2}$ hold in the sense{\rm :}
		\begin{align*}
				&\int_{\R_{+}^4}\Big(\mathcal{M}^{\varepsilon} \cdot \Psi_{t}
				+\frac{\mathcal{M}^{\varepsilon}}{\sqrt{\rho^{\varepsilon}}}
				\cdot\big(\frac{\mathcal{M}^{\varepsilon}}{\sqrt{\rho^{\varepsilon}}} \cdot \nabla\big) \Psi
				+P(\rho^{\varepsilon})
				\nabla\cdot\Psi\Big)\,\mathrm{d} \mathbf{x} \mathrm{d} t
				+\int_{\mathbb{R}^{3}} \mathcal{M}_{0}^{\varepsilon}(\mathbf{x}) \cdot \Psi(0, \mathbf{x})\,\mathrm{d} \mathbf{x} \\
				&=-\varepsilon \int_{\mathbb{R}_{+}^{4}}\Big(\frac{1}{2} \mathcal{M}^{\varepsilon} \cdot\big(\Delta\Psi
				+\nabla (\nabla\cdot\Psi)\big)
				+\frac{\mathcal{M}^{\varepsilon}}{\sqrt{\rho^{\varepsilon}}}\cdot\big(\nabla \sqrt{\rho^{\varepsilon}} \cdot \nabla\big)\Psi\Big)\,\mathrm{d} \mathbf{x} \mathrm{d} t\\
				&\quad -\varepsilon \int_{\mathbb{R}_{+}^{4}}
				\nabla \sqrt{\rho^{\varepsilon}} \cdot\big(\frac{\mathcal{M}^{\varepsilon}}{\sqrt{\rho^{\varepsilon}}} \cdot \nabla\big)
				\Psi\,\mathrm{d}\mathbf{x}\mathrm{d}t
+\int_{\mathbb{R}_{+}^{4}}\big(\rho^{\varepsilon} \nabla \Phi^{\varepsilon} \cdot \Psi\big)(t, \mathbf{x})\,\mathrm{d}\mathbf{x}.
		\end{align*}
		
		\item [\rm (\rmnum{4})] For any $t\geq 0$ and $\xi(\mathbf{x})\in C_{0}^1(\R^3)$,
		\begin{equation*}
			\int_{\R^3}\nabla \Phi^{\v}(t,\mathbf{x})\cdot \nabla\xi(\mathbf{x})\,\mathrm{d}\mathbf{x}=-\int_{\R^3}\rho^{\v}(t,\mathbf{x})\xi(\mathbf{x})\,\mathrm{d}\mathbf{x}.
		\end{equation*}
	\end{itemize}
\end{definition}
\noindent
Then we have
\smallskip
\begin{theorem}[Global existence of spherically symmetric solutions for CNSPEs]\label{thm1.2}
	Assume that the initial data function $(\rho_{0}^{\v},\mathcal{M}_{0}^{\v},\Phi_{0}^{\v})(\mathbf{x})$
	is given
	in \eqref{1.19}--\eqref{1.19-1} with $(\rho_{0}^{\v},m_{0}^{\v},\Phi_{0}^{\v})(r)$
	satisfying \eqref{1.27}--\eqref{1.27-2}.
	Then, for each fixed $\varepsilon \in (0, 1]$,
	there exists a global weak solution $(\rho^{\v},\mathcal{M}^{\v}, \Phi^{\v})(t,\mathbf{x})$
	of the Cauchy problem \eqref{1.8} and \eqref{1.19}--\eqref{1.19-1}
	in the sense of {\rm Definition \ref{def1.2}}
	with following spherically symmetric form{\rm :}
	\begin{equation}\label{1.22-2}
		(\rho^{\v},\mathcal{M}^{\v},\Phi^{\v})(t,\mathbf{x})
		=(\rho^{\v}(t,r),m^{\v}(t,r)\frac{\mathbf{x}}{r},\Phi^{\v}(t,r))\qquad \text{for }r=|\mathbf{x}|
	\end{equation}
	such that, for $t\geq 0$,
	\begin{align}\label{1.33}
&\int_{\mathbb{R}^{3}}\Big(\frac{1}{2}\Big|\frac{\mathcal{M}^{\varepsilon}}{\sqrt{\rho^{\varepsilon}}}\Big|^{2}+\rho^{\varepsilon} e(\rho^{\varepsilon})-\frac{1}{2}|\nabla\Phi^{\v}|^2\Big)\,\mathrm{d}\mathbf{x}
 \leq  \int_{\mathbb{R}^{3}}\Big(\frac{1}{2}\Big|\frac{\mathcal{M}_{0}^{\varepsilon}}{\sqrt{\rho_{0}^{\varepsilon}}}\Big|^{2}
		+\rho_{0}^{\varepsilon} e(\rho_{0}^{\varepsilon})-\frac{1}{2}|\nabla\Phi_{0}^{\v}|^2\Big)\,
		\mathrm{d}\mathbf{x}.
	\end{align}
	Furthermore, for $(\rho^{\varepsilon}, m^{\varepsilon}, \Phi^{\v})(t, r)$,
	there exists a measurable function $u^{\v}(t,r)$ with
	$$
	u^{\v}(t,r):=\frac{m^{\v}(t,r)}{\rho^{\v}(t,r)}\qquad {\it a.e.}\text{ on }\big\{(t,r)\,:\,\rho^{\v}(t,r)\neq 0\big\},
	$$
	and $u^{\v}(t,r):=0~{\it a.e.}\text{ on }\big\{(t,r)\,:\,\rho^{\v}(t,r)=0\text{ or }r=0\big\}$
	such that $m^{\v}(t,r)=(\rho^{\v}u^{\v})(t,r)$ {\it a.e.} on $\R_{+}^2:=\R_+\times\R_+$.
    Moreover, the following properties hold{\rm :}
	\begin{align}
		&({\rm \rmnum{1}})~ \int_{0}^{\infty} \rho^{\varepsilon}(t, r)\, r^{2}\mathrm{d}r=\int_{0}^{\infty} \rho_{0}^{\varepsilon}(r)\, r^{2}\mathrm{d} r=\frac{M}{\omega_{3}}
		\qquad \text {for } t \geq 0,\label{1.28}\\
		&({\rm \rmnum{2}})~ \int_{0}^{\infty} \eta^{*}(\rho^{\varepsilon}, m^{\varepsilon})(t, r)\, r^{2}\mathrm{d} r+\varepsilon \int_{\mathbb{R}_{+}^{2}}(\rho^{\varepsilon}|u^{\varepsilon}|^{2})(t, r)
		\,r^2\mathrm{d}r \mathrm{d} t + \|\nabla\Phi^{\v}\|_{L^2(\R^3)}\nonumber\\
		&\qquad +\|\Phi^{\v}\|_{L^{6}(\R^3)}+\int_{0}^{\infty}\Big(\int_{0}^{r} \rho^{\varepsilon}(t, z)\, z^{2}\mathrm{d}z \Big)
		\rho^{\varepsilon}(t, r)\, r\mathrm{d} r\leq C\left(M, E_{0}\right) \qquad\text {for } t \geq 0,\label{1.29}\\
		&({\rm \rmnum{3}})~ \sup_{t\in [0,T]}\varepsilon^{2} \int_{0}^{\infty}\left|\left(\sqrt{\rho^{\varepsilon}}\right)_{r}\right|^{2}\, r^{2}\mathrm{d} r+\varepsilon \int_{0}^{T} \int_{0}^{\infty}\frac{P'(\rho^{\v})}{\rho^{\v}}\left|\rho_{r}^{\varepsilon}\right|^{2}\, r^{2}\mathrm{d}r \mathrm{d}t\nonumber\\
		&\qquad\quad \leq C(M, E_{0}, T),\label{1.30}\\[1mm]
		&({\rm \rmnum{4}})~~ \int_{0}^{T} \int_{d}^{D} \rho^{\varepsilon}\left|u^{\varepsilon}\right|^{3}\,r^{2}\mathrm{d}r \mathrm{d} t \leq C(d, D, M, E_{0}, T),\label{1.31}\\
		&({\rm \rmnum{5}})~~ \int_{0}^{T} \int_{d}^{\infty}(\rho^{\v})^{\g_2+1}\, r^{2}\mathrm{d} r \mathrm{d} t
		\leq C(d, M, E_{0}, T),\label{1.32}
	\end{align}
	for any
	$T \in \R_{+}$ and
	interval $[d, D]\Subset (0, \infty)$,
	where $C(M, E_{0})$, $C(M, E_{0}, T)$, and $C(d, D, M, E_{0}, T)$
	are
positive constants independent of $\varepsilon$.
	In addition, for $\varepsilon\in (0,1]$,
	\begin{equation}\label{1.34}
		\partial_{t} \eta^{\psi}\left(\rho^{\varepsilon}, m^{\varepsilon}\right)+\partial_{r} q^{\psi}\left(\rho^{\varepsilon}, m^{\varepsilon}\right)
     \qquad \text { is compact in } W_{\operatorname{loc}}^{-1,p}(\mathbb{R}_{+}^{2})
	\end{equation}
	for any $p\in [1,2)$, where  $\psi(s)$ is any smooth function with compact support on $\mathbb{R}$.
\end{theorem}

\medskip
\begin{remark}
In this paper, we require the density-dependent viscosity coefficients $\mu(\rho)$ and $\lambda(\rho)$
to satisfy the BD entropy relation {\rm (}{\it cf.} {\rm \cite{Bresch-Desjardins-2002, Bresch-Desjardins-2004,
		Bresch-Desjardins-2007, Bresch-Desjardins-Lin-2003}}{\rm ):}
		\begin{align}\label{BD}
			\rho\mu'(\rho)=\mu(\rho)+\lambda(\rho),
		\end{align}
which is important for us to derive the estimate for the derivative of the density.
Under the physical restriction \eqref{1.8-1} and the BD entropy relation \eqref{BD}, $\lambda(\rho)$ cannot be a non-zero constant.
Since we focus mainly on the global existence of weak solutions for CEPEs by the vanishing viscosity limit of weak solutions of CNSPEs
which means the viscous terms will vanish eventually, we consider only the special case  $(\mu(\rho),\lambda(\rho))=(\rho,0)$ in the present paper,
which corresponds to  the well-known Saint-Venant model of shallow water.
	
Recently, in {\rm \cite{Bresch-Vasseur-Yu-2022,Ducomet-Necasova-Vasseur}},
the global existence of  weak solutions was established
for the compressible Navier-Stokes equations and CNSPEs for a class of general density-dependent viscous coefficients
satisfying the BD entropy relation, respectively.
Motivated by {\rm \cite{Bresch-Vasseur-Yu-2022,Ducomet-Necasova-Vasseur}},
it should be able to extend our results to a class of more general viscous coefficients.
However, for such general viscous coefficients $\mu(\rho)$ and $\lambda(\rho)$ satisfying the BD relation,
we have to check the uniform estimates of the solutions and the validity of vanishing viscosity limit $\v \rightarrow 0$
so that major modifications to our present paper are required, which is out of scope of this paper.
\end{remark}

\smallskip
Now we introduce the notion of finite-energy solutions of CEPEs \eqref{1.1}.

\smallskip
\begin{definition}\label{def1.1}
A measurable vector function $(\rho,\mathcal{M},\Phi)$ is said to be a finite-energy solution
of the Cauchy problem \eqref{1.1} and \eqref{1.2}--\eqref{1.2-1} provided that
\begin{itemize}
\item[\rm (\rmnum{1})]$\rho(t,\mathbf{x})\geq 0$ {\it a.e.},
  and $(\mathcal{M},\frac{\mathcal{M}}{\sqrt{\rho}})(t, \mathbf{x})=\mathbf{0}$
 {\it a.e.} on $\{(t,\mathbf{x})\in \R_{+}^{4} :\,\rho(t,\mathbf{x})=0\}$ {\rm (}vacuum{\rm )}.
\item[\rm (\rmnum{2})] For {\it a.e.} $t>0$, the total energy is finite{\rm :}
		\begin{equation}\label{1.16-1}
			\left\{\begin{aligned}
				&\int_{\mathbb{R}^{3}}\Big(\frac{1}{2}\Big|\frac{\mathcal{M}}{\sqrt{\rho}}\Big|^{2}
				+\rho e(\rho)+\frac{1}{2}|\nabla\Phi|^{2}\Big)(t, \mathbf{x}) \,\mathrm{d} \mathbf{x} \leq C(E_{0}, M), \\
				&\int_{\mathbb{R}^{3}}\Big(\frac{1}{2}\Big|\frac{\mathcal{M}}{\sqrt{\rho}}\Big|^{2}
				+\rho e(\rho)-\frac{1}{2}|\nabla \Phi|^{2}\Big)(t, \mathbf{x})\,\mathrm{d} \mathbf{x}
            \\				
            &\qquad
            \leq \int_{\mathbb{R}^{3}}\Big(\frac{1}{2}\Big|\frac{\mathcal{M}_{0}}{\sqrt{\rho_{0}}}\Big|^{2}
               +\rho_{0} e(\rho_{0})-\frac{1}{2}|\nabla \Phi_{0}|^{2}\Big)(\mathbf{x})\,\mathrm{d} \mathbf{x}.
			\end{aligned}\right.
		\end{equation}
\item[\rm (\rmnum{3})] For any $\zeta(t, \mathbf{x})\in C_{0}^{1}([0,\infty)\times \R^{3})$,
		\begin{equation}\label{1.17}
			\int_{\R_{+}^{4}}(\rho\zeta_{t}+\mathcal{M}\cdot \nabla\z)\,\mathrm{d}\mathbf{x}\mathrm{d}t+\int_{\R^3}(\rho_{0}\z)(0,\mathbf{x})\,\mathrm{d}\mathbf{x}=0.
		\end{equation}
\item[\rm (\rmnum{4})] For any $\Psi(t, \mathbf{x})=(\Psi_{1},\Psi_{2},\Psi_{3})(t,\mathbf{x})\in (C_{0}^{1}([0,\infty)\times \R^3))^{3}$,
		\begin{align}\label{1.18}
				&\int_{\R_{+}^{4}}\Big(\mathcal{M}\cdot \partial_{t}\Psi
				+\frac{\mathcal{M}}{\sqrt{\rho}}\cdot (\frac{\mathcal{M}}{\sqrt{\rho}}\cdot \nabla)\Psi+P(\rho)\,
				\nabla\cdot\Psi\Big)\,\mathrm{d}\mathbf{x}\mathrm{d}t+\int_{\R^3}\mathcal{M}_{0}(\mathbf{x})\cdot \Psi(0,\mathbf{x})\,\mathrm{d}\mathbf{x}\nonumber\\
				&
				=\int_{\mathbb{R}_{+}^{4}}(\rho \nabla\Phi \cdot \Psi)(t, \mathbf{x})\,\mathrm{d} \mathbf{x}.
		\end{align}
\item[\rm (\rmnum{5})] For any $\xi(\mathbf{x})\in C_{0}^{1}(\R^3)$,
		\begin{equation}\label{1.18-1}
			\int_{\mathbb{R}^{3}} \nabla \Phi(t, \mathbf{x}) \cdot \nabla\xi(\mathbf{x}) \,\mathrm{d}\mathbf{x}
			=- \int_{\mathbb{R}^{3}} \rho(t, \mathbf{x}) \xi(\mathbf{x})\,\mathrm{d}\mathbf{x}\qquad\,\, \text {for {\it a.e.} } t \geq 0.
		\end{equation}
	\end{itemize}
\end{definition}

\begin{remark}\label{rem1}
In the spherically symmetric form, {\rm Definition \ref{def1.1}} becomes the following:
A measurable vector function $(\rho,\mathcal{M},\Phi)(t,{\bf{x}})=(\rho(t,r),m(t,r)\frac{{\bf{x}}}{r},\Phi(t,r))$
is said to be a spherically symmetric finite-energy solution
of the Cauchy problem \eqref{1.1} and \eqref{1.11}--\eqref{1.11-1} provided that
\begin{itemize}
\item[\rm (\rmnum{1})] $\rho(t,r)\geq 0$ {\it a.e.}, and $(m,\frac{m}{\sqrt{\rho}})(t,r)=\mathbf{0}$ {\it a.e.}
   on $\{(t,r)\in \R_{+}^2\,:\,\rho(t,r)=0\}$ {\rm (}vacuum{\rm )}.
\item[\rm (\rmnum{2})] For {\it a.e.} $t>0$, the total energy is finite{\rm :}
			\begin{equation}\label{1.16-2}
				\left\{\begin{aligned}
					&\int_{0}^{\infty}\Big(\frac{1}{2}\Big|\frac{m}{\sqrt{\rho}}\Big|^{2}
					+\rho e(\rho)+\frac{1}{2}|\Phi_{r}|^{2}\Big)(t, r) \,r^2\mathrm{d} r \leq C(E_{0}, M), \\
					&\int_{0}^{\infty}\Big(\frac{1}{2}\Big|\frac{m}{\sqrt{\rho}}\Big|^{2}
					+\rho e(\rho)-\frac{1}{2}|\Phi_{r}|^{2}\Big)(t, r)\,r^2\mathrm{d}r\\
					&\qquad
\leq \int_{0}^{\infty}\Big(\frac{1}{2}\Big|\frac{m_{0}}{\sqrt{\rho_{0}}}\Big|^{2}+\rho_{0} e(\rho_{0})-\frac{1}{2}| \Phi_{0r}|^{2}\Big)(r)\,r^2\mathrm{d}r.
				\end{aligned}\right.
			\end{equation}
\item[\rm (\rmnum{3})] For any $\zeta(t, r)\in C_{0}^{1}([0,\infty)\times \R)$,
			\begin{align}\label{1.17-1}
				&\int_{\R_{+}^2}(\rho\zeta_{t}+m\zeta_{r})(t,r)\,r^2{\rm d}r{\rm d}t  +\int_{0}^{\infty}\rho_{0}(r)\zeta(0,r)\,r^2\mathrm{d}r=0.
			\end{align}
\item[\rm (\rmnum{4})] For any $\psi(t,r)\in C_{0}^{1}([0,\infty)\times \R)$ with $\psi(t,0)=0$ for all $t\geq 0$,
			\begin{align}\label{1.18-2}
				&\int_{\R_{+}^2}m(t,r)\psi_{t}(t,r)\,r^2{\rm d}r{\rm d}t +\int_{\R_{+}^2}\big(\frac{m^2}{\rho}\big)(t,r)\,\psi_{r}(t,r)\,r^2{\rm d}r{\rm d}t\nonumber\\
				&\quad  +\int_{\R_{+}^2}P(\rho(t,r))\,(\psi_{r}+\frac{2}{r}\psi)(t,r)\,r^2{\rm d}r{\rm d}t+\int_{0}^{\infty}m_{0}(r)\,\psi(0,r)\,r^2{\rm d}r\nonumber\\
				&=\int_{\R_{+}^2}(\rho\Phi_{r})(t,r)\,\psi(t,r)\,r^2\mathrm{d} r\mathrm{d}t.
			\end{align}
			\item[(\rmnum{5})] For any $\xi(r)\in C_{0}^{1}(\R)$ and {\it a.e. } $t\geq 0$,
			\begin{align}\label{1.18-5}
				&\int_{0}^{\infty}\Phi_{r}(t,r)\,\xi_{r}(r)\,r^2{\rm d}r=-\int_{0}^{\infty} \rho(t,r)\,\xi( r)\,r^2{\rm d}r.
			\end{align}
		\end{itemize}
\end{remark}

To establish the strong convergence of the inviscid
limit of
solutions $(\rho^{\varepsilon}, \mathcal{M}^{\varepsilon}, \Phi^{\v})(t,\mathbf{x})$
of CNSPEs \eqref{1.8} obtained in Theorem \ref{thm1.2} as $\v\to 0$,
we establish the following $L^p$ compensated compactness framework
for the 1-D Euler equations \eqref{2.10a} with general pressure law \eqref{1.3}--\eqref{1.5},
in which restriction $\gamma_{2}\in (\frac{6}{5}, \gamma_{1}]$ in \eqref{1.5}
can be relaxed to $\gamma_{2}\in (1, \gamma_1]$.

\smallskip

\begin{theorem}[$L^p$ compensated compactness framework]\label{thm1.3}
	Let
$$
(\rho^{\v},m^{\v})(t,r)=(\rho^\v, \rho^\v u^\v)(t,r)
$$
be a sequence of measurable functions
	with $\rho^{\v}\geq 0$ a.e. on $\R_{+}^2$
	satisfying the following two conditions{\rm :}
\begin{itemize}
\item[\rm (i)] For any $T>0$ and $K\Subset\mathbb{R}_+$, there exists $C(K,T)>0$ independent of $\v$ such that
$$
\int_0^T\int_K\big((\rho^\v)^{\g_2+1}+\rho^{\v}|u^{\v}|^3\big)\,\mathrm{d}r\mathrm{d}t\le C(K,T).
$$
\item[\rm (ii)] For any entropy pair $(\eta^\psi,q^\psi)$ defined in \eqref{2.21e}--\eqref{2.21q} with
any smooth function $\psi(s)$ of compact support on $\mathbb{R}$,	
$$
\partial_{t} \eta^{\psi}(\rho^{\varepsilon}, m^{\varepsilon})
+\partial_{r} q^{\psi}(\rho^{\varepsilon}, m^{\varepsilon})
\qquad\mbox{is compact in $ W_{\mathrm{loc}}^{-1,1}(\mathbb{R}_{+}^{2})$}.
$$
\end{itemize}
\noindent
Then there exists a subsequence $($still denoted$)$ $(\rho^{\v},m^{\v})(t,r)$ and a vector
	function $(\rho,m)(t,r)$ such that, as $\v\to 0$,
	\begin{equation}\label{1.35}
		\begin{aligned}
			&\rho^{\v}(t,r)\to \rho(t,r)~~ \mbox{in}~L^{q_1}_{\rm loc}(\R_{+}^2)\qquad\,\,\,\,\text{for }q_1\in [1,\g_2+1),\\
			& m^{\v}(t,r)\to m(t,r)~~ \mbox{in}~L^{q_2}_{\rm loc}(\R_{+}^2)\qquad \text{for }q_2\in [1,\frac{3(\g_2+1)}{\g_2+3}),
		\end{aligned}
	\end{equation}
	where $L_{\mathrm{loc}}^{p}(\mathbb{R}_{+}^{2})$ represents $L^{p}([0, T] \times K)$ for any $T>0$
	and compact set $K \Subset \mathbb{R}_+$.
\end{theorem}

\smallskip
Now, we are ready to state our main theorem.
\smallskip

\begin{theorem}[Global existence of finite-energy solutions]\label{thm1.1}
	Let the pressure function $P(\rho)$ satisfy \eqref{1.3}--\eqref{1.5}, and let the spherically symmetric initial
	data $(\rho_0,\mathcal{M}_0, \Phi_{0})(\mathbf{x})$ be given in \eqref{1.11}--\eqref{1.11-1}
	with $(\rho_{0},m_{0},\Phi_{0})(r)$ satisfying \eqref{1.14}--\eqref{1.14-1} and \eqref{1.15-1}.
	Assume that $\g_2>\frac{4}{3}$, or $M<M_{\rm c}$ as $\g_2\in (\frac{6}{5},\frac{4}{3}]$.
	Then there exists a global finite-energy solution $(\rho,\mathcal{M},\Phi)(t,\mathbf{x})$ of
	\eqref{1.1} and \eqref{1.11}--\eqref{1.11-1} with spherical symmetry form \eqref{1.10}
	in the sense of {\rm Definition \ref{def1.1}}.
\end{theorem}

\smallskip
\begin{remark}
	For the steady gaseous star problem, there is no white dwarf star if the total mass is larger
	than the so-called Chandrasekhar limit when $\g\in (\frac{6}{5},\frac{4}{3}]${\rm ;} see {\rm \cite{Chandrasekhar-1938}}.
	{\rm Theorem \ref{thm1.1}} requires similar restriction on the total mass when
	$\g_2\in (\frac{6}{5},\frac{4}{3}]$ for non-steady gaseous stars.
	Moreover, in view of \eqref{1.18-0}, for the non-steady white dwarf star, the critical mass is exactly the Chandrasekhar limit
	in the case that $P(\rho)=\kappa_2\rho^{\frac{4}{3}}$.
	It would be interesting to analyze whether
	the critical mass defined in \eqref{1.18-3}--\eqref{1.18-4} for $\gamma_2\in (\frac{6}{5},\frac{4}{3})$
	is optimal.
\end{remark}

\smallskip
\begin{remark}
	{\rm Theorem \ref{thm1.1}} can be extended to the $3$-D compressible Euler equations,
	{\it i.e.}, \eqref{1.1} with $\Phi=0$. Moreover, the inviscid
	limit from
	the compressible Navier-Stokes equations to
Euler equations
	with far-field vacuum can also be justified.
\end{remark}

\smallskip
\begin{remark}
	{\rm Theorem \ref{thm1.1}} also holds for the plasmas case, {\it i.e.},
	$k_{g}=-1$ in \eqref{1.1},
	by a similar proof. In this case, the restriction: $M<M_{\rm c}$ can be removed, and
	condition $\g_{2}>\frac{6}{5}$ can be relaxed to $\g_2>1$ if the additional
	assumption{\rm :} $\rho_{0}\in L^{\frac{6}{5}}(\R^3)$ is imposed. We omit the proof in this paper
	for brevity and, instead, refer the reader to {\rm \cite{Chen-He-Wang-Yuan-2021}}
	for details.
\end{remark}

\section{\,Properties of the General Pressure Law and Related Internal Energy}
In this section, we present some useful estimates
involving the general pressure $P(\r)$ with
\eqref{1.3}--\eqref{1.5}  and the corresponding internal energy $e(\rho)$,
which are used in the subsequent development.

Denote $c(\rho):=\sqrt{P'(\rho)}$ as the speed of sound, and
\begin{equation}\label{A.0}
	k(\rho):=\int_{0}^{\rho}\frac{\sqrt{P'(y)}}{y}\,\mathrm{d} y.
\end{equation}
By direct calculation, we can obtain the following asymptotic behaviors of $P(\rho)$, $e(\rho)$, and $k(\rho)$.

\smallskip
\begin{lemma}\label{lemA.1}
	Assume that $\rho_{*}$ given in \eqref{1.4} is small enough and $\rho^{*}$ given in \eqref{1.5}
	is large enough such that the following estimates hold{\rm :}
\begin{enumerate}
\item [\rm (\rmnum{1})] When $\rho\in (0,\rho_{*}]$,
		\begin{equation}\label{A.1-1}
			\left\{\begin{aligned}
				&\underline{\kappa}_{1}\rho^{\gamma_1}\leq P(\rho)\leq \bar{\kappa}_{1}\rho^{\g_1},\\
				&\underline{\kappa}_{1}\gamma_1\rho^{\gamma_1-1}\leq P'(\rho)\leq \bar{\kappa}_{1}\g_1\rho^{\g_1-1},\\
				&\underline{\kappa}_{1}\gamma_1(\gamma_1-1)\rho^{\gamma_1-2}\leq P''(\rho)\leq \bar{\kappa}_{1}\gamma_1(\gamma_1-1)\rho^{\g_1-2},
			\end{aligned}
			\right.
		\end{equation}
		and when $\rho\in [\rho^{*},\infty)$,
\begin{equation}\label{A.2-1}
\left\{\begin{aligned}
	&\underline{\kappa}_{2}\rho^{\gamma_2}\leq P(\rho)\leq  \bar{\kappa}_{2}\rho^{\g_2},\\
    &\underline{\kappa}_{2}\gamma_2\rho^{\gamma_2-1}\leq P'(\rho)\leq  \bar{\kappa}_{2}\g_2\rho^{\g_2-1},\\
    &\underline{\kappa}_{2}\gamma_2(\gamma_2-1)\rho^{\gamma_2-2}\leq P''(\rho)\leq  \bar{\kappa}_{2}\gamma_2(\gamma_2-1)\rho^{\g_2-2},
	\end{aligned}\right.
	\end{equation}
where we have denoted $\,\underline{\kappa}_{i}:=(1-\mathfrak{a}_0)\kappa_{i}$
and $\bar{\kappa}_{i}:=(1+\mathfrak{a}_0) \kappa_{i}$ with $\mathfrak{a}_0=\frac{3-\g_1}{2(\g_1+1)}$ and $i=1,2$.
		
\item[\rm (\rmnum{2})] For $e(\rho)$ and $k(\rho)$, there exists $C>0$ depending on
	$(\gamma_1, \gamma_2,  \k_1,\k_2, \rho_{*}, \rho^{*})$ such that
		\begin{align}
			&C^{-1}\rho^{\g_1-1}\leq e(\rho)\leq C\rho^{\gamma_1-1},\,\,
			C^{-1}\rho^{\g_1-2}\leq e'(\rho)\leq C\rho^{\g_1-2}
			\,\,\,\,\, \text{ for }\rho\in (0,\rho_{*}],\label{A.8-1}\\
			&C^{-1}\rho^{\g_2-1}\leq e(\rho)\leq C\rho^{\gamma_2-1},\,\,
			C^{-1}\rho^{\g_2-2}\leq e'(\rho)\leq C\rho^{\g_2-2}
			\,\,\,\,\, \text{ for }\rho\in [\rho^{*},\infty), \label{A.9-1}
		\end{align}
		and, for $i=0,1$,
		\begin{align}
			&\frac{\rho^{\theta_{1}-i}}{C}\leq k^{(i)}(\rho)\leq C\rho^{\theta_{1}-i},\,\,
			\frac{\rho^{\theta_{1}-2}}{C}\leq |k''(\rho)|\leq C\rho^{\theta_1-2}\,\,\,\,\, \text{for } \rho\in (0,\rho_{*}],\label{A.4-1}\\
			&\frac{\rho^{\theta_{2}-i}}{C}\leq k^{(i)}(\rho)\leq C\rho^{\theta_{2}-i},
			\,\, \frac{\rho^{\theta_{2}-2}}{C}\leq |k''(\rho)|\leq C\rho^{\theta_2-2} \,\,\,\,\,\text{for } \rho\in [\rho^{*},\infty),
			\label{A.5-1}
		\end{align}
		where $\t_{1}=\frac{\g_1-1}{2}$ and $\t_2=\frac{\g_2-1}{2}$.
	\end{enumerate}
\end{lemma}

\smallskip
It follows from \eqref{A.1-1}--\eqref{A.2-1} that
\begin{equation}\label{A.3}
	\frac{(3\g_1-1)(\g_1-1)}{\g_1+5}P'(\rho)\leq \rho P''(\rho)\leq \frac{(5+\g_1)(\g_1-1)}{3\g_1-1}P'(\rho)<2P'(\rho),
\end{equation}
when $\rho\in [0,\rho_{*}]\cup[\rho^{*},\infty)$. For later use, we denote
\begin{align}
\nu:=1-\frac{(3\gamma_{1}-1)(\gamma_{1}-1)}{2(5+\gamma_{1})}
	<1,
	\qquad d(\rho):=2+\frac{\rho k''(\rho)}{k'(\rho)}.\label{A.16}
\end{align}
Then it follows from \eqref{A.3} that
\begin{equation}\label{A.7}
	0<\Big\vert\frac{\rho k''(\rho)}{k'(\rho)}\Big\vert
	=1-\frac{\rho P''(\rho)}{2P'(\rho)}\leq \nu<1\qquad \text{for }\rho\in (0,\rho_{*}]\cup [\rho^{*},\infty).
\end{equation}

Motivated by \cite{Schrecker-Schulz-2020}, we have

\smallskip
\begin{lemma}\label{lemA.3}
	$0<d(\r)\leq C$ for all $\rho>0$,  and
	\begin{equation}\label{A.17-1}	
		\big|d(\r)-(1+\t_2)\big|\leq C\r^{-\epsilon}\qquad \mbox{for}\,\, \rho \gg1.
	\end{equation}
\end{lemma}

\noindent\textbf{Proof.} It follows from \eqref{1.3} that
$
d(\rho)=1+\frac{\rho P''(\rho)}{2P'(\rho)}>0.
$
Moreover, by \eqref{A.7}, it is direct to see that $d(\rho)$ is bounded.
Using \eqref{1.5}, we see that, for $\rho\geq \rho^{*}$,
$$
\left\{
\begin{aligned}
	&P'(\r)=\g_2\k_2\r^{\g_2-1}\big(1+\mathcal{P}_2(\r)+\r \mathcal{P}_2'(\r)\big),\\
	&P''(\r)=\g_2(\g_2-1)\k_2\r^{\g_2-2}\big(1+\mathcal{P}_2(\r)+3\r \mathcal{P}_2'(\r)+\r^2\mathcal{P}_2''(\r)\big).
\end{aligned}
\right.
$$
Then, for $\r\geq \max\{\r^{*},(8C^{*})^{1/\epsilon}\}$,
\begin{equation*}
	\big|d(\rho)-(1+\theta_2)\big|=\Big\vert\frac{\rho P''(\rho)}{2P'(\rho)}-\theta_2\Big\vert
	=\Big\vert\frac{\t_2\big(2\r \mathcal{P}_2'(\r)+\r^2\mathcal{P}_2''(\r)\big)}{1+P_2(\r)
      +3\r \mathcal{P}_2(\r)+\r^2\mathcal{P}_2''(\r)}\Big\vert\leq C(\t_2, C^{*})\r^{-\epsilon},
\end{equation*}
where we have used that
$|\mathcal{P}_2^{(j)}(\r)|\leq C^{*}\r^{-\epsilon-j}$ for $j=0,1,2,$ in the last inequality. $\hfill\square$

\smallskip
Hereafter, for simplicity of notation, we assume that \eqref{A.17-1} holds for $\rho\geq \rho^{*}$.
Furthermore, using $\eqref{1.5}$ and $e'(\rho)=\frac{P(\rho)}{\rho^2}$, we obtain that, for $\rho\geq \rho^{*}$,
\begin{equation}\label{ae1}
	e(\rho)
	=\frac{\kappa_{2}}{\g_2-1}\big(\rho^{\g_2-1}-(\rho^{*})^{\g_2-1}\big)
	+\k_2\int_{\rho^{*}}^{\rho}s^{\g_2-2}\mathcal{P}_{2}(s)\,\mathrm{d}s+\int_{0}^{\rho^{*}}\frac{P(s)}{s^2}\,\mathrm{d}s,
\end{equation}
which, with $e(0)=0$ and $|\mathcal{P}_{2}(\rho)|\leq C^{*}\rho^{-\epsilon}$, yields that,
for any parameter $\beta>0$,
\begin{equation*}
	\lim\limits_{\r\to 0}\rho^{\frac{\g_2-1}{5\g_2-6}}(\beta+e(\rho))^{-\frac{1}{5\g_2-6}}=0,
	\quad\,\,\lim\limits_{\r\to \infty}\rho^{\frac{\g_2-1}{5\g_2-6}}(\beta+e(\rho))^{-\frac{1}{5\g_2-6}}
	=(\frac{\k_2}{\g_2-1})^{-\frac{1}{5\g_2-6}}.
\end{equation*}
Then we see that
\begin{align}
	&C_{\max}(\beta):=\sup_{\rho\geq 0}\rho^{\frac{\g_2-1}{5\g_2-6}}(\beta+e(\rho))^{-\frac{1}{5\g_2-6}}\in [(\frac{\k_2}{\g_2-1})^{-\frac{1}{5\g_2-6}},\,\infty),\label{A.10-0}\\
	&\rho^{\frac{6(\g_2-1)}{5\g_2-6}}(\beta \rho+\rho e(\rho))^{-\frac{1}{5\g_2-6}}\leq C_{\max}(\beta)\rho\qquad\, \text{for }\rho >0.\label{A.10}
\end{align}

With a careful analysis of
$C_{\max}(\beta)$, we obtain some estimates
of $M_{\rm c}$ defined in \eqref{1.18-3}.

\smallskip
\begin{proposition}\label{prop3.1}
		Let $h(\rho)=P(\rho)\rho^{-1}-(\g_2-1)e(\rho)$,
		and let $\accentset{\sim}{M}_{\rm c}$ be the critical mass obtained in {\rm \cite[(2.8)] {Chen-He-Wang-Yuan-2021}}
		for the polytropic gases in \eqref{1.6-1} with $\gamma\in (\frac{6}{5},\frac{4}{3})$.
		Then $M_{\rm c}$ defined in \eqref{1.18-3}--\eqref{1.18-4} satisfies that
		$M_{\rm c}\leq \accentset{\sim}{M}_{\rm c}${\rm ;} in particular,
		$M_{\rm c}<\accentset{\sim}{M}_{\rm c}$ when
		$h'(\rho)> 0$ for all $\rho>0$.
		For example,
		\begin{equation}\label{A.14}
			P_{\delta}(\rho):=\int_{0}^{\rho^{\frac{1}{3}}}\frac{s^4}{\sqrt{\delta+s^{2+\epsilon_{0}}}}\,\mathrm{d}s\qquad\,\,
			\text{for $\delta>0$ and $\epsilon_{0}\in (0,\frac{4}{5})$}
		\end{equation}
		satisfies conditions \eqref{1.3}--\eqref{1.5}.
		If $M_{\rm c}(\delta)$ is the critical mass defined in  \eqref{1.18-3}--\eqref{1.18-4}
		for pressure $P_{\delta}(\rho)$, then $M_{\rm c}(\delta)<\accentset{\sim}{M}_{\rm c}$ for any $\delta>0$.	
\end{proposition}

\smallskip
\noindent
\textbf{Proof.}
	For $\g_2\in (\frac{6}{5},\frac{4}{3})$, it follows from \eqref{1.18-3}--\eqref{1.18-4}
	and \eqref{A.10-0} that, for any fixed $\beta>0$,
	\begin{align}
		M_{\rm c}(\beta)
		&= \Big(\frac{2}{9(\g_2-1)}(C_{\max}(\beta))^{\frac{5\g_2-6}{3(\g_2-1)}}\omega_{3}^{-\frac{4-3\g_2}{3(\g_2-1)}}
		\omega_{4}^{-\frac{2}{3}}\Big)^{-\frac{3(\g_2-1)}{5\g_2-6}}
		\Big(\frac{E_{0}+\omega_{3}^{-1}\beta M_{\rm c}(\beta)}{4-3\g_2}\Big)^{-\frac{4-3\g_2}{5\g_2-6}}\nonumber\\
		&\leq \Big(\frac{2}{9(\g_2-1)}\Big(\frac{\kappa_2}{\g_2-1}\Big)^{-\frac{1}{3(\g_2-1)}}\omega_{3}^{-\frac{4-3\g_2}{3(\g_2-1)}}\omega_{4}^{-\frac{2}{3}}\Big)^{-\frac{3(\g_2-1)}{5\g_2-6}}
		\Big(\frac{E_{0}}{4-3\g_2}\Big)^{-\frac{4-3\g_2}{5\g_2-6}}\nonumber\\
		&=\accentset{\sim}{M}_{\rm c},\label{A.15}
	\end{align}
	which yields that $M_{\rm c}\leq \accentset{\sim}{M}_{\rm c}$.

Let $g(\rho):=\rho^{\frac{\g_2-1}{5\g_2-6}}\big(\beta+e(\rho)\big)^{-\frac{1}{5\g_2-6}}$.
Then $C_{\max}(\beta)=\max_{\rho\geq 0}g(\rho)$.
Since $e'(\rho)=\frac{P(\rho)}{\rho^2}$, a direct calculation shows that
\begin{equation}\label{A.31}
	g'(\rho)=\frac{1}{5\g_2-6}\rho^{5-4\g_2}\big(\beta+e(\rho)\big)^{-\frac{5(\g_2-1)}{5\g_2-6}}
	\big((\g_2-1)\beta-h(\rho)\big).
\end{equation}

{
	If $h'(\rho)> 0$ for all $\rho> 0$, then $h(\rho)\geq h(0)=0$. Let $K_{0}:=\max_{\rho> 0}h(\rho)>0$.
	For $\beta$ small enough such that $0<\beta<\frac{K_{0}}{\g_2-1}$, there exists a unique
	point $\rho_{\beta}>0$ such that $g'(\rho_{\beta})=0$, {\it i.e.},
	\begin{equation}\label{A.30}
		h(\rho_{\beta})=(\g_2-1)\beta,
	\end{equation}
	and
	$
	C_{\max}(\beta)=g(\rho_{\beta})
	=(\g_2-1)^{\frac{1}{5\g_2-6}}
	\big(P(\rho_{\beta})\rho_{\beta}^{-\g_2}\big)^{-\frac{1}{5\g_2-6}}.
	$
	Moreover, it follows from \eqref{A.30} that $\lim_{\beta\to 0+}\rho_{\beta}=0$.
	Thus, we see from \eqref{1.4} that
	$$
	\lim\limits_{\beta\to 0+}C_{\max}(\beta)
	=(\g_2-1)^{\frac{1}{5\g_2-6}}
	\lim\limits_{\rho_{\beta}\to 0}\big(P(\rho_{\beta})\rho_{\beta}^{-\g_2}\big)^{-\frac{1}{5\g_2-6}}=\infty,
	$$
	which, with $\eqref{A.15}_{1}$, implies that $\lim_{\beta\to 0+}M_{\rm c}(\beta)=0$.}

{
	On the other hand, it follows directly from \eqref{A.10-0} and $\eqref{A.15}_{1}$
	that $\lim_{\beta\to \infty}M_{\rm c}(\beta)=0$.
	Therefore, the maximum
	value $M_{\rm c}$ of $M_{\rm c}(\beta)$ must be attained
	at some point $\beta_{0}\in (0,\infty)$ with
	$M_{\rm c}=M_{\rm c}(\beta_0)< \accentset{\sim}{M}_{\rm c}$ due to $\eqref{A.15}_{1}$.}

{
	For the pressure function $P_{\delta}(\rho)$ in \eqref{A.14}, it is direct to check that
	conditions \eqref{1.3}--\eqref{1.5} are satisfied
	and $\g_2=\frac{4}{3}-\frac{\epsilon_{0}}{6}\in (\f65,\f43)$.
	Let $e_{\delta}(\rho)$ be the corresponding internal energy with
	$e_{\delta}'(\rho)=\frac{P_{\delta}(\rho)}{\rho^2}$ and $e_{\delta}(0)=0$,
	and
	$h_{\delta}(\rho):=P_{\delta}(\rho)\rho^{-1}-(\g_2-1)e_{\delta}(\rho)$.
	It follows from a direct calculation that
	$$
	h_{\delta}'(\rho)=\rho^{-2}\big(-\g_2P_{\delta}(\rho)+\rho P_{\delta}'(\rho)\big)
	=:\rho^{-2}T_{\delta}(\rho),
	$$
	$T_{\delta}(0)=0$, and
	$
	T_{\delta}'(\rho)=-(\g_2-1)P_{\delta}'(\rho)+\rho P_{\delta}''(\rho)=\frac{2+\epsilon_{0}}{18}\delta\rho^{\frac{2}{3}}(\delta+\rho^{\frac{2+\epsilon_{0}}{3}})^{-\frac{3}{2}}>0
	$
	for $\rho>0$.
	Thus, $T_{\delta}(\rho)>0$ for $\rho>0$, which implies that $h_{\delta}'(\rho)> 0$ for $\rho>0$,
	so that $M_{\rm c}(\delta)<\accentset{\sim}{M}_{\rm c}$ for any $\delta>0$.$\hfill\square$}

\section{\,Entropy Analysis: Weak Entropy Pairs}
Compared with the polytropic gas case in \cite{Chen-He-Wang-Yuan-2021},
there is no explicit formula of the entropy kernel for the general pressure law
\eqref{1.3}--\eqref{1.5}
so that we have to analyze the entropy equation \eqref{1.25-1} carefully to obtain several desired estimates.

\subsection{\,A special entropy pair}
\begin{figure}
	\centering
	\includegraphics[width=0.8\linewidth]{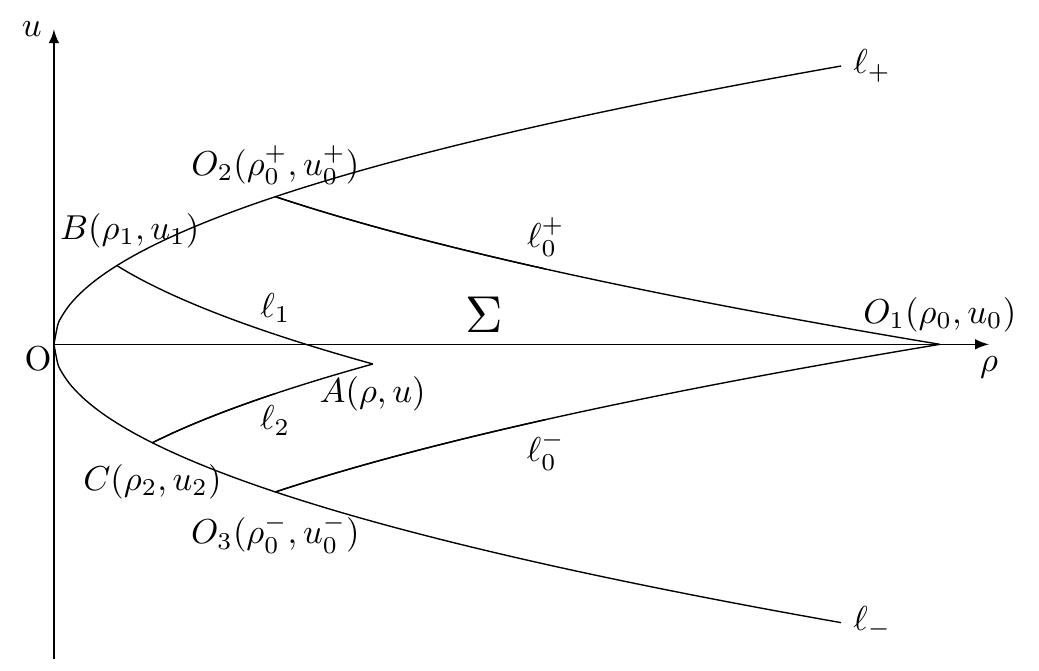}
	\caption{The schematic diagram of the characteristic curves of \eqref{4.6}}
	\label{fig-1}
\end{figure}

\FloatBarrier
In order to obtain the higher integrability of the velocity, we are going to construct a special entropy pair
such that $\rho|u|^3$ can be controlled by the entropy flux.
Indeed, such a special entropy $\hat{\eta}(\rho,u)$ is constructed as
\begin{equation}\label{4.1-0}
	\hat{\eta}(\rho,u)=\begin{cases}
		\frac{1}{2}\rho u^2+\rho e(\rho)\quad&\text{for }u\ge k(\rho),\\
		-\frac{1}{2}\rho u^2-\rho e(\rho)\,\,\, &\text{for }u\le -k(\rho),
	\end{cases}
\end{equation}
for $k(\rho)=\int_{0}^{\rho}\frac{\sqrt{P'(y)}}{y}\,\mathrm{d}y$
and, in the intermediate region $-k(\rho)\le u\le k(\rho)$,
$\hat{\eta}(\rho,u)$ is the unique solution of the Goursat problem
of the entropy equation \eqref{1.25-1}:
\begin{equation}\label{4.1}
	\left\{\begin{aligned}
		\dis&\eta_{\rho\rho}-k'(\rho)^2\eta_{uu}=0,\\
		\dis&\eta(\rho,u)\vert_{u=\pm k(\rho)}=\pm \big(\frac{1}{2}\rho u^2+\rho e(\rho)\big){\rm ;}
	\end{aligned}
	\right.
\end{equation}
see Fig. \ref{fig-1}.
Set
\begin{equation}\label{4.5}
	V_1=\frac{1}{2k'(\rho)}\big(\eta_{\rho}+k'(\rho)\eta_{u}\big),\qquad V_2=\frac{1}{2k'(\rho)}\big(\eta_{\rho}-k'(\rho)\eta_u\big).
\end{equation}
Then \eqref{4.1} can be rewritten as
\begin{equation}\label{4.6}
	\begin{cases}
		\dis\frac{\partial V_1}{\partial \rho}-k'(\rho)\frac{\partial V_1}{\partial u}
		=-\frac{k''(\rho)}{2k'(\rho)}(V_1+V_2),\\[4mm]
		\dis\frac{\partial  V_2}{\partial \rho}+k'(\rho)\frac{\partial V_2}{\partial u}=-\frac{k''(\rho)}{2k'(\rho)}(V_1+V_2).
	\end{cases}
\end{equation}
The corresponding characteristic boundary conditions become
\begin{equation}\label{4.6-1}
	\left\{\begin{aligned}
		&V_1\vert_{u=\pm k(\rho)}=\pm \frac{1}{2k'(\rho)}\big(\frac{1}{2}u^2+e(\rho)+\rho e'(\rho)\big)\pm \frac{1}{2}\rho u,\\
		&V_2\vert_{u=\pm k(\rho)}=\pm \frac{1}{2k'(\rho)}\big(\frac{1}{2}u^2+e(\rho)+\rho e'(\rho)\big)\mp \frac{1}{2}\rho u.
	\end{aligned}
	\right.
\end{equation}
Since $\frac{k''(\rho)}{k'(\rho)}$ has the singularity at vacuum $\rho=0$,
the Goursat problem \eqref{4.6}--\eqref{4.6-1} is singular, which requires a careful analysis.

It follows from \eqref{4.6} that there exist two characteristic curves originating from origin $O(0,0)$
in the $(\rho,u)$--plane:
\begin{equation}\label{4.2}
	\ell_{+}:=\{(\rho,u)\, :\, u=k(\rho)\},  \qquad \ell_{-}:=\{(\rho,u)\,:\ u=-k(\rho)\}.
\end{equation}
For any given point $O_1(\rho_0,u_0)$ with $u_0=0$, we can draw two backward characteristic curves
$\ell_{0}^{\pm}$
through $O_1(\rho_0,u_0)$; see Fig. \ref{fig-1}.
Let $O_2(\rho_{0}^{+},u_{0}^{+})$ be the intersection point of $\ell_{0}^{+}$ and $\ell_{+}$,
and let $O_3(\rho_{0}^{-},u_{0}^{-})$ be the intersection point of $\ell_{0}^{-}$ and $\ell_{-}$.
Let $\Sigma$ be the region surrounded by arc $\widehat{OO_2O_{1}O_{3}}$,
and let $\overline{\Sigma}$ be the closure of $\Sigma$.

\smallskip
\begin{lemma}\label{lem4.1}
	The Goursat problem \eqref{4.1} admits a unique solution $\hat{\eta}\in C^2(\R_{+} \times \R)$ such that
	\begin{itemize}
		\item[\rm (\rmnum{1})]
		$|\hat{\eta}(\rho,u)|\leq C\big(\rho|u|^2+\rho^{\gamma(\rho)}\big)$ for $(\rho,u)\in \R_{+}\times \R$,
		where $\gamma(\rho)=\gamma_1$ if $\rho\in [0,\rho_{*}]$ and $\gamma(\rho)=\gamma_2$ if $\rho\in (\rho_{*},\infty)$.
		
		\smallskip
		\item[\rm (\rmnum{2})] If $\hat{\eta}$ is regarded as a function of $(\r, u)$,
		\begin{align*}
			|\hat{\eta}_{\rho}(\rho,u)|\leq C\big(|u|^2+\rho^{2\theta(\rho)}\big),\,\,\,
           |\hat{\eta}_{u}(\rho,u)|\leq C\big(\rho|u|+\rho^{\theta(\rho)+1}\big)
			\quad\,\,\,\,\text{for }(\rho,u)\in \R_{+}\times \R,
		\end{align*}
		and, if $\hat{\eta}$ is regarded as a function of $(\r, m)$,
		\begin{align*}
		 |\hat{\eta}_{\rho}(\r,m)|\leq C\big(|u|^2+\rho^{2\t(\r)}\big),\,\,\,|\hat{\eta}_{m}(\r,m)|\leq C(|u|+\rho^{\t(\r)})
			\quad\,\,\,\,\text{for $(\rho,m)\in \R_{+}\times \R$},
         \end{align*}
        where $\theta(\rho):=\frac{\gamma(\rho)-1}{2}$.
		
		\smallskip
		\item[\rm (\rmnum{3})] If $\hat{\eta}_{m}$ is regarded as a function of $(\r,u)$,
		\begin{align*}
			\qquad |\hat{\eta}_{m\rho}(\rho,u)|\leq C\rho^{\theta(\rho)-1},\quad |\hat{\eta}_{mu}(\rho,u)|\leq C,
		\end{align*}
		and, if $\hat{\eta}_{m}$ is regarded as a function of $(\r,m)$,
		$$
		|\hat{\eta}_{m\rho}(\rho,m)|\leq C\rho^{\theta(\rho)-1},\quad |\hat{\eta}_{mm}(\rho,m)|\leq C\rho^{-1}.
		$$
		
		\item[\rm (\rmnum{4})]  If $\hat{q}$ is the corresponding entropy flux determined by \eqref{1.23},
		then $\hat{q}\in C^2(\R_{+}\times \R)$ and
		\begin{align*}
			&\hat{q}(\rho,u)=\frac{1}{2}\rho |u|^3\pm \rho u\big(e(\rho)+\rho e'(\rho)\big)\qquad\quad\text{for }\pm u\geq k(\rho),\\
			&|\hat{q}(\rho,u)|\leq C\r^{\g(\rho)+\t(\r)}\qquad\qquad\qquad\qquad\,\,\,\,\,\,\text{for }|u|<k(\rho),\\
			&\hat{q}(\rho,u)\geq \frac{1}{2}\rho |u|^3\qquad\qquad\qquad\qquad\qquad\quad\,\,\,\text{for } |u|\geq k(\rho),\\
			&|\hat{q}-u\hat{\eta}|\leq C\big(\rho^{\g(\rho)}|u|+\rho^{\g(\rho)+\t(\rho)}\big)\qquad\qquad\text{for } (\rho,u)\in \R_{+}\times \R.
		\end{align*}
	\end{itemize}
\end{lemma}

\noindent {\bf Proof.} To prove that \eqref{4.1} has a unique $C^2$--solution $\hat{\eta}$
in $\R_{+}\times \R$, it suffices to prove that \eqref{4.6}--\eqref{4.6-1} admits
a unique $C^1$--solution  $(V_1, V_2)$ in $\Sigma$ for any given point $O(\rho_0, u_0)$.
We use the Picard iteration and divide the proof into six steps.

\smallskip
1.  For any point $A(\rho,u)\in \Sigma$, there are two backward characteristic curves
through $A(\rho, u)$:
\begin{equation}\label{4.3}
	\begin{aligned}
		&\ell_1:=\big\{(s,\,u^{(1)}(s))\,:\,u^{(1)}(s)=-k(s)+u+k(\rho),\, 0< s\leq \rho\big\},\\
		&\ell_2:=\big\{(s,\,u^{(2)}(s))\,:\,u^{(2)}(s)=-k(\rho)+u+k(s),\, 0< s\leq \rho\big\}.
	\end{aligned}
\end{equation}
Let $B(\rho_1,u_1)$ be the intersection point of $\ell_{+}$ and $\ell_1$,
and let $C(\rho_2,u_2)$ be the intersection point of $\ell_{-}$ and $\ell_2$.
It follows from \eqref{4.2}--\eqref{4.3} that
\begin{equation}\label{4.4}
	u_{1}=k(\rho_1)=\frac{k(\rho)+u}{2},\qquad u_2=-k(\rho_2)=\frac{u-k(\rho)}{2}.
\end{equation}
Using \eqref{4.6} and integrating $V_1$ and $V_2$ along the characteristic curves
$\ell_{1}$ and $\ell_2$ respectively, we have
\begin{equation}\label{4.7}
	V_i(\rho,u)
	=V_i(\rho_i,u_i)-\int_{\rho_i}^{\rho}\frac{k''(s)}{2k'(s)}\,
	\sum\limits_{j=1}^2V_j(s,{u}^{(i)}(s))\,\mathrm{d}s\qquad \mbox{for $i=1,2$}.
\end{equation}
Denote $V_{i}^{(0)}(\rho,u):=V_{i}(\rho_{i},u_{i})$. It follows from \eqref{4.6-1} and \eqref{4.4} that
\begin{equation}\label{4.8}
	V_i^{(0)}(\rho,u)=(-1)^{i+1}\Big(\frac{1}{2}\rho_ik(\rho_i)+\frac{1}{4}\frac{k^2(\rho_i)}{k'(\rho_i)}+\frac{e(\rho_i)+\rho_ie'(\rho_i)}{2k'(\rho_i)}\Big)\qquad \mbox{for $i=1,2$}.
\end{equation}
We define the iterated scheme:
\begin{equation}\label{4.9}
	V_i^{(n+1)}(\rho,u):=V_i(\rho_i,u_i)-\int_{\rho_i}^{\rho}\frac{k''(s)}{2k'(s)}\, \sum\limits_{j=1}^{2}V_j^{(n)}(s,u^{(i)}(s))\, \mathrm{d}s\qquad \mbox{for $i=1,2$}.
\end{equation}
Then we obtain two sequences $\{V_{i}^{(n)}\}_{n=0}^{\infty}$ for $i=1,2$.
We now prove that $\{V_{i}^{(n)}\}_{n=0}^{\infty}$ are uniformly convergent
in $\overline{\Sigma}$, which is equivalent to proving
that
\begin{equation}\label{4.9-1}
	V_{i}^{(0)}(\rho,u)+\sum\limits_{n=1}^{\infty}\big(V_{i}^{(n)}-V_{i}^{(n-1)}\big)(\rho,u),
	\quad \mbox{$i=1,2$},
\end{equation}
are uniformly convergent in $\overline{\Sigma}$.

From Lemma \ref{lemA.1} and \eqref{4.8}, we know that $V_{i}^{(0)}, i=1,2$,
are continuous in $\overline{\Sigma}$ and there exists a constant $C_1>0$ depending only on $\rho_{*}$ and $\rho^{*}$ such that
\begin{align*}
	|V_i^{(0)}(\rho,u)|
	\leq
	\begin{cases}
		C_{1}\rho_{i}^{1+\theta_{1}}\,\,\, &\text{for }\rho_{i}\leq \rho_{*},\\
		C_{1}\rho_{i}^{1+\theta_{2}}\,\,\, &\text{for }\rho_{i}\geq \rho_{*},
	\end{cases}
	\qquad i=1,2,
\end{align*}
which, with the fact that $\rho_{i}\leq \rho$, yields that
\begin{align}\label{4.17}
	|V_i^{(0)}(\rho,u)|&\leq
	\begin{cases}
		C_{1}\rho^{1+\theta_{1}}\quad &\text{for }\rho\leq \rho_{*},\\
		\tilde{C}_{1}\rho^{1+M_{1}}\,\,\, &\text{for }\rho_{*}\leq \rho\leq \rho^{*},\\
		\hat{C}_{1}\rho^{1+\theta_2}\quad &\text{for }\rho\geq \rho^{*},
	\end{cases}
	\qquad i=1,2,
\end{align}
where $\tilde{C}_{1}\geq C_{1}(\rho_{*})^{\theta_{2}-M_{1}}$, $\hat{C}_{1}\geq C_{1}$,
and $M_{1}$ are positive constants to be chosen later.

It follows from \eqref{A.16} and \eqref{A.7} that there exist a constant $\nu<1$
and a constant $C_{0}\gg 1$ depending on $\rho_{*}$ and $\rho^{*}$ such that
\begin{equation}\label{4.18}
	\Big\vert\frac{k''(\rho)}{k'(\rho)}\Big\vert\leq
	\begin{cases}
		\nu\rho^{-1}\quad&\text{for $0<\rho\leq \rho_{*}$ and $\rho\geq \rho^{*}$},\\[0.5mm]
		C_{0}\rho^{-1}\,\,\,&\text{for $\rho_{*}< \rho < \rho^{*}$}.
	\end{cases}
\end{equation}

For the estimate of $|V_{i}^{(1)}-V_{i}^{(0)}|$, we divide it into six cases:

\smallskip
{\it Case 1. $\rho_{i}\leq \rho\leq \rho_{*}$}: It follows from \eqref{4.9}
and \eqref{4.17}--\eqref{4.18} that
\begin{equation}\label{4.10}
	\begin{aligned}
		\big|(V_i^{(1)}-V_{i}^{(0)})(\rho,u)\big|
		&\leq \int_{\rho_i}^{\rho} C_{1}\nu\,s^{\t_1}\,\mathrm{d}s
		\leq C_1
		\rho^{1+\theta_1}\,\ta_1,
	\end{aligned}
\end{equation}
where $\ta_1:=\frac{\nu}{1+\theta_1}\in (0,1)$.

\smallskip
{\it Case 2. $\rho_{i}\leq \rho_{*}\leq \rho\leq \rho^{*}$}: Then
\begin{align}\label{4.10-1}
	&\big|(V_i^{(1)}-V_{i}^{(0)})(\rho,u)\big|
	\leq \Big(\int_{\rho_{i}}^{\rho_{*}}+\int_{\rho_{*}}^{\rho}\Big)
	\Big\vert\frac{k''(s)}{2k'(s)}\Big\vert\,
	\sum\limits_{j=1}^2\big|V_{j}^{(0)}(s,u^{(i)}(s))\big|\,\mathrm{d}s
	\nonumber\\
	&\leq \int_{\rho_i}^{\rho_{*}}\frac{\nu}{2s}\, (2C_{1}s^{1+\t_1})\,\mathrm{d}s+\int_{\rho_{*}}^{\rho}\frac{C_{0}}{2s}\, (2\tilde{C}_{1}s^{1+M_{1}})\,\mathrm{d}s\nonumber\\[1mm]
	&\leq \tilde{C}_{1}\big(\rho^{1+M_{1}}-(\rho_{*})^{1+M_{1}}\big)\ta_{M_1}
	+C_{1}(\rho_{*})^{1+\theta_1}\ta_1   \leq \tilde{C}_{1}\rho^{1+M_{1}}\,\ta_{M_1},
\end{align}
where $\ta_{M_1}:=\frac{C_0}{1+M_{1}}$
and,  in the last inequality of \eqref{4.10-1},  we have chosen
\begin{equation}\label{Ch1}
	\tilde{C}_{1}\geq C_{1}(\rho_{*})^{\theta_1-M_{1}}\ta_1\ta_{M_1}^{-1}.
\end{equation}

\smallskip
{\it Case 3. $\rho_{*}\leq \rho_{i}\leq \rho\leq \rho^{*}$}:  It is direct to see that
\begin{align}\label{4.10-2}
	\big|(V_i^{(1)}-V_{i}^{(0)})(\rho,u)\big|
	\leq \int_{\rho_{i}}^{\rho}\tilde{C}_{1}C_{0}\,s^{M_1}\,{\rm d}s
	\leq \tilde{C}_{1}
	\rho^{1+M_{1}}\,\ta_{M_1}.
\end{align}

{\it Case 4. $\rho_{i}\leq \rho_{*}<\rho^{*}\leq \rho$}:  Then
\begin{align}\label{4.10-3}
	&\big|(V_i^{(1)}-V_{i}^{(0)})(\rho,u)\big|
	\leq \Big(\int_{\rho_{i}}^{\rho_{*}}+\int_{\rho_{*}}^{\rho^{*}}+\int_{\rho^{*}}^{\rho}\Big)
	\Big\vert\frac{k''(s)}{2k'(s)}\Big\vert\,
	\sum\limits_{j=1}^2\big|V_{j}^{(0)}(s,u^{(i)}(s))\big|\,\mathrm{d}s
	\nonumber\\
	&\leq \int_{\rho_{i}}^{\rho_{*}}\frac{\nu}{2s}\, (2C_{1}s^{1+\theta_{1}})\,\mathrm{d}s
	+\int_{\rho_{*}}^{\rho^{*}}\frac{C_{0}}{2s}\,(2\tilde{C}_{1}s^{1+M_{1}})\,\mathrm{d}s
	+\int_{\rho^{*}}^{\rho}\frac{\nu}{2s}\, (2\hat{C}_{1}s^{1+\theta_{2}})\,\mathrm{d}s\nonumber\\[1mm]
	&\leq \hat{C}_{1}\big(\rho^{1+\theta_2}-(\rho^{*})^{1+\theta_2}\big)\ta_2
	+\tilde{C}_{1}\big((\rho^{*})^{1+M_{1}}-
	(\rho_{*})^{1+M_{1}}\big)\ta_{M_1}
	+C_{1}(\rho_{*})^{1+\theta_{1}}\ta_1\nonumber\\
	&\leq \hat{C}_{1}\rho^{1+\theta_2}\,\ta_2,
\end{align}
where $\ta_2:=\frac{\nu}{1+\theta_2}\in (0,1)$ and, in the last inequality of \eqref{4.10-3},  we have used \eqref{Ch1} and chosen
\begin{equation}\label{Ch2}
	\hat{C}_{1}\geq \tilde{C}_{1}(\rho^{*})^{M_{1}-\theta_{2}}\ta_{M_1}\ta_2^{-1}.
\end{equation}

\smallskip
{\it Case 5. $\rho_{*}\leq \rho_{i}\leq \rho^{*}\leq \rho$}: It follows similarly that
\begin{align}\label{4.10-4}
	&\big|(V_i^{(1)}-V_{i}^{(0)})(\rho,u)\big|
	\leq \Big(\int_{\rho_{i}}^{\rho^{*}}+\int_{\rho^{*}}^{\rho}\Big)
	\Big\vert\frac{k''(s)}{2k'(s)}\Big\vert\,\sum\limits_{j=1}^2\big|V_{j}^{(0)}(s,u^{(i)}(s))\big|\,\mathrm{d}s
	\nonumber\\
	&\leq \hat{C}_{1}\big(\rho^{1+\theta_2}-
	(\rho^{*})^{1+\theta_2}\big)\ta_2+\tilde{C}_{1}(\rho^{*})^{1+M_{1}}\ta_{M_1}
	\leq \hat{C}_{1}\rho^{1+\theta_2}\,\ta_2,
\end{align}
where we have used \eqref{Ch2} in the last inequality of \eqref{4.10-4}.

\smallskip
{\it Case 6. $\rho^{*}\leq \rho_{i}\leq \rho$}: We see that
\begin{align}\label{4.10-5}
	\big|(V_i^{(1)}-V_{i}^{(0)})(\rho,u)\big|
	& \leq \int_{\rho_{i}}^{\rho} \hat{C}_{1}\nu \, s^{\theta_{2}}\,\mathrm{d}s
	\leq \hat{C}_{1}\rho^{1+\theta_2}\,\ta_2.
\end{align}
Combining \eqref{4.10}--\eqref{4.10-5}, we obtain
\begin{equation}\label{4.10-6}
	\big|(V_i^{(1)}-V_{i}^{(0)})(\rho,u)\big|\leq
	\begin{cases}
		C_{1}\rho^{1+\theta_{1}}\,\ta_1\quad&\text{for }\rho\leq \rho_{*},\\
		\tilde{C}_{1}\rho^{1+M_{1}}\,\ta_{M_1}\,\,\, &\text{for }\rho_{*}\leq \rho\leq \rho^{*},\\
		\hat{C}_{1}\rho^{1+\theta_2}\,\ta_2 \quad&\text{for }\rho\geq \rho^{*},
	\end{cases}
	\quad i=1,2,
\end{equation}
if \eqref{Ch1} and \eqref{Ch2} hold.

To utilize the induction arguments, we make the induction assumption for $n=k$:
\begin{equation}\label{4.11-1}
	\big|(V_i^{(k)}-V_{i}^{(k-1)})(\rho,u)\big|
	\leq
	\begin{cases}
		C_{1}\rho^{1+\theta_{1}}\,\ta_1^k\quad&\text{for }\rho\leq \rho_{*},\\
		\tilde{C}_{1}\rho^{1+M_{1}}\,\ta_{M_1}^k\,\,\,&\text{for }\rho_{*}\leq \rho\leq \rho^{*},\\
		\hat{C}_{1}\rho^{1+\theta_2}\,\ta_2^k \quad&\text{for }\rho\geq \rho^{*},
	\end{cases}
	\quad i=1,2.
\end{equation}
We now make the estimate for $n=k+1$.
To estimate $\big|V_{i}^{(k+1)}(\rho,u)-V_{i}^{(k)}(\rho,u)\big|$,
it suffices to consider the case: $\rho_{i}\leq \rho_{*}<\rho^{*}\leq \rho$, since
the other cases can be done by similar arguments as in \eqref{4.10}--\eqref{4.10-5}.
Noting \eqref{4.11-1} and $\rho_{i}\leq \rho_{*}<\rho^{*}\leq \rho$, and using similar arguments in \eqref{4.10-3},  we have
\begin{align*}
	&\big|(V_i^{(k+1)}-V_{k}^{(k)})(\rho,u)\big|\nonumber\\
	&\leq \Big(\int_{\rho_{i}}^{\rho_{*}}+\int_{\rho_{*}}^{\rho^{*}}+\int_{\rho^{*}}^{\rho}\Big)
	\Big\vert\frac{k''(s)}{2k'(s)}\Big\vert\,
	\sum\limits_{j=1}^2\big|V_{j}^{(k)}(s,u^{(i)}(s))-V_{j}^{(k-1)}(s,u^{(i)}(s))\big|\,\mathrm{d}s\nonumber\\
	&\leq \int_{\rho_{i}}^{\rho_{*}} \nu C_{1}\ta_1^ks^{\theta_{1}}\,\mathrm{d}s
	+\int_{\rho_{*}}^{\rho^{*}} C_{0}\tilde{C}_{1}\ta_{M_1}^ks^{M_{1}}\,\mathrm{d}s
	+\int_{\rho^{*}}^{\rho} \nu\hat{C}_{1}\ta_2^ks^{\theta_{2}}\,\mathrm{d}s\nonumber\\
	&\leq \hat{C}_{1}\big(\rho^{1+\theta_2}-(\rho^{*})^{1+\theta_2}\big)\ta_2^{k+1}
	+\tilde{C}_{1}\big((\rho^{*})^{1+M_{1}}-(\rho_{*})^{1+M_{1}}\big)\ta_{M_1}^{k+1}
	+C_{1}(\rho^{*})^{1+\theta_{1}}\,\ta_1^{k+1}\nonumber\\
	&\leq \hat{C}_{1}\rho^{1+\theta_2}\ta_2^{k+1},
\end{align*}
where we have chosen $\tilde{C}_1$ and $\hat{C}_1$ such that
\begin{equation}\label{Ch3}
	\tilde{C}_{1}\geq C_{1}(\rho_{*})^{\theta_{1}-M_{1}}\big(\ta_1\ta_{M_1}^{-1}\big)^{k+1},\qquad
	\hat{C}_{1}\geq \tilde{C}_{1}(\rho^{*})^{M_{1}-\theta_{2}} \big(\ta_{M_1}\ta_2^{-1}\big)^{k+1}.
\end{equation}
Therefore, under assumption \eqref{Ch3}, we obtain
\begin{align*}
	\big|(V_i^{(k+1)}-V_{i}^{(k)})(\rho,u)\big|\leq
	\begin{cases}
		C_{1}\rho^{1+\theta_{1}}\,\ta_1^{k+1}\quad  &\text{for }\rho\leq \rho_{*},\\
		\tilde{C}_{1}\rho^{1+M_{1}}\,\ta_{M_1}^{k+1}\,\,\, &\text{for }\rho_{*}\leq \rho\leq \rho^{*},\\
		\hat{C}_{1}\rho^{1+\theta_2}\,\ta_2^{k+1}\,\,\, &\text{for }\rho\geq \rho^{*}.
	\end{cases}
	\quad i=1,2.
\end{align*}

Recalling that $\theta_{1}\geq \theta_2$, we can take $C_{0}$ and $M_{1}$ large enough such that
\begin{equation}\label{Ch4}
	0<\ta_1\leq \ta_{M_1}\leq \ta_2<1.
\end{equation}
Combining \eqref{4.17} with \eqref{4.10-6}--\eqref{Ch4}, and taking
\begin{equation}\label{Ch5}
	\ta_{M_1}=\ta_2,\qquad \tilde{C}_{1}=C_{1}(\rho_{*})^{\theta_{2}-M_{1}},
	\qquad \hat{C}_{1}=\tilde{C}_{1}(\rho^{*})^{M_{1}-\theta_2},
\end{equation}
by induction, we conclude that, for any $n\geq 1$,
\begin{align}\label{4.12}
	\big|(V_{i}^{(n)}-V_{i}^{(n-1)})(\rho,u)\big|\leq
	\begin{cases}
		C_{1}\rho^{1+\theta_{1}}\,\ta_1^n\quad&\text{for }\rho\leq \rho_{*},\\
		\tilde{C}_{1}\rho^{1+M_{1}}\,\ta_{M_1}^n\,\,\, &\text{for }\rho_{*}\leq \rho\leq \rho^{*},\\
		\hat{C}_{1}\rho^{1+\theta_2}\,\ta_2^n\,\,\, &\text{for }\rho\geq \rho^{*},
	\end{cases}
	\quad\,i=1,2.
\end{align}

\smallskip
Noting that \eqref{Ch4} and $\r\leq \r_0$ for $(\rho,u)\in \overline{\Sigma}$, we have proved that
the two sequences in \eqref{4.9-1}, $i=1,2$, are uniformly convergent in $\overline{\Sigma}$ so that
sequence $\{(V_1^{(n)},  V_2^{(n)})\}$ is uniformly convergent in $\overline{\Sigma}$.
Let $(V_1, V_2)$ be the limit function of sequence $(V_{1}^{(n)},V_{2}^{(n)})$.
Noting the continuity and the uniform convergence of $(V_{1}^{(n)},V_{2}^{(n)})$,
$(V_{1},V_2)$ is continuous in $\overline{\Sigma}$.
Taking the limit: $n\to \infty$ in \eqref{4.9},
we conclude that $(V_{1},V_{2})$ is the continuous solution of \eqref{4.7}.

\smallskip
2.
It follows from \eqref{4.17} and \eqref{4.12} that,
for $(\rho, u)\in \{\rho\geq 0,\,|u|\leq k(\r)\}$ and $i=1,2$,
\begin{align}\label{4.20}
	\big|V_{i}(\rho,u)\big|
	&\leq  \big\vert V_{i}^{(0)} (\rho,u)\big\vert+\sum\limits_{n=1}^{\infty}\big\vert (V_{i}^{(n)}-V_{i}^{(n-1)})(\r,u)\big\vert\nonumber\\
	&\leq
	\begin{cases}
		C\rho^{1+\theta_{1}}\quad&\text{for }\rho\leq \rho_{*},\\
		C\rho^{1+M_{1}}\,\,\,&\text{for }\rho_{*}\leq \rho\leq \rho^{*},\\
		C\rho^{1+\theta_2}\quad&\text{for }\rho\geq \rho^{*}.
	\end{cases}
\end{align}
On the other hand, we see from \eqref{4.5} that, for $|u|\leq k(\rho)$,
\begin{equation}\label{4.21}
	\begin{aligned}
		&|\hat{\eta}_{\rho}(\rho,u)|=\left\vert k'(\rho)(V_1(\rho,u)+V_2(\rho,u))\right\vert\leq C\rho^{\g(\rho)-1},\\
&|\hat{\eta}_{u}(\rho,u)|= |V_1(\rho,u)-V_2(\rho,u)|\leq C\rho^{1+\theta(\rho)}.
	\end{aligned}
\end{equation}
Hence, for $|u|\leq k(\r)$, it holds that
\begin{equation}\label{4.22}
\begin{aligned}
	&|\hat{\eta}(\rho,u)|\leq \int_{\bar{\rho}}^{\r} |\hat{\eta}_{\rho}(s,u)|\,\mathrm{d}s+|\hat{\eta}(\bar{\rho},u)|\leq C\rho^{\g(\rho)},\\
	&|\hat{\eta}_{m}(\rho,m)|=|\r^{-1}\hat{\eta}_{u}(\r,u)|\leq C\rho^{\theta(\rho)},
\end{aligned}
\end{equation}
where $(\bar{\rho},u)$ is the point satisfying $k(\bar{\rho})=|u|$,
and we have used the boundary data in \eqref{4.1}.

\smallskip
3.
We now show that $V_1$ and $V_2$ have continuous first-order derivatives with respect to $(\rho, u)$.
Using \eqref{4.3}--\eqref{4.4} and Lemma \ref{lemA.1}, we have
\begin{equation}\label{4.12-6}
	\frac{\partial u^{(i)}(s)}{\partial u}=1,\qquad
	\Big\vert\frac{\partial \rho_{i}}{\partial u}\Big\vert= \frac{1}{2 |k'(\rho_{i})|}
	\leq C_{2}\rho_{i}^{1-\t(\rho_i)}.
\end{equation}
Applying $\partial_u$ to \eqref{4.9} and using \eqref{4.12-6} yield that, for $i=1,2$,
\begin{align}\label{4.12-2}
	\frac{\partial V_{i}^{(n+1)}}{\partial u}(\r,u)
	&=\frac{\partial V_{i}(\rho_{i},u_{i})}{\partial u}
	+ \frac{k''(\rho_{i})}{2k'(\rho_{i})}\,
	\sum\limits_{j=1}^{2}V_{j}^{(n)}(\rho_{i},u_{i})\,\frac{\partial \rho_{i}}{\partial u}\nonumber\\
	&\quad -\int_{\rho_{i}}^{\rho}\frac{k''(s)}{2k'(s)}\,\sum\limits_{j=1}^{2}\frac{\partial V_{j}^{(n)}}{\partial u}(s, u^{(i)}(s))\,\mathrm{d}s.
\end{align}
It follows from \eqref{4.8}, \eqref{4.12-6}, Lemma \ref{lemA.1}, and a direct calculation that, for $i=1,2$,
\begin{align}\label{4.12-1}
	\Big|\frac{\partial V_{i}^{(0)}}{\partial u}(\rho,u)\Big|
	&=\Big\vert\frac{{\rm{d}} V_{i}(\rho_i,u_i)}{{\rm{d}} \rho_i}\frac{\partial \rho_{i}}{\partial u}\Big\vert \leq
	C_{2}\rho_i
	\leq
	\begin{cases}
		\bar{C}_2\rho\qquad\,&\text{for }\rho\leq \rho_{*},\\
		\tilde{C}_2\rho^{1+M_2}\,\,\, &\text{for }\rho_{*}\leq \rho\leq \rho^{*},\\
		\hat{C}_2\rho\qquad\,&\text{for }\rho\geq \rho^{*},
	\end{cases}
\end{align}
where $C_{2}$ is chosen to
be a common, fixed, and large enough constant in \eqref{4.12-6} and \eqref{4.12-1}
depending only on $\rho_{*}$ and $\rho^{*}$,
and $\bar{C}_2 \geq C_2,\tilde{C}_2\geq C_2(\rho_{*})^{-M_2}, \hat{C}_2\geq C_2$, and $M_2$
are some large positive constants to be chosen later.

To estimate $\big|(\frac{\partial V_{i}^{(1)}}{\partial u})(\rho,u)-(\frac{\partial V_{i}^{(0)}}{\partial u})(\rho,u)\big|$,
we divide it into six cases:

\smallskip
{\it Case 1. $\rho_{i}\leq \rho\leq \rho_{*}$}: It follows from \eqref{4.17}--\eqref{4.18}
and \eqref{4.12-2}--\eqref{4.12-1} that
\begin{align}\label{4.12-4}
	\Big|\frac{\partial (V_{i}^{(1)}-
		V_{i}^{(0)})}{\partial u}(\rho,u)\Big|
	&\leq \int_{\rho_{i}}^{\rho}\Big\vert\frac{k''(s)}{2k'(s)}
	\Big\vert\, \sum\limits_{j=1}^2\Big\vert\Big(\frac{\partial V_{j}^{(0)}}{\partial u}\Big)(s, u^{(i)}(s))\Big\vert\,\mathrm{d}s\nonumber\\
	&\quad + \Big\vert\frac{k''(\rho_{i})}{2k'(\rho_{i})}\Big\vert\, \sum\limits_{j=1}^2\big\vert V_{j}^{(0)}(\rho_{i},u_{i})\big\vert\,
	\Big\vert\frac{\partial \rho_{i}}{\partial u}\Big\vert\nonumber\\
	&\leq \int_{\rho_{i}}^{\rho}\frac{\nu}{2s}\,(2\bar{C}_2s)\,\mathrm{d}s+\frac{\nu}{2\rho_{i}}\,(2C_{1}\rho_{i}^{1+\theta_{1}})\,(C_{2}\rho_{i}^{1-\theta_1})
	\nonumber\\	
	& =\bar{C}_{2}\nu\rho-\bar{C}_{2}\nu\rho_{1}+C_{1}C_{2}\nu\rho_{i}\leq \bar{C}_{2}\nu\rho,
\end{align}
where, in the last inequality of \eqref{4.12-4}, we have chosen
\begin{equation}\label{Ch6}
	\bar{C}_{2}\geq C_{1}C_{2}.
\end{equation}

\smallskip
{\it Case 2. $\r_{i}\leq \rho_{*}\leq \rho\leq \rho^{*}$}:  Then, similarly, we have
\begin{align}\label{4.12-8}
	&\Big|\frac{\partial (V_{i}^{(1)}-
		V_{i}^{(0)})}{\partial u}(\rho,u)\Big|\nonumber\\
	&\leq \Big(\int_{\rho_{i}}^{\rho_{*}}+\int_{\rho_{*}}^{\rho}\Big)
	\Big\vert\frac{k''(s)}{2k'(s)}\Big\vert\,  \sum\limits_{j=1}^2\Big\vert\Big(\frac{\partial V_{j}^{(0)}}{\partial u}\Big)(s, u^{(i)}(s))\Big\vert\,\mathrm{d}s
	\nonumber\\
	&\quad + \Big\vert\frac{k''(\rho_{i})}{2k'(\rho_{i})}\Big\vert\,\sum\limits_{j=1}^2\big\vert V_{j}^{(0)}(\rho_{i},u_{i})\big\vert\, \Big\vert\frac{\partial \rho_{i}}{\partial u}\Big\vert
	\nonumber\\
	&\leq \tilde{C}_{2}\big(\rho^{1+M_{2}}-
	(\rho_{*})^{1+M_2}\big)\ta_{M_2}
	+\bar{C}_{2}\nu(\rho_{*}-
	\rho_{i})+C_{1}C_{2}\nu\rho_{i}
	\leq \tilde{C}_{2}\rho\, \ta_{M_2},
\end{align}
where $\ta_{M_2}:=\frac{C_{0}}{1+M_{2}}$
and,
in the last inequality of \eqref{4.12-8}, we have used \eqref{Ch6} and chosen
\begin{equation}\label{Ch7}
	\tilde{C}_{2}\geq \bar{C}_2(\rho_{*})^{-M_2}\,\nu \ta_{M_2}^{-1}.
\end{equation}

\smallskip
{\it Case 3. $\rho_{*}\leq \rho_{i}\leq \rho\leq \rho^{*}$}: It follows that
\begin{align}\label{4.12-11}
	&\Big|\frac{\partial (V_{i}^{(1)}-
		V_{i}^{(0)})}{\partial u}(\rho,u)\Big|\nonumber\\
	&\leq \int_{\rho_{i}}^{\rho}\Big\vert\frac{k''(s)}{2k'(s)}\Big\vert\,
      \sum\limits_{j=1}^2\Big\vert\Big(\frac{\partial V_{j}^{(0)}}{\partial u}\Big)(s, u^{(i)}(s))\Big\vert\,\mathrm{d}s
       + \Big\vert\frac{k''(\rho_{i})}{2k'(\rho_{i})}\Big\vert\, \sum\limits_{j=1}^2\big\vert V_{j}^{(0)}(\rho_{i},u_{i})\big\vert\,
       \Big\vert\frac{\partial \rho_{i}}{\partial u}\Big\vert\nonumber\\
	&\leq \tilde{C}_{2}\big(\rho^{1+M_{2}}-
	\rho_{i}^{1+M_2}\big)\ta_{M_2}+\tilde{C_1}C_2C_{0}\rho_{i}^{1+M_{1}-\theta_2}
    \leq \tilde{C}_{2}\rho\,\ta_{M_2},
\end{align}
where,  in the last inequality of \eqref{4.12-11}, we have chosen
\begin{equation}\label{Ch8}
	M_{1}\geq M_2+\theta_2,\qquad \tilde{C}_2\geq \tilde{C}_1C_2(1+M_2)(\rho^{*})^{M_{1}-M_2-\theta_2}.
\end{equation}

\smallskip
{\it Case 4. $\rho_{i}\leq \rho_{*}<\rho^{*}\leq \rho$}:  For this case, similarly, we have
\begin{align}\label{4.12-13}
	&\Big|\frac{\partial (V_{i}^{(1)}-
		V_{i}^{(0)})}{\partial u}(\rho,u)\Big|\nonumber\\
	&\leq \Big(\int_{\rho_{i}}^{\rho_{*}}+\int_{\rho_{*}}^{\rho^{*}}+\int_{\rho^{*}}^{\rho}\Big)
	\Big\vert\frac{k''(s)}{2k'(s)}\Big\vert\,  \sum\limits_{j=1}^2\Big\vert\Big(\frac{\partial V_{j}^{(0)}}{\partial u}\Big)(s, u^{(i)}(s))\Big\vert\,\mathrm{d}s
	\nonumber\\
	&\quad + \Big\vert\frac{k''(\rho_{i})}{2k'(\rho_{i})}\Big\vert\,\sum\limits_{j=1}^2\big\vert V_{j}^{(0)}(\rho_{i},u_{i})\big\vert\, \Big\vert\frac{\partial \rho_{i}}{\partial u}\Big\vert\nonumber\\
	&\leq\hat{C}_2(\rho-
	\rho^{*})\nu+ \tilde{C}_{2}
	\big((\rho^{*})^{1+M_2}-
	(\rho_{*})^{1+M_2}\big) \ta_{M_2}
	+\bar{C}_{2}(\rho_{*}-
	\rho_{i})\nu +C_{1}C_2\rho_{i}\nu\nonumber\\
	&\leq \hat{C}_2\rho\, \nu,
\end{align}
where,  in the last inequality of \eqref{4.12-13}, we have used \eqref{Ch6} and \eqref{Ch7} and chosen
\begin{equation}\label{Ch9}
	\hat{C}_2\geq \tilde{C}_2(\rho^{*})^{M_2}\,\nu^{-1}\ta_{M_2}.
\end{equation}

\smallskip
{\it Case 5. $\rho_{*}\leq \rho_{i}\leq \rho^{*}\leq \rho$}: Then
\begin{align}\label{4.12-14}
	&\Big|\frac{\partial (V_{i}^{(1)}-
		V_{i}^{(0)})}{\partial u}(\rho,u)\Big|\nonumber\\
	&\leq \Big(\int_{\rho_{i}}^{\rho^{*}}+\int_{\rho^{*}}^{\rho}\Big)
	\Big\vert\frac{k''(s)}{2k'(s)}\Big\vert\,  \sum\limits_{j=1}^2\Big\vert\Big(\frac{\partial V_{j}^{(0)}}{\partial u}\Big)(s, u^{(i)}(s))\Big\vert\,\mathrm{d}s\nonumber\\
	&\quad + \Big\vert\frac{k''(\rho_{i})}{2k'(\rho_{i})}\Big\vert\,\sum\limits_{j=1}^2\big\vert V_{j}^{(0)}(\rho_{i},u_{i})\big\vert\, \Big\vert\frac{\partial \rho_{i}}{\partial u}\Big\vert\nonumber\\
	&\leq \hat{C}_2(\rho -
	\rho^{*})\nu
	+\tilde{C}_{2}\big((\rho^{*})^{1+M_2} -
	\rho_{i}^{1+M_2}\big)\ta_{M_2}+\tilde{C}_{1}C_{2}C_{0}\rho_{i}^{1+M_1-\theta_2}
	\leq \hat{C}_2 \rho\,\nu,
\end{align}
where we have used \eqref{Ch8} and \eqref{Ch9} in the last inequality of \eqref{4.12-14}.

\smallskip
{\it Case 6. $\rho^{*}\leq \rho_{i}\leq \rho$}: It follows similarly that
\begin{align}\label{4.12-15}
&\Big|\frac{\partial (V_{i}^{(1)}-
		V_{i}^{(0)})}{\partial u}(\rho,u)\Big|\nonumber\\
&\leq \int_{\rho_{i}}^{\rho}\Big\vert\frac{k''(s)}{2k'(s)}\Big\vert\,
	\sum\limits_{j=1}^2\Big\vert\Big(\frac{\partial V_{j}^{(0)}}{\partial u}\Big)(s, u^{(i)}(s))\Big\vert\,\mathrm{d}s
+ \Big\vert\frac{k''(\rho_{i})}{2k'(\rho_{i})}\Big\vert\, \sum\limits_{j=1}^2\big\vert V_{j}^{(0)}(\rho_{i},u_{i})\big\vert\,
	\Big\vert\frac{\partial \rho_{i}}{\partial u}\Big\vert\nonumber\\
	&\le \hat{C}_2\rho\nu -\hat{C}_2 \rho_{i}\nu+\hat{C}_{1}C_{2}\rho_i\nu
	\leq \hat{C}_2\rho\,\nu,
\end{align}
where, in the last inequality of \eqref{4.12-15},  we have chosen
\begin{equation}\label{Ch10}
	\hat{C}_2\geq \hat{C}_1C_2.
\end{equation}

\smallskip
Combining \eqref{4.12-4}--\eqref{Ch10}, we conclude that, for $i=1,2$,
\begin{align}\label{4.12-22}
	\Big|\frac{\partial (V_{i}^{(1)}-
		V_{i}^{(0)})}{\partial u}(\rho,u)\Big|\leq
	\begin{cases}
		\bar{C}_{2}\rho\,\nu \quad  &\text{for }\rho\leq \rho_{*},\\
		\tilde{C}_{2}\rho^{1+M_{2}}\,\ta_{M_2}\,\,\, &\text{for }\rho_{*}\leq \rho\leq \rho_{*},\\
		\hat{C}_{2}\rho\,\nu \quad &\text{for }\rho\geq \rho_{*},
	\end{cases}
\end{align}
provided \eqref{Ch6}, \eqref{Ch7}, \eqref{Ch8}, \eqref{Ch9}, and \eqref{Ch10} hold.

\smallskip
To use the induction arguments, we make the induction assumption for $n=k$: For $i=1,2$,
\begin{align}\label{4.12-7}
	\Big|\frac{\partial (V_{i}^{(k)}-
		V_{i}^{(k-1)})}{\partial u}(\rho,u)\Big|\leq
	\begin{cases}
		\bar{C}_{2}\rho\,\nu^{k}\quad&\text{for }\rho\leq \rho_{*},\\
		\tilde{C}_{2}\rho^{1+M_{2}}\,\ta_{M_2}\,\,\, &\text{for }\rho_{*}\leq \rho\leq \rho_{*},\\
		\hat{C}_{2}\rho\,\nu^{k}\quad&\text{for }\rho\geq \rho_{*}.
	\end{cases}
\end{align}
To estimate $\vert\frac{\partial (V_{i}^{(k+1)}-
	V_{i}^{(k)})}{\partial u}(\rho,u)\vert$,
it suffices to consider the case: $\rho_{i}\leq \rho_{*}<\rho^{*}\leq \rho$
for simplicity of presentation, since the other cases can be estimated
by similar arguments in \eqref{4.12-4}--\eqref{Ch10}.
In fact, for the case:  $\rho_{i}\leq \rho_{*}<\rho^{*}\leq \rho$,
it follows from \eqref{4.12} and \eqref{4.12-7} that
\begin{align*}
	&\Big\vert\frac{\partial (V_{i}^{(k+1)}-
		V_{i}^{(k)})}{\partial u}(\rho,u)\Big\vert\nonumber\\
	&\leq \Big(\int_{\rho_{i}}^{\rho_{*}}+\int_{\rho_{*}}^{\rho^{*}}+\int_{\rho^{*}}^{\rho}\Big)
	\Big\vert\frac{k''(s)}{2k'(s)}\Big\vert\,\sum\limits_{j=1}^{2}
	\Big\vert\frac{\partial (V_{j}^{(k)}-
		V_{j}^{(k-1)})}{\partial u}(s,u^{(i)}(s))\Big\vert\mathrm{d}s
	\nonumber\\
	&\quad\,\,
	+\Big\vert\frac{k''(\rho_{i})}{2k'(\rho_{i})}\Big\vert\,\Big\vert\frac{\partial \rho_{i}}{\partial u}\Big\vert\,
	\sum\limits_{j=1}^2\big\vert (V_{j}^{(k)}-V_{j}^{(k-1)})(\rho_{i},u_{i})\big\vert\nonumber\\
	&\leq \int_{\rho_{i}}^{\rho_{*}} \bar{C}_2\nu^{k+1}\,\mathrm{d}s
	+\int_{\rho_{*}}^{\rho^{*}} C_{0}\tilde{C}_2\ta_{M_2}^ks^{+M_2}\,\mathrm{d}s
	+\int_{\rho^{*}}^{\rho} \hat{C}_2\nu^{k+1}\,\mathrm{d}s
	+ C_{2}\nu \rho_{i}^{-\theta_1}\, C_{1}\ta_{M_2}^k\rho_{i}^{1+\theta_1}\nonumber\\
	&\leq\hat{C}_{2}(\rho-
	\rho^{*})\,\nu^{k+1}  +\tilde{C}_{2} \big((\rho^{*})^{1+M_2}-
	(\rho_{*})^{1+M_2}\big)\ta_{M_2}^{k+1}
	+\bar{C}_{2}(\rho_{*}-
	\rho_{i})\nu^{k+1} +C_{1}C_{2}\rho_i \nu^{k+1}\nonumber\\
	&
	\leq \hat{C}_{2}\rho\, \nu^{k+1},\nonumber
\end{align*}
where we have chosen
\begin{align}
		&M_{1}\geq M_{2}+\theta_2,\qquad \bar{C}_{2}\geq C_{1}C_{2}, \nonumber\\
		&\tilde{C}_2\geq \max\Big\{\tilde{C}_{1}C_{2}(1+M_2)(\rho^{*})^{M_{1}-M_2-\theta_2}\Big(\frac{1+M_{2}}{1+M_1}\Big)^{k},
		\,\bar{C}_2\frac{\big(\nu \ta_{M_2}^{-1}\big)^{k+1}}{(\rho_{*})^{M_2}}\Big\}, \label{Ch11}\\
		&\hat{C}_2\geq \max\Big\{\hat{C}_{1}C_2,\,\tilde{C}_2(\rho^{*})^{M_2} \big(\nu^{-1}\ta_{M_2}\big)^{k+1}\Big\}.\nonumber
\end{align}
Thus, under assumption \eqref{Ch11}, we conclude that, for $i=1,2$,
\begin{align*}
	\Big|\frac{\partial (V_{i}^{(k+1)}-
		V_{i}^{(k)})}{\partial u}(\rho,u)\Big|\leq
	\begin{cases}
		\bar{C}_{2}\rho\,\nu^{k+1}\quad &\text{for }\rho\leq \rho_{*},\\
		\tilde{C}_{2}\rho^{1+M_{2}}\,\ta_{M_2}^{k+1}\,\,\, &\text{for }\rho_{*}\leq \rho\leq \rho_{*},\\
		\hat{C}_{2}\rho\,\nu^{k+1}\quad &\text{for }\rho\geq \rho_{*}.
	\end{cases}
\end{align*}

Combining \eqref{Ch5} with \eqref{4.12-22}--\eqref{Ch11} and taking
\begin{align*}
	&\ta_{M_2}=\nu,\qquad
	\bar{C}_2=\max\{C_{2},\,C_{1}C_2\},\\
	&\tilde{C}_{2}= \max\big\{\bar{C}_2(\rho_{*})^{-M_2},\tilde{C}_1C_2(1+M_2)(\rho^{*})^{M_{1}-M_2-\theta_2}\big\},
	\quad
	\hat{C}_{2}=\max\big\{\tilde{C}_2(\rho^{*})^{M_2},\hat{C}_1C_2\big\},
\end{align*}
we have proved that, for any $n\ge 1$ and $i=1,2$,
\begin{align}\label{4.12-9}
	\Big|\frac{\partial (V_{i}^{(n)}-
		V_{i}^{(n-1)})}{\partial u}(\rho,u)\Big|\leq
	\begin{cases}
		\bar{C}_{2}\rho\,\nu^{n}\quad &\text{for }\rho\leq \rho_{*},\\
		\tilde{C}_{2}\rho^{1+M_{2}}\,\nu^{n}\,\,\, &\text{for }\rho_{*}\leq \rho\leq \rho_{*},\\
		\hat{C}_{2}\rho\, \nu^{n}\quad\,&\text{for }\rho\geq \rho_{*}.
	\end{cases}
\end{align}

Noting that $\nu<1$ and $\rho\leq \rho_0$ for $(\rho,u)\in \overline{\Sigma}$,
we know that $\big\{\frac{\partial V_{i}^{(n)}}{\partial u}\big\}$ is uniformly convergent
in $\overline{\Sigma}$. It is direct to check that the limit function is $\frac{\partial V_{i}}{\partial u}$.
Due to the continuity and uniform convergence of $\{\frac{\partial V_{i}^{(n)}}{\partial u}\}$,
it is clear that $\frac{\partial V_{i}}{\partial u}$ is continuous in $\overline{\Sigma}$.

On the other hand, it follows from \eqref{4.7} that
\begin{equation*}
	\left\{\begin{aligned}
		&\frac{\partial V_{1}^{(n)}}{\partial \r}=k'(\r)\frac{\partial V_{1}^{(n)}}{\partial u}
		-\frac{k''(\rho)}{2k'(\rho)}\big(V_{1}^{(n-1)}+V_{2}^{(n-1)}\big),\\
		&\frac{\partial V_{2}^{(n)}}{\partial \rho}=-k'(\r)\frac{\partial V_{2}^{(n)}}{\partial u}-\frac{k''(\rho)}{2k'(\rho)}\big(V_{1}^{(n-1)}+V_{2}^{(n-1)}\big),
	\end{aligned}
	\right.
\end{equation*}
which, with \eqref{4.18}, \eqref{4.12}, and \eqref{4.12-9}, yields that, for $k\geq 0$ and $i=1,2$,
\begin{align}\label{4.7-2}
	&\Big\vert\frac{\partial (V_{i}^{(n)}-
		V_{i}^{(n-1)})}{\partial \rho}(\rho,u)\Big\vert\nonumber\\
	&\leq k'(\rho)\Big\vert\frac{\partial(V_{i}^{(n)}-
		V_{i}^{(n-1)})}{\partial u}(\r,u)\Big\vert
	+\Big\vert \frac{k''(\r)}{2k'(\r)}\Big\vert
	\, \sum\limits_{j=1}^2\big\vert (V_{j}^{(n)}-V_{j}^{(n-1)})(\rho,u)\big\vert\nonumber\\
	&\leq
	\begin{cases}
		C\rho^{\theta_1}\,\nu^{n}\qquad&\text{for }\rho\leq \rho_{*},\\
		C\rho^{M_1}\,\nu^{n}\,\,\, &\text{for }\rho_{*}\leq \rho\leq \rho_{*},\\
		C\rho^{\theta_2}\,\nu^{n}\qquad&\text{for }\rho\geq \rho_{*},
	\end{cases}
\end{align}
for some large constant $C>0$.
Thus, $\frac{\partial V_{i}^{(n)}}{\partial \rho}$ converges uniformly to $\frac{\partial V_{i}}{\partial \rho}$ in $\overline{\Sigma}$. It is clear that  $\frac{\partial V_{i}}{\partial \r}$ is continuous.
Therefore, $(V_{1}(\rho,u), V_{2}(\r,u))$ is a $C^1$--solution of the Goursat problem \eqref{4.6}--\eqref{4.6-1},
which implies that $\hat{\eta}$ is a $C^2$--solution of \eqref{4.1}.

\smallskip
4.
From \eqref{4.12-1} and \eqref{4.12-9}, we obtain that, for $i=1,2$,
\begin{equation}\label{4.12-10}
	\Big\vert\frac{\partial V_{i}}{\partial u}(\rho,u)\Big\vert
	\leq \Big\vert\frac{\partial V_{i}^{(0)}}{\partial u}(\rho,u)\Big\vert +\sum\limits_{n=1}^{\infty}\Big\vert
	\frac{\partial (V_{i}^{(n)}-V_{i}^{(n-1)})}{\partial u}(\r,u)\Big\vert
	\leq  C\rho
\end{equation}
for $\rho\geq 0$ and $|u|\leq k(\rho)$.
Similarly, using \eqref{4.7-2}, we see that, for $i=1,2$,
\begin{equation*}
	\Big\vert\frac{\partial V_{i}}{\partial \rho}\Big\vert\leq C\rho^{\t(\rho)}\qquad
	\text{for $\rho\geq 0$ and $|u|\leq k(\rho)$}.
\end{equation*}
Therefore, for $|u|\leq k(\rho)$,
it follows from \eqref{4.5} and \eqref{4.12-10} that
\begin{equation*}
	\begin{aligned}
		&|\hat{\eta}_{uu}(\rho,u)|=|\partial_{u}V_1(\rho,u)-\partial_{u}V_2(\rho,u)|\leq C\rho,\\
		&|\hat{\eta}_{\rho u}(\rho,u)|=|k'(\rho)(\partial_{u}V_1(\rho,u)+\partial_{u}V_2(\rho,u))|\leq C\rho^{\theta(\rho)}.
	\end{aligned}
\end{equation*}
If $\hat{\eta}_{m}$ is regarded as a function of $(\rho,u)$, we have
\begin{equation*}
	|\hat{\eta}_{m\rho }(\rho,u)|\leq C\rho^{\theta(\rho)-1},\quad  |\hat{\eta}_{mu}(\rho,u)|\leq C
	\qquad\,\, \text{ for }|u|\leq k(\rho).
\end{equation*}
If $\hat{\eta}_{m}$ is regarded as a function of $(\rho,m)$, we see that, for $|u|\leq k(\rho)$,
$$
\begin{aligned}
&|\hat{\eta}_{m\rho}(\rho,m)|=|\hat{\eta}_{m\rho}(\rho,u)+u\hat{\eta}_{mu}(\rho,u)|\leq C\rho^{\theta(\rho)-1},\\
&|\hat{\eta}_{mm}(\rho,m)|=|\rho^{-1}\hat{\eta}_{mu}(\rho,u)|\leq C\rho^{-1}.
\end{aligned}
$$

\smallskip
5.
We now prove the uniqueness of $\hat{\eta}$, which is equivalent to
the uniqueness of solutions of \eqref{4.6}--\eqref{4.6-1}
in the class of $C^1$--solutions satisfying \eqref{4.20}.
Suppose that there exist two $C^{1}$ solutions $(V_{1},V_{2})$ and $(\tilde{V}_{1},\tilde{V}_{2})$
of \eqref{4.6}--\eqref{4.6-1} satisfying the uniform estimate \eqref{4.20}.
Then it follows from \eqref{4.7} that, for $i=1,2$,
\begin{equation}\label{U1}
	V_{i}(\rho,u)-\tilde{V}_{i}(\rho,u)=-\int_{\rho_{i}}^{\rho}\frac{k''(s)}{2k'(s)}\,\sum\limits_{j=1}^2\big(V_{j}(s,u^{(i)}(s))-\tilde{V}_{j}(s,u^{(i)}(s))\big)\,\mathrm{d}s.
\end{equation}
Applying the uniform estimates \eqref{4.20} and similar arguments as in \eqref{4.12}--\eqref{U1} yields
\begin{align*}
	\max_{(\rho,u)\in \R_{+}\times \R\atop |u|\leq k(\rho)}\big\vert V_{i}(\rho,u)-\tilde{V}_{i}(\rho,u)\big\vert\leq\left\{
	\begin{aligned}
		&C\rho^{1+\theta_{1}}\,\ta_1^n\quad\,\,\,\,\,\,\,\text{for }\rho\leq \rho_{*},\\
		&C\rho^{1+M_{1}}\,\ta_{M_1}^n \quad\,\text{for }\rho_{*}\leq \rho\leq \rho^{*},\\
		&C\rho^{1+\theta_2}\,\ta_2^n\qquad\,\text{for }\rho\geq \rho^{*},
	\end{aligned}
	\right.
\end{align*}
for any $n\ge 0$, where $C \gg 1$ is independent of $n$.
Taking $n\to \infty$, we obtain that
$V_{i}(\rho,u)\equiv \tilde{V}_{i}(\rho,u)$ for $|u|\leq k(\rho)$ which,
with \eqref{4.5} and $\hat{\eta}(0,u)\equiv 0$, yields the uniqueness of $\hat{\eta}$.

\smallskip
6.
We now estimate the entropy flux $\hat{q}$.
It follows from \eqref{1.23} that, for all entropy pairs,
\begin{equation}\label{4.14}
	q_{\rho}=u\eta_{\rho}+\rho k'(\rho)^2\eta_u,\qquad q_{u}=\rho\eta_{\rho}+u\eta_{u}.
\end{equation}
Then there exists an entropy flux $\hat{q}(\rho,u)\in C^2(\R_{+}\times \R)$
corresponding to the special entropy $\hat{\eta}$:
\begin{equation*}
	\hat{q}(\rho,u)=\frac{1}{2}\rho |u|^3\pm\rho u\big(e(\rho)+\rho e'(\rho)\big)\qquad\, \text{for } \pm u\geq k(\rho).
\end{equation*}
It follows from \eqref{4.21} and \eqref{4.14} that
$|\hat{q}_{\rho}(\rho,u)|=|u\hat{\e}_{\r}+\r k'(\r)^2\hat{\e}_{u}|\leq C\r^{\g(\r)+\t(\r)-1}$
for $|u|\leq k(\rho)$,
which implies
\begin{equation}\label{4.27-2}
	|\hat{q}(\r,u)| = \Big|\int_{\bar{\rho}}^{\rho}\hat{q}_{\rho}(s,u) \,\mathrm{d}s+\hat{q}(\bar{\rho},u) \Big|\leq C\r^{\g(\r)+\t(\r)}\qquad\, \text{for } |u|\leq k(\r),
\end{equation}
where $(\bar{\rho},u)$ is the point satisfying $k(\bar{\rho})=|u|$.

\smallskip
For $|u|\leq k(\rho)$, it follows from \eqref{4.22} and \eqref{4.27-2} that
$|\hat{q}-u\hat{\eta}|\leq |\hat{q}|+|u||\hat{\eta}|\leq C\rho^{\g(\rho)+\t(\rho)}$.
In region $\{(\rho,u)\,:\,|u|\geq k(\rho)\}$,
it is direct to check that all the estimates in Lemma \ref{lem4.1} hold
by using \eqref{4.1-0}.
Therefore, the proof of Lemma \ref{lem4.1} is now complete. $\hfill\square$
\FloatBarrier

\subsection{\,Estimates of the weak entropy pairs}
In order to show the compactness of the weak entropy dissipation measures below,
we now derive some estimates of the weak entropy pairs.
To achieve this, from \eqref{2.21e}--\eqref{2.21q}, it requires
to analyze the entropy kernel and entropy flux kernel, respectively.

The entropy kernel $\chi=\chi(\r,u,s)$ is a fundamental solution of the entropy equation \eqref{1.25-1}:
\begin{equation}\label{6.1}
	\left\{\begin{aligned}
		\dis&\chi_{\r\r}-\frac{P'(\r)}{\r^2}\chi_{uu}=0,\\
		\dis&\chi\vert_{\r=0}=0,\quad \chi_{\r}\vert_{\r=0}=\delta_{u=s}.
	\end{aligned}
	\right.
\end{equation}
As pointed out in \cite{Chen-LeFloch-2000} that equation \eqref{6.1}
is invariant under the Galilean transformation, which implies that
$\chi(\rho,u,s)=\chi(\r,u-s,0)=\chi(\r,0,s-u)$.
For simplicity, we write it as  $\chi(\r,u,s)=\chi(\r,u-s)$ below
when no confusion arises.

\smallskip
The corresponding entropy flux kernel $\sigma(\r,u,s)$ satisfies the Cauchy problem
for $\sigma-u\chi$:
\begin{equation}\label{6.2}
	\left\{\begin{aligned}
		\dis&(\sigma-u\chi)_{\r\r}-\frac{P'(\r)}{\r^2}(\sigma-u\chi)_{uu}=\frac{P''(\r)}{\r}\chi_{u},\\
		\dis&(\sigma-u\chi)\vert_{\r=0}=0,\quad  (\sigma-u\chi)_{\r}\vert_{\r=0}=0.
	\end{aligned}
	\right.
\end{equation}
We recall from \cite{Chen-LeFloch-2000} that $\sigma-u\chi$
is also Galilean invariant.
From \eqref{1.3}--\eqref{1.5}, $P(\rho)$
satisfies all the conditions in \cite{Chen-LeFloch-2000,Chen-LeFloch-2003}.

For later use, we introduce the definition of fractional
derivatives ({\it cf.} \cite{Chen-1986, Chen-LeFloch-2000, Lions-Perthame-Souganidis-1996}).
For any real $\alpha>0$, the fractional derivative $\partial_{s}^{\alpha}f$ of a function $f=f(s)$ is
$$
\partial_{s}^{\alpha}f(s)=\Gamma(-\alpha)f*[s]_{+}^{-\alpha-1},
$$
where $\Gamma(x)$ is the Gamma function and the convolution should be understood
in the sense of distributions.
The following formula:
$$
\partial_{s}^{\alpha}(sg(s))=s\partial_{s}^{\alpha+1}g+(\alpha+1)\partial_{s}^{\alpha}g
$$
holds for fractional derivatives.
We now present two useful lemmas
for the entropy kernel $\chi(\rho ,u)$ and the entropy flux kernel $\sigma(\rho, u)$ when $\rho$ is bounded.

\smallskip
\begin{lemma}[{\cite[Theorems 2.1--2.2]{Chen-LeFloch-2000}}]\label{thm6.1}
	The entropy kernel $\chi(\rho,u)$ admits the expansion{\rm :}
	\begin{equation}\label{6.3}
		\chi(\rho,u)=a_{1}(\rho)G_{\lambda_1}(\rho,u)+a_2(\rho)G_{\lambda_1+1}(\rho,u)+g_1(\rho,u)\qquad
		\text{for $\r\in [0,\infty)$},
	\end{equation}
	where $k(\r)=\int_{0}^{\r}\frac{\sqrt{P'(y)}}{y}\,\mathrm{d}y$ and
\begin{equation}\label{6.5-1}
	\begin{aligned}
			&	G_{\lambda_1}(\rho,u)=[k(\rho)^2-u^2]_{+}^{\lambda_1},\qquad
			\lambda_1=\frac{3-\g_1}{2(\g_1-1)}>0,\\
			&a_1(\r)=M_{\lambda_1}k(\r)^{-\lambda_1}k'(\r)^{-\frac{1}{2}},\quad
			M_{\lambda_1}
			=\left(\frac{2\lambda_1}{\sqrt{2\lambda_1+1}}\int_{-1}^{1}(1-z^2)^{\lambda_1}\,\mathrm{d}z\right)^{-1},\\
			&a_2(\r)=-\frac{1}{4(\lambda_1+1)}k(\r)^{-\lambda_1-1}k'(\r)^{-\frac{1}{2}}
			\int_{0}^{\r}k(s)^{\lambda_1}k'(s)^{-\frac{1}{2}}a_{1}''(s)\,\mathrm{d} s.
\end{aligned}
\end{equation}
	Moreover, $\operatorname{supp}\chi(\rho,u)\subset \{(\rho,u)\,:\, |u|\leq k(\rho)\}$,
	and $\chi(\rho,u)>0$
	in $\{(\rho,u)\, : \, |u|<k(\rho)\}$. The remainder term $g_1(\rho,\cdot)$ and its fractional
    derivative $\partial_{u}^{\lambda_1+1}g_1(\rho,\cdot)$ are H\"{o}lder continuous.
	Furthermore, for any fixed $\r_{\max}>0$, there exists
	$C(\rho_{\max})>0$
	depending only on $\r_{\rm max}$ such that
	\begin{equation}\label{6.4}
		|g_1(\r,u-s)|\leq C(\rho_{\max})[k(\r)^2-(u-s)^2]_{+}^{\lambda_1+\alpha_0+1},
	\end{equation}
	for any $0\leq \r\leq \r_{\max}$ and some $\alpha_0\in (0,1)$.
	In addition, for any $0\leq \r\leq \r_{\max}$,
	\begin{equation}\label{6.5}
		|a_1(\r)|+\r^{1-2\theta_1}|a_1'(\r)|+\r^{2-2\t_1}|a_1''(\r)|
		+
		|a_2(\r)|+\r|a_2'(\r)|+\r^2|a_2''(\r)|\leq C(\r_{\max}).
	\end{equation}
\end{lemma}

\noindent
\textbf{Proof.}
Since \eqref{6.3}--\eqref{6.4} have been derived in \cite[Theorem 2.2]{Chen-LeFloch-2000},
it suffices to prove \eqref{6.5}.
From \eqref{A.4-1}, we find that
$|a_{1}(\r)|\leq C(\rho_{\max})$ for $0\leq \rho\leq \rho_{\max}$.
For $|a_1'(\r)|$, a direct calculation shows that
\begin{equation*}
	a_1'(\r)=-\lambda_1M_{\lambda_1}k(\rho)^{-\lambda_1-1}k'(\r)^{\frac{1}{2}}
	-\frac{1}{2}M_{\lambda_1}k(\r)^{-\lambda_1}k'(\r)^{-\frac{3}{2}}k''(\r).
\end{equation*}	
It follows from \eqref{1.4} that $k(\r)=C_{1}\rho^{\t_1}\big(1+O(\rho^{2\theta_1})\big)$
as $\rho\in [0,\rho_{\max}]$
for some constant $C_1>0$ that may depend on $\k_1$ and $\g_1$.
Then, by direct calculation, we observe that the term involving $\rho^{-1}$ in $a_1'(\r)$ vanishes
so that $|a_1'(\r)|\leq C(\rho_{\max})\rho^{2\theta_1-1}$.
Similarly, we obtain that $|a_1''(\r)|\leq C(\rho_{\max})\rho^{2\theta_1-2}$.
Finally, using $\eqref{6.5-1}_3$, we can obtain the estimates
of $a_2(\rho)$ in $\eqref{6.5}$ by a direct calculation.
This completes the proof.
$\hfill\square$

\smallskip
\begin{lemma}[{\cite[Theorem 2.3]{Chen-LeFloch-2000}}]\label{thm6.2}
	The entropy flux kernel $\sigma(\rho,u)$ admits the expansion
	\begin{equation*}
		(\sigma-u\chi)(\rho,u)
		=-u\big(b_{1}(\rho)G_{\lambda_1}(\rho,u)+b_2(\rho)G_{\lambda_1+1}(\rho,u)\big)
		+g_2(\rho,u)\qquad \text{for $\r\in [0,\infty)$},
	\end{equation*}
	where
	\begin{equation}\label{6.8-1}
		\begin{aligned}
			&b_1(\r)=M_{\lambda_1}\r k(\r)^{-\lambda_1-1}k'(\r)^{\frac{1}{2}}>0,\\
			&\begin{aligned}
				b_2(\rho)&=-\frac{1}{4(\lambda_1+1)}\rho k'(\rho)^{\frac{1}{2}}k(\rho)^{-(\lambda_1+2)}
				\int_{0}^{\rho}k(s)^{\lambda_1}k'(s)^{-\frac{1}{2}}a_1''(s)\,\mathrm{d} s
				\\&\quad -\frac{1}{4(\lambda_1+1)}k(\rho)^{-(\lambda_1+2)}k'(\rho)^{-\frac{1}{2}}
				\int_{0}^{\rho}k(s)^{\lambda_1+1}k'(s)^{-\frac{1}{2}}b_1''(s)\,\mathrm{d}s
				\\&\quad +\frac{1}{4(\lambda_1+1)}k(\rho)^{-(\lambda_1+2)}k'(\rho)^{-\frac{1}{2}}
				\int_{0}^{\rho} s k(s)^{\lambda_1}k'(s)^{\frac{1}{2}}a_1''(s)\,\mathrm{d}s.
			\end{aligned}
		\end{aligned}
	\end{equation}
	The remainder term $g_2(\rho,\cdot)$ and its fractional derivative $\partial_{u}^{\lambda_1+1}g_2(\rho,\cdot)$
	are H\"{o}lder continuous.
	Moreover, for any fixed $\r_{\max}>0$, there exists $C(\rho_{\max})>0$
	depending only on $\r_{max}$ such that
	\begin{equation*}
		|g_2(\r,u-s)|\leq C(\rho_{\max})[k(\r)^2-(u-s)^2]_{+}^{\lambda_1+\alpha_0+1},
	\end{equation*}
	for any $0\leq \r\leq \r_{\max}$ and some $\alpha_0\in (0,1)$.
	Furthermore, similar to the proof of \eqref{6.5}, for any $0\leq \r\leq \r_{\max}$,
	\begin{equation}\label{6.8}
		|b_1(\r)|+\rho^{1-2\t_1}|b_1'(\r)|+\r^{2-2\t_1}|b_2''(\r)|
		+|b_2(\r)|+\rho|b_2'(\r)|+\rho^2|b_2''(\r)|\leq C(\rho_{\max}).
	\end{equation}
\end{lemma}

\begin{remark}
	In {\rm \cite[Theorem 2.2]{Chen-LeFloch-2000}}, it is proved
	that $a_2(\rho)$ and $b_{2}(\rho)$ satisfy $|a_2(\rho)|+|b_2(\rho)|\leq C\rho k(\rho)^{-2}$
	for the pressure law given in {\rm \cite[(2.1)]{Chen-LeFloch-2000}}. In this paper, we have improved them
	to be $\eqref{6.5}$ and $\eqref{6.8}$ under conditions \eqref{1.3}--\eqref{1.5}.
\end{remark}

\smallskip
For later use, we recall a useful representation formula for $\chi(\rho,u)$.

\smallskip
\begin{lemma}[First representation formula, {\cite[Lemma 3.4]{Schrecker-Schulz-2020}}]\label{lem6.3}
	Given any $(\rho,u)$ with $|u|\leq k(\r)$ and $0\leq \rho_0<\r$,
	\begin{equation*}
		\begin{aligned}
			\chi(\rho, u)=&\,\frac{1}{2(\rho-\rho_0) k^{\prime}(\rho)}
			\int_{\rho_0}^{\rho} k^{\prime}(s)\, \tilde{d}(s)\big(\chi(s, u+k(\rho)-k(s))+ \chi(s, u-k(\rho)+k(s))\big) \,\mathrm{d}s \\
			&\,-\frac{1}{2(\rho-\rho_0) k^{\prime}(\rho)} \int_{-(k(\rho)-k(\rho_0))}^{k(\rho)-k(\rho_0)} \chi(\rho_0, u-s) \,\mathrm{d}s,
		\end{aligned}
	\end{equation*}
	where
	$\,
	\tilde{d}(\r):=2+(\r-\rho_0)\frac{k''(\r)}{k'(\r)}.
	$
\end{lemma}

\smallskip
\begin{remark}
	In the statement of {\rm \cite[Lemma 3.4]{Schrecker-Schulz-2020}}, $\rho_{0}$ is positive.
	However, the proof of {\rm \cite[Lemma 3.4]{Schrecker-Schulz-2020}} is also valid
	for $\rho_{0}=0$ without modification{\rm ;} see also
	{\rm \cite[(3.38)]{Chen-LeFloch-2000}}.
\end{remark}

\smallskip
Given any $\psi\in C_{0}^2(\R)$, a regular weak entropy pair $(\e^{\psi},\,q^{\psi})$ can be given by
\begin{align}\label{6.10}
	\e^{\psi}(\r,u)=\int_{\R} \psi(s)\,\chi(\r,u,s)\, \mathrm{d}s,\qquad
	q^{\psi}(\r,u)=\int_{\R}\psi(s)\,\sigma(\r,u,s) \,\mathrm{d}s.
\end{align}
It follows from \eqref{6.1} that
\begin{equation}\label{6.1-1}
	\left\{\begin{aligned}
		\dis&\eta_{\r\r}^{\psi}-k'(\r)^2\eta_{uu}^{\psi}=0,\\
		\dis&\eta^{\psi}\vert_{\r=0}=0,\quad \eta_{\r}^{\psi}\vert_{\r=0}=\psi(u).
	\end{aligned}
	\right.
\end{equation}
Using Lemmas \ref{thm6.1}--\ref{lem6.3}, we
obtain the following lemma
for the weak entropy pair $(\eta^{\psi},q^{\psi})$.

\smallskip
\begin{lemma}\label{lem6.2}
	For any weak entropy $(\e^{\psi},q^{\psi})$ defined in \eqref{6.10},
	there exists a constant $C_{\psi}>0$ depending only on $\rho^{*}$ and $\psi$
	such that, for all $\rho\in [0,2\rho^{*}]$,
	\begin{equation*}
		|\eta^{\psi}(\r,u)|+
		|q^{\psi}(\rho,u)|\leq C_{\psi}\rho.
	\end{equation*}
	If $\eta^{\psi}$ is regarded as a function of $(\rho,m)$, then
	\begin{align*}
		|\eta_{m}^{\psi}(\rho,m)|+|\r\e_{mm}^{\psi}(\rho,m)|\leq C_{\psi},\qquad
		|\eta_{\r}^{\psi}(\r,m)|\leq C_{\psi}\big(1+\r^{\t_1}\big).
	\end{align*}
	Moreover, if $\eta_m^{\psi}$ is regarded as a function of $(\r,u)$, then
	\begin{equation*}
		|\eta_{mu}^{\psi}(\r,u)|+|\rho^{1-\t_1}\eta_{m\rho}^{\psi}(\rho,u)|\leq C_{\psi}.
	\end{equation*}
\end{lemma}

\noindent\textbf{Proof.} All the estimates can be found in \cite[Lemma 3.8]{Schrecker-Schulz-2019}
or \cite[Lemma 4.13]{Schrecker-Schulz-2020} except the estimate of $\eta_{\r}^{\psi}(\rho,m)$.
In fact, applying Lemma \ref{lem6.3} to \eqref{6.1-1} and using
$\dis d(\r):=2+\frac{\r k''(\r)}{k'(\r)}$, we have
\begin{align}
	\eta^{\psi}(\rho, u)=
	&\, \frac{1}{2\rho k^{\prime}(\rho)} \int_{0}^{\rho} k^{\prime}(s)\, d(s)\,
	\eta^{\psi}(s, u+k(\rho)-k(s))\, \mathrm{d}s \nonumber\\
	&\,+\frac{1}{2\rho k^{\prime}(\rho)} \int_{0}^{\rho} k^{\prime}(s)\, d(s)\,
	\eta^{\psi}(s, u-k(\rho)+k(s)) \,\mathrm{d}s
	:= I_1+I_2.\label{6.21-1}
\end{align}

\smallskip
We regard $\eta^{\psi}$ as a function of $(\r,m)$. Then we have
\begin{equation}\label{6.21-4}
	\partial_{\r}\eta^{\psi}(\r,m)=\partial_{\rho}\eta^{\psi}(\rho,u)-\frac{u}{\r}\partial_{u}\eta^{\psi}(\rho,u).
\end{equation}
Without loss of generality, we assume that $\operatorname{supp}\psi\subset [-L,L]$ for some $L>0$.
Then a direct calculation shows that $\eta^{\psi}(\rho,u)=0$ if $|u|\geq k(\rho)+L$.
Noticing $\eta_{u}^{\psi}(\r,u)=\r\eta_{m}^{\psi}(\rho,m)$, we have
\begin{equation}\label{6.21-5}
	\Big\vert\frac{u}{\r}\pa_{u}\eta^{\psi}(\rho,u)\Big\vert\leq |u|\,|\eta_{m}^{\psi}(\rho,m)|
	\leq C_{\psi}(1+\rho^{\t_1})\qquad \mbox{for }0\leq \rho\leq 2\rho^{*}.
\end{equation}
Thus, it suffices to calculate $\pa_{\r}\eta^{\psi}(\rho,u)$. It follows from \eqref{6.21-1} that
$\partial_{\r}\eta^{\psi}(\r,u)=\partial_{\r}I_1+\partial_{\r}I_2$. A direct calculation shows that
\begin{align}
	\partial_{\r}I_1
	&=\frac{1}{2}\big(-\rho^{-2}(k'(\r))^{-1}-\r^{-1}(k'(\rho))^{-2}k''(\rho)\big)
	\int_{0}^{\rho} k^{\prime}(s)\, d(s)\, \eta^{\psi}(s, u+k(\rho)-k(s)) \, \mathrm{d}s\nonumber\\
	&\quad +\frac{1}{2\rho}\int_{0}^{\r}k'(s)\,d(s)\,\eta_{u}^{\psi}(s,u+k(\r)-k(s)) \,\mathrm{d}s + \frac{1}{2\rho}d(\r)\eta^{\psi}(\r,u).\label{6.21-2}
\end{align}
Using \eqref{A.4-1} and Lemma \ref{lemA.3}, we obtain that
\begin{equation*}
	|\partial_{\r}I_1|
\le C_\psi (1+\rho^{\theta_1})
	\qquad\,\, \mbox{for $0\leq \rho\leq 2\rho^{*}$}.
\end{equation*}
Similarly, we obtain that
$|\partial_{\r}I_2|\leq C_{\psi}(1+\rho^{\t_1})$.
Thus, we conclude that
$
|\partial_{\rho}\eta^{\psi}(\rho,u)|\leq |\partial_{\r}I_1|+|\partial_{\r}I_2|\leq C_{\psi}(1+\rho^{\t_1}),
$
which, with \eqref{6.21-4}--\eqref{6.21-5}, implies that
$
|\partial_{\r}\eta^{\psi}(\r,m)|\leq C_\psi(1+\rho^{\t_1}).
$
$\hfill\square$

\smallskip
We notice that all the above estimates for the weak entropy pairs
in Lemmas \ref{thm6.1}--\ref{lem6.2} hold when the density is bounded.
To establish the $L^p$-compensated compactness framework, we need the entropy pair estimates
when the density is large, namely $\r\geq \r^{*}$.
From now on in this subsection,
we use the representation formula of Lemma \ref{lem6.3} to estimate $(\e^{\psi},q^{\psi})$
in the large density region $\rho\geq \rho^{*}$.

\smallskip
\begin{lemma}\label{lem6.2-1}
	There exists a positive constant $C>0$ depending only on $\r^{*}$ such that
	\begin{equation*}
		\|\chi(\r,\cdot)\|_{L_{u}^{\infty}}\leq C\r\qquad \text{for $\r\geq \r^{*}$}.
	\end{equation*}
\end{lemma}

\noindent{\bf Proof.}
For $\rho\geq \rho^{*}$, $\chi(\rho,u)$ satisfies
\begin{equation*}
	\left\{\begin{aligned}
		&\dis\chi_{\rho\rho}-k'(\rho)^2\chi_{uu}=0,\\
		&\dis\chi(\rho,u)\vert_{\rho=\rho^{*}}=\chi(\rho^{*},u),\quad
		\chi_{\rho}(\rho,u)\vert_{\rho=\rho^{*}}=\chi_{\rho}(\rho^{*},u).
	\end{aligned}\right.
\end{equation*}
where $\chi(\r^{*},u)$ and $\chi_{\r}(\r^{*},u)$ are given in Lemma \ref{thm6.1}.
Then, applying Lemma \ref{lem6.3},
we obtain that, for $\rho> \rho^{*}$,
\begin{align*}
	\|k'(\rho)\chi(\rho,\cdot)\|_{L_{u}^{\infty}}
	&\leq \frac{1}{\rho-\rho^{*}}\int_{\rho^{*}}^{\rho}d_{*}(s)\|k'(s)\chi(s,\cdot)\|_{L_{u}^{\infty}}\,\mathrm{d}s
	\nonumber\\
	&\quad
	+\frac{1}{2(\rho-\rho^{*})}\int_{-(k(\rho)-k(\rho^{*}))}^{k(\rho)-k(\rho^{*})}\|\chi(\rho^{*},\cdot)\|_{L_{u}^{\infty}}
\,\mathrm{d}s\nonumber\\
	&\leq \frac{1}{\rho-\rho^{*}}\int_{\rho^{*}}^{\rho}d_{*}(s)\|k'(s)\chi(s,\cdot)\|_{L_{u}^{\infty}}\,\mathrm{d}s+C\rho^{\t_2-1},\label{6.14}
\end{align*}
where
$d_{*}(\rho):=2+(\rho-\rho^{*})\frac{k''(\rho)}{k'(\rho)}$.
By a similar proof to that for Lemma \ref{lemB.2},
we have
\begin{equation*}
	\|k'(\rho)\chi(\rho,\cdot)\|_{L_{u}^{\infty}}\leq C\rho^{\t_{2}}\qquad \text{for $\rho\geq 2\rho^{*}$},
\end{equation*}
which, with \eqref{A.5-1}, yields that
$\|\chi(\rho,\cdot)\|_{L_{u}^{\infty}}\leq C\rho$
for $\rho\geq 2\rho^{*}$.
For $\rho^{*}\leq \rho \leq 2\rho^{*}$, it follows from Lemma \ref{thm6.1} that
$\|\chi(\rho,\cdot)\|_{L_{u}^{\infty}}\leq C\leq C\rho$.
$\hfill\square$

\medskip
\begin{lemma}\label{lem6.4}
	Let $\rho\geq \rho^{*}$ and $\psi\in C_{0}^2(\mathbb{R})$.
	Then,  in the $(\rho,u)$--coordinates,
	\begin{equation*}
		|\eta^{\psi}(\rho,u)|+|\eta_{u}^{\psi}(\rho,u)|+|\eta_{uu}^{\psi}(\rho,u)|\leq C_{\psi}\rho,
		\quad |\e_{\r}^{\psi}(\r,u)|+\rho^{1-\t_2}|\eta_{\rho\rho}^{\psi}(\r,u)|\leq C_{\psi}\rho^{\t_2}.
	\end{equation*}
In the $(\rho,m)$--coordinates,
		$\,|\e^{\psi}_{\rho}(\rho,m)|+	\rho^{\t_2}|\eta_{m}^{\psi}(\rho,m)|+\rho^{1+\t_2}|\e_{mm}^{\psi}(\rho,m)|\leq C_{\psi}\r^{\t_2}$.

\smallskip
\noindent
If $\e^{\psi}_{m}(\rho,m)$ is regarded as a function of $(\rho, u)$, then
	$\,|\eta_{mu}^{\psi}|+\rho^{1-\t_2}|\e_{m\rho}^{\psi}|\leq C_{\psi}$.

\smallskip
\noindent
	All the above constants $C_{\psi}>0$ depend only on $\|\psi\|_{C^2}$ and $\operatorname{supp}\psi$.
\end{lemma}

\smallskip
\noindent{\bf Proof.} We divide the proof into five steps.

\smallskip
1. Using \eqref{6.10} and Lemma \ref{lem6.2-1},
we obtain that, for $\r\geq \r^{*}$,
\begin{equation}\label{6.23}
	|{\eta}^{\psi}(\r,u)|+|{\eta}_{u}^{\psi}(\r,u)|+|{\eta}_{uu}^{\psi}(\r,u)|
	\leq \|{\chi}(\rho,\cdot)\|_{L^{\infty}(\mathbb{R})}\|(\psi,\psi',\psi'')\|_{L^1(\mathbb{R})}
	\leq  C_{\psi}\rho.
\end{equation}

2. For the estimate of $\eta_{\rho}^{\psi}(\rho,u)$, the proof is similar to Lemma \ref{lem6.2}. Indeed,
$\eta^{\psi}$ satisfies
\begin{equation}\label{6.25}
	\left\{\begin{aligned}
		\dis&\eta_{\rho\rho}^{\psi}-k'(\rho)^2\eta_{uu}^{\psi}=0,\\
		\dis&\eta^{\psi}(\rho,u)\vert_{\rho=\rho^{*}}=\eta^{\psi}(\rho^{*},u),\quad
		\eta_{\rho}^{\psi}(\rho,u)\vert_{\rho=\rho^{*}}=\eta_{\r}^{\psi}(\rho^{*},u).
	\end{aligned}\right.
\end{equation}
It follows from  \eqref{6.25} and Lemma \ref{lem6.3} that
\begin{align}
	\eta^{\psi}(\rho, u)
	&=\frac{1}{2(\rho-\rho^{*}) k^{\prime}(\rho) }\int_{\rho^{*}}^{\rho}
	d_{*}(s)k^{\prime}(s)\eta^{\psi}(s, u+k(\rho)-k(s))\,\mathrm{d}s\nonumber\\
	&\quad -\frac{1}{2(\rho-\rho^{*}) k^{\prime}(\rho) }\int_{\rho^{*}}^{\rho}
	d_{*}(s)k^{\prime}(s)\eta^{\psi}(s, u-k(\rho)+k(s))\,\mathrm{d}s\nonumber\\
	&\quad-\frac{1}{2(\rho-\rho^{*}) k^{\prime}\left(\rho\right) }\int_{-(k(\rho)-k(\rho^{*}))}^{k(\rho)-k(\rho^{*})}\eta^{\psi}(\rho^{*},u-s)\,\mathrm{d}s,\label{6.27}
\end{align}	
where $d_{*}(\r)=2+(\r-\r^{*})\frac{k''(\r)}{k'(\r)}$ and
$0<d_{*}(\rho)\leq 3$
for $\rho\geq \rho^{*}$ from \eqref{A.7}.
Then, following the similar arguments as in the proof of Lemma \ref{lem6.2},
we can obtain that $|\eta_{\rho}^{\psi}(\rho,u)|\leq C_{\psi}\r^{\t_2}$ for $\rho\geq 2\rho^{*}$.
Moreover, from Lemma \ref{lem6.2},
$|\eta_{\rho}^{\psi}(\rho,u)|\leq C_{\psi}\leq C_{\psi}\r^{\t_2}$ for $\rho \in [\rho^{*},2\rho^{*}]$ so that
\begin{equation}\label{6.27-1}
	|\eta_{\rho}^{\psi}(\rho,u)|\leq C_{\psi}\r^{\t_2}\qquad \text{for }\rho\geq \rho^{*}.
\end{equation}

3.
For $\eta_{\r\r}^{\psi}(\r,u)$, it follows from \eqref{6.23}--\eqref{6.25} that
\begin{equation*}
	|\eta_{\rho\rho}^{\psi}(\r,u)|\leq |k'(\rho)^2|\,|\eta_{uu}^{\psi}(\r,u)|
	\leq  C_{\psi}\rho^{2\theta_2-1}=C_{\psi}\rho^{\gamma_2-2}\qquad \text{ for }\rho\geq \rho^{*}.
\end{equation*}

\smallskip
4.
In the $(\rho,m)$--coordinates, it is clear that
\begin{align*}
	\e_{m}^{\psi}(\r,m)=\r^{-1}\e_{u}^{\psi}(\r,u),\quad
	\e_{mm}^{\psi}(\r,m)=\r^{-2}\e_{uu}^{\psi}(\r,u),\quad
	\e_{\rho}^{\psi}(\rho,m)=\eta_{\rho}^{\psi}(\rho,u)-\frac{m}{\rho^2}\eta_{u}^{\psi}(\rho,u).
\end{align*}
On the other hand,  if $\eta_{m}^{\psi}$ is regarded as a function $(\r,u)$, it is direct to obtain
$$
\e_{mu}^{\psi}(\r, u)=\partial_{u}\big(\r^{-1}\e_{u}^{\psi}(\r,u)\big)=\r^{-1}\e_{uu}^{\psi}(\r,u).
$$
Thus, using \eqref{6.23} and \eqref{6.27-1},
\begin{align*}
	&|\e_{m}^{\psi}(\r,m)|+ |\e_{mu}^{\psi}(\r,u)|+\r |\e_{mm}^{\psi}(\r,m)|\leq C_{\psi}\r^{-1},\\
	&|\e_{\r}^{\psi}(\r,m)|\leq |\eta_{\rho}^{\psi}(\rho,u)|+\frac{|m|}{\rho^2}\big|\eta_{u}^{\psi}(\rho,u)\big|
	\leq C_{\psi}\r^{\t_2}+C_{\psi}(L+k(\rho))\leq C_{\psi}\rho^{\t_2},
\end{align*}
where we have used that $\operatorname{supp}\psi\subset [-L,L]$
and $\eta_{u}^{\psi}(\r,u)=\r\eta_{m}^{\psi}(\rho,m)$.

\smallskip
5.
For the estimates of $\eta_{m\r}^{\psi}(\r,u)=\partial_{\r}\eta_{m}^{\psi}(\r,u)$, it follows from \eqref{6.27} that
\begin{align}
	\eta_m^{\psi}(\rho,m)=\frac{1}{\rho}\eta_u^{\psi}(\rho,u)
&=\frac{1}{2\rho(\rho-\rho^{*})k'(\rho)}\int_{\rho^{*}}^{\rho}d_{*}(s)\,k'(s)\,\eta_{u}^{\psi}(s,u+k(\rho)-k(s))\,\mathrm{d}s\nonumber
	\\&\quad+\frac{1}{2\rho(\rho-\rho^{*})k'(\rho)}\int_{\rho^{*}}^{\rho}d_{*}(s)\,k'(s)\,\eta_{u}^{\psi}(s,u-k(\rho)+k(s))\,\mathrm{d}s\nonumber
	\\&\quad-\frac{1}{2\rho(\rho-\rho^{*})k'(\rho)}\int_{-(k(\rho)-k(\rho^{*}))}^{k(\rho)-k(\rho^{*})}\eta_{u}^{\psi}(\rho^{*},u-s)\,\mathrm{d}s
	\nonumber\\&:=J_1+J_2+J_3,\label{6.28}
\end{align}
and $\pa_{\r}\eta_{m}^{\psi}(\rho,u)=\partial_{\rho}J_1+\partial_{\rho}J_2+\partial_{\rho}J_3$.

\smallskip
A direct calculation shows that
\begin{align}
	\partial_{\rho}J_1&=\partial_{\rho}\Big(\frac{1}{2\rho(\rho-\rho^{*})k'(\rho)}\Big)
	\int_{\rho^{*}}^{\rho}d_{*}(s)\,k'(s)\,\eta_{u}^{\psi}(s,u+k(\rho)-k(s))\,\mathrm{d}s\nonumber\\
	&\quad +\frac{1}{2\rho(\rho-\rho^{*})}\int_{\rho^{*}}^{\rho}d_{*}(s)\,k'(s)\,\eta_{uu}^{\psi}(s,u+k(\rho)-k(s))\,\mathrm{d}s\nonumber\\
	&\quad
	+\frac{1}{2\rho(\rho-\rho^{*})}d_{*}(\rho)\,\eta_{u}^{\psi}(\rho,u)\nonumber\\
	&:=J_{1,1}+J_{1,2}+J_{1,3},\label{6.29}\\
\partial_{\rho}J_2&=\partial_{\rho}\Big(\frac{1}{2\rho(\rho-\rho^{*})k'(\rho)}\Big)\int_{\rho^{*}}^{\rho}d_{*}(s)\,k'(s)\,\eta_{u}^{\psi}(s,u-k(\rho)+k(s))\,\mathrm{d}s
	\nonumber\\
	&\quad -\frac{1}{2\rho(\rho-\rho^{*})}\int_{\rho^{*}}^{\rho}d_{*}(s)\,k'(s)\,\eta_{uu}^{\psi}(s,u-k(\rho)+k(s))\,\mathrm{d}s\nonumber\\
	&\quad +\frac{1}{2\rho(\rho-\rho^{*})}d_{*}(\rho)\eta_{u}^{\psi}(\rho,u)\nonumber\\[1mm]
	&:=J_{2,1}+J_{2,2}+J_{2,3}.\label{6.29-1}
\end{align}

Clearly, we have
\begin{equation*}
	\Big\vert\partial_{\rho}\Big(\frac{1}{2\rho(\rho-\rho^{*})k'(\rho)}\Big)\Big\vert
	=\frac{1}{2}\Big\vert\frac{\rho-\rho^*}{2\rho^{2}\,(\rho-\rho^{*})^2\,k'(\rho)}
	+\frac{k''(\rho)}{\rho\,(\rho-\rho^{*})\,k'(\rho)^2}\Big\vert
	\leq\frac{C}{(\rho-\rho^{*})^2\rho^{\theta_2}},
\end{equation*}
which, with \eqref{6.23} and $0<d_{*}(\rho)\leq 3$ for $\rho\geq \rho^{*}$,
yields
\begin{align}\label{6.32}
		|J_{1,1}+J_{2,1}|&=\Big\vert\partial_{\rho}\Big(\frac{1}{2\rho\,(\rho-\rho^{*})\,k'(\rho)}\Big)
		\int_{\rho^{*}}^{\rho}k'(s)\,d_{*}(s)\,\eta_{u}^{\psi}(s,u+k(\rho)-k(s))\,\mathrm{d}s\Big\vert\nonumber
		\\&\quad+\Big\vert\partial_{\rho}\Big(\frac{1}{2\rho\,(\rho-\rho^{*})\,k'(\rho)}\Big)
		\int_{\rho^{*}}^{\rho}k'(s)\,d_{*}(s)\,\eta_{u}^{\psi}(s,u-k(\rho)+k(s))\,\mathrm{d}s\Big\vert\nonumber
		\\&\leq \frac{C_{\psi}}{(\rho-\rho^{*})^2\,\rho^{\theta_2}}\int_{\rho^{*}}^{\rho}s^{\theta_2}\,\mathrm{d}s\leq \frac{C_{\psi}}{\rho-\rho^{*}}.
\end{align}

\smallskip
It follows from \eqref{6.23} and  $0<d_{*}(\rho)\leq 3$ for $\rho\geq \rho^{*}$
that
\begin{equation}\label{6.30}
	\begin{aligned}
		|J_{1,2}|+|J_{2,2}|&=\Big\vert\frac{1}{2\rho(\rho-\rho^{*})}\int_{\rho^{*}}^{\rho}d_{*}(s)\,k'(s)\,\eta_{uu}^{\psi}(s,u+k(\rho)-k(s))\,\mathrm{d}s\Big\vert
		\\&\quad+\Big\vert\frac{1}{2\rho(\rho-\rho^{*})}\int_{\rho^{*}}^{\rho}d_{*}(s)\,k'(s)\,\eta_{uu}^{\psi}(s,u-k(\rho)+k(s))\,\mathrm{d}s\Big\vert
		\\&\leq \frac{C_{\psi}}{\rho(\rho-\rho^{*})}\int_{\rho^{*}}^{\rho}s^{\theta_2}\,\mathrm{d}s
		\leq \frac{C_{\psi}}{\rho(\rho-\rho^{*})}\r^{\t_2}(\r-\r^{*})\leq C_{\psi}\r^{\t_2-1}.
	\end{aligned}
\end{equation}

For $J_{1,3}+J_{2,3}$, it is direct to see that
\begin{equation}\label{6.33}
	|J_{1,3}+J_{2,3}|\leq \Big\vert\frac{1}{\rho(\rho-\rho^{*})}\,d_{*}(\rho)\,\eta_{u}^{\psi}(\rho,u)\Big\vert\leq \frac{C_{\psi}}{\rho-\rho^{*}}.
\end{equation}

\smallskip
For $\partial_{\rho}J_3$, we notice that
\begin{equation*}
	\begin{aligned}
		\partial_{\rho}J_3&=-\partial_{\rho}\Big(\frac{1}{2\rho\,(\rho-\rho^{*})\,k'(\rho)}\Big)
		\int_{-(k(\rho)-k(\rho^{*}))}^{k(\rho)-k(\rho^{*})}\eta_{u}^{\psi}(\rho^{*},u-s)\,\mathrm{d}s
		\\&\quad\, -\frac{1}{2\rho\,(\rho-\rho^{*})}\,\big(\eta_{u}^{\psi}(\rho^{*},u-k(\rho)+k(\rho^{*}))
		+\eta_{u}^{\psi}(\rho^{*},u+k(\rho)-k(\rho^{*}))\big),
	\end{aligned}
\end{equation*}
which, with $0<\theta_2\leq 1$, yields
\begin{align}\label{6.35}
	|\partial_{\rho}J_3|
	\leq \frac{C_{\psi}\rho^{*}}{(\rho-\rho^{*})^2}\,|k(\rho)-k(\rho^{*})|+\frac{C_{\psi}\rho^{*}}{\rho\,(\rho-\rho^{*})}
	\leq \frac{C_{\psi}}{\rho-\rho^{*}}(1+\rho^{\theta_2-1})\leq \frac{C_{\psi}}{\rho-\rho^{*}}.
\end{align}
Combining \eqref{6.29}--\eqref{6.35} with \eqref{6.28} yields
that
$|\eta_{m\r}^{\psi}(\r,u)|\leq C_{\psi}\r^{\t_2-1}$ for $\r\geq 2\r^{*}$.
For $\r^{*}\leq \r\leq 2\r^{*}$, it follows from  Lemma \ref{lem6.2} that
$
|\eta_{m\r}^{\psi}(\r,u)|\leq C_{\psi}\r^{\t_1-1}\leq C_{\psi}
$
for $\r^{*}\leq \r\leq 2\r^{*}$.
Thus, we conclude that
$|\eta_{m\rho}^{\psi}(\r,u)|\leq C_{\psi}\rho^{\theta_2-1}$
for $\rho\geq \rho^{*}$.
$\hfill\square$

\smallskip
We now estimate $q^{\psi}$ for $\rho\geq \rho^{*}$. It follows from \eqref{6.2} that $h:=\sigma-u\chi$ satisfies
\begin{equation*}
	\left\{\begin{aligned}
		\dis &h_{\rho\rho}-k'(\rho)^2h_{uu}=\frac{P''(\rho)}{\rho}\chi_{u},\\
		\dis &h(\rho^{*},u)=(\sigma-u\chi)(\rho^{*},u),\quad 	
		h_{\rho}(\rho^{*},u)=(\sigma-u\chi)_{\rho}(\rho^{*},u),
	\end{aligned}\right.
\end{equation*}
where $(\sigma-u\chi)(\rho^{*},u)$ and $(\sigma-u\chi)_{\rho}(\rho^{*},u)$ are given by Lemma \ref{thm6.2}.
Similar to Lemma \ref{lem6.3}, we have the following representation formula for $h$.

\smallskip
\begin{lemma}[Second representation formula {\cite[Lemmas 3.4 and 3.9]{Schrecker-Schulz-2020}}]\label{lem6.4-0}
	For any $(\rho,u)$ with $|u|\leq k(\r)$ and $\r>\rho^{*}$,
	\begin{align}\label{6.37}
		h(\rho, u)
		=&\,\frac{1}{2(\rho-\rho^{*}) k^{\prime}(\rho)} \int_{\rho^{*}}^{\rho} k^{\prime}(s) d_{*}(s)
		\big(h(s, u+k(\rho)-k(s))+
		h(s, u-k(\rho)+k(s))\big) \,\mathrm{d}s
		\nonumber\\
		&+\frac{1}{2(\rho-\rho^{*})k'(\rho)}\int_{\rho^{*}}^{\rho}(s-\rho^{*})\frac{P''(s)}{s}
		\chi(s,u+k(\rho)-k(s))\,\mathrm{d}s
		\nonumber\\
		&-\frac{1}{2(\rho-\rho^{*})k'(\rho)}\int_{\rho^{*}}^{\rho}(s-\rho^{*})\frac{P''(s)}{s}\chi(s,u-k(\rho)+k(s))\,\mathrm{d}s\nonumber\\
		&-\frac{1}{2(\rho-\rho^{*}) k^{\prime}(\rho)} \int_{-(k(\rho)-k(\rho^{*}))}^{k(\rho)-k(\rho^{*})} h(\rho^{*}, u-s) \,\mathrm{d}s,
	\end{align}
	where  $d_{*}(\rho)=2+(\rho-\rho^{*})\frac{k''(\rho)}{k'(\rho)}$.
\end{lemma}

\smallskip
\begin{lemma}\label{lem6.4-1}
	There exists a constant $C>0$ depending only on $\r^{*}$ such that
	\begin{equation*}
		\|(\sigma-u\chi)(\r,u)\|_{L_{u}^{\infty}}\leq C\r^{1+\t_2}
		\qquad\,\, \text{for }\r\geq \r^{*}.
	\end{equation*}
\end{lemma}

\noindent{\bf Proof.}
It follows from \eqref{A.2-1}, \eqref{6.37}, and Lemma \ref{lem6.2-1} that
\begin{align*}
		\|k'(\rho)h(\rho,\cdot)\|_{L_{u}^{\infty}}
		&\leq \frac{1}{\rho-\rho^{*}}\int_{\rho^{*}}^{\rho}d_{*}(s)\|k'(s)h(s,\cdot)\|_{L_{u}^{\infty}}\,\mathrm{d}s\nonumber\\
		&\quad +\frac{C}{\rho-\rho^{*}}\int_{\rho^{*}}^{\rho}s^{\g_2-1}\,\mathrm{d}s+C\rho^{\t_2-1}\\
		&\leq  \frac{1}{\rho-\rho^{*}}\int_{\rho^{*}}^{\rho}d_{*}(s)\|k'(s)h(s,\cdot)\|_{L_{u}^{\infty}}\,\mathrm{d}s+C\r^{2\t_2},
\end{align*}
which, with \eqref{A.5-1} and a similar proof to that for Lemma \ref{lemB.2},
yields
\begin{equation*}
	\|k'(\rho)h(\rho,\cdot)\|_{L_{u}^{\infty}}\leq C\rho^{2\t_2}
	\,\Longrightarrow \,\|h(\r,\cdot)\|_{L_{u}^{\infty}}\leq C\rho^{1+\t_2}\qquad \text{for }\rho\geq 2\rho^{*}.
\end{equation*}
For $\rho^{*}\leq \rho\leq 2\rho^{*}$, it follow from Lemma \ref{thm6.2} that
$\|h(\r,\cdot)\|_{L_{u}^{\infty}}\leq C\leq C\rho^{1+\t_2}$.
$\hfill\square$

\smallskip
\begin{lemma}\label{lem6.5}
	For $\rho\geq \rho^{*}$ and $\psi\in C_{0}^2(\R)$,
	\begin{equation}\label{6.43}
		|q^{\psi}(\r,u)|\leq C_{\psi}\r^{1+\theta_2}.
	\end{equation}
\end{lemma}

\noindent{\bf Proof.} Recall that
\begin{equation}\label{6.44}
	\begin{aligned}
		q^{\psi}(\r,u)
		&=\int_{\R}\big(\sigma(\r,u,s)-u\chi(\r,u-s)\big)\psi(s)\,\mathrm{d}s+u\int_{\R}\chi(\r,u-s)\psi(s)\,\mathrm{d}s\\
		&:=h^{\psi}(\r,u)+u\,\e^{\psi}(\r,u).
	\end{aligned}
\end{equation}
It follows from Lemma \ref{lem6.4-1}
that
\begin{equation}\label{6.45}
	|h^{\psi}(\r,u)|\leq C\|(\sigma-u\chi)(\r,\cdot)\|_{L_{u}^{\infty}(\R)}\|\psi\|_{L^1(\R)}\leq C_{\psi}\rho^{1+\theta_2}.
\end{equation}
Since there exists
$L>0$ such that $\operatorname{supp}\psi\subset [-L,L]$,
then it follows from Lemma \ref{lem6.4} that
$|u\e^{\psi}(\r,u)|\leq (k(\r)+L)|\e^{\psi}(\r,u)|\leq C_{\psi}\r^{1+\t_2}$
for $\rho\geq \rho^{*}$,
which,  with \eqref{6.44}--\eqref{6.45}, yields \eqref{6.43}.
$\hfill\square$

\subsection{\,Singularities of the entropy kernel and the entropy flux kernel}
As indicated in \cite{Chen-LeFloch-2000, Schrecker-Schulz-2019, Schrecker-Schulz-2020},
understanding the singularities of the entropy kernel and the entropy flux kernel is essential for the reduction
of the Young measure. Thus, it requires some detailed estimates of the singularities of the entropy kernel
and the entropy flux kernel.
The arguments in this subsection are similar to \cite[\S 6]{Schrecker-Schulz-2020},
the main difference is that a more subtle Gr\"{o}nwall inequality (see Lemma \ref{lemB.2}) is needed
to obtain the desired estimates of the singularities.

\smallskip
\begin{lemma}\label{lem7.1}
	For $\r\geq \r^{*}$, the coefficient functions $a_1(\r)$ and $a_2(\r)$
	and the remainder term $g_1(\r,u)$ in {\rm Lemma \ref{thm6.1}} satisfy
	\begin{align*}
		&|a_1(\rho)|+\rho^{\t_2}|a_2(\rho)|\leq C\rho^{\frac{1}{2}-\frac{\t_2}{2\t_1}},
		\qquad\,\,
		\|g_1(\r,u)\|_{L_{u}^{\infty}(\R)}
		\leq
		\begin{cases}
			C\rho \quad &\text{if $\t_2<\t_1$},\\
			C\rho\ln \r
			\,\,\, &\text{if $\t_2=\t_1$},
		\end{cases}
		\\
		&\|\partial_ug_1(\rho,u)\|_{L_{u}^{\infty}(\R)}
		+ \|\partial_u^{\lambda_1+1}g_1(\rho,u)\|_{L_{u}^{\infty}(\R)}
		+\|\partial_u^{\lambda_1+1+\alpha_0}g_1(\rho,u)\|_{L_{u}^{\infty}(\R)}\leq C\rho,
	\end{align*}
	where $\alpha_0\in (0,1)$ is the H\"{o}lder exponent.
\end{lemma}

\smallskip
\noindent{\bf Proof.} We divide the proof into five steps.

\smallskip
1.
It follows from \eqref{6.5} that, for $\rho\in [0,\rho_{*}]$,
\begin{equation}\label{7.3}
	|a_1(\rho)|+\rho^{1-2\theta_1}|a_1'(\rho)|+\rho^{2-2\theta_1}|a_1''(\rho)|
	+|a_2(\rho)|+\rho|a_2'(\rho)|+\rho^{2}|a_2''(\rho)|\leq C.
\end{equation}
For $\rho\geq \rho^{*}$, it follows from \eqref{A.5-1}, \eqref{6.5-1}, and a direct calculation that
\begin{align}
	&|a_1(\rho)|+\r| a_1'(\r)|+\r^2| a_1''(\r)|\leq C\rho^{\frac{1}{2}-\frac{\t_2}{2\t_1}},\label{7.6}\\
	&|a_{2}(\rho)|=\Big\vert\frac{1}{4\lambda_1+1}k(\rho)^{-\lambda_1-1}k'(\rho)^{-\frac{1}{2}}
	\int_{0}^{\rho}k(s)^{\lambda_1}k'(s)^{-\frac{1}{2}}a_1''(s)\,\mathrm{d}s\Big\vert
	\nonumber\\
	&\qquad\quad \leq C\rho^{\frac{1}{2}-\frac{\t_2}{2\t_1}-\t_2}\Big(\rho_{*}^{\t_1}+1
	+\int_{\rho^{*}}^{\rho}s^{-1-\t_2}\,\mathrm{d}s\Big)
	\leq C\rho^{\frac{1}{2}-\t_2-\frac{\theta_2}{2\t_1}}.\label{7.8}
\end{align}
Moreover, calculating the derivatives explicitly, we conclude
$$
|a_{2}(\rho)|+\r|a_2'(\rho)|+\r^2|a_2''(\rho)|\leq C\rho^{\frac{1}{2}-\t_2-\frac{\theta_2}{2\t_1}}
\qquad\,\, \text{for }\rho \geq \rho^{*}.
$$

\smallskip
2. For the remainder term $g_1(\r,u)$, it follows from \cite[Proof of Theorem 2.1]{Chen-LeFloch-2000} that
\begin{equation*}
	\left\{\begin{aligned}
		\dis&\partial_{\rho\rho}g_1(\rho,u)-k'(\rho)^2\partial_{uu}g_1(\rho,u)=A(\rho)k(\rho)^{-1}f_{\lambda_1+1}(\frac{u}{k(\rho)}),
		\\\dis &g_1(0,\cdot)=0,\quad \partial_{\rho}g_1(0,\cdot)=0,
	\end{aligned}
	\right.
\end{equation*}
where $f_{\lambda_1}(y)=[1-y^2]_{+}^{\lambda_1}$ and $A(\rho)=-a_2''(\rho)k(\rho)^{2\lambda_1+3}$.
By \eqref{7.3},
$A(\rho)\sim O(\rho^{-1+2\theta_1})$ as $\rho\to 0$ and $|A(\rho)|\leq \rho^{-\frac{3}{2}+\theta_2(1+\frac{1}{2\t_1})}$ for $\rho\geq \rho^{*}$.
Similar to Lemma \ref{lem6.3}, we have the following representation formula for $g_1(\rho,u)$:
\begin{align}
	k'(\rho)g_1(\rho,u)
	&=\frac{1}{2\rho}\int_{0}^{\rho}d(s)\,k'(s)\,
	\big(g_1(s,u+k(\rho)-k(s))+g_1(s,u-k(\rho)+k(s))\big)\,\mathrm{d}s\nonumber\\
	&\quad\, +\frac{1}{2\rho}\int_{0}^{\rho}s\,A(s)\,k(s)^{-1}\,
	\Big(\int_{u-k(\rho)+k(s)}^{u+k(\rho)-k(s)}f_{\lambda_1+1}(\frac{y}{k(s)})\,\mathrm{d}y\Big)\,\mathrm{d}s,\label{7.12}
\end{align}
where $d(\rho)=2+\rho\frac{k''(\rho)}{k'(\rho)}$ satisfying \eqref{A.17-1}.
Since $\rho A(\rho)k(\rho)^{-1}\sim O(\rho^{\theta_1})$ as $\rho \to 0$,
the second integral is well-defined.
Then it follows directly from  \eqref{A.4-1}--\eqref{A.5-1} and \eqref{7.12} that
\begin{align}
	&\|k'(\rho)g_1(\rho,\cdot)\|_{L_{u}^{\infty}(\R)}\nonumber\\
	&\leq \frac{1}{\rho}\int_{0}^{\rho}d(s)\|k'(s)g_1(s,\cdot)\|_{L_{u}^{\infty}(\R)}\,\mathrm{d}s+\frac{2k(\rho)}{\rho}
	\Big(\int_{0}^{\rho_{*}}+\int_{\r_{*}}^{\r^{*}}+\int_{\r^{*}}^{\r}\Big)|sA(s)k(s)^{-1}|\,\mathrm{d}s\nonumber
	\\&\leq \frac{1}{\rho}\int_{0}^{\rho}d(s)\|k'(s)g_1(s,\cdot)\|_{L_{u}^{\infty}(\R)}\,\mathrm{d}s
	+C\r^{\t_2-1}\big(\r_{*}^{1+\t_1}+1+\r^{\frac{1}{2}+\frac{\t_2}{2\t_1}}\big)\nonumber\\
	&\leq \frac{1}{\rho}\int_{0}^{\rho}d(s)\|k'(s)g_1(s,\cdot)\|_{L_{u}^{\infty}(\R)}\,\mathrm{d}s
	+C\rho^{-\frac{1}{2}+\frac{\t_2}{2\t_1}+\t_2}\qquad\,\, \text{for }\rho\geq \rho^{*}.\label{7.12-1}
\end{align}
Since $g_1(\r,u)$ is H\"{o}lder continuous and $\operatorname{supp}g_1(\r,\cdot)\subset [-k(\r),k(\r)]$,
it follows from \eqref{A.4-1}--\eqref{A.5-1} that
$\|k'(\rho)g_1(\rho,\cdot)\|_{L_{u}^{\infty}(\R)}$ is locally integrable with respect to $\rho\in [0,\infty)$.
Applying Lemma \ref{lemB.2} to \eqref{7.12-1}, we obtain that, for $\r\geq \r^{*}$,
\begin{equation*}
	\|k'(\rho)g_1(\r,\cdot)\|_{L_{u}^{\infty}(\R)}\leq \left\{\begin{aligned}
		&C\rho^{\t_2}\qquad\quad\;\;\text{if $\t_2<\t_1$},\\
		&C\rho^{\t_2}\ln \r\qquad \text{if $\t_2=\t_1$},
	\end{aligned}
	\right.
\end{equation*}
which, with \eqref{A.5-1}, yields that, for $\r\geq \r^{*}$,
\begin{equation*}
	\|g_1(\r,\cdot)\|_{L_{u}^{\infty}(\R)}\leq \left\{\begin{aligned}
		&C\rho\qquad\quad\;\;\text{if }\t_2<\t_1,\\
		&C\rho\ln \r\qquad \text{if }\t_2=\t_1.
	\end{aligned}
	\right.
\end{equation*}	

3. Applying $\partial_{u}$ to \eqref{7.12}, we have
\begin{align}\label{7.13}
	&k'(\rho)\partial_ug_1(\rho,u)\nonumber\\
	&=\frac{1}{2\rho}\int_{0}^{\rho}d(s)k'(s)\Big(\partial_ug_1(s,u+k(\rho)-k(s))
	+\partial_ug_1(s,u-k(\rho)+k(s))\Big)\,\mathrm{d}s
	\nonumber\\
	&\quad  +\frac{1}{2\rho}\int_{0}^{\rho}sA(s)k(s)^{-1}f_{\lambda_1+1}(\frac{u+k(\rho)-k(s)}{k(s)})
	\,\mathrm{d}s\nonumber\\
	&\quad-\frac{1}{2\rho}\int_{0}^{\rho}sA(s)k(s)^{-1}f_{\lambda_1+1}(\frac{u-k(\rho)+k(s)}{k(s)})\,\mathrm{d}s.
\end{align}
Since $|f_{\lambda_1+1}(s)|\leq 1$, by similar arguments as in Step 2, we can obtain
\begin{equation*}
	\|\partial_{u}g_1(\r,\cdot)\|_{L_{u}^{\infty}(\R)}\leq C\r\qquad\,\, \text{for }\r\geq \r^{*}.
\end{equation*}	

4. Applying the fractional derivative $\partial_{u}^{\lambda_1}$ to \eqref{7.13}, we have
\begin{align*}
	&k'(\rho)\partial_u^{\lambda_1+1}g_1(\rho,u)\\
	&=\frac{1}{2\rho}\int_{0}^{\rho}d(s)k'(s)
	\Big((\partial_u^{\lambda_1+1}g_1)(s,u+k(\rho)-k(s))
	+(\partial_u^{\lambda_1+1}g_1)(s,u-k(\rho)+k(s))\Big)\,\mathrm{d}s
	\\&\quad +\frac{1}{2\rho}\int_{0}^{\rho}sA(s)k(s)^{-1-\lambda_1}
	(\partial_{u}^{\lambda_1}f_{\lambda_1+1})(\frac{u+k(\rho)-k(s)}{k(s)})\,\mathrm{d}s\nonumber\\
	&\quad -\frac{1}{2\rho}\int_{0}^{\rho}sA(s)k(s)^{-1-\lambda_1}(\partial_{u}^{\lambda_1}f_{\lambda_1+1})(\frac{u-k(\rho)+k(s)}{k(s)})\,\mathrm{d}s,
\end{align*}
where we have taken into account the homogeneity of the factional derivative in the last term.
Using the Fourier transform relation as in \cite[(I.26)--(I.27)]{Lions-Perthame-Souganidis-1996}, we can obtain
\begin{equation*}
	\big|\mathscr{F}\big((\partial_{u}^{\lambda_1}f_{\lambda_1+1})(u)\big)(\xi)\big|
	= C_{\lambda_1+1}|\xi|^{-\frac{3}{2}}\,\vert J_{\lambda_1+\frac{3}{2}}(|\xi|)\vert\leq \frac{\tilde{C}_{\lambda_1+1}}{1+\xi^2}
\end{equation*}
for some positive constants $C_{\lambda_1+1}$ and $\tilde{C}_{\lambda_1+1}$ depending only on $\lambda_1+1$,
where we have used the asymptotic relations for the first kind of Bessel functions $J_{\lambda_1+\frac{3}{2}}(|\xi|)$
to obtain the final inequality.
Since $(1+|\xi|^2)^{-1}$ is integrable, applying the Fourier inversion theorem,
we see that $(\partial_{u}^{\lambda_1}f_{\lambda_1+1})(u)$ is uniformly bounded.
Hence, by similar arguments as in Step 2, we have
\begin{equation*}
	\|\big(\partial_u^{\lambda_1+1}g_1\big)(\rho,\cdot)\|_{L_{u}^{\infty}(\mathbb{R})}\leq C\rho
	\qquad\,\, \text{for $\rho\geq\rho^{*}$}.
\end{equation*}

5. By Lemma \ref{thm6.1}, we assume that $\alpha_0\in (0,1)$ is the H\"{o}lder exponent of $(\partial_{u}^{\lambda_1+1})g_1(\r,u)$.
Then, applying the fractional derivative $\partial_{u}^{\lambda_1}$ to \eqref{7.13}, we have
\begin{align*}\label{7.21-1}
	&k'(\rho)\partial_u^{\lambda_1+1+\alpha_0}g_1(\rho,u)\nonumber\\	
	&=\frac{1}{2\rho}\int_{0}^{\rho}d(s)k'(s)(\partial_u^{\lambda_1+1+\alpha_0}g_1)(s,u+k(\rho)-k(s))\,\mathrm{d}s\nonumber\\
	&\quad +\frac{1}{2\rho}\int_{0}^{\rho}d(s)k'(s)(\partial_u^{\lambda_1+1+\alpha_0}g_1)(s,u-k(\rho)+k(s))\,\mathrm{d}s\nonumber\\&\quad +\frac{1}{2\rho}\int_{0}^{\rho}sA(s)k(s)^{-1-\lambda_1-\alpha_0}
	(\partial_{u}^{\lambda_1+\alpha_0}f_{\lambda_1+1})(\frac{u+k(\rho)-k(s)}{k(s)})\,\mathrm{d}s\nonumber\\
	&\quad -\frac{1}{2\rho}\int_{0}^{\rho}sA(s)k(s)^{-1-\lambda_1-\alpha_0}(\partial_{u}^{\lambda_1+\alpha_0}f_{\lambda_1+1})(\frac{u-k(\rho)+k(s)}{k(s)})\,\mathrm{d}s.
\end{align*}
Noting
\begin{equation*}
	\big|\mathscr{F}\big((\partial_{u}^{\lambda_1+\alpha_0}f_{\lambda_1+1})(u)\big)(\xi)\big|
	=C_{\lambda_1+1}|\xi|^{-\frac{3}{2}+\alpha_0}\big|J_{\lambda_1+\frac{3}{2}}(|\xi|)\big|
	\leq \frac{\tilde{C}_{\lambda_1+1}}{1+\xi^{2-\alpha_0}},
\end{equation*}
and using the Fourier inversion theorem, we find that
$(\partial_{u}^{\lambda_1+\alpha_0}f_{\lambda_1+1})(u)$ is uniformly bounded.
By similar arguments as in Step 2, we obtain that
$\|\big(\partial_u^{\lambda_1+1+\alpha_0}g_1\big)(\rho,\cdot)\|_{L_{u}^{\infty}(\mathbb{R})}\leq C\rho$
for $\rho\geq\rho^{*}$.
This completes the proof.
$\hfill\square$

\smallskip
From {\rm Lemmas \ref{thm6.1}} and {\rm \ref{lem7.1}}, we conclude

\smallskip
\begin{corollary}\label{cor7.1}
	$\chi(\rho, \cdot)$ is H\"{o}lder continuous and
	$$
	\|\chi(\r,\cdot)\|_{C_{u}^{\tilde{\alpha}}}\leq C(1+\rho|\ln \r|)\qquad\,\,
	\text{for $\,\tilde{\alpha}\in (0,\min\{\lambda_1,1\}]\text{ and }\rho\geq 0$}.
	$$
\end{corollary}

\begin{lemma}\label{lem7.2}
	For $\rho\geq \rho^{*}$, the coefficient functions $b_{1}(\r)$ and $b_2(\r)$
	and the remainder term $g_2(\r,u)$ in {\rm Lemma \ref{thm6.2}} satisfy
	\begin{equation*}
		\begin{aligned}
			&|b_1(\rho)|
			+\rho^{\t_2}|b_2(\rho)|\leq C\rho^{\frac{1}{2}-\frac{\t_2}{2\t_1}},\qquad\,\,
			\|g_2(\r,\cdot)\|_{L_{u}^{\infty}(\R)}\leq
			\begin{cases}
				C\rho^{1+\t_2}\quad &\text{if $\t_2<\t_1$},\\
				C\rho^{1+\t_2}\ln \r\,\,\, &\text{if $\t_2=\t_1$},
			\end{cases}
			\\[1mm]
			&\|\partial_{u}g_2(\r,\cdot)\|_{L_{u}^{\infty}(\R)}+\|\big(\partial_u^{\lambda_1+1}g_2\big)(\rho,\cdot)\|_{L_{u}^{\infty}(\mathbb{R})}
			+\|\big(\partial_u^{\lambda_1+1+\alpha_0}g_2\big)(\rho,\cdot)\|_{L_{u}^{\infty}(\mathbb{R})}\leq C\rho^{1+\theta_2},
		\end{aligned}	
	\end{equation*}
	where $\alpha_0\in (0,1)$ is the H\"{o}lder exponent.
\end{lemma}

\smallskip
\noindent\textbf{Proof.} We divide the proof into five steps.

\smallskip
1.
It follows from \eqref{6.8} that, for $\rho\in [0,\rho_{*}]$,
\begin{equation}\label{7.26}
	|b_1(\rho)|+\rho^{1-2\theta_1}|b_1'(\rho)|+\rho^{2-2\theta_1}|b_1''(\rho)|
	+|b_2(\rho)|+\rho| b_2'(\rho)|+\rho^{2}|b_2''(\rho)|\leq C.
\end{equation}
From \eqref{6.8-1} and  \eqref{A.5-1}, we have
\begin{equation}\label{7.27}
	|b_1(\rho)|+|\rho b_1'(\rho)|+|\rho^{2}b_1''(\rho)|\leq C\rho^{\frac{1}{2}-\frac{\t_2}{2\t_1}}
	\qquad\,\, \mbox{for $\rho\geq \rho^{*}$}.
\end{equation}
Using \eqref{7.3}--\eqref{7.6} and \eqref{7.26}--\eqref{7.27}, we obtain that, for  $\rho\geq \rho^{*}$,
\begin{align*}
	&\Big\vert\int_{0}^{\rho}k(s)^{\lambda_1}k'(s)^{-\frac{1}{2}}a_1''(s)\,\mathrm{d}s\Big\vert
	\leq C\int_{0}^{\rho^{*}}s^{-1+\theta_1}\mathrm{d}s
	+C\int_{\rho_{*}}^{\rho} s^{-1-\t_2}\,\mathrm{d}s\leq C,\\
	&\Big\vert\int_{0}^{\rho}k(s)^{\lambda_1+1}k'(s)^{-\frac{1}{2}}b_1''(s)\,\mathrm{d}s\Big\vert
	+\Big\vert\int_{0}^{\rho}s k(s)^{-\lambda_1}k'(s)^{\frac{1}{2}}a_1''(s)\,\mathrm{d}s\Big\vert
	\leq C\ln \rho,
\end{align*}
which, with \eqref{6.8-1}, yields that
$|b_2(\rho)|\leq C\rho^{\frac{1}{2}-\t_2(1+\frac{1}{2\t_1})}$
for $\r\geq \rho^{*}$.
Moreover, by calculating the derivatives explicitly, we obtain
\begin{equation}\label{7.29-1}
	\r|b_2'(\rho)|+\r^2|b_2''(\r)|\leq C\rho^{\frac{1}{2}-\t_2(1+\frac{1}{2\t_1})}\qquad\,\, \text{for $\r\geq \rho^{*}$}.
\end{equation}

2. For the remainder term $g_2(\rho,u)$, recalling from \cite[Proof of Theorem 2.2]{Chen-LeFloch-2000}, $g_2$ satisfies
\begin{equation*}
	\left\{\begin{aligned}
		\dis&\partial_{\rho\rho}g_2(\rho,u)-k'(\rho)^2\partial_{uu}g_2(\rho,u)
		=ub_2''(\rho)k(\rho)^{2\lambda_1+2}f_{\lambda_1+1}(\frac{u}{k(\rho)})
		+\frac{P''(\rho)}{\rho}\partial_{u}g_1(\rho,u),\\
		\dis &g_2(0,u)=0,\quad \partial_{\rho}g_2(0,u)=0,
	\end{aligned}
	\right.
\end{equation*}
where $f_{\lambda_1}(y)=[1-y^2]_{+}^{\lambda_1}$.
Similar to the arguments for Lemma \ref{lem6.4-0}, we obtain
\begin{align}
	&k'(\rho)g_2(\rho,u)\nonumber\\
	&=\frac{1}{2\rho}
	\int_{0}^{\rho}d(s)k'(s)\Big(g_2(s,u+k(\rho)-k(s))+g_2(s,u-k(\rho)+k(s))\Big)\,\mathrm{d}s
	\nonumber\\
	&\quad +\frac{1}{2\rho}\int_{0}^{\rho}sb_2''(s)k(s)^{2\lambda_1+2}
	\Big(\int_{u-k(\rho)+k(s)}^{u+k(\rho)-k(s)}y f_{\lambda_1+1}(\frac{y}{k(s)})\,\mathrm{d}y\Big)\,\mathrm{d}s\nonumber\\
	&\quad +\frac{1}{2\rho}\int_{0}^{\rho}P''(s)\Big(g_1(s,u+k(\rho)-k(s))-g_{1}(s,u-k(\rho)+k(s))\Big)\,\mathrm{d}s,\label{7.31}
\end{align}
which yields
\begin{align}
	&\|k'(\rho)g_2(\rho,\cdot)\|_{L_{u}^{\infty}(\R)}\nonumber\\
	&=\frac{1}{\rho}\int_{0}^{\rho}d(s)\|k'(s)g_2(s,\cdot)\|_{L_{u}^{\infty}(\R)}\,\mathrm{d}s
	+\frac{C}{\rho}\int_{0}^{\rho}s|b_2''(s)|k(s)^{2\lambda_1+4}\,\mathrm{d}s\nonumber\\
	&\quad +\frac{C}{\rho}\int_{0}^{\rho}P''(s)\|g_1(s,\cdot)\|_{L_{u}^{\infty}(\R)}\,\mathrm{d}s\nonumber\\
	&:=\frac{1}{\rho}\int_{0}^{\rho}d(s)\|k'(s)g_2(s,\cdot)\|_{L_{u}^{\infty}(\R)}\,\mathrm{d}s+I_1+I_2.\label{7.31-1}
\end{align}
It follows from  \eqref{A.1-1}--\eqref{A.2-1}, \eqref{A.4-1}--\eqref{A.5-1},  \eqref{7.26}, \eqref{7.29-1},
and Lemma \ref{lem7.1} that, for $\rho\geq \rho^{*}$
\begin{align}
	|I_1|&=\Big\vert\frac{C}{\rho}\Big(\int_{0}^{\r_{*}}+\int_{\r_{*}}^{\r^{*}}+\int_{\r^{*}}^{\r}\Big)
	s|b_2''(s)|k(s)^{2\lambda_1+4}\,\mathrm{d}s\Big\vert\nonumber\\
	&\leq \frac{C}{\rho}\Big(1+\int_{\r^{*}}^{\r}s^{-\frac{1}{2}+2\t_2+\frac{\t_1}{2\t_2}}\,\mathrm{d}s\Big)
	\leq C\r^{-\frac{1}{2}+2\t_2+\frac{\t_1}{2\t_2}}\leq C\r^{2\t_2},\label{7.31-2}\\
	|I_2|&=\Big\vert\frac{C}{\rho}\Big(\int_{0}^{\r_{*}}+\int_{\r_{*}}^{\r^{*}}+\int_{\r^{*}}^{\r}\Big)
	P''(s)\|g_1(s,\cdot)\|_{L_{u}^{\infty}(\R)}\,\mathrm{d}s\Big\vert\nonumber\\
	&\leq \frac{C}{\rho}\Big(1+\int_{\r^{*}}^{\r}s^{2\t_2-1}\|g_1(s,\cdot)\|_{L_{u}^{\infty}(\R)}\,\mathrm{d}s\Big)
	\leq \begin{cases}
		C\r^{2\t_2}\quad&\text{if }\t_2<\t_1,\\
		C\r^{2\t_2}\ln \r\,\,\, &\text{if }\t_2=\t_1.
	\end{cases}\label{7.31-2-1}
\end{align}
Substituting \eqref{7.31-2}--\eqref{7.31-2-1} into \eqref{7.31-1},
applying Lemma \ref{lemB.2} and Corollary \ref{corB.1}, and using \eqref{A.5-1},
we obtain that, for $\rho\geq \rho^{*}$,
\begin{equation*}
	\|g_2(\r,\cdot)\|_{L_{u}^{\infty}(\R)}
	\leq \begin{cases}
		C\rho^{1+\t_2}\quad &\text{if $\t_2<\t_1$},\\
		C\rho^{1+\t_2}\ln \r\,\,\, &\text{if $\t_2=\t_1$}.
	\end{cases}
\end{equation*}

3. Denoting $\widetilde{f}(s)=sf_{\lambda_1+1}(s)$ and applying $\partial_{u}$ to \eqref{7.31}, we have
\begin{align}\label{7.31-5}	
	&k'(\rho)\partial_{u}g_2(\rho,u)\nonumber\\
	&=\frac{1}{2\rho}\int_{0}^{\rho}d(s)k'(s)\Big(\partial_{u}g_2(s,u+k(\rho)-k(s))
	+\partial_{u}g_2(s,u-k(\rho)+k(s))\Big)\,\mathrm{d}s\nonumber\\
	&\quad +\frac{1}{2\rho}\int_{0}^{\rho}sb_2''(s)k(s)^{2\lambda_1+3}
	\Big(\widetilde{f}(\frac{u+k(\r)-k(s)}{k(s)})
	-\widetilde{f}(\frac{u-k(\r)+k(s)}{k(s)})\Big)\,\mathrm{d}s\nonumber\\
	&\quad +\frac{1}{2\rho}\int_{0}^{\rho}P''(s)
	\Big(\partial_{u}g_1(s,u+k(\rho)-k(s))-\partial_{u}g_{1}(s,u-k(\rho)+k(s))\Big)\,\mathrm{d}s.
\end{align}
Since $|\widetilde{f}(s)|\leq 1$, by similar arguments as in Step 2, we obtain
\begin{equation*}
	\|\partial_{u}g_2(\r,\cdot)\|_{L_{u}^{\infty}(\R)}\leq C\rho^{1+\t_2}\qquad\, \text{for $\rho\geq \rho^{*}$}.
\end{equation*}

4.
Applying $\partial_{u}^{\lambda_1}$ to \eqref{7.31-5}, we have
\begin{align*}
	&k'(\rho)\partial_{u}^{\lambda_1+1}g_2(\rho,u)\\
	&=\frac{1}{2\rho}\int_{0}^{\rho}d(s)k'(s)\Big((\partial_{u}^{\lambda_1+1}g_2)(s,u+k(\rho)-k(s))
	+(\partial_{u}^{\lambda_1+1}g_2)(s,u-k(\rho)+k(s))\Big)\,\mathrm{d}s
	\\&\quad +\frac{1}{2\rho}\int_{0}^{\rho}sb_2''(s)k(s)^{\lambda_1+3}
	\Big(\widetilde{f}^{(\lambda_1)}(\frac{u+k(\rho)-k(s)}{k(s)})
	-\widetilde{f}^{(\lambda_1)}(\frac{u-k(\rho)+k(s)}{k(s)})\Big)\,\mathrm{d}s
	\\&\quad +\frac{1}{2\rho}\int_{0}^{\rho}P''(s)\Big((\partial_{u}^{\lambda_1+1}g_1)(s,u+k(\rho)-k(s))
	-(\partial_{u}^{\lambda_1+1}g_{1})(s,u-k(\rho)+k(s))\Big)\,\mathrm{d}s,
\end{align*}
where $\widetilde{f}^{(\lambda_1)}:=\partial_{s}^{\lambda_1}\widetilde{f}(s)$.
Since $\widetilde{f}^{(\lambda_1)}$ is uniformly bounded, similar arguments as in Step 2 yield
\begin{equation*}
	\|\big(\partial_{u}^{\lambda_1+1}g_2\big)(\rho,\cdot)\|_{L_{u}^{\infty}(\mathbb{R})}
	\leq C\rho^{1+\theta_2}\qquad\, \text{for $\rho \geq \rho^{*}$}.
\end{equation*}

\smallskip
5. Applying $\partial_{u}^{\lambda_1+\alpha_0}$ to \eqref{7.31-5}, we have
\begin{align*}
	&k'(\rho)\partial_{u}^{\lambda_1+1+\alpha_0}g_2(\rho,u)\\
	&=\frac{1}{2\rho}\int_{0}^{\rho}d(s)k'(s)(\partial_{u}^{\lambda_1+1+\alpha_0}g_2)(s,u+k(\rho)-k(s))\,\mathrm{d}s\\
	&\quad +\frac{1}{2\rho}\int_{0}^{\rho}d(s)k'(s)(\partial_{u}^{\lambda_1+1+\alpha_0}g_2)(s,u-k(\rho)+k(s))\,\mathrm{d}s
	\\&\quad +\frac{1}{2\rho}\int_{0}^{\rho}sb_2''(s)k(s)^{\lambda_1+3-\alpha_0}
	\Big(\widetilde{f}^{(\lambda_1)}(\frac{u+k(\rho)-k(s)}{k(s)})
	-\widetilde{f}^{(\lambda_1)}(\frac{u-k(\rho)+k(s)}{k(s)})\Big)\,\mathrm{d}s
	\\&\quad
	+\frac{1}{2\rho}\int_{0}^{\rho}P''(s)(\partial_{u}^{\lambda_1+1+\alpha_0})g_1(s,u+k(\rho)-k(s))
	\,\mathrm{d}s\\
	&\quad -\frac{1}{2\rho}\int_{0}^{\rho}P''(s)(\partial_{u}^{\lambda_1+1+\alpha_0})g_{1}(s,u-k(\rho)+k(s))\,\mathrm{d}s.
\end{align*}
Noting that $\widetilde{f}^{(\lambda_1+\alpha_0)}(s)$ is uniformly bounded,  by similar arguments as in Step 2,
we have
\begin{equation*}
	\|\partial_{u}^{\lambda_1+1+\alpha_0}g_2(\rho,\cdot)
	\|_{L_{u}^{\infty}(\mathbb{R})}\leq C\rho^{1+\theta_2}
	\qquad \text{for $\rho \geq \rho^{*}$}.
\end{equation*}
This completes the proof.
$\hfill\square$

\smallskip
The following lemma provides the explicit singularities of $\chi(\r,u-s)$ and $(\sigma-u\chi)(\r,u-s)$.

\smallskip
\begin{lemma}\label{lem7.3}
	The fractional derivatives $\partial_{u}^{\lambda_1+1}\chi$ and $\partial_{u}^{\lambda_1+1}(\sigma-u\chi)$
	admit the expansions{\rm :}
	\begin{align}
		&\partial_{s}^{\lambda_1+1}\chi(\rho,u-s)\nonumber\\
		&=\sum\limits_{\pm}\Big(A_{1,\pm}(\rho)\,\delta(s-u\pm k(\rho))+A_{2,\pm}(\rho)\,H(s-u\pm k(\rho))\Big)
		\nonumber\\
		&\quad +\sum\limits_{\pm}\Big(A_{3,\pm}(\rho)\,PV(s-u\pm k(\rho))+A_{4,\pm}(\rho)\, Ci(s-u\pm k(\rho))\Big)
		\nonumber		\\
		&\quad +r_{\chi}(\rho,u-s),\label{7.40}\\[2mm]
		&\partial_{s}^{\lambda_1+1}(\sigma-u\chi)(\rho,u-s)\nonumber\\
		&=\sum\limits_{\pm}(s-u)\Big(B_{1,\pm}(\rho)\,\delta(s-u\pm k(\rho))
		+B_{2,\pm}(\rho)\,H(s-u\pm k(\rho))\Big)\nonumber\\
		&\quad +\sum\limits_{\pm}(s-u)\Big(B_{3,\pm}(\rho)\,PV(s-u\pm k(\rho))+B_{4,\pm}\,Ci(s-u\pm k(\rho))\Big)
		\nonumber\\
		&\quad +\sum\limits_{\pm}\Big(B_{5,\pm}(\rho)\,H(s-u\pm k(\rho))+B_{6,\pm}(\rho)\,Ci(s-u\pm k(\rho))\Big)
		\nonumber\\&
		\quad +r_{\sigma}(\rho,u-s),\label{7.41}
	\end{align}
	where $\delta$ is the Dirac measure, $H$ is the Heaviside function, $PV$ is the principle value distribution,
	and $Ci$ is the Cosine integral{\rm :}
	$$
	Ci(s):=-\int_{|s|}^{\infty}\frac{\cos y}{y}\,\mathrm{d}y
	=\log |s|+\int_{0}^{|s|}\frac{\cos y-1}{y}\,\mathrm{d}y+C_0
	\qquad \mbox{for $s\in \R$}
	$$
	for some constant $C_0>0$.
	The remainder terms $r_{\chi}$ and $r_{\sigma}$ are H\"{o}lder continuous functions.
	Moreover, there exists a positive constant $C=C(\gamma_1,\gamma_2,\rho_{*},\rho^{*})$ such that, for  $\rho\geq \rho^{*}$,
	\begin{align*}
		&\sum\limits_{j=1,\pm}^{4}|A_{j,\pm}(\rho)|+\sum\limits_{j=1,\pm}^6|B_{j,\pm}|
		\leq C\rho^{\frac{1}{2}-\frac{\t_2}{2}},
		\\ &\|r_{\chi}(\rho,\cdot)\|_{C^{\alpha_1}(\mathbb{R})}
		\leq C\rho,\quad \|r_{\sigma}(\rho,\cdot)\|_{C^{\alpha_1}(\mathbb{R})}\leq C\rho^{1+\t_2},
	\end{align*}
	where $\alpha_1\in (0,\alpha_0]$ is the common H\"{o}lder exponent of $r_{\chi}$ and $r_{\sigma}$.
\end{lemma}

\smallskip
\noindent\textbf{Proof. }
From \cite[Lemma 6.4]{Schrecker-Schulz-2020}, we obtain \eqref{7.40}--\eqref{7.41},
where the coefficients are given by
\begin{align*}
	&A_{1, \pm}(\rho)=a_{1}(\rho) k(\rho)^{\lambda_1} A_{1}^{\lambda_{1}},
	\quad\,\,\,\, A_{2, \pm}(\rho)=\pm a_{1}(\rho) k(\rho)^{\lambda_1-1} A_{3}^{\lambda_1}+ a_{2}(\rho) k(\rho)^{\lambda_1+1} A_{1}^{\lambda_1+1},\\
	&A_{3, \pm}(\rho)=\pm a_{1}(\rho) k(\rho)^{\lambda_1} A_{2}^{\lambda_1},
	\quad 	A_{4, \pm}(\rho)=\pm a_{1}(\rho) k(\rho)^{\lambda_1-1} A_{4}^{\lambda_1} \pm a_{2}(\rho) k(\rho)^{\lambda_1+1} A_{2}^{\lambda_1+1},\\[1mm]
	&r_{\chi}(\rho,u-s)= a_1(\rho)k(\rho)^{\lambda_1-1}\tilde{q}(\frac{s-u}{k(\rho)})
	+a_2(\rho)k(\rho)^{\lambda_1+1}\tilde{r}(\frac{s-u}{k(\rho)})-A_{4}^{\lambda_1} k(\rho)^{\lambda_1-1}(\mathrm{log}k(\rho))^2\\
	&\qquad\qquad\qquad\,
	+\partial_{s}^{\lambda_1+1}g_1(\rho,u-s),
\end{align*}
where $A_{i}^{\lambda_1}, i=1,\cdots,4$, are constants depending only on $\lambda_1$,
and $\tilde{r}$ and $\tilde{q}$ are uniformly bounded H\"{o}lder continuous functions.
Thus, using Lemma \ref{lem7.1}, we see that, for $\rho\geq \rho^{*}$,
\begin{align*}
	&|A_{i,\pm}(\rho)|\leq C\rho^{\frac{1}{2}-\frac{\t_2}{2\t_1}}\rho^{\theta_2(\frac{1}{2\t_1}-\frac{1}{2})}
	\leq C\rho^{\frac{1}{2}-\frac{1}{2}\t_2}\qquad\,\, \text{for $\,i=1,3$},\\
	&|A_{j,\pm}(\rho)|\leq C\rho^{\frac{1}{2}-\frac{3\t_2}{2}}+C\rho^{\frac{1}{2}-\frac{\t_2}{2\t_1}-\t_2}\rho^{\t_2(\frac{1}{2\t_1}+\frac{1}{2})}
	\leq C\rho^{\frac{1}{2}-\frac{\t_2}{2}}
	\qquad\,\, \text{for $\,j=2,4$}.\\
	&\|r_{\chi}(\rho,\cdot)\|_{C^{\alpha}(\R)}:=\|r_{\chi}(\rho,\cdot)\|_{L^{\infty}(\R)}+[r_{\chi}(\rho,\cdot)]_{C^{\alpha}(\R)}\\
	&\leq C\big(\rho^{\frac{1}{2}-\frac{\t_2}{2}}+\rho^{\frac{\t_2}{2\t_1}-\frac{3}{2}\t_2}|\ln\rho|^2+\rho\big)
	\leq C\rho.
\end{align*}
Similarly, we have
\begin{align*}
	&B_{1,\pm}(\rho)=b_1(\rho)k(\rho)^{\lambda_1}A_{1}^{\lambda_1},\quad
	B_{2,\pm}(\rho)=\pm b_{1}(\rho)k(\rho)^{\lambda_1-1}A_3^{\lambda_1}+b_2(\rho)k(\rho)^{\lambda_1+1}A_{1}^{\lambda_1+1},\\
	&B_{3,\pm}(\rho)=\pm b_1(\rho)k(\rho)^{\lambda_1}A_{2}^{\lambda_1},\quad
	B_{4,\pm}(\rho)=\pm b_1(\rho)k(\rho)^{\lambda_1-1}A_{4}^{\lambda_1}\pm b_2(\rho)k(\rho)^{\lambda_1+1}A_{2}^{\lambda_1+1},\\
	&B_{5,\pm}(\rho)=(\lambda_1+1)b_1(\rho)k(\rho)^{\lambda_1}A_{1}^{\lambda_1},\quad
	B_{6,\pm}(\rho)=\pm (\lambda_1+1)b_1(\rho)k(\rho)^{\lambda_1}A_{2}^{\lambda_1},\\
	&r_{\sigma}(\rho,u-s)\\
	&=(s-u)\Big(b_1(\rho)k(\rho)^{\lambda_1-1}\big(-A_4^{\lambda_1}(\mathrm{log}k(\rho))^2+\tilde{q}(\frac{s-u}{k(\rho)})\big)
	+b_2(\rho)k(\rho)^{\lambda_1+1}\tilde{r}(\frac{s-u}{k(\rho)})\Big)
	\\
	&\quad  +(\lambda_1+1)\Big(b_1(\rho)k(\rho)^{\lambda_1}\tilde{r}(\frac{s-u}{k(\rho)})
	+b_2(\rho)k(\rho)^{\lambda_1+2}\tilde{\ell}(\frac{s-u}{k(\rho)})\Big)
	+\partial_{s}^{\lambda_1+1}g_2(\rho,u-s),
\end{align*}
where $\tilde{\ell}$ is also a uniformly bounded H\"{o}lder continuous function.
Using Lemma \ref{lem7.2}, we conclude that, for $\rho\geq \rho^{*}$,
\begin{align*}
	&|B_{i,\pm}(\rho)|\leq C\rho^{\frac{1}{2}-\frac{\t_2}{2\t_1}}\rho^{\theta_2(\frac{1}{2\t_1}-\frac{1}{2})}
	\leq C\rho^{\frac{1}{2}-\frac{\t_2}{2}}\qquad \text{for }i=1,3,5,6,\\
	&|B_{j,\pm}(\rho)|\leq C\rho^{\frac{1}{2}-\frac{3\t_2}{2}}+C\rho^{\frac{1}{2}-(\frac{1}{2\t_1}+1)\t_2}\rho^{\theta_2(\frac{1}{2\t_1}+\frac{1}{2})}
	\leq C\rho^{\frac{1}{2}-\frac{\t_2}{2}}\qquad \text{for }j=2,4,\\
	&\|r_{\sigma}(\rho,\cdot)\|_{C^{\alpha}(\R)}:=\|r_{\sigma}(\rho,\cdot)\|_{L^{\infty}(\R)}+[r_{\sigma}(\rho,\cdot)]_{C^{\alpha}(\R)}
	\leq C\big(\rho^{\frac{1}{2}-\frac{\t_2}{2}}|\ln \rho|^2+\rho^{1+\t_2}\big)
	\leq C\rho^{1+\t_2}.
\end{align*}
This completes the proof.
$\hfill\square$

\section{\,Uniform Estimates of Approximate Solutions}\label{UE}
As in \cite{Chen-He-Wang-Yuan-2021}, we construct the approximate solutions via the following approximate
free boundary problem for CNSPEs:
\begin{equation}\label{2.1}
	\left\{
	\begin{aligned}
		&\rho_{t}+(\rho u)_r+\frac{2}{r}\rho u=0,\\
		&(\rho u)_t+(\rho u^2+P(\rho))_{r}+\frac{2}{r}\rho u^2+\frac{\r}{r^{2}}\int_{a}^{r}\r(t,y)\,y^{2}\mathrm{d}y
       =\v\Big(\rho(u_r+\frac{2}{r}u)\Big)_{r}
		-\frac{2\v}{r}\rho_{r}u,
	\end{aligned}
	\right.
\end{equation}
for $(t,r)\in \Omega_{T}$ with
\begin{equation}\label{2.2}
	\Omega_{T}=\{(t,r)\in [0,\infty)\times \R\,:\,a\leq r\leq b(t),\;0\leq t\leq T\},
\end{equation}
where $\{r=b(t)\,:\,0\leq t\leq T\}$ is a free boundary determined by
\begin{equation}\label{2.3}
	b'(t)=u(t,b(t))\quad \text{for $t>0$},\qquad\,\,\,	b(0)=b,
\end{equation}
and $a=b^{-1}$ with $b\gg 1$. On the free boundary $r=b(t)$, we impose the stress-free boundary condition:
\begin{equation}\label{2.4}
	\big(P(\rho)-\v\rho(u_{r}+\frac{2}{r}u)\big)(t,b(t))=0\qquad \text{for }t>0.
\end{equation}
On the fixed boundary $r=a=b^{-1}$, we impose the Dirichlet boundary condition:
\begin{equation}\label{2.5}
	u(t,r)\vert_{r=a}=0\qquad \text{for }t>0.
\end{equation}
The initial condition is
\begin{equation}\label{2.6}
	(\rho,\rho u)\vert_{t=0}=(\rho_{0}^{\v,b},\rho_{0}^{\v,b}u_{0}^{\v,b})\qquad\, \text{for }r\in [a,b].
\end{equation}

\subsection{\,Basic estimates}\label{BE}
Denote
\begin{align}
	&E_{0}^{\varepsilon, b}:=\omega_3\int_{a}^{b} \rho_{0}^{\varepsilon, b}\Big(\frac{1}{2}\big|u_{0}^{\varepsilon, b}\big|^{2}
	+e(\rho_{0}^{\varepsilon, b})\Big) \,r^{2}\mathrm{d}r,
	\qquad E_{1}^{\varepsilon, b}:=\omega_{3}\varepsilon^{2}
	\int_{a}^{b}\Big|\big(\sqrt{\rho_{0}^{\varepsilon, b}}\big)_{r}\Big|^{2}\, r^{2}\mathrm{d} r.\nonumber
\end{align}
For given total energy $E_{0}^{\v,b}>0$,
the critical mass $M_{\rm c}^{\v,b}$ is defined in \eqref{1.18-0}--\eqref{1.18-4} by replacing $E_{0}$ with $E_{0}^{\v,b}$.

For the approximate initial data $(\rho_{0}^{\v},m_{0}^{\var})$ imposed in \eqref{1.19}
satisfying \eqref{1.27}--\eqref{1.27-2}, using similar arguments
in \cite[Appendix A]{Chen-He-Wang-Yuan-2021}, we can construct a sequence of smooth
functions $(\rho_{0}^{\v,b},u_{0}^{\v,b})$ defined on $[a,b]$, which is compatible with the boundary
conditions \eqref{2.4}--\eqref{2.5}, such that
\begin{itemize}
\item [(\rmnum{1})] There exists a constant $C_{\v,b}>0$ depending on $(\v,b)$ so that,
for all $\v\in (0,1]$ and $b>1$,
	\begin{equation}\label{2.7-4}
		0<C_{\v,b}^{-1}\leq \rho_{0}^{\v,b}(r)\leq C_{\v,b}<\infty.
	\end{equation}
\item [(\rmnum{2})] For all $\v\in (0,1]$ and $b>1$,
	\begin{align}
		&\int_{a}^{b}\rho_{0}^{\v,b}(r)r^{2}dr=\frac{M}{\omega_3},
		\qquad E_{0}^{\v,b}\leq C(1+E_{0}),\qquad  E_{1}^{\v,b}\leq C(1+M)\v,\label{2.9}\\
		&\rho_{0}^{\v,b}(b)\cong b^{-(3-\alpha)}\qquad\,\, \text{with }\alpha:=\min \{\frac{1}{2},\frac{3(\gamma_1-1)}{\gamma_1}\}.
		\label{2.9b}
	\end{align}
\item [(\rmnum{3})] For each fixed $\v\in (0,1]$, as $b\to \infty$, $(E_{0}^{\v,b},E_{1}^{\v,b})\rightarrow (E_0^{\v},E_{1}^{\v})$ and
	\begin{equation}\label{2.8}
		(\rho_{0}^{\v,b},\rho_{0}^{\v,b}u_{0}^{\v,b})\longrightarrow (\rho_0^{\v},m_{0}^{\v})
		\qquad \text{in $L^{\tilde{q}}([a,b];r^{2}dr)\times L^1([a,b];r^{2}dr)$ for $\tilde{q}\in\{1,\g_2\}$}.
	\end{equation}
\item[(\rmnum{4})] For each fixed $\v\in (0,\v_0]$, there exists a large constant $\mathcal{B}(\v)>0$ such that
	\begin{equation}\label{2.9-1}
		M<M_{\rm c}^{\v,b}\qquad \text{for }b\geq \mathcal{B}(\v)\text{ and }\g_2\in (\frac{6}{5},\frac{4}{3}],
	\end{equation}
	where
	$M_{\rm c}^{\v,b}$ is defined in \eqref{1.18-0}--\eqref{1.18-4} by replacing $E_{0}$ with $E_{0}^{\v,b}$.
\end{itemize}

\smallskip
We point out that \eqref{2.9b} is important for us to close the BD-type entropy estimate
in Lemma \ref{lem2.2} and to obtain the higher integrability of the density in Lemma \ref{lem2.4} below.

\smallskip
Once the free boundary problem \eqref{2.1}--\eqref{2.6} is solved,
we define the potential function $\Phi$ to be the solution of the Poisson equation:
\begin{equation*}
	\Delta \Phi=\r\mathbf{I}_{\Omega_{t}},\qquad \lim\limits_{|\mathbf{x}|\to \infty}\Phi(\mathbf{x})=0,
\end{equation*}
with $\Omega_{t}:=\{\mathbf{x}\in \R^3\,:\,a\leq |\mathbf{x}|\leq b(t)\}$,
for which $\rho$ has been extended to be zero outside $\Omega_{t}$.
In fact, we can show that $\Phi(t, \mathbf{x})=\Phi(t, r)$ with
\begin{equation}\label{2.9-3}
	\Phi_{r}(t, r)=\left\{\begin{aligned}
		0\qquad \qquad \qquad\quad &\quad\text { for } 0 \leq r \leq a, \\
		\frac{1}{r^{2}} \int_{a}^{r} \rho(t, y) \, y^{2}\mathrm{d}y & \quad\text { for } a \leq r \leq b(t), \\
		\frac{M}{\omega_{3}} \frac{1}{r^{2}}\qquad \qquad\quad&\quad\text{ for } r \geq b(t),
	\end{aligned}\right.
\end{equation}
so that $\Phi(t,r)$ can be recovered by integrating \eqref{2.9-3}.

In this section, parameters $(\v,b)$ are fixed with $\v\in (0,\v_0]$ and
$b\geq \max\{\rho_{*}^{-\frac{\g_1}{3}},\mathcal{B}(\v)\}$ such that \eqref{2.9-1} holds
and $\rho_{0}^{\v,b}(b)\leq \rho_{*}$.
The global existence of smooth solutions of our approximate problem \eqref{2.1}--\eqref{2.6}
whose initial data satisfy \eqref{2.7-4}--\eqref{2.9-1} and pressure satisfies \eqref{1.3}--\eqref{1.5}
can be obtained by using similar arguments in \cite[\S 3]{Duan-Li-2015} with $\g_2\in (\frac{4}{3},\infty)$,
or with $ \g_2\in (\frac{6}{5},\frac{4}{3}]$ and $M<M_{\rm c}^{\v,b}(\g_2)$,
so the details are omitted here for simplicity.

\smallskip
Noting that the upper and lower bounds of $\rho^{\v,b}$ in \cite{Duan-Li-2015} depend on
parameters $(\v, b)$, we now establish some uniform estimates, independent of $b$,
such that the limit: $b\to \infty$ can be taken to obtain the global weak solutions
of problem \eqref{1.8} and \eqref{1.19}--\eqref{1.19-1} in \S 6 below as approximate
solutions of problem \eqref{1.1} and \eqref{1.11}--\eqref{1.11-1}.
Throughout this section, we drop the superscript in both the approximate solutions
$(\rho^{\v,b},u^{\v,b})(r)$ and the approximate initial data $(\rho_{0}^{\v,b},u_{0}^{\v,b})$
for simplicity.

\smallskip
For smooth solutions, it is convenient to analyze \eqref{2.1}--\eqref{2.6}
in the Lagrangian coordinates. It follows from \eqref{2.3} that
\begin{equation*}
	\frac{\mathrm{d}}{\mathrm{d}t}\int_{a}^{b(t)}\rho(t,r)\,r^{2}\mathrm{d}r
	=(\rho u)(t,b(t))b(t)^{2}-\int_{a}^{b(t)}(\rho u r^{2})_{r}(t,r)\,\mathrm{d}r=0,
\end{equation*}
which implies
\begin{equation}\label{2.11}
	\int_{a}^{b(t)}\rho(t,r)\,r^{2}\mathrm{d}r
	=\int_{a}^{b}\rho_{0}(r)\,r^{2}\mathrm{d}r=\frac{M}{\omega_3}  \qquad\, \mbox{for all $t\geq 0$}.
\end{equation}
For $r\in [a,b(t)]$ and $t\in [0,T]$, the Lagrangian coordinates $(\tau, x)$ are defined by
$$
\tau=t, \qquad x(t,r)=\int_{a}^{r}\rho(t,y)\,y^{2}\mathrm{d}y,
$$
which translate $[0,T]\times [a,b(t)]$ into a fixed domain $[0,T]\times [0,\frac{M}{\omega_3}]$.
By direct calculation, we see that
$\nabla_{(t,r)}x=(-\rho u r^{2},\rho r^{2})$,
$\nabla_{(t,r)}\tau=(1,0)$,
$\nabla_{(\tau,x)}r=(u,\rho^{-1}r^{-2})$,
and $\nabla_{(\tau,x)}t=(1,0)$.
In the Lagrangian coordinates, the initial-boundary value problem \eqref{2.1}--\eqref{2.6} becomes
\begin{equation}\label{2.13}
	\left\{
	\begin{aligned}
		&\rho_{\tau}+\rho^2(r^{2}u)_{x}=0,\\
		&u_{\tau}+r^{2}P_{x}=-\frac{x}{r^2}+\v r^{2}(\rho^2(r^{2}u)_{x})_{x}-2\v r\rho_{x} u,
	\end{aligned}
	\right.
\end{equation}
for $(\tau,x)\in [0,T]\times [0,\frac{M}{\omega_3}]$, and
\begin{equation}\label{2.14}
	u(\tau, 0)=0,\quad (P-\v\rho^2(r^{2}u)_{x})(\tau,\frac{M}{\omega_{3}})=0\qquad\,\, \text{for }\tau\in [0,T],
\end{equation}
where $r=r(\tau,x)$ is defined by
$\,\frac{\mathrm{d}}{\mathrm{d}\tau}r(\tau,x)=u(\tau,x)$
for $(\tau,x)\in [0,T]\times [0,\frac{M}{\omega_3}]$,
and the fixed boundary $x=\frac{M}{\omega_3}$ corresponds to
the free boundary: $b(\tau)=r(\tau,\frac{M}{\omega_3})$ in the Eulerian coordinates.

\smallskip
\begin{lemma}[Basic energy estimate]\label{lem2.1}
	The smooth solution $(\rho, u)(t,r)$ of problem \eqref{2.1}--\eqref{2.6} satisfies
	\begin{equation*}
		\begin{aligned}
			&\int_{a}^{b(t)}\Big(\frac{1}{2}\rho u^2+\rho e(\rho)\Big)\,r^{2}\mathrm{d}r
			-\frac{1}{2}\int_{a}^{\infty}\frac{1}{r^{2}}\Big(\int_{a}^{r}\rho(t,z)\,z^{2}\mathrm{d}z\Big)^2\,\mathrm{d}r\\			
			&\quad+\v\int_{0}^{t}\int_{a}^{b(s)}\Big(\rho u_{r}^2+2\frac{\rho u^2}{r^2}\Big)\,r^{2}\mathrm{d}r\mathrm{d}s
			+2\v\int_{0}^{t}(\rho u^2)(s,b(s))b(s)\,\mathrm{d}s\\
			&=\int_{a}^{b}\Big(\frac{1}{2}\rho_{0}u_{0}^2+\rho_{0}e(\rho_0)\Big)\,r^{2}\mathrm{d}r
			-\frac{1}{2}\int_{a}^{\infty}\frac{1}{r^{2}}\Big(\int_{a}^{r}\rho_{0}(t,z)z^{2}dz\Big)^2\,\mathrm{d}r,
		\end{aligned}
	\end{equation*}
	where $\rho(t, r)$ has been understood to be $0$ for $r\in[0, a]\cup (b(t),\infty)$ in the second term
	of the left-hand side {\rm (LHS)} and the second term of the right-hand side {\rm (RHS)}.
	In particular, there exists a positive constant $C(E_{0},M)$ depending only on the total initial energy $E_{0}$ and initial-mass $M$
	such that the following estimates hold for the two separate cases{\rm :}\\
	Case 1. $\dis \g_2\in (\frac{6}{5},\frac{4}{3}]$ and $M<M_{\rm c}^{\v,b}${\rm :}
	\begin{align}
		&\int_{a}^{b(t)} \rho\Big(\frac{1}{2}u^{2}+e(\rho)\Big)\,r^{2}\mathrm{d} r
+\varepsilon\int_{0}^{t}\int_{a}^{b(s)}\rho\Big(u_{r}^{2}+\frac{2u^{2}}{r^{2}}\Big)(t,r)\,r^{2}\mathrm{d}r\mathrm{d}s\nonumber\\
		&\quad +2 \varepsilon \int_{0}^{t}(\rho u^{2})(s, b(s))b(s)\, \mathrm{d} s\leq C(E_{0},M).\label{2.16-2}
	\end{align}
	Case 2. $\dis \g_2>\frac{4}{3}${\rm :}
	\begin{align}
		&\int_{a}^{b(t)} \frac{1}{2} \rho\Big(u^{2}+e(\rho)\Big)\, r^{2}\mathrm{d} r
		+\varepsilon \int_{0}^{t} \int_{a}^{b(s)} \rho\Big(u_{r}^{2}+\frac{2u^{2}}{r^{2}}\Big)(t, r)\, r^{2}\mathrm{d} r \mathrm{d} s\nonumber\\
		&\quad +2 \varepsilon \int_{0}^{t}(\rho u^{2})(s,b(s)) b(s)\, \mathrm{d} s \leq C(E_{0},M).\label{2.16-3}
	\end{align}
\end{lemma}

\noindent
{\bf Proof.} We divide the proof into three steps.
\smallskip

1.
Using \eqref{1.15} and similar calculations as in the proof \cite[Lemma 3.1]{Chen-He-Wang-Yuan-2021}, we have
\begin{equation}\label{2.21-1}
	\begin{aligned}
		&\int_{a}^{b(t)} \rho\Big(\frac{1}{2} u^{2}+e(\rho)\Big)\, r^{2}\mathrm{d} r
		- \int_{a}^{b(t)}\Big(\int_{a}^{r} \rho(t, z)\, z^{2}\mathrm{d} z\Big) \rho\, r\mathrm{d} r \\
		&\quad+\varepsilon \int_{0}^{t} \int_{a}^{b(s)}
		\Big(\rho u_{r}^{2}+2 \rho \frac{u^{2}}{r^{2}}\Big)\,r^{2}\mathrm{d} r \mathrm{d} s
		+2 \varepsilon \int_{0}^{t}(\rho u^{2})(s, b(s))\, b(s)\mathrm{d} s \\
		&=\int_{a}^{b} \rho_{0}\Big(\frac{1}{2} u_{0}^{2}+e\left(\rho_{0}\right)\Big)\, r^{2}\mathrm{d} r
		-\int_{a}^{b}\Big(\int_{a}^{r} \rho_{0}(z) z^{2} \,\mathrm{d} z\Big) \rho_{0}(r)\, r\mathrm{d} r.
	\end{aligned}
\end{equation}
\smallskip

2. We now control the second term on the LHS of \eqref{2.21-1} and the second term on the RHS of \eqref{2.21-1} to close the estimates.
By similar calculations as in \cite[Lemma 3.1]{Chen-He-Wang-Yuan-2021}, one can obtain
\begin{equation}\label{2.21-4}
	\int_{a}^{b(t)}\Big(\int_{a}^{r} \rho\, z^{2}\mathrm{d} z\Big) \rho\, r \mathrm{d} r
	=\frac{1}{2 \omega_{3}}\|\nabla \Phi\|_{L^{2}\left(\mathbb{R}^{3}\right)}^{2}
	=\frac{1}{2}\int_{a}^{\infty}\frac{1}{r^{2}}\Big(\int_{a}^{r}\rho\, z^{2}\mathrm{d}z\Big)^{2} \,\mathrm{d}r,
\end{equation}
where we have understood $\rho$ to be zero for $r\in [0,a)\cup (b(t),\infty)$ in \eqref{2.21-4}.

\smallskip
{3. Now we use the internal energy to control the gravitational potential term.
	First, we obtain from \eqref{ae1} that there exist two constants $C_{1},C_{2}>0$ depending only on $\rho^{*}$ such that
	\begin{align*}
		\big|\rho e(\rho)-\frac{\kappa_2\rho^{\gamma_2}}{\gamma_2-1}\big|
         \leq C_{1}\rho^{\max\{\gamma_{2}-\epsilon,0\}}\,\,\,\, \text{for }\rho\geq \rho^{*},
		\quad\,\,\,
    \big|\rho e(\rho)-\frac{\kappa_2\rho^{\gamma_2}}{\gamma_2-1}\big|\leq C_{2}\rho^{\gamma_2}\,\,\,\, \text{for }\rho\leq \rho^{*}.
	\end{align*}
	Thus, we have
	\begin{align}\label{ae2}
		&\Big\vert\int_{a}^{b(t)}\big(\rho e(\rho)-\frac{\kappa_2}{\gamma_2-1}\rho^{\gamma_{2}}\big)\,r^2{\rm d}r\Big\vert\nonumber\\
		&=\int_{\rho(t,r)\geq K}\Big\vert\rho e(\rho)-\frac{\kappa_2}{\gamma_2-1}\int_{a}^{b(t)}\rho^{\gamma_{2}}\Big\vert\,r^2{\rm d}r
		+\int_{\rho(t,r)\leq K}\Big\vert\rho e(\rho)-\frac{\kappa_2}{\gamma_2-1}\int_{a}^{b(t)}\rho^{\gamma_{2}}\Big\vert\,r^2{\rm d}r\nonumber\\
		&\leq C_{1}K^{-\min\{\gamma_2,\epsilon\}}\int_{a}^{b(t)}\rho^{\gamma_{2}}\,r^2{\rm d}r+C_{2}\omega_{3}^{-1}K^{\gamma_{2}-1}\,M,
	\end{align}
	where $K>\rho^{*}$ is some large constant to be chosen later.}

{Multiplying \eqref{2.9-1} by $\Phi$ and integrating by parts yield
	\begin{equation}\label{2.21-5}
		\|\nabla \Phi\|_{L^{2}(\mathbb{R}^{3})}^{2} \leq\|\Phi\|_{L^{6}(\mathbb{R}^{3})}\|\rho\|_{L^{\frac{6}{5}}(\Omega_{t})}
		\leq \sqrt{A_{3}}\|\nabla \Phi\|_{L^{2}(\mathbb{R}^{3})}\|\rho\|_{L^{\frac{6}{5}}(\Omega_{t})},
	\end{equation}
	where we have used the positive constant $A_3:=\frac{4}{3}\omega_{4}^{-\frac{2}{3}}>0$ that is
	the sharp constant for the Sobolev inequality in $\R^3$ (see Lemma \ref{lemC.0}). Then it follows from \eqref{2.21-4} and \eqref{2.21-5} that
	\begin{align}\label{2.21-6}
		&\int_{a}^{b(t)}\Big(\int_{a}^{r} \rho\, z^{2}\mathrm{d} z\Big) \rho\, r \mathrm{d} r=\frac{1}{2 \omega_{3}}\|\nabla \Phi\|_{L^{2}(\mathbb{R}^{3})}^{2}\leq \frac{2}{3\omega_{3}}\omega_{4}^{-\frac{2}{3}}\|\rho\|_{L^{\frac{6}{5}}(\Omega_{t})}^2\nonumber\\
		&\leq \frac{2}{3\omega_{3}}\omega_{4}^{-\frac{2}{3}}\Big(\int_{\Omega_{t}}\rho^{\frac{6(\g_2-1)}{5\g_2-6}}
		\big(\beta \rho+\rho e(\rho)\big)^{-\frac{1}{5\g_2-6}}\,\mathrm{d}\mathbf{x}\Big)^{\frac{5\g_2-6}{3(\g_2-1)}}
		\Big(\int_{\Omega_{t}}\big(\beta \rho +\rho e(\rho)\big)\,\mathrm{d}\mathbf{x}\Big)^{\frac{1}{3(\g_2-1)}}\nonumber\\
		&\leq \frac{2}{3}\omega_{4}^{-\frac{2}{3}}\omega_{3}^{\frac{4-3\gamma_2}{3(\g_2-1)}}
		\Big(\int_{\Omega_{t}}C_{\max}(\beta)\rho\,\mathrm{d}\mathbf{x}\Big)^{\frac{5\g_2-6}{3(\g_2-1)}}
		\Big(\int_{a}^{b(t)}\big(\beta\rho+\rho e(\rho)\big)\,r^2\mathrm{d}r\Big)^{\frac{1}{3(\g_2-1)}}\nonumber\\
		&=B_{\beta}M^{\frac{5\g_2-6}{3(\g_2-1)}}\Big(\int_{a}^{b(t)}\big(\beta \rho +\rho e(\rho)\big)\,r^2\mathrm{d}r\Big)^{\frac{1}{3(\g_2-1)}},
	\end{align}
	where $B_{\beta}$ is the constant defined in \eqref{1.18-4}.}

{When $\gamma_2>\frac{4}{3}$, {\it i.e.}, $\frac{1}{3(\gamma_2-1)}<1$, it follows from \eqref{2.21-6}
	by taking $\beta=1$ that
	\begin{align}
		&\int_{a}^{b(t)} \rho e(\rho)\, r^{2}\mathrm{d} r
		- \int_{a}^{b(t)}\Big(\int_{a}^{r} \rho\, z^{2}\mathrm{d} z\Big) \rho\, r\mathrm{d} r \nonumber\\
		&\geq \int_{a}^{b(t)} \rho e(\rho)\, r^{2}\mathrm{d} r
		-B_{1}M^{\frac{5 \gamma_2-6}{3(\gamma_2-1)}}\Big(\left(\omega_{3}^{-1} M\right)^{\frac{1}{3(\g_2-1)}}+\Big(\int_{a}^{b(t)}\r e(\r)\,r^{2}\mathrm{d}r\Big)^{\frac{1}{3(\g_2-1)}}\Big)\nonumber\\
		&\geq \frac{1}{2} \int_{a}^{b(t)} \rho e(\rho)\, r^{2}\mathrm{d} r-C(M),\label{2.21-9}
	\end{align}
	which, with \eqref{2.21-1}, yields \eqref{2.16-3}.}

\smallskip
{When $\gamma_{2}=\frac{4}{3}$, {\it i.e.}, $\frac{1}{3(\gamma_2-1)}=1$. It has been proved in \cite[Theorem 3.1]{Cheng-Cheng-Lin}
	that there exists an optimal constant $C_{\min}=6\kappa_2M_{\rm ch}^{-\frac{3}{2}}$ such that
	\begin{align}\label{pe1}
		\int_{a}^{b(t)}\Big(\int_{a}^{r}\, \rho\, z^{2}\mathrm{d} z\Big) \rho\, r \mathrm{d} r=\frac{1}{2 \omega_{3}}\|\nabla \Phi\|_{L^{2}(\mathbb{R}^{3})}^{2}
		&\leq  \frac{C_{\min}}{2\omega_{3}}\|\rho\|_{L^{1}(\Omega_{t})}^{\frac{2}{3}}\|\rho\|_{L^{\frac{4}{3}}(\Omega_{t})}^{\frac{4}{3}}\nonumber\\
		&=\frac{C_{\min}}{2}M^{\frac{2}{3}}\int_{a}^{b(t)}\rho^{\frac{4}{3}}\,r^2{\rm d}r,
	\end{align}
	which, with \eqref{ae2}, yields
	\begin{align}
		&\int_{a}^{b(t)} \rho e(\rho)\, r^{2}\mathrm{d} r
		- \int_{a}^{b(t)}\Big(\int_{a}^{r} \rho\, z^{2}\mathrm{d} z\Big) \rho\, r\mathrm{d} r\nonumber\\
		&\geq \Big(3\kappa_2-\frac{C_{\min}}{2}M^{\frac{2}{3}}-C_{1}K^{-\min\{\gamma_2,\epsilon\}}\Big)\int_{a}^{b(t)}\rho^{\frac{4}{3}}\,r^2{\rm d}r-C(M,K).\label{2.21-10}
	\end{align}
	Since $M<M_{\rm ch}$, we can always choose $K>\rho^{*}$ large enough such that
	$$
	3\kappa_2-\frac{C_{\min}}{2}M^{\frac{2}{3}}-C_{1}K^{-\min\{\gamma_2,\epsilon\}}>0.
	$$
	Then one can deduce \eqref{2.16-2} for $\gamma_2=\frac{4}{3}$ from \eqref{2.21-1}, \eqref{2.21-10}, and the fact that $\rho^{\gamma_2}\geq C\rho e(\rho)$.}

\smallskip
{When $\g_{2}\in (\frac{6}{5},\frac{4}{3})$, we define
	\begin{equation*}
		F(s;\beta)=s-B_{\beta}M^{\frac{5\gamma_2-6}{3(\gamma_2-1)}}\left(\omega_{3}^{-1}\beta M+s\right)^{\frac{1}{3(\g_2-1)}}\qquad\,\, \text{for }s\geq 0\text{ and any fixed }\beta>0.
	\end{equation*}
	A direct calculation shows that
	\begin{equation*}
		\left\{\begin{aligned}
			&\frac{\mathrm{d}F(s;\beta)}{\mathrm{d}s}=1-\frac{1}{3(\gamma_2-1)} B_{\beta}M^{\frac{5 \gamma_2-6}{3(\gamma_2-1)}}
			\big(\omega_{3}^{-1}\beta M+s\big)^{\frac{4-3 \gamma_2}{3(\gamma_2-1)}},\\
			&\frac{\mathrm{d}^2F(s;\beta)}{\mathrm{d}s^2}=-\frac{4-3 \gamma_2}{9(\gamma_2-1)^{2}} B_{\beta}
			M^{\frac{5\gamma_2-6}{3(\gamma_2-1)}}
			\big(\omega_{3}^{-1}\beta M+s\big)^{\frac{7-6\gamma_2}{3(\gamma_2-1)}},
		\end{aligned}\right.
	\end{equation*}
	which yields that $\frac{\mathrm{d}^2F(s;\beta)}{\mathrm{d}s^2}<0$ for  $s>0$ since $\gamma_2<\frac{4}{3}$.
	Thus, $F(s;\beta)$ is concave with respect to $s>0$.
	We denote
	\begin{equation}\label{2.21-13}
		s_{*}(\beta)=\Big(\frac{ B_{\beta}}{3(\gamma_2-1)}\Big)^{-\frac{3(\gamma_2-1)}{4-3 \gamma_2}} M^{-\frac{5 \gamma_2-6}{4-3\gamma_2}}-\omega_{3}^{-1}\beta M,
	\end{equation}
	which is the critical point of $F(s)$ satisfying $\frac{\mathrm{d}F(s;\beta)}{\mathrm{d}s}(s_{*}(\beta))=0$.
	The maximum of $F(s;\beta)$ with respect to $s>0$ is
	\begin{equation}\label{2.21-14}
		F(s_{*}(\beta);\beta)=(4-3\gamma_2)\Big(\frac{B_{\beta}}{3(\gamma_2-1)}\Big)^{-\frac{3(\gamma_2-1)}{4-3 \gamma_2}} M^{-\frac{5 \gamma_2-6}{4-3 \gamma_2}}-\omega_{3}^{-1}\beta M.
\end{equation}}

{It follows from the definition of $M_{\rm c}^{\v,b}$ that, if $M<M_{\rm c}^{\v,b}$,
	there exists $\beta_{0}>0$ such that $M<M_{\rm c}^{\v,b}(\beta_0)$. Then, from \eqref{2.21-13}--\eqref{2.21-14}, we have
	\begin{align}
		& F(s_{*}(\beta_{0});\beta_{0})>\frac{E_{0}^{\v,b}}{\omega_{3}},\label{2.21-16}\\
		&s_{*}(\beta_{0})>\Big(\frac{ B_{ \beta_0}}{3(\gamma_2-1)}\Big)^{-\frac{3(\gamma_2-1)}{4-3 \gamma_2}}
		\big(M_{\rm c}^{\varepsilon, b}(\beta_{0})\big)^{-\frac{5 \gamma_2-6}{4-3 \gamma_2}}
		-\omega_{3}^{-1}\beta_{0}M_{\rm c}^{\v,b}(\beta_{0})\nonumber\\
		&\qquad\,\,\,=\frac{1}{4-3\gamma_2} \big(E_{0}^{\varepsilon, b}+\omega_{3}^{-1}\beta_{0}M_{\rm c}^{\v,b}(\beta_{0})\big)
		-\omega_{3}^{-1}\beta_{0}M_{\rm c}^{\v,b}(\beta_{0})> \frac{E_{0}^{\v,b}}{\omega_{3}},\label{2.21-17}
	\end{align}
	where we have used that $\frac{1}{4-3\g_2}>\frac{5}{2}>1$
	for $\g_2\in (\frac{6}{5},\frac{4}{3})$.
	Then, combining \eqref{2.21-1} and \eqref{2.21-6} with \eqref{2.21-16}--\eqref{2.21-17}, we obtain
	\begin{align}\label{2.21-18}
		\begin{aligned}
		&F(\int_{a}^{b(t)} \rho e(\rho)\, r^{2}\mathrm{d} r;\beta_{0}) \leq \frac{E_{0}^{\v,b}}{\omega_{3}}<F(s_{*}(\beta_0);\beta_{0}),\\
		&\int_{a}^{b} \big(\rho_{0}e(\rho_{0})\big)(r)\, r^{2}\mathrm{d} r \leq \frac{E_{0}^{\v,b}}{\omega_{3}}<s_{*}(\beta_0).
		\end{aligned}
	\end{align}
	Hence, due to the continuity of $ \int_{a}^{b(t)}\big(\rho e(\rho)\big)(t,r)\,r^{2}\mathrm{d}r$
	with respect to $t$, the strict inequality:
	\begin{equation}\label{2.21-15}
		\int_{a}^{b(t)}\big(\rho e(\rho)\big)(t,r)\,r^{2}\mathrm{d} r<s_{*}(\beta_0)
	\end{equation}
	must hold. Otherwise, there exists $t_0>0$ such that
	$\int_{a}^{b(t_0)}\big(\rho e(\rho)\big)(t_{0},r)\,r^{2}\mathrm{d}r=s_{*}(\beta_0)$, which yields
	$$
	F(\int_{a}^{b(t_{0})} (\rho e(\rho))(t_{0},r)\, r^{2}\mathrm{d} r;\beta_{0})=F(s_{*}(\beta_0);\beta_{0})>\frac{E_{0}^{\v,b}}{\omega_{3}}.
	$$
	This contradicts \eqref{2.21-18}.
	Thus, we prove \eqref{2.21-15} under condition \eqref{2.9-1}.}

\smallskip
Therefore, under condition \eqref{2.9-1}, it follows from \eqref{2.21-13} and \eqref{2.21-15} that
	\begin{align}\label{2.21-20}
		&F(\int_{a}^{b(t)}\rho e(\rho)\,r^{2}\mathrm{d}r;\beta_{0})\nonumber\\
    &\geq \int_{a}^{b(t)}\rho e(\rho)\,r^{2}\mathrm{d}r
    \nonumber\\		&\quad
		-B_{\beta_{0}}M^{\frac{5\g_2-6}{3(\g_2-1)}}
		\big(s_{*}(\beta_0)+\omega_{3}^{-1}\beta_{0}M\big)^{\frac{4-3\g_2}{3(\g_2-1)}}
		\Big(\int_{a}^{b(t)}\rho e(\rho)\,r^{2}\mathrm{d}r+\omega_{3}^{-1}\beta_{0}M\Big)
		\nonumber\\&=(4-3\g_2)\int_{a}^{b(t)}\rho e(\rho)\,r^{2}\mathrm{d}r-3(\g_2-1)\omega_{3}^{-1}\beta_{0}M^{\frac{5}{3}}.
	\end{align}
	Combining \eqref{2.21-1} and \eqref{2.21-10} with \eqref{2.21-20}, we conclude \eqref{2.16-2}.
	$\hfill\square$

\smallskip
\begin{corollary}\label{cor2.1}
	Under the assumptions of {\rm Lemma \ref{lem2.1}} and noting \eqref{A.9-1},
	\begin{equation*}
		\int_{a}^{b(t)}\rho^{\g_2}(t,r)\,r^{2}\mathrm{d}r
		\leq C
		\int_{a}^{b(t)}\big(\rho+\r e(\r)\big)(t,r)\,r^{2}\mathrm{d}r\leq C(M,E_0)
		\qquad \text{ for $t\geq 0$}.
	\end{equation*}	
\end{corollary}

\begin{corollary}\label{cor2.2}
	Under the assumptions of {\rm Lemma \ref{lem2.1}},
	it follows from \eqref{2.9-3}, \eqref{2.16-2}--\eqref{2.16-3}, and \eqref{2.21-4} that,
    for $t\ge 0$ and $r\ge 0$,
	\begin{align*}
		&\left|r^{2} \Phi_{r}(t, r)\right| \leq \frac{M}{\omega_{3}},\\
		&\int_{a}^{b(t)}\Big(\int_{a}^{r} \rho(t, y)\, y^{2}\mathrm{d} y\Big) \rho(t, r)\, r\mathrm{d} r
		+\|\Phi(t)\|_{L^{6}(\mathbb{R}^{3})}+\|\nabla \Phi(t)\|_{L^{2}(\mathbb{R}^{3})} \leq C(M, E_{0}).
\end{align*}
\end{corollary}

For later use, we analyze the boundary value of density $\rho$.
Using $\eqref{2.13}_1$ and \eqref{2.14}, we have
\begin{equation}\label{2.22}
	\rho_{\tau}(\tau,\frac{M}{\omega_3})=-\f{1}{\v}\,P(\tau,\frac{M}{\omega_3})\leq 0,
\end{equation}
which yields that $\rho(\tau, \frac{M}{\omega_3})\leq \rho_{0}(\frac{M}{\omega_3})$.
In the Eulerian coordinates, it is equivalent to
\begin{equation}\label{2.24}
	\rho(t,b(t))\leq \rho_0(b).
\end{equation}
Moreover, noting \eqref{2.9} and $b\geq (\rho_{*})^{-\g_1/3}$,
we see that $\rho(t,b(t))\leq \rho_0(b)\leq \rho_{*}$ for all $t\geq 0$.
From  $\eqref{A.1-1}_1$ and \eqref{2.22}, there exists a positive constant $\tilde{C}$
depending only on $(\g_1, \k_1)$ such that
	$\rho_{\tau}(\tau,\frac{M}{\omega_3})=-\frac{1}{\v}\, P(\tau,\frac{M}{\omega_3})
	\geq -\frac{\tilde{C}}{\v}\big(\rho(\tau,\frac{M}{\omega_3})\big)^{\g_1}$,
which implies
\begin{equation*}
	\rho(\tau,\frac{M}{\omega_3})\geq \rho_{0}(\frac{M}{\omega_3})
	\Big(1+\frac{\tilde{C}(\g_1-1)}{\v}\big(\rho_{0}(\frac{M}{\omega_3})\big)^{\g_1-1}\tau\Big)^{-\frac{1}{\g_1-1}}.
\end{equation*}
Therefore, in the Eulerian coordinates,
\begin{equation}\label{2.27}
	\rho(t,b(t))\geq \rho_{0}(b)\Big(1+\frac{\tilde{C}(\g_1-1)}{\v} (\rho_{0}(b))^{\g_1-1}t\Big)^{-\frac{1}{\g_1-1}}
	\qquad \mbox{for $t\geq 0$}.
\end{equation}

\begin{lemma}[BD-type entropy estimate]\label{lem2.2}
	Under the conditions of {\rm Lemma \ref{lem2.1}}, for any given $T>0$,
	\begin{align}
		&\v^2\int_{a}^{b(t)}\big|\left(\sqrt{\rho}\right)_{r}\big|^2\,r^{2}\mathrm{d}r
		+\v\int_{0}^{t}\int_{a}^{b(s)}\frac{P'(\rho)}{\rho}|\rho_{r}|^2\,r^{2}\mathrm{d}r\mathrm{d}s
		+\frac{1}{3}\,P(\rho(t,b(t)))\,b(t)^3
		\nonumber\\
		&\quad +\frac{1}{3\v}\int_{0}^t \big(P(\rho)P'(\rho)\big)(s,b(s))\,b(s)^3\,\mathrm{d}s
		\leq C(E_0,M,T)\qquad\,\, \mbox{for $t\in [0,T]$}.\label{2.28}
	\end{align}
\end{lemma}

\noindent{\bf Proof.} We divide the proof into three steps.

\smallskip
1. Using \eqref{1.15} and similar calculations as in the proof \cite[Lemma 3.3]{Chen-He-Wang-Yuan-2021}, we have
\begin{align}
	&\int_{a}^{b(t)}\Big(\frac{1}{2}\big(u+\v \frac{\rho_{r}}{\rho}\big)^2\rho+\rho e(\rho)\Big)
	\,r^2\mathrm{d}r
	-\int_{a}^{b(t)}\Big(\int_{a}^{r}\rho(t,y)\,y^{2}\mathrm{d}y\Big)\rho\, r\mathrm{d}r\nonumber\\
	&\quad +\v\int_{0}^{t}\int_{a}^{b(s)}\frac{P'(\rho)}{\rho}\rho_{r}^2\,r^{2}\mathrm{d}r\mathrm{d}s
	+\frac{1}{3}P(\rho(t,b(t)))\,b(t)^3+\frac{1}{3\v}\int_{0}^{t}\big(P(\rho)P'(\rho)\big)(s,b(s))\,b(s)^3\mathrm{d}s \nonumber
	\\&=\int_{a}^{b}\Big(\frac{1}{2}\big(u_0+\v \frac{\rho_{0,r}}{\rho}\big)^2+e(\rho_0)\Big)\rho_{0}
	\,r^{2}\mathrm{d}r
	-\int_{a}^{b}\Big(\int_{a}^{r}\rho_{0}(y)\,y^{2}\mathrm{d}y\Big)\rho_0(r)\,r\mathrm{d}r\nonumber\\
	&\quad+\frac{1}{3}P(\rho_0\left(b\right))b^3
	+\v\int_{0}^{t}\int_{a}^{b(s)}\rho^2\,r^{2}\mathrm{d}r\mathrm{d}s
	-\frac{M\v}{\omega_3}\int_{0}^{t}\rho(s,b(s))\,\mathrm{d}s,\nonumber
\end{align}	
which, with Lemma \ref{lem2.1}, yields
\begin{align}
	&\v^2\int_{a}^{b(t)}|\left(\sqrt{\rho}\right)_{r}|^2\,r^{2}\mathrm{d}r
	+\v\int_{0}^{t}\int_{a}^{b(s)}\frac{P'(\rho)}{\rho}|\rho_{r}|^2\,r^{2}\mathrm{d}r\mathrm{d}s\nonumber\\
	&\quad+\frac{1}{3}P(\rho(t,b(t)))\,b(t)^3+\frac{1}{3\v}\int_{0}^t\big(P(\rho)P'(\rho)\big)(s,b(s))\,b(s)^3\mathrm{d}s \nonumber	\\
	&\leq C(E_0,M)+\frac{1}{3}P(\rho_0(b))b^3
	+\v\int_{0}^{t}\int_{a}^{b(s)}\rho^2\,r^{2}\mathrm{d}r\mathrm{d}s
	-\frac{M\v}{\omega_3}\int_{0}^{t}\rho(s,b(s))\,\mathrm{d}s.\label{2.37}
\end{align}

2. For the second term on the RHS of \eqref{2.37}, it follows from \eqref{2.9} and $\eqref{A.1-1}_1$ that
\begin{equation}\label{2.37-1}
	\frac{1}{3}P(\rho_0(b))b^3\leq C.
\end{equation}
For the last term on the RHS of \eqref{2.37},  using \eqref{2.24}, we have
\begin{equation}\label{2.37-2}
	\Big|\frac{M\varepsilon }{\omega_{3}} \int_{0}^{t} \rho(s, b(s))\,\mathrm{d}s\Big| \leq C(M) \rho_{0}(b) T \leq C(M, T).
\end{equation}

3. To close the estimates, we need to control the third term on the RHS of \eqref{2.37}, that is,
$$
\v\int_{0}^{t}\int_{a}^{b(s)}\rho^2\,r^{2}\mathrm{d}r\mathrm{d}s=\frac{\v}{\omega_3}\int_{0}^{t}\|\rho(s,\cdot)\|_{L^2(\Omega_{s})}^2\,\mathrm{d}s.
$$
We divide the estimate of the above term into the following two cases:

\smallskip
Case 1. For $\g_2\geq 2$, it follows from Corollary \ref{cor2.1} that
\begin{equation}\label{2.37-3}
	\v\int_{0}^{t}\int_{a}^{b(s)}\rho^2\,r^{2}\mathrm{d}r\mathrm{d}s
	\leq \v\int_{0}^{t}\int_{a}^{b(s)}(\rho+\rho^{\g_2})\,r^{2}\mathrm{d}r\mathrm{d}s\leq C(E_0,M,T).
\end{equation}

Case 2. For $\dis \g_2\in (\frac{6}{5},2)$, then $3\g_2>2$. A direct calculation shows that
\begin{equation}\label{2.37-4}
	\v\int_{0}^{t}\int_{a}^{b(s)}\rho^2\mathbf{I}_{\{\rho\leq 2\rho^{*}\}}\,r^2\mathrm{d}r\mathrm{d}s
	\leq 2\v\rho^{*}\int_{0}^{t}\int_{a}^{b(s)}\rho\, r^{2}\mathrm{d}r\mathrm{d}s\leq C(M,\rho^{*}).
\end{equation}
Denote $\sqrt{F(\rho)}:=\int_{0}^{\rho}\sqrt{\frac{P'(s)}{s}}\,\mathrm{d}s$. Then it follows from  $\eqref{A.2-1}_2$  that
\begin{align*}
	\sqrt{F(\rho)}&\geq \big(1-2^{-\frac{\g_2}{2}}\big)
	\frac{2\sqrt{(1-\mathfrak{a}_0)\k_2\g_2}}{\g_2}\rho^{\frac{\g_2}{2}}
	:=C(\g_2)^{-\frac{\g_2}{2\overline{\vartheta}}}\rho^{\frac{\g_2}{2}}\qquad \text{for }\rho\in [2\rho^{*},\infty),
\end{align*}
which, with  Corollary \ref{cor2.1}, implies that, for
$\overline{\vartheta}=\frac{3(2-\g_2)}{4}$,
\begin{align}\label{2.37-5}
	\|\rho \mathbf{I}_{\{\rho\geq 2\rho^{*}\}}\|_{L^2(\Omega_{t})}
	&\leq \|\rho \mathbf{I}_{\{\rho\geq 2\rho^{*}\}}\|_{L^{3\g_2}(\Omega_{t})}^{\overline{\vartheta}}
	\|\rho \mathbf{I}_{\{\rho\geq 2\rho^{*}\}}\|_{L^{\g_2}(\Omega_{t})}^{1-\overline{\vartheta}}\nonumber\\
	&\leq C(\g_2)\|\sqrt{F(\rho)}\|_{L^{6}(\Omega_{t})}^{\frac{2\overline{\vartheta}}{\g_2}}\|\rho \|_{L^{\g_2}(\Omega_{t})}^{1-\overline{\vartheta}}.
\end{align}

For $B_{R}(\textbf{0})\subset \R^3$, the following Sobolev's inequality  holds:
\begin{equation}\label{2.37-6}
	\|f\|_{L^{6}(B_{R}(\mathbf{0}))}
	\leq C\big(\|\nabla f\|_{L^{2}(B_{R}(\mathbf{0}))}+ R^{-1}\|f\|_{L^{2}(B_{R}(\mathbf{0}))}\big).
\end{equation}
It follows from \eqref{2.11} and Corollary \ref{cor2.1} that
\begin{equation*}\label{2.37-7}
	\begin{aligned}
		\frac{M}{\omega_{3}} &=\int_{a}^{b(t)} \rho(t, r)\, r^{2}\mathrm{d} r
		\leq\Big(\int_{a}^{b(t)} \rho^{\gamma_2}\,r^{2} \mathrm{d} r\Big)^{\frac{1}{\gamma_2}}
		\Big(\int_{a}^{b(t)} r^{2} \mathrm{d} r\Big)^{1-\frac{1}{\gamma_2}}
		\nonumber\\
		&\leq C
		b(t)^{\frac{3(\gamma_2-1)}{\gamma_2}}
		\Big(\int_{a}^{b(t)} \rho^{\gamma_2}\, r^{2} \mathrm{d} r\Big)^{\frac{1}{\gamma_2}},
	\end{aligned}
\end{equation*}
which yields
\begin{equation}\label{2.37-8}
	b(t)^{-1}
	\leq C M^{-\frac{\gamma_2}{3(\gamma_2-1)}}\Big(\int_{a}^{b(t)} \rho^{\gamma_2}\, r^{2}\mathrm{d} r\Big)^{\frac{1}{3(\gamma_2-1)}}
	\leq C\left(M, E_{0}\right).
\end{equation}

Using $\eqref{A.1-1}_2$--$\eqref{A.2-1}_2$ leads to $F(\rho)\leq C(\rho+\rho^{\g_2})$,
which, with \eqref{2.37-6}--\eqref{2.37-8} and Corollary \ref{cor2.1},
implies
\begin{align}
	\big\|\sqrt{F(\rho)}\big\|_{L^{6}(\Omega_{t})}
	&\leq C\left(\big\|\nabla (\sqrt{F(\rho)})\big\|_{L^2(\Omega_{t})}
	+b(t)^{-1}\big\|\sqrt{F(\rho)}\big\|_{L^2(\Omega_{t})}\right)\nonumber
	\\&\leq C\Big(\int_{a}^{b(t)}\frac{P'(\rho)}{\rho}|\rho_{r}|^2\,r^{2}\mathrm{d}r\Big)^{\frac{1}{2}}
	+C(M,E_0)\Big(\int_a^{b(t)}F(\rho)\,r^{2}\mathrm{d}r\mathrm{d}t\Big)^{\frac{1}{2}}\nonumber\\
	&\leq C(M,E_0)\Big(1+\big(\int_{a}^{b(t)}\frac{P'(\rho)}{\rho}|\rho_{r}|^2\,r^{2}\mathrm{d}r\big)^{\frac{1}{2}}\Big).\label{2.37-9}
\end{align}
Substituting \eqref{2.37-9} into \eqref{2.37-5}, we obtain
\begin{align}
	\v\int_{0}^{t}\int_{a}^{b(s)}\rho^2r^{2}\mathbf{I}_{\{\rho\geq 2\rho^{*}\}}\,\mathrm{d}r\mathrm{d}s
	&\leq C(M,E_0,T)\v\Big(1+\big(\int_{0}^{t}\int_{a}^{b(s)}\frac{P'(\rho)}{\rho}|\rho_{r}|^2\,r^{2}
	\mathrm{d}r\mathrm{d}s\big)^{\frac{2\overline{\vartheta}}{\g_2}}\Big)\nonumber
	\\&\leq C(M,E_0,T)+\frac{\v}{2}\int_{0}^{t}\int_{a}^{b(s)}\frac{P'(\rho)}{\rho}|\rho_{r}|^2\,r^{2}\mathrm{d}r\mathrm{d}s,\label{2.37-10}
\end{align}
where we have used
$\frac{2\bar{\vartheta}}{\g_2}\in (0,1)$ for $\g_2>\frac{6}{5}$.
Finally, substituting \eqref{2.37-1}--\eqref{2.37-4} and \eqref{2.37-10} into \eqref{2.37},
we conclude \eqref{2.28}.
$\hfill\Box$

\smallskip
In order to take the limit: $b\to \infty$, we need to make sure that domain $\Omega_{T}$
can be expanded to $[0,T]\times \R_{+}$ for fixed $\v>0$: $\lim\limits_{b\to \infty}b(t)=\infty$.

\begin{lemma}[Expanding of domain $\Omega_{T}$]\label{lem2.3}
	Given $T>0$ and $\v\in (0,\v_0]$, there exists
	$C_1(M,E_0,T,\v)>0$ such that, if $b\geq C_1(M,E_0,T,\v)$,
	\begin{equation}\label{2.38}
		b(t)\geq \frac{1}{2}b\qquad\,\, \text{for $t\in [0,T]$}.
	\end{equation}
\end{lemma}

\noindent{\bf Proof.}
Noting $b(0)=b$
and the continuity of $b(t)$,
we first make the {\it a priori} assumption:
\begin{equation}\label{2.39}
	b(t)\geq \frac{1}{2}b.
\end{equation}
Integrating \eqref{2.3} over $[0,t]$ yields
\begin{equation}\label{2.40}
	b(t)=b+\int_{0}^tu(s,b(s))\,\mathrm{d}s.
\end{equation}
It follows from \eqref{2.27}, \eqref{2.39}, and Lemma \ref{2.1}   that
\begin{align}
	\int_{0}^{t}|u(s,b(s))|\,\mathrm{d}s
	&\leq \frac{C}{\sqrt{\v}}\Big(\int_{0}^{t}\v(\rho u^2 r)(s,b(s))\,\mathrm{d}s\Big)^{\frac{1}{2}}
	\Big(\int_{0}^{t}\frac{1}{\rho(s,b(s))b(s)}\,\mathrm{d}s\Big)^{\frac{1}{2}}\nonumber\\
	&\leq C(M,E_0)\v^{-\frac{1}{2}}\Big(\int_{0}^{t}\frac{(1+\tilde{C}(\g_1-1)\v^{-1}\rho_{*}^{\g_1-1}s)^{\frac{1}{\g_1-1}}}{\rho_0(b)b}
	\,\mathrm{d}s\Big)^{\frac{1}{2}}\nonumber\\
	&\leq C(M,E_0,T,\rho_{*},\g_1,\g_2,\v)\rho_{0}(b)^{-\frac{1}{2}}b^{-\frac{1}{2}}.\label{2.41}
\end{align}
We take
$
C_{1}(M,E_0,T,\v):= \max\big\{\rho_{*}^{-\frac{\g_1}{3}}, (4C(M,E_0,T,\rho_{*},\g_1,\g_2,\v))^{\frac{2}{\alpha}}, \mathcal{B}(\v)\big\},
$
which, with \eqref{2.9} and \eqref{2.41}, implies that\\
\begin{equation}\label{2.42}
	\Big\vert\int_{0}^{t}u(s,b(s))\,\mathrm{d}s\Big\vert\leq \int_{0}^{t}|u(s,b(s))|\,\mathrm{d}s\leq  \frac{1}{4}b,
\end{equation}
provided that $b\geq C_1(M,E_0,T,\v)$. Combining \eqref{2.42} with \eqref{2.40}, we have
\begin{equation}\label{2.43}
	b(t)\geq \frac{3}{4}b.
\end{equation}
Thus, we have closed the {\it a priori} assumption \eqref{2.39}.
Finally, using \eqref{2.43} and the continuity argument, we can conclude \eqref{2.38}.
$\hfill\Box$

\subsection{\,Higher integrability of the density and the velocity}
As implied in \cite{Chen-Perepelitsa-2010}, the higher integrabilities of
the density and the velocity are important for the $L^p$ compensated compactness framework.
However, for the general pressure law,
due to the lack of an explicit formula for the entropy kernel,
for the special entropy pair $(\eta^{\psi},q^{\psi})$ by taking the test function
$\psi=\f12 s|s|$ in \eqref{2.21e}--\eqref{2.21q},
we can not obtain that $q^{\psi} \gtrsim \rho |u|^3 + \rho^{\gamma+\theta}$ in general.
To derive the higher integrability of the velocity,
we use the special entropy pair constructed in Lemma \ref{lem4.1},
at the cost of the higher integrability of the density over domain $[0,T]\times [d, b(t)]$
for some $d>0$. Since $b(t)\to \infty$ as $b\to \infty$,
we indeed need the higher integrability of the density on the unbounded domain.
We point out that this is different from the case of \cite{Chen-He-Wang-Yuan-2021}
in which only the higher integrability on the bounded domain $[0,T]\times [d,D]$ for any given $0<d<D<\infty$
is needed.

\smallskip
\begin{lemma}[Higher integrability on the density]\label{lem2.4}
	Let $(\rho,u)$ be a smooth solution of \eqref{2.1}--\eqref{2.6}.
	Then, under the assumption of {\rm Lemma \ref{lem2.1}}, for any given $d>2b^{-1}>0$,
	\begin{equation}\label{2.44}
		\int_{0}^{T}\int_{d}^{b(t)}\rho P(\rho)\,r^{2}\mathrm{d}r\mathrm{d}t\leq C(d,M,E_0,T).
	\end{equation}
\end{lemma}

\noindent{\bf Proof.}
Let $\omega(r)$ be a smooth function with $\operatorname{supp}\omega\subset (\frac{d}{2},\infty)$
and $\omega(r)=1$ for $r\in [d,\infty)$. Multiplying $\eqref{2.1}_2$ by $w(y)y^{2}$, we have
\begin{align}
	&(y^{2}\rho u \omega)_{t}+(y^{2}\rho u^2 \omega)_{y}+(y^{2}P(\rho)\omega)_{y}
	-\omega_{y}\big(y^{2}\rho u^2+y^2P(\rho)\big)
	+\r\omega\int_{a}^{r}\r\, z^{2}\mathrm{d}z\nonumber\\
	&=2yP(\rho)\omega +\v(y^{2}\rho u_{y}\omega)_{y}-\v\omega_{y}y^{2}\rho u_{y}-2\v\rho u\,\omega.\label{2.45}
\end{align}	
Integrating \eqref{2.45} with respect to $y$ from $\frac{d}{2}$ to $r$ and then multiplying
the equation
by $\rho(t,r)$ yield
\begin{align}
	&r^{2}\rho(t,r)P(\rho(t,r))\omega(r)\nonumber\\
	&=-\rho\frac{\mathrm{d}}{\mathrm{d}t}\int_{\frac{d}{2}}^{r}\rho u\,\omega\, y^2\mathrm{d}y
	-r^{2}\rho^2 u^2\omega(r)+\rho\int_{\frac{d}{2}}^{r}\omega_{y}\rho u\,y^2\mathrm{d}y\nonumber\\
	&\quad
	+\rho\int_{\frac{d}{2}}^{r}\omega_{y}P(\rho)\,y^2\mathrm{d}y+2\rho\int_{\frac{d}{2}}^{r}P(\rho)\,\omega\,y\mathrm{d}y
    -\rho\int_{\frac{d}{2}}^{r}\rho\omega\Big(\int_{a}^{y}\rho\, z^{2}\mathrm{d}z\Big)\,\mathrm{d}r
	\nonumber\\
	&\quad+\v r^{2}\rho^2 u_{r}\omega(r)-\v\rho\int_{\frac{d}{2}}^{r}\omega_y\rho u_{y}\,y^2\mathrm{d}y
	-2\v\rho\int_{\frac{d}{2}}^{r}\rho u\,\omega\, \mathrm{d}y.\label{2.46}
\end{align}
Using $\eqref{2.1}_1$, we have
\begin{align*}
		\rho\frac{\mathrm{d}}{\mathrm{d}t}\int_{\frac{d}{2}}^{r}\rho u\omega\,y^2\mathrm{d}y
		&=\Big(\rho\int_{\frac{d}{2}}^{r}\rho u\omega\,y^2\mathrm{d}y\Big)_{t}
		+\Big(\rho u\int_{\frac{d}{2}}^{r}\rho u\omega \,y^2\mathrm{d}y\Big)_{r}
		\nonumber\\
		&\quad -\rho^2u^2\omega(r)r^2+\frac{2}{r}\rho u\int_{\frac{d}{2}}^{r}\rho u\omega\,y^2\mathrm{d}y,
\end{align*}
which, with \eqref{2.46}, yields that
\begin{align}
	&r^{2}\rho(t,r)P(\rho(t,r))\omega(r)\nonumber\\
	&=-\Big(\rho\int_{\frac{d}{2}}^{r}\rho u\omega\,y^2\mathrm{d}y\Big)_{t}
	-\Big(\rho u\int_{\frac{d}{2}}^{r}\rho u\omega\,y^2\mathrm{d}y\Big)_{r}
	-\frac{2}{r}\rho u\int_{\frac{d}{2}}^{r}\rho u\omega\,y^2\mathrm{d}y
	+\rho\int_{\frac{d}{2}}^{r}\omega_{y}\rho u^2\,y^2\mathrm{d}y\nonumber\\
	&\quad +\rho\int_{\frac{d}{2}}^{r}\omega_{y}P(\rho)\,y^2\mathrm{d}y
	+2\rho\int_{\frac{d}{2}}^{r}P(\rho)\omega\,y\mathrm{d}y
	+\v \rho^2 u_{r}\omega(r)r^2-\v\rho\int_{\frac{d}{2}}^{r}\rho u_{y}\omega_y\,y^{2}\mathrm{d}y\nonumber\\
	&\quad -2\v\rho\int_{\frac{d}{2}}^{r}\rho u\omega\,\mathrm{d}y-\rho\int_{\frac{d}{2}}^{r}\rho\omega\,
	\Big(\int_{a}^{y}\rho\, z^{2}\mathrm{d}z\Big)\,\mathrm{d}y.\label{2.48}
\end{align}
Multiplying \eqref{2.48} by $\omega(r)$ leads to
\begin{align}
	&r^{2}\rho(t,r)P(\rho(t,r))\omega^2(r)\nonumber\\
	&=-\Big(\rho\omega(r)\int_{\frac{d}{2}}^{r}\rho u\omega(y)\,y^2\mathrm{d}y\Big)_{t}
	-\Big(\rho u\omega(r)\int_{\frac{d}{2}}^{r}\rho u\omega(y)\,y^2\mathrm{d}y\Big)_{r}
	+\omega_{r}\rho u\int_{\frac{d}{2}}^{r}\rho u\omega\,y^2\mathrm{d}y
	\nonumber\\	
	&\quad -\frac{2}{r}\rho u\omega(r)\int_{\frac{d}{2}}^{r}\rho u\omega\,y^2\mathrm{d}y
	+\rho\omega(r)\int_{\frac{d}{2}}^{r}\rho u\omega_{y}\,y^2\mathrm{d}y
	+\rho\omega(r)\int_{\frac{d}{2}}^{r}P(\rho)\,\omega_{y}\,y^2\mathrm{d}y
	\nonumber\\
	&\quad +2\rho\omega(r)\int_{\frac{d}{2}}^{r}P(\rho)\omega\,y\mathrm{d}y
	-\rho\omega(r)\int_{\frac{d}{2}}^{r}\rho\omega\Big(\int_{a}^{y}\rho\, z^{2}\mathrm{d}z\Big)\,\mathrm{d}y
	-\v\rho\omega(r)\int_{\frac{d}{2}}^{r}\rho u_{y}\,\omega_{y}\,y^2\mathrm{d}y
	\nonumber\\		&\quad
	-2\v\rho\omega(r)\int_{\frac{d}{2}}^{r}\rho u\,\omega\,\mathrm{d}y
	+\v r^{2}\rho^2 u_{r}\omega^2(r)
	:= \sum\limits_{i=1}^{11}I_{i}.\label{2.49}
\end{align}

Using Lemma \ref{lem2.1} and \eqref{2.11}, we have
$$
\begin{aligned}
	&\Big|\int_{\frac{d}{2}}^{r}\big((\rho u+P(\rho))\omega(y)+\v\rho u_{y}\omega_y\big)\,y^2\mathrm{d}y\Big|\\
	&\leq C\int_{a}^{b(t)}\big(\rho u^2+\rho+\rho^{\gamma_2}\big)\omega(y)\,y^2\mathrm{d}y
	+\v\int_{a}^{b(t)}\rho (u_{y}^2+1)|\omega_y|\,y^{2}\mathrm{d}y
	\leq C(M,E_0,\|\omega\|_{C^1}),
\end{aligned}
$$
which yields
\begin{align}
	&\Big|\int_{0}^{T}\int_{\frac{d}{2}}^{b(t)}I_{i}\,\mathrm{d}r\mathrm{d}t\Big|\leq C(M,E_0,T,\|\omega\|_{C^1})(d^{-2}+d^{-4})
	\quad\,\, \text{for $i=3,4,\cdots,10$}, \label{2.54}\\
	&\Big|\int_{0}^{T}\int_{\frac{d}{2}}^{b(t)}I_1\,\mathrm{d}r\mathrm{d}t\Big|
	=\Big| \int_{\frac{d}{2}}^{b(t)}\rho(T,r)\omega(r)
	\Big(\int_{\frac{d}{2}}^{r}y^{2}\rho(T,y)u(T,y)\omega(y)\,\mathrm{d}y\Big)\mathrm{d}r\Big|\nonumber\\
	&\qquad\qquad\qquad\qquad\qquad +\Big|\int_{\frac{d}{2}}^{b(t)}\rho(0,r)\omega(r)
	\Big(\int_{\frac{d}{2}}^{r}y^{2}\rho(0,y)u(0,y)\omega(y)\,\mathrm{d}y\Big)\mathrm{d}r\Big|\nonumber\\
	&\qquad\qquad\qquad\qquad\quad
	\leq C(M,E_0,\|\omega\|_{C^1})d^{-2}. \label{2.53}
\end{align}

For $I_2$, using  \eqref{2.9}, \eqref{2.24}, \eqref{2.42}, and $b\gg 1$, we have
\begin{align}\label{2.56-2}
	\Big|\int_{0}^{T}\int_{\frac{d}{2}}^{b(t)}I_{2}\,\mathrm{d}r\mathrm{d}t\Big|
	&=\Big|\int_{0}^{T}\int_{\frac{d}{2}}^{b(t)}\Big(\rho u\omega(r)\int_{\frac{d}{2}}^{r}\rho u\omega(y)\,y^2\mathrm{d}y\Big)_{r}\,\mathrm{d}r\mathrm{d}t\Big|\nonumber\\
	&\leq \Big|\int_{0}^{T}(\rho u)(t,b(t))\Big(\int_{\frac{d}{2}}^{b(t)}\rho u\omega(y)\,y^2\mathrm{d}y\Big)\,\mathrm{d}t\Big|\nonumber\\
	&\leq C(E_0,M)b^{-3+\frac{1}{2}}b\leq C(E_0,M).
\end{align}

For $I_{11}$, we obtain
\begin{align}\label{2.57}
	\Big| \int_{0}^{T}\int_{\frac{d}{2}}^{b(t)}I_{11}\,\mathrm{d}r\mathrm{d}t\Big|
	&=\v\Big\vert \int_{0}^{T}\int_{\frac{d}{2}}^{b(t)}\rho^2 u_{r}\omega^2\,r^2\mathrm{d}r\mathrm{d}t\Big\vert
	\nonumber\\
	&\leq \v\int_{0}^{T}\int_{\frac{d}{2}}^{b(t)}\rho^3\omega^2\,r^{2}\mathrm{d}r\mathrm{d}t+C(M,E_0,T,\|\omega\|_{C^1}).
\end{align}	
We divide the estimate of
$\int_{0}^{T}\int_{\frac{d}{2}}^{b(t)}\v\rho^3\omega^2\,r^{2}\mathrm{d}r\mathrm{d}t$ into two cases:

\smallskip
{\it Case 1. $\g_2\in (\f65,2)$}:  For
$t\in [0,T]$,
denoting
$
A(t):=\{r\in [\frac{d}{2},b(t)]\,:\, \rho(t,r)\geq \rho^{*}\},
$
then it follows from \eqref{2.11} that $|A(t)|\leq C(d,\rho^{*})M$.
For any $r\in A(t)$, let $r_{0}$ be the closest point to $r$ so that
$\rho(t,r_{0})=\rho^{*}$ with $|r-r_{0}|\leq |A(t)|\leq C(d,\rho^{*})M$.
Then, for any smooth function $f(\rho)$,
\begin{align}
	\sup_{r\in A(t)}\big(f(\rho(t,r))\omega^{2}(r)\big)&\leq f(\rho(t,r_0))\omega^2(r_0)
	+\Big|\int_{r_0}^{r}\partial_y\big(f(\rho(t,y))\omega^{2}(y)\big)\,\mathrm{d}y\Big|\nonumber\\
	&\leq C(\|\omega\|_{C^1})|f(\rho^{*})|
	+\int_{A(t)}\big|\partial_y\big(f(\rho(t,y))\omega^{2}(y)\big)\big|\,\mathrm{d}y.\nonumber
\end{align}
Recalling  \eqref{A.2-1} and \eqref{A.9-1},
we notice that $P(\rho)\cong \r^{\g_2}$ and $e(\rho)\cong \r^{\g_2-1}$ for any $r\in A(t)$.
Then
\begin{align}\label{2.58}
	&\v\int_{0}^{T}\int_{\frac{d}{2}}^{b(t)}\rho^3\omega^2\,r^{2}\mathrm{d}r\mathrm{d}t\nonumber\\
	&= \v\int_{0}^{T}\int_{\frac{d}{2}}^{b(t)}\rho^3{\bf I}_{\{\rho\leq \rho^{\ast}\}}\omega^2\,r^2\mathrm{d}r\mathrm{d}t
	+ \v\int_{0}^{T}\int_{\frac{d}{2}}^{b(t)}\rho^3
	{\bf I}_{\{\rho\geq \rho^{\ast}\}}\omega^2\,r^{2}\mathrm{d}r\mathrm{d}t\nonumber\\
	&\leq C(M,E_0,\rho^{*},T)+ C(M,E_0)\, \v\int_{0}^{T}\Big(\int_{\f{d}{2}}^{b(t)}\rho e(\rho)\, r^2 dr\Big)
	\sup_{r\in A(t)}  \Big(\frac{\rho^2}{e(\rho)}\omega^2\Big)\mathrm{d}t\nonumber\\
	&\leq C(M,E_0,\rho^{*},T)+ C(M,E_0)\, \v\int_{0}^{T}\int_{A(t)}\Big\vert\Big(\frac{\rho^2}{e(\rho)}\omega^2\Big)_{r}(t,r)\Big\vert
	\,\mathrm{d}r\mathrm{d}t\nonumber\\
	&\leq  C(M,E_0)\,\v\int_{0}^{T}\int_{A(t)}\Big(\big(\frac{2\rho}{e(\rho)}-\frac{P(\rho)}{e(\rho)^2}\big)
	|\rho_{r}|\omega^2+\frac{\rho^2}{e(\rho)}\omega|\omega_{r}|\Big)\,\mathrm{d}r\mathrm{d}t\nonumber\\
	&\quad+C(M,E_0,\rho^{*},T).
\end{align}
A direct calculation shows that
\begin{align}
	&\int_{0}^{T}\int_{A(t)}\v\Big(\frac{2\rho}{e(\rho)}-\frac{P(\rho)}{e(\rho)^2}\Big)
	|\rho_{r}|\,\omega^2\,\mathrm{d}r\mathrm{d}t\nonumber\\
	&\leq \int_{0}^{T}\int_{A(t)}\v\frac{P'(\rho)}{\rho}|\rho_{r}|^2\omega^2\,r^{2}\mathrm{d}r\mathrm{d}t
     +\int_{0}^{T}\int_{A(t)}\v\Big(\frac{2\rho}{e(\rho)}-\frac{P(\rho)}{e(\rho)^2}\Big)^2\frac{\rho}{P'(\rho)}
      omega^2\,r^{-2}\,\mathrm{d}r\mathrm{d}t
	\nonumber\\
	&\leq C(M,E_0,T)+\int_{0}^{T}\int_{A(t)}\v\rho^{6-3\g_2}\omega^2\,r^{-2}\,\mathrm{d}r\mathrm{d}t\nonumber\\
	&\leq C(M,E_0,T)+\v C(M,E_{0})^{-1}\int_{0}^{T}\int_{A(t)}\rho^3\omega^2\,r^{2}\,\mathrm{d}r\mathrm{d}t\nonumber\\
	&\quad
	+C(M,E_0)\v\int_{0}^{T}\int_{A(t)}(r^{2})^{-\frac{3-\g_2}{\g_2-1}}\omega^2\,\mathrm{d}r\mathrm{d}t\nonumber\\
	&\leq C(M,E_0,T)+\v C(M,E_0)^{-1}\int_{0}^{T}\int_{\frac{d}{2}}^{b(t)}\rho^3\omega^2\,r^{2}\mathrm{d}r\mathrm{d}t.\label{2.59}\\
	&\int_{0}^{T}\int_{A(t)}\v\frac{\rho^2}{e(\rho)}\omega|\omega_{r}|\,\mathrm{d}r\mathrm{d}t
	\leq \int_{0}^{T}\Big(\v\sup_{r\in  A(t)}(\rho\omega)(t,r)\int_{A(t)}\frac{\rho}{e(\rho)}
            |\omega_{r}|\,\mathrm{d}r\Big)\,\mathrm{d}t\nonumber\\
	&\leq C(\rho^{*})\int_{0}^{T}\Big(\v\sup_{r\in A_{t}}(\rho\omega)(t,r)\int_{A(t)}\rho|\omega_{r}|\,\mathrm{d}r\Big)\,\mathrm{d}t\nonumber\\
	&\leq C(\rho^{*},M,\|\omega\|_{C^1})d^{-2}\int_{0}^{T}\,\v\sup_{r\in A(t)}(\rho\omega)(t,r)\,\mathrm{d}t\nonumber\\
	&\leq C(\rho^{*},M,\|\omega\|_{C^1},T)d^{-2}+C(\rho^{*},M,\|\omega\|_{C^1})d^{-2}
	\int_{0}^{T}\int_{A(t)}\v\big(|\rho_{r}|\omega+\rho|\omega_{r}|\big)\,\mathrm{d}r\mathrm{d}t\nonumber\\
	&\leq  C(\rho^{*},M,\|\omega\|_{C^1},T)d^{-2}
	\Big(1+
	\int_{0}^{T}\int_{A(t)}\v\Big(\frac{P'(\rho)}{\rho}|\rho_{r}|^2\omega
	+\rho|\omega_{r}|+\rho^{2-\g_2}\omega\Big)\,\mathrm{d}r\mathrm{d}t\Big)
	\nonumber\\
	&\leq C(\rho^{*},M,\|\omega\|_{C^1},T)d^{-2}+C(\rho^{*},M,E_0,T,\|\omega\|_{C^1})d^{-4}.\label{2.60}
\end{align}
Combining  \eqref{2.58}--\eqref{2.60}, we obtain that, for $\g_2\in (\f65,2)$,
\begin{equation}\label{2.61}
	\v\int_{0}^{T}\int_{\frac{d}{2}}^{b(t)}\rho^3\omega^2\,r^{2}\mathrm{d}r\mathrm{d}t
	\leq C(M,E_0,\rho^{*},T,\|\omega\|_{C^1}) (1+d^{-4}).
\end{equation}

\medskip
{\it Case 2. $\g_2\in [2,3)$}:  Using \eqref{2.11} and the same argument as for \eqref{2.60}, we have
\begin{align}\label{2.62}
	\v\int_{0}^{T}\int_{\frac{d}{2}}^{b(t)}\rho^3\omega^2\,r^{2}\mathrm{d}r\mathrm{d}t
	&\leq C(M)\int_{0}^{T}\v\sup_{r\in [\frac{d}{2},b(t)]}(\rho^2\omega)(t,r)\,\mathrm{d}t\nonumber\\[-1mm]
	&\leq C(M,\rho^{*},\|\omega\|_{C^1},T)+C(M)\int_{0}^{T}\v\sup_{r\in A(t)}(\rho^2\omega)(t,r)\,\mathrm{d}t\nonumber\\
	&\leq C(M,\rho^{*},\|\omega\|_{C^1},T)+C(M,E_0,\rho^{*},\|\omega\|_{C^1},T)d^{-2}.
\end{align}

Finally, integrating \eqref{2.49} over $[0,T]\times [\frac{d}{2},b(t)]$ and using \eqref{2.54}--\eqref{2.57}
and \eqref{2.61}--\eqref{2.62}, we conclude \eqref{2.44}.
$\hfill\Box$

\smallskip
\begin{corollary}\label{cor2.6}
	It follows from \eqref{A.2-1} and {\rm Lemma \ref{lem2.4}} that
	\begin{align}\label{2.64-1}
		\int_{0}^{T}\int_{d}^{b(t)}\rho^{\g_2+1}(t,r)\,r^{2}\mathrm{d}r\mathrm{d}t
		&\leq C\int_{0}^{T}\int_{d}^{b(t)}\big(\rho+\rho P(\rho)\big)(t,r)\,r^{2}\mathrm{d}r\mathrm{d}t\nonumber\\
		&\leq C(d, M, E_0,T).
	\end{align}
\end{corollary}

In order to use the $L^p$ compensated compactness framework,
we still need to obtain the higher integrability of the velocity (see \cite{Chen-Perepelitsa-2010}).
With the help of Lemma \ref{lem2.4}, we use the special entropy pair constructed
in Lemma \ref{lem4.1} to achieve this.

\smallskip
\begin{lemma}[Higher integrability of the velocity]\label{lem4.3}
	Let $(\rho,u)$ be the smooth solution of \eqref{2.1}--\eqref{2.6}.
	Then, under the assumption of {\rm Lemma \ref{lem2.1}},
	\begin{equation*}
		\int_{0}^{T}\int_{d}^{D}(\rho |u|^3)(t,r)\,r^{2}\mathrm{d}r\mathrm{d}t\leq C(d, D, \rho^{*}, M, E_0, T)
		\qquad\,\mbox{for any $(d, D)\Subset [a,b(t)]$}.
	\end{equation*}
\end{lemma}

\noindent{\bf Proof.}
Considering $\eqref{2.1}_1\times \hat{\eta}_{\rho}r^{2}+\eqref{2.1}_2\times\hat{\eta}_{m}r^{2}$, we can  obtain
\begin{align}
	&(\hat{\eta}r^{2})_{t}+(\hat{q}r^{2})_{r}+2r\big(-\hat{q}+\rho u\hat{\eta}_{\r}+\rho u^2\hat{\eta}_m\big)
	\nonumber\\
	&=\v\,r^{2}\Big((\rho u_{r})_{r}+2\rho\big(\frac{u}{r}\big)_{r}\Big)\hat{\eta}_{m}
	-\r\int_{a}^{r}\rho\, y^{2}\mathrm{d}y\,\hat{\eta}_{m}.\label{4.29}
\end{align}
Using \eqref{2.3}, a direct calculation yields
\begin{equation}\label{4.30}
	\frac{\mathrm{d}}{\mathrm{d}t}\int_{r}^{b(t)}\hat{\eta}\,y^{2}\mathrm{d}y
	=(u\hat{\eta})(t,b(t))\,b(t)^{2}+\int_{r}^{b(t)}\partial_t\hat{\eta}(t,y)\,y^{2}\mathrm{d}y.
\end{equation}
Integrating \eqref{4.29} over $[r,b(t))$ and using \eqref{4.30}, we have
\begin{align}
	\hat{q}(t,r)\,r^{2}
	&=-\v\int_{r}^{b(t)}\hat{\eta}_{m}(t,y)(\rho u_{y}y^{2})_{y}\,\mathrm{d}y
	+2\v\int_{r}^{b(t)}\hat{\eta}_{m}(t,y)\,\rho u\,\mathrm{d}y\nonumber\\
	&\quad\, +\Big(\int_{r}^{b(t)}\hat{\eta}(t,y)\,y^{2}\mathrm{d}y\Big)_{t}
	+(\hat{q}-u\hat{\eta})(t,b(t))\,b(t)^{2}
	\nonumber\\
	&\quad\, +2\int_{r}^{b(t)}\big(-\hat{q}+\rho u\hat{\eta}_{\rho}+\rho u^2\hat{\eta}_{m}\big)\,y\mathrm{d}y
	+\int_{r}^{b(t)}\Big(\int_{a}^{y}\rho\, z^2 {\rm d}z\Big)\rho\, \hat{\eta}_{m}\,\mathrm{d}y.\label{4.31}
\end{align}

We now control the terms on the RHS of \eqref{4.31}.
For the third term on the RHS of \eqref{4.31},  it follows from $\eqref{2.9b}$, \eqref{2.24}, \eqref{2.37-8},
and Lemmas \ref{lem4.1}, \ref{lem2.1}, and \ref{lem2.2}--\ref{lem2.3} that
\begin{align}\label{4.32}
	&\int_{0}^{T}\big|(\hat{q}-u\hat{\eta})(t,b(t))\big|b(t)^{2}\,\mathrm{d}t\nonumber\\
	&\leq C\int_{0}^{T}\left((\rho(t,b(t)))^{\g_1+\t_1}+(\rho^{\g_1}|u|)(t,b(t))\right) b(t)^{2}\,\mathrm{d}t \nonumber\\
	&\leq C\Big(\int_{0}^{T}\v(\rho|u|^2)(t,b(t))b(t)\,\mathrm{d}t\Big)^{\frac{1}{2}}
	\Big(\int_{0}^{T}\frac{1}{\v}(\rho(t,b(t)))^{2\g_1-1}b(t)^3\,\mathrm{d}t\Big)^{\frac{1}{2}} \nonumber\\
	&\quad + C(M,E_0,T)\int_{0}^{T} (\rho(t,b(t)))^{\t_1} b(t)^{-1}\,\mathrm{d}t
	\leq C(M,E_0,T).
\end{align}
For the first term on the RHS of \eqref{4.31}, integrating by parts yields
\begin{align}\label{4.33}
	\v\int_{r}^{b(t)}\hat{\eta}_{m}(t,y)\,(\rho u_{y}y^{2})_{y}\,\mathrm{d}y
	&=\v\hat{\eta}_{m}(t,b(t))\,(\rho u_{r})(t,b(t))\,b(t)^{2}-\v\hat{\eta}_{m}(t,r)\,(\rho u_{r})(t,r)\,r^{2}\nonumber\\
	&\quad -\v\int_{r}^{b(t)}\rho u_{y}\big(\hat{\eta}_{mu}u_{y}+\hat{\eta}_{m\rho}\rho_{y}\big)\,y^2\mathrm{d}y.
\end{align}
It follows from \eqref{2.4} and Lemma \ref{lem4.1} that
\begin{align*}\label{4.34}
	&\vert\v\hat{\eta}_{m}(t,b(t))\,(\rho u_{r})(t,b(t))\,b(t)^{2}\vert\nonumber\\
	&=\Big\vert\hat{\eta}_{m}(t,b(t))\Big(\v\rho\big(u_{r}+\frac{2}{r}u\big)(t,b(t))
	-2\v b(t)^{-1}\,(\rho u)(t,b(t))\Big)b(t)^{2}\Big\vert\nonumber\\
	&=\Big\vert\hat{\eta}_{m}(t,b(t))\, p(\rho)(t,b(t))\, b(t)^{2}-2\v\hat{\eta}_{m}(t,b(t))\,(\rho u)(t,b(t))\,b(t)\Big\vert\nonumber\\
	&\leq C\big((\rho^{\g_1}|u|)(t,b(t))+(\rho(t,b(t)))^{\g_1+\t_1}\big)b(t)^{2}
	+C\v\big((\rho |u|^2)(t,b(t))+(\rho(t,b(t)))^{\g_1}\big)b(t).
\end{align*}
which,
with similar arguments as in \eqref{4.32}, yields
\begin{equation}\label{4.35}
	\int_{0}^{T}\left\vert\v \hat{\eta}_{m}(t,b(t))\, (\rho u_{r})(t,b(t))\, b^{2}(t)\right\vert \mathrm{d}t\leq C(M, E_0,T).
\end{equation}
Hence, using  \eqref{2.28}, \eqref{4.33}--\eqref{4.35}, and Lemma \ref{lem4.1}, we have
\begin{align}
	&\int_{0}^{T}\int_{d}^{D}\Big\vert\v\int_{r}^{b(t)}\hat{\eta}_m(\rho u_yy^{2})_{y}\,\mathrm{d}y\Big\vert\, \mathrm{d}r\mathrm{d}t\nonumber\\
	&\leq \int_{0}^{T}\int_{d}^{D}\vert\v\hat{\eta}_{m}(t,b(t))\,(\rho u_{r})(t,b(t))\vert\,b(t)^{2}\,\mathrm{d}r\mathrm{d}t
	+\int_{0}^{T}\int_{d}^{D}\vert\v\hat{\eta}_{m}(t,r)\,(\rho u_{r})(t,r)\vert\,r^2\mathrm{d}r\mathrm{d}t\nonumber\\
	&\quad +\v\int_{0}^{T}\int_{d}^{D}\Big\vert\int_{r}^{b(t)}\r|u_{y}|\big(\hat{\eta}_{mu}|u_{y}|
     +\hat{\eta}_{m\r}|\r_{y}|\big)\,y^2\mathrm{d}y\Big\vert\,\mathrm{d}r\mathrm{d}t \nonumber\\
	&\leq C(D,M,E_0,T)+C\int_{0}^{T}\int_{d}^{D}\v|\rho u_{r}|\big(|u|+\rho^{\t(\rho)}\big)\,r^{2}\mathrm{d}r\mathrm{d}t\nonumber\\
	&\quad +C(D)\int_{0}^{T}\int_{d}^{D}\int_{r}^{b(t)}
	\var\rho|u_{y}|\big(|u_{y}|+\rho^{\t(\rho)-1}|\rho_{y}|\big)\,y^{2}\mathrm{d}y\mathrm{d}r\mathrm{d}t\nonumber\\
	&\leq C(D,M,E_0,T)+C\int_{0}^{T}\int_{d}^{D}\v\big(\rho|u|^2+\rho|u_{r}|^2+\rho^{\g(\rho)}\big)\,r^{2}\mathrm{d}r\mathrm{d}t\nonumber\\
	&\quad +C(D)\int_{0}^{T}\int_{d}^{b(t)}\v\rho|u_{y}|^2\,y^{2}\mathrm{d}y\mathrm{d}t
	+C(D)\int_{0}^{T}\int_{d}^{b(t)}\v\rho^{\g(\rho)-2}|\rho_{y}|^2\,y^{2}\mathrm{d}y\mathrm{d}t\nonumber\\
	&\leq C(D, M, E_0, T).\label{4.36}
\end{align}
For the second term, third term, and sixth term on the RHS of \eqref{4.31}, using \eqref{A.8-1}--\eqref{A.9-1}
and Lemmas \ref{lem4.1} and \ref{lem2.1},  we obtain
\begin{align}\label{4.37}
	&\Big\vert \int_{0}^{T}\int_{d}^{D}\Big(2\v\int_{r}^{b(t)}\hat{\eta}_{m}(t,y)\rho u\,\mathrm{d}y\Big)
	\,\mathrm{d}r\mathrm{d}t\Big\vert\nonumber\\
	&\leq C(d)\int_{0}^{T}\int_{d}^{D}\int_{r}^{b(t)}\v\big(|u|+\rho^{\t(\rho)}\big)\,\rho |u|
	\,y^2\mathrm{d}y\mathrm{d}r\mathrm{d}t\nonumber\\
	&\leq C(d, D)\int_{0}^{T}\int_{d}^{b(t)}\v\big(\rho|u|^2+\rho+\rho e(\rho)\big)\,y^{2}\mathrm{d}y\mathrm{d}t
	\leq C(d, D, M, E_0, T),
\end{align}
\begin{align}
	&\Big\vert\int_{0}^{T}\int_{d}^{D}\Big(\int_{r}^{b(t)}\hat{\eta}(t,y)\,y^{2}\mathrm{d}y\Big)_{t}\,\mathrm{d}r\mathrm{d}t\Big\vert\nonumber\\
	&\leq \int_{d}^{D}\int_{r}^{b(t)} |\hat{\eta}(T,y)|\, y^{2}\mathrm{d}y\mathrm{d}r
	+ \int_{d}^{D}\int_{r}^{b} |\hat{\eta}(0,y)|\, y^{2}\mathrm{d}y\mathrm{d}r\nonumber\\
	&\leq C\sup_{t\in [0,T]}\int_{d}^{D}\int_{r}^{b(t)} \big(\rho^{\g(\rho)}+\rho |u|^2\big)\, y^{2}\mathrm{d}y\mathrm{d}r\nonumber\\
	&\leq C\sup_{t\in [0,T]}\int_{d}^{D}\int_{r}^{b(t)}
	\big(\rho e(\rho)+\rho+\rho |u|^2\big)\,y^{2}\mathrm{d}y\mathrm{d}r
	\leq C(D, M, E_0, T),\label{4.37-1}
\end{align}
\begin{align}\label{4.38-1}
	&\Big\vert\int_{0}^T\int_{d}^{D}\Big(\int_{r}^{b(t)}\Big(\int_{a}^{y}\rho\, z^{2}\mathrm{d}z\Big)
	\rho\, \hat{\eta}_{m}\,\mathrm{d}y\Big)\,\mathrm{d}r\mathrm{d}t\Big\vert\nonumber\\
	&\leq \int_{0}^{T}\int_{d}^{D}\Big\vert\int_{r}^{b(t)}\Big(\int_{a}^{y}\rho\, z^{2}\mathrm{d}z\Big)\rho\, \hat{\eta}_{m}\,\mathrm{d}y\Big\vert\,\mathrm{d}r\mathrm{d}t\nonumber\\
	&\leq C(d,D,M)\int_{0}^{T}\int_{d}^{b(t)}\rho \big(|u|+\rho^{\g(\rho)}\big)\,r^2\mathrm{d}r\mathrm{d}t\nonumber\\
	&\leq C(d,D,M)\int_{0}^{T}\int_{d}^{b(t)}\big(\rho|u|^2+\rho+\rho e(\rho)\big)\,r^{2}\mathrm{d}r\mathrm{d}t\leq C(d, D, M, E_0, T).
\end{align}

For the fifth term on the RHS of \eqref{4.31}, we note from  \eqref{A.4-1}--\eqref{A.5-1} and Lemma \ref{lem4.1} that
\begin{align}
	&-\hat{q}+\rho u\hat{\eta}_{\rho}+\rho u^2\hat{\eta}_{m}=0\qquad\,\, \mbox{if $|u|\geq k(\rho)$},\label{4.39}\\
	&\, \vert-\hat{q}+\rho u\hat{\eta}_{\rho}+\rho u^2\hat{\eta}_{m}\vert
	\leq C\rho^{\g(\rho)+\t(\rho)}\leq C\big(\rho+\rho^{\g_2+\t_2}\big)\qquad\,\, \mbox{if $|u|\leq k(\rho)$}.\label{4.40}
\end{align}
Then it follows from \eqref{4.39}--\eqref{4.40} and  Corollary \ref{cor2.6} that
\begin{align}\label{4.41}
	&
	\int_{0}^{T}\int_{d}^{D}\Big|\int_{r}^{b(t)}\big(-\hat{q}+\rho u\hat{\eta}_{\rho}+\rho u^2\hat{\eta}_{m}\big)\,y\mathrm{d}y\Big| \,\mathrm{d}r\mathrm{d}t
	\nonumber\\
	&\leq C(d, D)\int_{0}^{T}\int_{d}^{b(t)} \big(\rho+\rho^{\g_2+\t_2}\big)\,y^{2}\mathrm{d}y\mathrm{d}t
	\leq C(d, D, M, E_0,T),
\end{align}
where we have used $\theta_2\in (0,1)$ since $\gamma_2\in(\f65,3)$.
Combining \eqref{4.31}--\eqref{4.32}, \eqref{4.36}--\eqref{4.38-1}, and \eqref{4.41},
we obtain that
$\int_{0}^{T}\int_{d}^{D}\hat{q}\,r^{2}\mathrm{d}r\mathrm{d}t\leq C(d, D, M, E_0,T)$,
which, along with \eqref{2.64-1} and Lemma \ref{lem4.1}, gives
\begin{align}\label{4.43}
	&\int_{0}^{T}\int_{[d, D]\cap \{r: |u|\geq k(\rho)\}}\rho |u|^3\,r^{2}\mathrm{d}r\mathrm{d}t\nonumber\\
	&\leq 2\int_{0}^{T}\int_{[d, D]\cap \{r: |u|\geq k(\rho)\}}\hat{q}\,r^{2}\mathrm{d}r\mathrm{d}t\nonumber\\
	&=2\int_{0}^{T}\int_{d }^{D}\hat{q}\,r^{2}\mathrm{d}r\mathrm{d}t
	-2\int_{0}^{T}\int_{[d, D]\cap \{r: |u|<k(\rho)\}}\hat{q}\, r^{2}\mathrm{d}r\mathrm{d}t\nonumber\\
	&\leq C(d, D, M, E_0, T)+C\int_{0}^{T}\int_{d}^{D}(\rho+\rho^{\g_2+1})\,r^{2}\mathrm{d}r\mathrm{d}t\nonumber\\
	&\leq C(d, D, M, E_0, T).
\end{align}
On the other hand, we have
\begin{align}\label{4.44}
	&\int_{0}^{T}\int_{[d, D]\cap \{r:\,|u|\leq  k(\rho)\}}\rho |u|^3\,r^{2}\mathrm{d}r\mathrm{d}t\nonumber\\
	&\leq C\int_{0}^{T}\int_{d}^{D}\rho^{\g(\rho)+\t(\rho)}\,r^{2}\mathrm{d}r\mathrm{d}t
	\nonumber\\	&
	\leq C\int_{0}^{T}\int_{d}^{D}\big(\rho+\rho p(\rho)\big)\,r^{2}\mathrm{d}r\mathrm{d}t
	\leq C(M, E_0,T).
\end{align}
Combining \eqref{4.43} with \eqref{4.44}, we obtain that
$\int_{0}^{T}\int_{d}^{D}\rho |u|^3\,r^{2}\mathrm{d}r\mathrm{d}t\leq C(d, D, M, E_0, T)$.
This completes the proof of Lemma \ref{lem4.3}. $\hfill\square$

\section{\,Existence of Global Weak Solutions of CNSPEs}
In this section, for fixed $\v>0$, we take the limit: $b \rightarrow \infty$
to obtain the global existence of solutions of the Cauchy problem for \eqref{1.8}.
Meanwhile, some uniform estimates in Theorem \ref{thm1.2} are obtained.
To take the limit, some careful attention is required,
since the weak solutions may involve the vacuum.
We use similar compactness arguments as in \cite{Chen-He-Wang-Yuan-2021, Chen-Wang-2020}
to handle the limit:  $b \rightarrow \infty$.
Throughout this section, we denote the smooth solutions of \eqref{2.1}--\eqref{2.6}
as $(\rho^{\varepsilon, b}, u^{\varepsilon, b})$ for simplicity.

First of all, we extend our solutions $(\rho^{\varepsilon, b}, u^{\varepsilon, b})$
to be zero on $([0, T] \times[0, \infty)) \backslash \Omega_{T}$.
It follows from Lemma \ref{lem2.3} that
\begin{equation}\label{5.1}
	\lim _{b \rightarrow \infty} \min_{t \in[0, T]} b(t)=\infty,
\end{equation}
which implies that domain $[0, T] \times[a, b(t)]$ expands to $[0, T] \times(0, \infty)$
as $b \rightarrow \infty$.
That is, for any
set $K \Subset (0, \infty)$, when $b \gg 1, K \Subset (a, b(t))$ for all $t \in[0, T]$.
Now we define
\begin{equation*}
	(\rho^{\v,b},\mathcal{M}^{\v,b},\Phi^{\v,b})(t,\mathbf{x})
	:=(\rho^{\v,b}(t,r),m^{\v,b}(t,r)\frac{\mathbf{x}}{r},\Phi^{\v,b}(t,r))
	\qquad \text{for }r=|\mathbf{x}|,
\end{equation*}
where $m^{\v,b}:=\rho^{\v,b}u^{\v,b}$.
Then it is direct to check that the corresponding functions
$(\rho^{\v,b},\mathcal{M}^{\v,b}, \Phi^{\v,b})(t,\mathbf{x})$ are classical solutions
of
\begin{equation*}
	\left\{
	\begin{aligned}
		&\partial_{t}\rho^{\v,b}+
		\nabla\cdot\mathcal{M}^{\v,b}=0,\\
		&\partial_{t}\mathcal{M}^{\v,b}+
		\nabla\cdot\Big(\frac{\mathcal{M}^{\v,b}\otimes \mathcal{M}^{\v,b}}{\rho^{\v,b}}\Big)+\nabla P(\rho^{\v,b})+\rho^{\v,b}\nabla\Phi^{\v,b}=\v
		\nabla\cdot\Big(\rho^{\v,b}D\big(\frac{\mathcal{M}^{\v,b}}{\rho^{\v,b}}\big)\Big),\\
		&\Delta \Phi^{\v,b}=\rho^{\v,b},
	\end{aligned}
	\right.
\end{equation*}
for $(t,\mathbf{x})\in [0,\infty)\times \Omega_{t}$ with $\mathcal{M}^{\v,b}\vert_{\partial B_{a}(\mathbf{0})}=0$.

Based on the estimates obtained in \S \ref{UE},
by the same arguments as in \cite[\S 4]{Chen-He-Wang-Yuan-2021},
we have

\smallskip
\begin{lemma}\label{lem5.1}
	For fixed $\v>0$, as $b\to\infty$ {\rm (}up to a subsequence{\rm )}, there exists a vector function $(\rho^{\v}, m^\v)(t,r))$ such that
	\begin{enumerate}
		\item[\rm (i)]
		$(\sqrt{\rho^{\v,b}},\rho^{\v,b})\to(\sqrt{\rho^{\v}}, \rho^{\v})$ {\it a.e.}  and strongly in
		$C(0,T;L_{\mathrm{loc}}^p)$
		for any $p\in [1,\infty)$, where $L_{\mathrm{loc}}^p$ denotes $L^p(K)$
		for any compact set $K\Subset (0,\infty)$. In particular, $\rho^{\v}\geq 0$ {\it a.e.} on $\R_{+}^2$.
		
		\smallskip
		\item[\rm (ii)]  The pressure function sequence $P(\rho^{\varepsilon, b})$ is uniformly bounded
		in $L^{\infty}(0, T ; L_{\mathrm {loc}}^{p}(\R))$ for all $p \in[1, \infty]$,
		and
		$$
		P(\rho^{\varepsilon, b}) \longrightarrow P(\rho^{\varepsilon}) \quad
         \text { strongly in } L^{p}(0, T ; L_{\mathrm{loc }}^{p}(\R))\qquad \text {for } p\in[1, \infty).
		$$	
		
		\item[\rm (iii)]
		The momentum function sequence $m^{\varepsilon, b}$
		converges strongly in $L^{2}(0, T ; L_{\mathrm{loc}}^{p}(\R))$
		to
		$m^{\varepsilon}$ for all $p \in[1, \infty)$.
		In particular,
		\begin{equation*}
			m^{\varepsilon, b}(t,r)=(\rho^{\varepsilon, b} u^{\varepsilon, b})(t,r) \longrightarrow m^{\varepsilon}(t, r)
			\qquad \text { {\it a.e.} in }[0, T] \times(0, \infty).
		\end{equation*}
		
		\item[\rm (iv)]
		$m^{\varepsilon}(t, r)=0$ {\it a.e.} on $\{(t, r)\,:\,\rho^{\varepsilon}(t, r)=0\}$.
		Furthermore, there exists a function $u^{\varepsilon}(t, r)$ such that
		$m^{\varepsilon}(t, r)=\rho^{\varepsilon}(t, r) u^{\varepsilon}(t, r)$ {\it a.e.},
		$u^{\varepsilon}(t, r)=0$ {\it a.e.} on $\{(t, r)\,:\,\rho^{\varepsilon}(t, r)=0\}$, and
		\begin{equation*}
			\begin{aligned}
				&m^{\varepsilon, b} \longrightarrow m^{\varepsilon}=\rho^{\varepsilon} u^{\varepsilon} \qquad
				\text { strongly in $L^{2}\left(0, T ; L_{\mathrm{loc}}^{p}(\R)\right)$ for $p \in[1, \infty)$}, \\
				&\frac{m^{\varepsilon, b}}{\sqrt{\rho^{\varepsilon, b}}} \longrightarrow \frac{m^{\varepsilon}}{\sqrt{\rho^{\varepsilon}}}
				=\sqrt{\rho^{\varepsilon}} u^{\varepsilon} \qquad \text { strongly in $L^{2}(0, T ; L_{\mathrm{loc}}^{2}(\R))$}.
			\end{aligned}
		\end{equation*}
	\end{enumerate}
\end{lemma}

Let $(\rho^{\v},m^{\v})(t,r)$ be the limit function obtained above. Using \eqref{2.11}, \eqref{5.1},
Lemmas \ref{lem2.1}--\ref{lem2.4}, \ref{lem4.3}, and \ref{lem5.1},
Corollaries \ref{cor2.1}--\ref{cor2.2} and \ref{cor2.6}, Fatou's lemma, and the lower semicontinuity,
we conclude the proof of \eqref{1.29}--\eqref{1.32}.

Now we show the convergence of the potential functions $\Phi^{\v,b}$.
Using the similar arguments as in \cite[Lemma 4.6]{Chen-He-Wang-Yuan-2021}, we have

\smallskip
\begin{lemma}\label{lem5.3-1}
	For fixed $\v>0$, there exists a function $\Phi^{\v}(t,\mathbf{x})=\Phi^{\v}(t,r)$ such that, as $b\to \infty$
	$($up to a subsequence$)$,
	\begin{align}
		&\Phi^{\v,b}
		\,{\rightharpoonup}\, \Phi^{\v}\quad
		\text{weak-star in $L^{\infty}(0,T;H_{\mathrm{loc}}^{1}(\R^3))$
			and weakly in $L^2{(0,T;H_{\mathrm{loc}}^{1}(\R^3))}$},\label{5.34-1}\\
		&\Phi_{r}^{\v,b}(t,r)r^{2}\,\to\,\Phi_{r}^{\v}(t,r)r^{2}=\int_{0}^{r}\rho^{\v}(t,y)\,y^{2}\mathrm{d}y
		\qquad \text{in $C_{\mathrm{loc}}([0,T]\times [0,\infty))$},\label{5.34-2}\\
		&\|\Phi^{\varepsilon}(t)\|_{L^{6}(\mathbb{R}^{3})}+\|\nabla \Phi^{\varepsilon}(t)\|_{L^{2}(\mathbb{R}^{3})}
		\leq C(M, E_{0}) \qquad \text {for $t \geq 0$}.\label{5.34-3}
	\end{align}
	Moreover, since $\gamma_2>\f65$,
	\begin{equation}\label{5.34-4}
		\int_{0}^{\infty}\big|\big(\Phi_{r}^{\varepsilon, b}-\Phi_{r}^{\varepsilon}\big)(t, r)\big|^{2}
		\,r^{2}\mathrm{d} r \rightarrow 0 \qquad \text{ as $b \rightarrow \infty$ $($up to a subsequence$)$}.
	\end{equation}
\end{lemma}

Using \eqref{5.34-4}, Fatou's lemma, and Lemmas \ref{lem2.1} and \ref{lem5.1},
we obtain the following energy inequality:
\begin{align}
		&\int_{0}^{\infty}\Big(\frac{1}{2}\Big|\frac{m^{\varepsilon}}{\sqrt{\rho^{\varepsilon}}}\Big|^{2}
		+\rho^{\varepsilon} e(\rho^{\varepsilon})\Big)(t, r) \,r^{2} \mathrm{d} r
		-\frac{1}{2} \int_{0}^{\infty}|\Phi^{\varepsilon}(t, r)|^{2} \,r^{2}\mathrm{d} r\nonumber \\
		&\leq \int_{0}^{\infty}\Big(\frac{1}{2}\Big|\frac{m_{0}^{\varepsilon}}{\sqrt{\rho_{0}^{\varepsilon}}}\Big|^{2}
		+\rho_{0}^{\varepsilon} e(\rho_{0}^{\varepsilon})\Big)(r)\, r^{2}\mathrm{d} r
		-\frac{1}{2} \int_{0}^{\infty}|\Phi_{0}^{\varepsilon}(r)|^{2}\,r^{2}\mathrm{d} r.\label{5.34-11}
\end{align}
We denote
$$
(\rho^{\v},\mathcal{M}^{\v}, \Phi^{\v})(t,\mathbf{x}):=(\rho^{\v}(t,r),m^{\v}(t,r)\frac{\mathbf{x}}{r},\Phi^{\v}(t,r)).
$$
Then \eqref{1.33} follows directly from \eqref{5.34-11}.
Moreover, we can prove that $(\rho^{\v},\mathcal{M}^{\v}, \Phi^{\v})$ is a global weak solution of
the Cauchy problem  \eqref{1.8} and \eqref{1.19}--\eqref{1.19-1} in the sense of Definition \ref{def1.2}.
In fact, by the same arguments in \cite[Remark 4.7 and Lemmas 4.9--4.11]{Chen-He-Wang-Yuan-2021},
we have

\smallskip
\begin{lemma}\label{lem5.4}
	Let $0 \leq t_{1}<t_{2} \leq T$, and let $\zeta(t, \mathbf{x}) \in C_{0}^{1}([0, T] \times \mathbb{R}^{3})$
	be any smooth function with compact support. Then
	\begin{equation}\label{5.35}
		\int_{\mathbb{R}^{3}} \rho^{\varepsilon}(t_{2}, \mathbf{x}) \zeta(t_{2}, \mathbf{x})\,\mathrm{d} \mathbf{x}
		=\int_{\mathbb{R}^{3}} \rho^{\varepsilon}(t_{1},\mathbf{x}) \zeta(t_{1},\mathbf{x})\, \mathrm{d} \mathbf{x}
		+\int_{t_{1}}^{t_{2}} \int_{\mathbb{R}^{3}}(\rho^{\varepsilon} \zeta_{t}
		+\mathcal{M}^{\varepsilon} \cdot \nabla \zeta)\,\mathrm{d} \mathbf{x}\mathrm{d} t.
	\end{equation}
	Moreover, \eqref{1.28} holds, and the total mass is conserved{\rm :}
	\begin{equation}\label{5.36}
		\int_{\mathbb{R}^{3}} \rho^{\varepsilon}(t, \mathbf{x})\, \mathrm{d} \mathbf{x}
		=\int_{\mathbb{R}^{3}} \rho_{0}^{\varepsilon}(\mathbf{x})\, \mathrm{d} \mathbf{x}=M \qquad\, \text {for $t \geq 0$}.
	\end{equation}
\end{lemma}

\begin{lemma}\label{lem5.5}
	Let $\Psi(t, \mathbf{x}) \in(C_{0}^{2}([0, T] \times \mathbb{R}^{3}))^{3}$ be any smooth function
	with compact support so that $\Psi(T, \mathbf{x})=0$. Then
	\begin{align}
		&\int_{\mathbb{R}_{+}^{4}}\Big\{\mathcal{M}^{\varepsilon} \cdot \partial_{t} \Psi
		+\frac{\mathcal{M}^{\varepsilon}}{\sqrt{\rho^{\varepsilon}}} \cdot\big(\frac{\mathcal{M}^{\varepsilon}}{\sqrt{\rho^{\varepsilon}}} \cdot \nabla\big) \Psi+P(\rho^{\varepsilon})
		\nabla\cdot\Psi-\r^{\v}\nabla\Phi^{\v}\cdot \Psi\Big\}\,\mathrm{d} \mathbf{x} \mathrm{d} t
		\nonumber\\
		&\quad +\int_{\mathbb{R}^{3}} \mathcal{M}_{0}^{\varepsilon}(\mathbf{x}) \cdot \Psi(0, \mathbf{x}) \,\mathrm{d} \mathbf{x}
		\nonumber\\
		&=-\varepsilon \int_{\mathbb{R}^{4}_+}\Big\{\frac{1}{2} \mathcal{M}^{\varepsilon} \cdot\big(\Delta \Psi
		+\nabla (\nabla\cdot \Psi)\big)+\frac{\mathcal{M}^{\varepsilon}}{\sqrt{\rho^{\varepsilon}}} \cdot\big(\nabla \sqrt{\rho^{\varepsilon}} \cdot \nabla\big) \Psi
		\nonumber\\
		&\qquad \qquad\quad +\nabla \sqrt{\rho^{\varepsilon}} \cdot\big(\frac{\mathcal{M}^{\varepsilon}}{\sqrt{\rho^{\varepsilon}}} \cdot \nabla\big) \Psi\Big\}\,\mathrm{d} \mathbf{x} \mathrm{d} t
		\nonumber\\
		&=\sqrt{\varepsilon} \int_{\mathbb{R}_{+}^{4}} \sqrt{\rho^{\varepsilon}}\Big\{V^{\varepsilon} \frac{\mathbf{x} \otimes \mathbf{x}}{r^{2}}+\frac{\sqrt{\varepsilon}}{r} \frac{m^{\varepsilon}}{\sqrt{\rho^{\varepsilon}}}\Big(I_{3 \times 3}-\frac{\mathbf{x} \otimes \mathbf{x}}{r^{2}}\Big)\Big\}: \nabla \Psi\, \mathrm{d} \mathbf{x} \mathrm{d} t
		\label{5.53}
	\end{align}
	with $V^{\varepsilon}(t, \mathbf{x}) \in L^{2}(0, T ; L^{2}(\mathbb{R}^{3}))$ as a function satisfying
	\begin{equation*}
		\int_{0}^{T} \int_{\mathbb{R}^{3}}\left|V^{\varepsilon}(t, \mathbf{x})\right|^{2} \,\mathrm{d} \mathbf{x}\mathrm{d} t
		\leq C(E_{0}, M),
	\end{equation*}
	where $C\left(E_{0}, M\right)>0$ is a constant independent of $T>0$.
\end{lemma}

\smallskip
\begin{lemma}\label{lem5.5-1}
	It follows from \eqref{5.34-2} that $\Phi^{\v}$ satisfies Poisson's equation in the classical sense
	except for the origin{\rm :}
	$\Delta\Phi^{\v}=\rho^{\v}(t,\mathbf{x})$ for $(t,\mathbf{x})\in [0,\infty)\times (\R^{3}\backslash \{\mathbf{0}\})$.
	Moreover, for any smooth function $\xi(\mathbf{x})\in C_{0}^{1}(\R^3)$ with compact support,
	\begin{equation}\label{5.80}
		\int_{\mathbb{R}^{3}} \nabla \Phi^{\varepsilon}(t, \mathbf{x}) \cdot \nabla \xi(\mathbf{x})\, \mathrm{d} \mathbf{x}
		=- \int_{\mathbb{R}^{3}} \rho^{\varepsilon}(t, \mathbf{x})\, \xi(\mathbf{x}) \,\mathrm{d} \mathbf{x} \qquad \text {for } t \geq 0.
	\end{equation}
\end{lemma}

\section{$\,W_{\mathrm{loc}}^{-1,p}$--Compactness of Weak Entropy Dissipation  Measures}
In this section, using
the estimates of the weak entropy pairs obtained in Lemmas \ref{lem6.2}, \ref{lem6.4}, and \ref{lem6.5}, we establish the compactness
of weak entropy dissipation measures:
$\pa_{t}\e^{\psi}(\rho^\v,m^\v)+\pa_{r}q^{\psi}(\r^\v,m^\v)$ for each weak entropy pair $(\e^{\psi},q^{\psi})$.
Unfortunately, we fail to obtain the same $H_{\mathrm{loc}}^{-1}$-compactness as in \cite{Chen-He-Wang-Yuan-2021, Chen-Wang-2020}, since we only obtain that
$q^{\v}$ is uniformly bounded in $L_{\mathrm{loc}}^{2}$ from
Lemma \ref{lem6.5} and Corollary \ref{cor2.6}.
Instead, using similar arguments as in \cite[Lemma 4.2]{Chen-He-Wang-Yuan-2021},
we can obtain the compactness in $W_{\mathrm{loc}}^{-1,p}$ for any $p\in [1,2)$.

\smallskip
\begin{lemma}[Compactness of the entropy dissipation measures]\label{lem6.6}
	Let $(\e^{\psi},q^{\psi})$ be a weak entropy pair defined in
	\eqref{6.10} for any smooth and compactly supported function $\psi(s)$ on $\R$. Then, for $\v \in (0,\v_0]$,
	\begin{equation}\label{6.48}
		\partial_{t}\eta^{\psi}(\rho^{\v},m^{\v})+\partial_{r}q^{\psi}(\rho^{\v},m^{\v})\qquad \text{is compact in }W_{\mathrm{loc}}^{-1,p}(\R_{+}^2) \quad\mbox{for any $p\in [1,2)$}.
	\end{equation}
\end{lemma}

\noindent{\bf Proof.}
To establish \eqref{6.48}, we first need to study the equation:
$\partial_{t}\eta^{\psi}(\rho^{\v},m^{\v})+\partial_{r}q(\rho^{\v},m^{\v})$ in the distributional sense,
which is more complicated than that in \cite{Chen-Perepelitsa-2010, Chen-Perepelitsa-2015}.
For simplicity, we denote $(\e^{\v,b},q^{\v,b})=(\e^{\psi}(\rho^{\v,b},m^{\v,b}),q^{\psi}(\r^{\v,b},m^{\v,b}))$
and $(\e^{\v},q^{\v})=(\e^{\psi}(\rho^{\v},m^{\v}),q^{\psi}(\rho^{\v},m^{\v}))$.
We divide it into four steps.

\smallskip
1. Considering $\eqref{2.1}_1\times \eta_{\rho}^{\v,b}+\eqref{2.1}_2\times \eta_{m}^{\v,b}$, we obtain
\begin{align}
	&\partial_{t} \eta(\rho^{\varepsilon, b}, m^{\varepsilon, b})+\partial_{r} q(\rho^{\varepsilon, b}, m^{\varepsilon, b})
	+\frac{2}{r} m^{\varepsilon, b}\big(\eta_{\rho}^{\varepsilon, b}+u^{\varepsilon, b} \eta_{m}^{\varepsilon, b}\big)
	\nonumber\\
	&=-\eta_{m}^{\varepsilon, b} \frac{\rho^{\varepsilon, b}}{r^{2}}
	\int_{0}^{r} \rho^{\varepsilon, b}(t, y)\,y^{2}\mathrm{d}y
	+\varepsilon \eta_{m}^{\varepsilon, b}\Big\{\big(\rho^{\varepsilon, b}(u_{r}^{\varepsilon, b}
	+\frac{2}{r} u^{\varepsilon, b})\big)_{r}-\frac{2}{r} \rho_{r}^{\varepsilon, b} u^{\varepsilon, b}\Big\},\label{6.49}
\end{align}
where $\rho^{\v,b}$ is understood to be zero in domain $[0,T]\times [0,a)$ so that
$\int_{a}^{r}\r^{\v,b}(t,z)\,z^{2}\mathrm{d}z$ can be written as
$\int_{0}^{r}\r^{\v,b}(t,z)\,z^{2}\mathrm{d}z$ in the potential term.
Let $\phi(t, r) \in C_{0}^{\infty}\left(\mathbb{R}_{+}^{2}\right)$ and $b \gg 1$
so that $\operatorname{supp}\phi(t, \cdot) \in(a, b(t))$.
Multiplying \eqref{6.49} by $\phi$ and integrating by parts yield
\begin{align}
	&\int_{\mathbb{R}_{+}^{2}}(\partial_{t} \eta^{\varepsilon, b}+\partial_{r} q^{\varepsilon, b}) \phi \,\mathrm{d} r \mathrm{d} t\nonumber\\
	&=-\int_{\mathbb{R}_{+}^{2}} \frac{2}{r} m^{\varepsilon, b}(\eta_{\rho}^{\varepsilon, b}
	+u^{\varepsilon, b} \eta_{m}^{\varepsilon, b}) \phi \,\mathrm{d} r \mathrm{d} t
	-\varepsilon \int_{\mathbb{R}_{+}^{2}} \rho^{\varepsilon, b}(\eta_{m}^{\varepsilon, b})_{r}
	\big(u_{r}^{\varepsilon, b}+\frac{2}{r} u^{\varepsilon, b}\big) \phi \,\mathrm{d} r\mathrm{d} t \nonumber\\
	&\quad\, -\varepsilon \int_{\mathbb{R}_{+}^{2}} \rho^{\varepsilon, b} \eta_{m}^{\varepsilon, b}
	\big(u_{r}^{\varepsilon, b}+\frac{2}{r} u^{\varepsilon, b}\big) \phi_{r} \,\mathrm{d} r \mathrm{d} t
	-\varepsilon \int_{\mathbb{R}_{+}^{2}} \frac{2}{r} \eta_{m}^{\varepsilon, b}
	\rho_{r}^{\varepsilon, b} u^{\varepsilon, b} \phi\, \mathrm{d} r \mathrm{d} t
	\nonumber\\		
	&\quad\,
	- \int_{\mathbb{R}_{+}^{2}} \frac{\rho^{\varepsilon, b}}{r^{2}}\eta_{m}^{\varepsilon, b}
	\Big(\int_{0}^{r} \rho^{\varepsilon, b}(t, y)\,y^{2}\mathrm{d} y\Big) \phi\, \mathrm{d} r \mathrm{d} t
	:=\sum_{i=1}^{5} I_{i}^{\varepsilon, b}.\label{6.50}
\end{align}

2. From Lemmas  \ref{thm6.1} and \ref{lem5.1},
it is clear to see that
\begin{equation}\label{6.51}
	\e^{\v,b}\longrightarrow \e^{\v}\qquad \text{{\it a.e.} in }\{(t,r)\,:\rho^{\v}\neq 0\}\quad \text{as }b\to \infty.
\end{equation}
In $\{(t,r)\,:\,\rho^{\v}(t,r)=0\}$, it follows from Lemmas \ref{lem6.2} and \ref{lem6.4} that
\begin{equation}\label{6.52}
	|\e^{\v,b}|\leq C\rho^{\v,b}\longrightarrow 0=\e^{\v}\qquad \text{as $b\to \infty$}.
\end{equation}
Combining \eqref{6.51}--\eqref{6.52}, we obtain
\begin{equation}\label{6.53}
	\e^{\v,b}\longrightarrow \e^{\v}\qquad\, \text{{\it a.e.} as $b\to \infty$}.
\end{equation}
Similarly, it follows from Lemmas \ref{thm6.2}, \ref{lem6.5}, and \ref{lem5.1}
that
\begin{equation}\label{6.54}
	q^{\v,b}\longrightarrow q^{\v}\qquad \text{{\it a.e.} as $b\to \infty$}.
\end{equation}

For $\gamma_2 \in(1,3)$ and any subset
$K\Subset (0,\infty)$, it follows from  Lemmas \ref{lem6.2},
\ref{lem6.4}, \ref{lem6.5},  and Corollary \ref{cor2.6} that
\begin{align*}
	\int_{0}^{T} \int_{K}\big(|\eta^{\varepsilon, b}|^{\g_2+1}
	+|q^{\varepsilon, b}|^{2}\big)\,\mathrm{d} r \mathrm{d} t
	&\leq C_{\psi}(K)\int_{0}^{T} \int_{K}
	\big(1+ |\rho^{\varepsilon, b}|^{\gamma_2+1}\Big)\, \mathrm{d} r \mathrm{d} t\nonumber\\
	&\leq C_{\psi}(K, M, E_{0}, T),\nonumber
\end{align*}
which implies that
$(\eta^{\varepsilon, b}, q^{\varepsilon, b})$
is uniformly bounded in $L_{\mathrm{loc}}^{2}(\mathbb{R}_{+}^{2})$.
This, with \eqref{6.53}--\eqref{6.54}, yields that, up to a subsequence,
\begin{equation*}
	(\eta^{\varepsilon, b}, q^{\varepsilon, b}) \rightharpoonup (\eta^{\varepsilon}, q^{\varepsilon})
	\qquad \text {in $L_{\mathrm{loc}}^{2}(\mathbb{R}_{+}^{2})$ as $b \rightarrow \infty$}.
\end{equation*}
Thus, for any $\phi \in C_{0}^{1}(\mathbb{R}_{+}^{2})$,
as $b \rightarrow \infty$ (up to a subsequence),
\begin{align}\label{6.58}
	\int_{\mathbb{R}_{+}^{2}}\big(\partial_{t} \eta^{\varepsilon, b}+\partial_{r} q^{\varepsilon, b}\big) \phi \,\mathrm{d} r \mathrm{d} t
	&=-\int_{\mathbb{R}_{+}^{2}}\big(\eta^{\varepsilon, b} \partial_{t} \phi+q^{\varepsilon, b} \partial_{r} \phi\big) \,\mathrm{d} r \mathrm{d} t\nonumber\\
	&\longrightarrow-\int_{\mathbb{R}_{+}^{2}}\big(\eta^{\varepsilon} \partial_{t} \phi+q^{\varepsilon} \partial_{r} \phi\big) \,\mathrm{d} r \mathrm{d} t.
\end{align}
Furthermore, $(\eta^{\varepsilon}, q^{\varepsilon})$ is uniformly bounded in $L_{\mathrm{loc }}^{2}(\mathbb{R}_{+}^{2})$, which implies that
\begin{equation}\label{6.59}
	\partial_{t} \eta^{\varepsilon}+\partial_{r} q^{\varepsilon} \quad \text { is uniformly bounded in } \varepsilon>0 \text { in } W_{\text {loc }}^{-1, 2}(\mathbb{R}_{+}^{2}).
\end{equation}
Since $q^{\v,b}$ is only uniformly bounded in $L_{\mathrm{loc}}^{2}(\R_{+}^2)$
in view of  Lemma \ref{lem6.5} and Corollary \ref{cor2.6}, we cannot conclude that
$\pa_{t}\eta^{\v}+\pa_{r}q^{\v}$ is uniformly bounded in $\v>0$ in $W_{\mathrm{loc}}^{-1,p}(\R_{+}^2)$ with $p>2$,
which is different from \cite{Chen-He-Wang-Yuan-2021}.

\smallskip
3. For the terms on the RHS of \eqref{6.50}, using Lemmas \ref{lem6.2}, \ref{lem6.4}, and \ref{lem6.5},
and similar calculations as in \cite[Lemma 5.11]{Chen-He-Wang-Yuan-2021}, we obtain that
\begin{equation}\label{6.62}
	\begin{aligned}
		&I_{1}^{\varepsilon, b} \longrightarrow-\int_{\mathbb{R}_{+}^{2}} \frac{2}{r} m^{\varepsilon}(\eta_{\rho}^{\varepsilon}
		+u^{\varepsilon} \eta_{m}^{\varepsilon}) \phi \,\mathrm{d} r \mathrm{d} t
		\qquad \text {up to a subsequence as $b \rightarrow \infty$}, \\
		&\int_{0}^{T} \int_{K}\frac{2}{r}\big| m^{\varepsilon}(\eta_{\rho}^{\varepsilon}
		+u^{\varepsilon} \eta_{m}^{\varepsilon})\big|^{\frac{7}{6}}\, \mathrm{d} r \mathrm{d} t \leq C_{\psi}(K, M, E_{0}, T),
	\end{aligned}
\end{equation}
and there exist local bounded Radon measures $(\mu_{1}^{\varepsilon}, \mu_{2}^{\varepsilon}, \mu_{3}^{\varepsilon})$
on $\mathbb{R}_{+}^{2}$ such that, as $b \rightarrow \infty$ (up to a subsequence),
\begin{align*}
	-(\varepsilon \rho^{\varepsilon, b}(\eta_{m}^{\varepsilon, b})_{r}(u_{r}^{\varepsilon, b}+\frac{2}{r} u^{\varepsilon, b}),
	\,\frac{2\v}{r} \eta_{m}^{\varepsilon, b}\rho_{r}^{\varepsilon, b} u^{\varepsilon, b},\,
	\kappa \eta_{m}^{\varepsilon, b} \frac{\rho^{\varepsilon, b}}{r^{n-1}} \int_{0}^{r} \rho^{\varepsilon, b}(t, z)
	\,z^{2}\mathrm{d} z) \rightharpoonup (\mu_{1}^{\varepsilon}, \mu_{2}^{\varepsilon},  \mu_{3}^{\varepsilon}).
\end{align*}
In addition, for $i=1,2,3$,
\begin{equation}\label{6.65}
	\mu_{i}^{\varepsilon}((0, T) \times V) \leq C_{\psi}(K, T, E_{0})
	\qquad\,\, \mbox{for any open subset $V \subset K$}.
\end{equation}
Then, up to a subsequence, we have
\begin{equation}\label{6.66}
	I_{2}^{\varepsilon, b}+I_{4}^{\varepsilon, b}+I_{5}^{\v,b} \longrightarrow\langle\,\mu_{1}^{\varepsilon}
	+\mu_{2}^{\varepsilon}+\mu_{3}^{\v},\,\phi\rangle \qquad \text { as $b \rightarrow \infty$}.
\end{equation}
Moreover, there exists a function $f^{\varepsilon}$ such that, as $b \rightarrow \infty$ (up to a subsequence),
\begin{equation}\label{6.68}
	\begin{aligned}
		&-\sqrt{\varepsilon} \rho^{\varepsilon, b} \eta_{m}^{\varepsilon, b}\big(u_{r}^{\varepsilon, b}
		+\frac{2}{r} u^{\varepsilon, b}\big) \rightharpoonup f^{\varepsilon} \qquad
		\text { weakly in $L_{\mathrm{loc}}^{\frac{4}{3}}(\mathbb{R}_{+}^{2})$},\\
		&\int_{0}^{T} \int_{K}\left|f^{\varepsilon}\right|^{\frac{4}{3}} \,\mathrm{d} r \mathrm{d} t \leq C_{\psi}(K, M, E_{0}, T).
	\end{aligned}
\end{equation}
Then it follows from \eqref{6.68} that
\begin{equation}\label{6.69}
	I_{3}^{\varepsilon, b} \longrightarrow \sqrt{\varepsilon} \int_{0}^{T} \int_{K} f^{\varepsilon} \phi_{r} \,\mathrm{d} r \mathrm{d} t
	\qquad \text { as $b \rightarrow \infty$ (up to a subsequence)}.
\end{equation}

4. Taking $b \rightarrow \infty$ (up to a subsequence) on both sides of \eqref{6.50},
it follows from \eqref{6.58}, \eqref{6.62}, \eqref{6.66}, and \eqref{6.69} that
\begin{equation*}
	\partial_{t} \eta^{\varepsilon}+\partial_{r} q^{\varepsilon}=-\frac{2}{r} \rho^{\varepsilon} u^{\varepsilon}(\eta_{\rho}^{\varepsilon}+u^{\varepsilon} \eta_{m}^{\varepsilon})+\mu_{1}^{\varepsilon}+\mu_{2}^{\varepsilon}+\mu_{3}^{\varepsilon}-\sqrt{\varepsilon} f_{r}^{\varepsilon}
\end{equation*}
in the distributional sense. It follows from \eqref{6.62}--\eqref{6.65} that
\begin{equation}\label{6.71}
	-\frac{2}{r} \rho^{\varepsilon} u^{\varepsilon}\left(\eta_{\rho}^{\varepsilon}+u^{\varepsilon} \eta_{m}^{\varepsilon}\right)+\mu_{1}^{\varepsilon}+\mu_{2}^{\varepsilon}+\mu_{3}^{\varepsilon}
\end{equation}
is a locally uniformly bounded Radon measure sequence.
From \eqref{6.68}, we know that
\begin{equation}\label{6.72}
	\sqrt{\varepsilon} f_{r}^{\varepsilon} \longrightarrow 0 \qquad
	\text { in $W_{\operatorname{loc}}^{-1, \frac{4}{3}}(\mathbb{R}_{+}^{2})$ as $\varepsilon \rightarrow 0$}.
\end{equation}
Then it follows from \eqref{6.71}--\eqref{6.72} that the sequence:
\begin{equation}\label{6.73}
	\partial_{t} \eta^{\varepsilon}+\partial_{r} q^{\varepsilon} \qquad
	\text{ is confined in a compact subset of $W_{\text {loc }}^{-1, p_{2}}(\mathbb{R}_{+}^{2})$
		for some $p_{2} \in(1,2)$},
\end{equation}
which also implies that
\begin{equation}\label{6.73-1}
	\pa_{t}\eta^{\v}+\pa_{r}q^{\v}\qquad \text{is compact in
		$W_{\mathrm{loc}}^{-1,p}(\R_{+}^2)$ with $1\leq p\leq p_{2}$}.
\end{equation}

On the other hand, the interpolation compactness theorem ({\it cf.} \cite{Chen-1986, Ding-Chen-Luo}) indicates that,
for $p_{2}>1, p_{1} \in\left(p_{2}, \infty\right]$, and $p_{0} \in\left[p_{2}, p_{1}\right)$,
\begin{align*}
	&\big(\text{compact set of }W_{\mathrm{loc}}^{-1, p_{2}}(\mathbb{R}_{+}^{2})\big)
	\cap \big(\text{bounded set of }W_{\mathrm{loc}}^{-1, p_{1}}(\mathbb{R}_{+}^{2})\big)
	\nonumber\\
	&\quad \subset \big(\text{compact set of }W_{\mathrm{loc}}^{-1, p_{0}}(\mathbb{R}_{+}^{2})\big),
\end{align*}
which is a generalization of Murat's lemma in \cite{Murat-1978, Tartar-1979}.
Combining this theorem for $1<p_{2}<2$ and $p_{1}=2$ with the facts in \eqref{6.59} and \eqref{6.73}, we conclude that
\begin{equation}\label{6.73-2}
	\pa_{t}\eta^{\v}+\pa_{r}q^{\v}\qquad
	\text{is compact in $W_{\mathrm{loc}}^{-1,p}(\R_{+}^2)$ with $p_2\leq p<2$}.
\end{equation}
Combining \eqref{6.73-2} with \eqref{6.73-1}, we conclude \eqref{6.48}.
$\hfill\square$

\section{\,$L^p$ Compensated Compactness Framework}
In this section, with the help of our understanding of the singularities of the entropy kernel
and entropy flux kernel obtained in Lemma \ref{lem7.3}, we now establish the
$L^p$ compensated compactness framework and
complete the proof of Theorem \ref{thm1.3}.
The key ingredient is to prove the reduction of the Young measure.
The arguments are similar to
\cite[\S 4]{Schrecker-Schulz-2019} and \cite[\S 7]{Schrecker-Schulz-2020},
based on \cite{Chen-Perepelitsa-2010,Chen-LeFloch-2000},
so we only sketch the proof for self-containedness.

We denote the upper half-plane by
$$
\mathbb{H}:=\{(\rho,u)\in \mathbb{R}^2\,:\, \rho>0\}
$$
and consider the following subset of continuous functions:
\begin{equation*}
	\overline{C}(\mathbb{H}):=\left\{\phi \in C(\overline{\mathbb{H}})\, :
	\,\, \begin{array}{l}
		\phi(\rho,u) \,\text{ is constant on the vacuum states }\, \{\rho=0\} \, \text { and} \\
		\text{the map: } (\rho,u)\, \mapsto\, \lim\limits _{s \rightarrow \infty} \phi(s \rho, s u)
		\text { belongs to } C(\mathbb{S}^{1} \cap \overline{\mathbb{H}})
	\end{array}\right\},
\end{equation*}
where $\mathbb{S}^{1} \subset \mathbb{R}^{2}$ is the unit circle.
Since $\overline{C}(\mathbb{H})$ is a complete sub-ring of the continuous functions on $\mathbb{H}$ containing the constant functions,
there exists a compactification $\overline{\mathcal{H}}$ of $\mathbb{H}$ such that $C(\overline{\mathcal{H}})$ is isometrically
isomorphic to $\overline{C}(\mathbb{H})$ ({\it cf}. \cite{LeFloch-Westdickenberg-2007}),
written $C(\overline{\mathcal{H}})\cong \overline{C}(\mathbb{H})$.
The topology of $\overline{\mathcal{H}}$ is the weak-star topology induced by $C(\overline{\mathcal{H}})$, {\it i.e.},
a sequence $\{v_{n}\}_{n \in \mathbb{N}}$ in $\overline{\mathcal{H}}$ converges to $v \in \overline{\mathcal{H}}$
if $|\varphi(v_{n})-\varphi(v)| \rightarrow 0$ for all $\varphi \in C(\overline{\mathcal{H}})$, which
is separable and metrizable ({\it cf}. \cite{LeFloch-Westdickenberg-2007}).
Denote by $V$ the weak-star closure of the vacuum states $\{(\rho,u) \in \R^2\,:\, \rho=0\}$ and
define $\mathcal{H}:=\mathbb{H}\cup V$. In view of the functions that lie in $\overline{C}(\mathbb{H})$,
the topology of $\overline{\mathcal{H}}$ does not distinguish points in
$V$.
Since $\overline{\mathcal{H}}$ is homeomorphic to a compact metric space,
we may apply the fundamental theorem of Young measures
in Alberti-M\"{u}ller \cite[Theorem 2.4]{Alberti-Muller-2001}.

\smallskip
\begin{lemma}[{\cite[Theorem 2.4]{Alberti-Muller-2001}}]\label{lem7.4}
	Given any sequence of measurable functions $\left(\rho^{\v}, \rho^{\v}u^{\v}\right): \mathbb{R}_{+}^{2} \rightarrow \overline{\mathbb{H}}$,
	there exists a subsequence $($still denoted$)$ $(\r^{\v},\rho^{\v}u^{\v})$
	generating a Young measure $\nu_{(t, r)} \in \operatorname{Prob}(\overline{\mathcal{H}})$ in the sense that,
	for any $\phi \in \overline{C}(\mathbb{H})$,
	\begin{equation*}
		\phi(\rho^{\v}(t, r), u^{\v}(t, r))\, \overset{*}{\rightharpoonup}\,
		\int_{\overline{\mathcal{H}}} \iota(\phi)(\r,u) \,\mathrm{d} \nu_{(t, r)}(\r, u) \qquad \text { in $L^{\infty}(\mathbb{R}_{+}^{2})$},
	\end{equation*}
	where $\iota:\overline{C}(\mathbb{H})\to C(\overline{\mathcal{H}})$ is an isometrically isomorphism.
	Moreover, sequence $(\r^{\v},\rho^{\v}u^{\v})$ converges to $(\r, \r u):\,\R_{+}^2\to \overline{\mathcal{H}}$
	{\it a.e.} if and only if
	$$
	\nu_{(t,r)}(\rho,u)=\delta_{(\rho(t,r),m(t,r))}  \qquad \text{a.e. $(t,r)\in \R_{+}^2$},
	$$
	in the phase coordinates $(\rho,m)$ with $m=\rho u$.
\end{lemma}

From now on, we often use the same letter $\nu_{(t,r)}$ for an element of
$\big(\overline{C}(\mathbb{H})\big)^{*}$ or $\big(C(\overline{\mathcal{H}})\big)^{*}$,
and use the same letter for $\iota(\phi)$ and $\phi$ for simplicity,
when no confusion arises.

The following lemma shows that the Young measure $\nu_{(t,r)}$, generated by the sequence of
measurable approximate solutions $(\rho^{\v},\rho^{\v}u^{\v})$ satisfying the assumptions of Theorem \ref{thm1.3},
is only supported on the interior of $\mathcal{H}$.
Moreover, the Young measure $\nu_{(t,r)}$ can be extended to a larger class of test functions
than just $\overline{C}(\mathbb{H})$. This is proved in \cite[Proposition 5.1]{Chen-Perepelitsa-2010};
also see
\cite[Proposition 2.3]{LeFloch-Westdickenberg-2007}.

\smallskip
\begin{lemma}[{\cite[Proposition 5.1]{Chen-Perepelitsa-2010}}]\label{lem7.6}
	The following statements hold{\rm :}
	\begin{itemize}
		\item[\rm (\rmnum{1})] For the Young measure $\nu_{(t,r)}$ generated by the sequence of measurable approximate
		solutions $(\rho^{\v},\rho^{\v}u^{\v})$ satisfying the assumptions in {\rm Theorem \ref{thm1.3}},
		\begin{equation*}
			(t,r)\mapsto \int_{\mathbb{H}}(\rho^{\g_2+1}+\rho|u|^3)\,\mathrm{d}\nu_{(t,r)}(\rho,u)\in L_{\mathrm{loc}}^1(\mathbb{R}_{+}^2)  .
		\end{equation*}
		\item[\rm (\rmnum{2})] Let $\phi(\r,u)$ be a function such that
		\begin{itemize}
			\item[\rm (a)] $\phi$ is continuous on $\overline{\mathbb{H}}$ and  $\phi=0$ on $\partial \mathbb{H}$,
			\item[\rm (b)] there exists a constant $\mathfrak{a}>0$ such that $\operatorname{supp}\phi\subset \{u+k(\rho)\geq -\mathfrak{a},u-k(\rho)\leq \mathfrak{a}\}$,
			\item[\rm (c)] $|\phi(\rho,u)|\leq \rho^{\beta(\gamma_2+1)}$ for all $(\rho,u)$ with large $\rho$ and some $\beta\in (0,1)$.
		\end{itemize}
		\smallskip
		Then $\phi$ is $\nu_{(t,r)}$--integrable for $(t,r)\in \mathbb{R}_{+}^2$ {\it a.e.} and
		\begin{equation*}
			\phi(\rho^{\v}(t,r), u^{\varepsilon}(t,r))\rightharpoonup
			\int_{\mathbb{H}}\phi(\r,u) \,\mathrm{d}\nu_{(t,r)}(\r,u)\qquad \text{in }\;L_{\mathrm{loc}}^1(\mathbb{R}_{+}^2).
		\end{equation*}
		\item[\rm (\rmnum{3})] $\nu_{(t,r)}\in \operatorname{Prob}(\mathcal{H})$ for $(t,r)\in \R_{+}^2$ a.e., that is,
		$\,
		\nu_{(t,r)}\big(\overline{\mathcal{H}}\backslash (\mathbb{H}\cup V)\big)=0.
		$
	\end{itemize}	
\end{lemma}

We now prove the commutation relation. Since we only have the $W_{\mathrm{loc}}^{-1,p}$--compactness
of the entropy dissipation measures for $p\in [1,2)$,
the classical div-curl lemma in \cite{Murat-1978} fails to obtain
the commutation relation.
Thus, we adopt an improved version of the div-curl lemma.

\smallskip
\begin{lemma}[{\cite[Theorem]{Conti-Dolzmann-Muller-2011}}]\label{lem7.7}
	Let $\Omega\subset \R^n$ be an open bounded set, and $p,q\in (1,\infty)$ with $\frac{1}{p}+\frac{1}{q}=1$.
	Let $\mathbf{v}^{\var}$ and $\mathbf{w}^{\var}$ are sequences of vector fields such that
	\begin{equation*}
		\mathbf{v}^{\var}\rightharpoonup \mathbf{v} \,\, \text{ in }L^p(\Omega;\mathbb{R}^n),
		\quad \mathbf{w}^{\var}\rightharpoonup \mathbf{w} \,\, \text{ in }L^q(\Omega;\R^n)\qquad \text{ as }\var\to 0.
	\end{equation*}
	Suppose that $\mathbf{v}^{\v}\cdot \mathbf{w}^{\v}$ is equi-integrable uniformly in $\v$, and
	\begin{equation*}
		\begin{split}
			&
			\operatorname{div}
			\mathbf{v}^{\var} \quad\,\,\, \text{ is $($pre-$)$compact in }W^{-1,1}(\Omega),\\
			&\operatorname{curl} \mathbf{w}^{\var} \quad \text{ is $($pre-$)$compact in }W^{-1,1}(\Omega;\R^{n\times n}).
		\end{split}
	\end{equation*}
	Then
	$\mathbf{v}^{\var}\cdot \mathbf{w}^{\var}\,\rightharpoonup \,\mathbf{v}\cdot \mathbf{w}$
	in $\mathcal{D}^{\prime}(\Omega)$.
\end{lemma}

\smallskip
\begin{lemma}[Commutation relation]\label{lem7.8}
	Let $\{(\rho^{\var},\rho^{\v}u^{\v})\}_{\v>0}$ be the measurable approximate solutions
	satisfying the assumptions of {\rm Theorem \ref{thm1.3}},
	and let $\nu_{(t,r)}$ be a Young measure generated by
	the family $\{(\rho^{\varepsilon},\rho^{\v}u^{\v})\}_{\v>0}$ in {\rm Lemma \ref{lem7.6}}. Then
	\begin{equation}\label{7.64}
		\overline{\chi(s_1)\sigma(s_2)-\chi(s_2)\sigma(s_1)}=\overline{\chi(s_1)}\;\overline{\sigma(s_2)}-\overline{\chi(s_2)}\;\overline{\sigma(s_1)}
	\end{equation}
	for all $s_1,s_2\in \mathbb{R}$, where
	$\dis \overline{f}:=\int f\, \mathrm{d}\nu_{(t,r)}$,
	$\chi(s_i)=\chi(\cdot, \cdot-s_i)$, and $\sigma(s_i)=\sigma(\cdot,\cdot-s_i)$.
\end{lemma}

\noindent\textbf{Proof.}
For any $\psi \in C_0^2(\R)$, it follows from Lemmas \ref{lem6.2}, \ref{lem6.4}, and \ref{lem6.5} that
\begin{equation}\label{7.65}
	|\eta^{\psi}(\rho,m)|\leq C_{\psi}\rho,\qquad |q^{\psi}(\rho,m)|\leq  C_{\psi}\big(\r+\rho^{1+\theta_2}\big).
\end{equation}
It is clear that the support of $(\eta^{\psi}, q^{\psi})$ is contained
in $\left\{k(\rho)+u \geq -L, u-k(\r) \leq L\right\}$ for some $L>0$
depending only on $\operatorname{supp}\,\psi$.
For any $\psi_{1}, \psi_{2} \in C_{0}^{2}(\mathbb{R})$, we consider the sequences of vector fields:
\begin{equation*}
	\mathbf{v}^{\varepsilon}=(\eta^{\psi_{1}}(\rho^{\varepsilon}, \rho^{\varepsilon} u^{\varepsilon}),
	q^{\psi_{1}}(\rho^{\varepsilon}, \rho^{\varepsilon} u^{\varepsilon})),
	\qquad \mathbf{w}^{\varepsilon}=(q^{\psi_{2}}(\rho^{\varepsilon}, \rho^{\varepsilon} u^{\varepsilon}),
	-\eta^{\psi_{2}}(\rho^{\varepsilon}, \rho^{\varepsilon} u^{\varepsilon})).
\end{equation*}
Noting $\rho^{\v}\in L_{\mathrm{loc}}^{1+\g_2}(\R_{+}^2)$ and \eqref{7.65},
we see that both $\mathbf{v}^{\varepsilon}$ and $\mathbf{w}^{\varepsilon}$ are uniformly bounded
sequences in $L_{\mathrm{loc}}^2(\R_{+}^2)$.
Moreover, by Lemma \ref{lem7.6} and the uniqueness of weak limits, we obtain
\begin{equation*}
	\mathbf{v}^{\varepsilon} \rightharpoonup(\overline{\eta^{\psi_{1}}}, \overline{q^{\psi_{1}}})\,\,\,
	\text { in $L_{\mathrm{loc}}^{2}(\R_{+}^2)$},
	\qquad \mathbf{w}^{\varepsilon} \rightharpoonup (\overline{q^{\psi_{2}}},-\overline{\eta^{\psi_{2}}})
	\,\,\,\text { in $L_{\mathrm{loc}}^{2}(\R_{+}^2)$}.
\end{equation*}
By direct calculation, we see that
\begin{equation*}
	\begin{aligned}
		&\operatorname{div}
		\mathbf{v}^{\varepsilon}
		=\partial_{t}\eta^{\psi_1}(\rho^{\var},\rho^{\var}u^{\var})+\partial_{r}q^{\psi_1}(\rho^{\var},\rho^{\var}u^{\var})
		\quad\,\,\,\text{are compact in }W_{\mathrm{loc}}^{-1,1}(\R_{+}^2),
		\\&\operatorname{curl} \mathbf{w}^{\varepsilon}=\partial_{t}\eta^{\psi_2}(\rho^{\var},\rho^{\var}u^{\var})
		+\partial_{r}q^{\psi_2}(\rho^{\var},\rho^{\var}u^{\var})\quad \text{are compact in }W_{\mathrm{loc}}^{-1,1}(\R_{+}^2).
	\end{aligned}
\end{equation*}
Using \eqref{7.65}, we obtain that
$
\vert\mathbf{v}^{\v}\cdot \mathbf{w}^{\v}\vert\leq C\big((\rho^{\v})^2+(\r^{\v})^{2+\t_2}\big)
$
for $\rho>0$
which, with \eqref{1.32}, yields that $\mathbf{v}^{\v}\cdot \mathbf{w}^{\v}\in L_{\mathrm{loc}}^{\frac{1+\gamma_2}{2+\t_2}}(\R_{+}^2)$ uniformly in $\v$.
Thus, $\mathbf{v}^{\v}\cdot \mathbf{w}^{\v}$ is equi-integrable uniformly in $\v$ since $\frac{1+\gamma_2}{2+\t_2}>1$. It follows from Lemma \ref{lem7.7} that
\begin{equation}\label{7.70}
	\mathbf{v}^{\varepsilon} \cdot \mathbf{w}^{\varepsilon}\, \rightarrow \, \overline{\eta^{\psi_{1}}}\;\overline{q^{\psi_{2}}}-\overline{\eta^{\psi_{2}}}\; \overline{q^{\psi_{1}}} \qquad
	\text {in the sense of distributions in
		$\R_{+}^2$}.
\end{equation}

On the other hand, using \eqref{7.65} and Lemma \ref{lem7.6}, we find that
\begin{equation*}
	\mathbf{v}^{\varepsilon} \cdot \mathbf{w}^{\varepsilon}
	\,\rightarrow\, \overline{\eta^{\psi_{1}} q^{\psi_{2}}-\eta^{\psi_{2}} q^{\psi_{1}}}
	\qquad \text { in $L_{\mathrm{loc}}^{1}(\R_{+}^2)$},
\end{equation*}
which, with \eqref{7.70}, yields that
\begin{equation}\label{7.72}
	\overline{\eta^{\psi_{1}} q^{\psi_{2}}-\eta^{\psi_{2}} q^{\psi_{1}}}=\overline{\eta^{\psi_{1}}}
	\;\overline{q^{\psi_{2}}}-\overline{\eta^{\psi_{2}}} \;\overline{q^{\psi_{1}}}.
\end{equation}
It follows from \eqref{7.72} and Fubini's theorem that
\begin{align*}
	\int_{\mathbb{R}^2}
	\Big(\overline{\chi(s_{1}) \sigma(s_{2})-\chi(s_{2}) \sigma(s_{1})}
- \overline{\chi(s_{1})}\;\overline{\sigma(s_{2})}+\overline{\chi(s_{2})}\;\overline{\sigma(s_{1})}\Big)
	\psi_{1}(s_{1}) \psi_{2}(s_{2})  \,\mathrm{d} s_{1} \mathrm{d} s_{2}=0.
\end{align*}
Since $\psi_{1}, \psi_{2} \in C_{0}^{2}(\mathbb{R})$ are arbitrary, we conclude
\begin{equation*}
	\overline{\chi(s_{1})\sigma(s_{2})-\chi(s_{2}) \sigma(s_{1})}=\overline{\chi(s_{1})}
	\;\overline{\sigma(s_{2})}-\overline{\chi(s_{2})} \;\overline{\sigma(s_{1})}\qquad \mbox{for any $s_1, s_2\in \R$}.
\end{equation*}
This completes the proof.
$\hfill\square$

\smallskip
\begin{theorem}[Reduction of the Young measure]\label{thm7.1}
	Let $\nu_{(t,r)}\in \mathrm{Prob}(\mathcal{H})$ be the Young
	measure generated
	by sequence $\{(\rho^{\varepsilon},\rho^{\v}u^{\v})\}_{\v>0}$ in {\rm Lemma \ref{lem7.6}}.
	Then either $\nu_{(t,r)}$ is contained in $V$ or the support of $\nu_{(t,r)}$ is a single point in $\mathbb{H}$.
\end{theorem}

\smallskip
\noindent\textbf{Proof.} The proof is similar
to \cite{Chen-LeFloch-2000, Lions-Perthame-Souganidis-1996, Schrecker-Schulz-2019, Schrecker-Schulz-2020}.
Since the estimates of the entropy kernel and entropy flux kernel
are different, we sketch the proof for self-containedness.

Taking $s_1,s_2,s_3\in \R$ and multiplying \eqref{7.64} by $\overline{\chi(s_3)}$, one obtains
\begin{equation*}
	\overline{\chi(s_3)}\;\overline{\chi(s_1)\,\sigma(s_2)-\chi(s_2)\,\sigma(s_1)}
	=\overline{\chi(s_3)}\;\overline{\chi(s_1)}\;\overline{\sigma(s_2)}-\overline{\chi(s_3)}\;\overline{\chi(s_2)}\;\overline{\sigma(s_1)}.
\end{equation*}
Cyclically permuting index $s_j$ and adding the resultant equations together, we have
\begin{align*}
	&\overline{\chi(s_{1})}\; \overline{\chi(s_{2}) \sigma(s_{3})-\chi(s_{3}) \sigma(s_{2})}\nonumber\\
	&= \overline{\chi(s_{3})}\; \overline{\chi(s_{2}) \sigma(s_{1})-\chi(s_{1}) \sigma(s_{2})} -\overline{\chi(s_{2})}\; \overline{\chi(s_{3}) \sigma(s_{1})-\chi(s_{1}) \sigma(s_{3})} .
\end{align*}
Applying the fractional derivative operators $P_{2}:=\partial_{s_{2}}^{\lambda_1+1}$
and $P_{3}:=\partial_{s_{3}}^{\lambda_1+1}$ in the sense of distributions to obtain
\begin{align}\label{7.78}
		&\overline{\chi(s_1)}\;\overline{P_2\chi(s_2)\,P_3\sigma(s_3)-P_3\chi(s_3)\,P_2\sigma(s_2)}\nonumber\\
		&=\overline{P_3\chi(s_3)}\;\overline{P_2\chi(s_2)\,\sigma(s_1)-\chi(s_1)\,P_2\sigma(s_2)}\nonumber\\
		&\quad -\overline{P_2\chi(s_2)}\;\overline{P_3\chi(s_3)\,\sigma(s_1)-\chi(s_1)\,P_3\sigma(s_3)},
\end{align}
where, for example, distribution $\overline{P_2\chi(s_2)}$ is defined by
$$
\langle \overline{P_2\chi(s_2)},\psi\rangle=-\int_{\R}\overline{\partial_{s_2}^{\lambda_1}\chi(s_2)}
\,\psi'(s_2)\,\mathrm{d}s_2\qquad\, \text{for $\psi\in \mathcal{D}(\R)$}.
$$

We take two standard but different functions $\phi_{2}, \phi_{3} \in C_{0}^{\infty}(-1,1)$ such that
$\dis\int_{\mathbb{R}} \phi_{j}(s_{j})\,\mathrm{d} s_{j}=1$
with $\phi_{j} \geq 0$ for $j=2,3$.
For $\tau>0$, denote $\phi_{j}^{\tau}(s_{j}):=\frac{1}{\tau}\phi_{j}(\frac{s_{j}}{\tau})$.
As indicated in \cite{LeFloch-Westdickenberg-2007}, we can always choose $\phi_2$ and $\phi_3$ such that
\begin{equation}\label{7.78-1}
	Y(\phi_{2}, \phi_{3})=\int_{-\infty}^{\infty} \int_{-\infty}^{s_{2}} \big(\phi_{2}(s_{2}) \phi_{3}(s_{3})
	-\phi_{2}(s_{3}) \phi_{3}(s_{2})\big) \,\mathrm{d} s_{3} \mathrm{d}s_{2}>0.
\end{equation}
Multiplying \eqref{7.78} by $\phi_{2}^{\tau}(s_1-s_2)\phi_3^{\tau}(s_1-s_3)$
and integrating the resultant equation with respect to $(s_2, s_3)$ yield
\begin{align}\label{7.79}
	&\overline{\chi(s_1)}\;\overline{P_2\chi_{2}^{\tau}\,P_3\sigma_3^{\tau}-P_3\chi_3^{\tau}\,P_2\sigma_{2}^{\tau}}\nonumber\\
	&=\overline{P_3\chi_{3}^{\tau}}\;\overline{P_2\chi_2^{\tau}\,\sigma_1-\chi_1\,P_2\sigma_2^{\tau}} -\overline{P_2\chi_2^{\tau}}\;\overline{P_3\chi_3^{\tau}\,\sigma_1-\chi_1\,P_3\sigma_3^{\tau}},
\end{align}
where we have used the notion:
$
\overline{P_{j} \chi_{j}^{\tau}}=\overline{P_{j} \chi_{j}} * \phi_{j}^{\tau}(s_{1})
=\int_{\R} \overline{\partial_{s_{j}}^{\lambda} \chi(s_{j})} \frac{1}{\tau^{2}} \phi_{j}^{\prime}(\frac{s_{1}-s_{j}}{\tau})
\,\mathrm{d}s_{j}$
for $j=2,3$.

Multiplying \eqref{7.79} by $\psi(s_1)\in \mathcal{D}(\R)$,
integrating the resultant equation with respect to $s_1$,
then taking limit $\tau\to 0$ and applying Lemmas \ref{lem7.11}--\ref{lem7.12} below,
we obtain
\begin{equation}\label{7.110}
	Y(\phi_{2}, \phi_{3})
	\int_{\mathcal{H}} Z(\rho) \sum_{\pm}(K^{\pm})^{2}\, \overline{\chi(u \pm k(\rho))}\, \psi(u \pm k(\rho)) \,\mathrm{d} \nu_{(t,r)}(\rho, u)=0.
\end{equation}
Noting that $Z(\rho)>0$ for $\rho>0$ from Lemma \ref{lem7.10} below,
$Y(\phi_2,\phi_3)> 0$ from \eqref{7.78-1}, and $\psi(s)$ is an arbitrary test function,
we deduce from \eqref{7.110} that
\begin{equation}\label{7.111}
	\int_{\mathcal{H}} Z(\rho)\,\overline{\chi(u\pm k(\rho))}\,\mathrm{d} \nu_{(t,r)}(\rho, u)=0.
\end{equation}
We define $\mathbb{S}=\{s\in \R\,:\,\overline{\chi(s)}>0\}$.
It follows from \cite{Schrecker-Schulz-2019} that $\mathbb{S}$ admits the representation:
$$
\mathbb{S}=\left\{s\in \R\,:\, u-k(\rho)<s<u+k(\rho)\,\,
\text{ with $(\rho,u)\in \operatorname{supp}\nu_{(t,r)}$}\right\}.
$$

For the case: $\mathbb{S}=\emptyset$, it is clear that $\overline{\chi(s)}=0$ for all $s\in \R$
so that $\operatorname{supp}\nu_{(t,r)}\subset V$,
since $\chi(s)>0$ for all $\rho>0$ and $s\in (u-k(\rho), u+k(\rho))$.

\newcommand*{\medcup}{\mathbin{\scalebox{1.5}{\ensuremath{\cup}}}}
For the case: $\mathbb{S}\neq \emptyset$, it follows from \eqref{7.78-2} below
that $s\mapsto \overline{\chi(s)}$ is a continuous map. Then $\mathbb{S}$ is an open set
so that $\mathbb{S}$ is at most a countable union of open intervals. Thus, we may write
$$
\mathbb{S}=\underset{k}{\medcup}(\zeta_k,\, \xi_k)
$$
for at most countably many numbers $\zeta_k:=u_{k}-k(\rho_{k})$ and
$\xi_k:=u_{k}+k(\rho_{k})$ with $(\rho_{k},u_{k})\in \mathrm{supp}\nu_{(t,r)}$
in the extended real line $\R\cup\{\pm\infty\}$
such that $\zeta_k<\xi_k\leq \zeta_{k+1}$ for all $k$.
For later use, we denote the Riemann invariants $z(\rho,u):=u-k(\rho) $ and $w(\rho,u):=u+k(\rho)$.
Thus, noting that $\mathrm{supp}\,\chi(s)=\{(\rho,u)\,:\,z(\rho,u)\leq s\leq w(\rho,u)\}$,
we obtain
$$
\operatorname{supp}\nu_{(t,r)}\subset \underset{k}{\medcup}\{(\rho,u)\in
\mathbb{H}\,:\, \zeta_k\leq z(\rho,u)<w(\rho,u)\leq \xi_{k}\}\cup V.
$$

If $\zeta_{k}$ and $\xi_{k}$ are both finite,
due to the fact that $k(\rho)$ is a strictly monotone increasing and unbounded function of $\rho$,
it is clear that
$\left\{(\rho, u) \,:\, \zeta_{k} \leq z(\rho, u)\leq w(\rho, u) \leq \xi_{k}\right\}$ is bounded.
Now we deduce from \eqref{7.111} that
$$
\operatorname{supp} \nu_{(t,r)} \cap \{(\rho, u)
\in \mathbb{H}\,:\, \zeta_{k}< z(\rho,u)<w(\rho, u)<\xi_{k}\}
=\emptyset\qquad \text {for all } k.
$$
Thus, the support of measure $\nu_{(t,r)}$ must be contained in the vacuum set $V$
and at most a countable union of points $P_k(\rho_{k},\,u_{k})$:
$$
\operatorname{supp} \nu_{(t,r)} \subset V \cup \big(\underset{\{k: \zeta_{k}, \xi_{k}
	\text{ are finite}\}}{\medcup}P_{k}(\rho_{k}, u_{k})\big).
$$
Therefore, we may write
\begin{equation}\label{7.113}
	\nu_{(t,r)}=\nu_{V}+\sum_{k} \alpha_{k} \delta_{P_{k}} \qquad\, \mbox{for all $\alpha_{k} \in[0,1]$}
\end{equation}
with measure $\nu_{V}$ supported on the vacuum set $V$.
For later use, we denote
$$
\chi(P_{k},s):=\chi(\rho_{k}, u_{k}, s), \quad \sigma(P_{k}, s):=\sigma(\rho_{k}, u_{k}, s)
\qquad\,\, \text{for $s\in \R$}.
$$
We claim that, if $\chi(P_{k},s)>0$, then $\chi(P_{k'},s)=0$ for all $k\neq k'$.
Indeed, recall that $\operatorname{supp}\chi(s)=\{(\r,u)\,:\,z(\rho,u)\leq s\leq w(\rho,u)\}$
and that $\chi(\rho, u, s)>0$ if and only if $z(\rho,u)<s<w(\rho,u)$.
If $\chi(P_{k},s)>0$, then $\zeta_{k}<s<\xi_{k}$.
If, in addition, $\chi(P_{k'},s)>0$ for some $k\neq k'$,
it must hold that $\zeta_{k'}< s < \xi_{k'}$.
However, since $\xi_{k-1}\leq\zeta_{k}< \xi_{k}\leq\zeta_{k+1}$,
this is impossible for any $P_{k'}$ with $k' \neq k$.

Thus, taking $s_{1}, s_{2} \in \mathbb{R}$ such that $\chi(P_{k},s_{1})\chi(P_{k},s_{2})>0$,
we deduce from the commutation relation
\eqref{7.64} and \eqref{7.113} that
$$
(\alpha_{k}-\alpha_{k}^{2})
\big(\chi(P_{k}, s_{1}) \sigma(P_{k}, s_{2})-\chi(P_{k}, s_{2}) \sigma(P_{k}, s_{1})\big)=0.
$$
Now, choosing $s_{1}$ and $s_{2}$ such that the second factor in this expression is non-zero,
we obtain that $\alpha_k=0$ or $1$ for all $k$. This completes the proof. $\hfill\square$

\smallskip
Combining Theorem \ref{thm7.1} with Lemma \ref{lem7.4},
we conclude that $(\rho^{\v},m^{\v})$ converges to $(\rho,m)$ almost everywhere.
Moreover, noting that
$|m^{\v}|^{\frac{3(\gamma_{2}+1)}{\gamma_{2}+3}}\leq C\big((\rho^{\v})^{\g_{2}+1}+\rho^{\v}|u^{\v}|^3\big)$
for any $T,d,D>0$, we have
$$
\int_{0}^{T}\int_{d}^{D}|m^{\v}|^{\frac{3(\gamma_{2}+1)}{\gamma_{2}+3}}\,\mathrm{d}r\mathrm{d}t
\leq C\int_{0}^{T}\int_{d}^{D}\big((\rho^{\v})^{\g_{2}+1}+\rho^{\v}|u^{\v}|^3\big)\,\mathrm{d}r\mathrm{d}t\leq C(d,D,T),
$$
which implies that $m^{\v}$ is uniformly bounded in $L_{\mathrm{loc}}^{\frac{3(\gamma_{2}+1)}{\gamma_{2}+3}}(\R_{+}^2)$
with respect to $\v$. This implies that \eqref{1.35} holds.
Therefore, the proof of Theorem \ref{thm1.3} is complete.

Now, we are going to prove the auxiliary lemmas, Lemmas \ref{lem7.11}--\ref{lem7.12}, which
are used in the proof of Theorem \ref{thm7.1}.
We first recall two useful lemmas in \cite{Chen-LeFloch-2000, LeFloch-Westdickenberg-2007}.

\smallskip
\begin{lemma}[{\cite[Lemmas 3.8--3.9]{LeFloch-Westdickenberg-2007}}]\label{lem7.9}
	Let $\mathfrak{R} \in C_{\mathrm{loc }}^{\alpha}(\mathbb{R})$ be a H\"older continuous function
	for some $\alpha \in(0,1)$, and let $g \in C_{0}^{\alpha}(\mathbb{R})$ be a H\"{o}lder continuous function
	with compact support. Assume $L_{0}>2$ such that $\operatorname{supp} g \subset B_{L_0-2}(0)$.
\begin{itemize}
\item[\rm (\rmnum{1})] For
		any pair of distributions $T_{2}, T_{3} \in \mathcal{D}^{\prime}(\mathbb{R})$
		from the following collection{\rm :}
		$$
		(T_{2}, T_{3})=(\delta, Q_{3}), \,\, (PV, Q_{3}), \,\, (Q_{2}, Q_{3})
		$$
		with $Q_{2}, Q_{3} \in\{H, Ci, \mathfrak{R}\}$,
		there exists a constant $C>0$ independent of $(\rho, u)$ such that
		\begin{align*}
			\sup\limits_{\tau \in(0,1)}\Big\vert\int_{-\infty}^{\infty} g(s_{1})\Big\{\big(&T_{2}(s_{2}-u \pm k(\rho)) T_{3}(s_{3}-u \pm k(\rho))\\
			&-T_{2}(s_{3}-u\pm k(\rho))T_{3}(s_{2}-u\pm k(\rho))\big) * \phi_{2}^{\tau} * \phi_{3}^{\tau}\Big\}
			\left(s_{1}\right)\,\mathrm{d} s_{1}\Big\vert\\
			\leq C\|g\|_{C^{ \alpha}(\mathbb{R})}\big(1+&\|\mathfrak{R}\|_{C^{ \alpha}(\overline{B_{L_0}(0)})}\big)^{2}.
		\end{align*}
\item[\rm (\rmnum{2})] For
		any pair of distributions from
		\begin{equation*}
			(T_{2}, T_{3})=(\delta, \delta), \,\,(PV, PV),\,\, (Q_{2}, Q_{3}),\,\,
			(\delta, PV), \,\, (PV, Q_{3}),
		\end{equation*}
with $Q_{2}, Q_{3} \in\{H, Ci, \mathfrak{R}\}$,
there exists $C>0$ independent of $(\rho,u)$ such that
		\begin{align*}
			&\sup\limits_{\tau \in(0,1)} \Big\vert \int_{-\infty}^{\infty}\Big\{\big((s_{2}-s_{3}) T_{2}(s_{2}
			-u \pm k(\rho))T_{3}(s_{3}-u \pm k(\rho))\big)* \phi_{2}^{\tau} * \phi_{3}^{\tau}\Big\}(s_{1}) \,\mathrm{d} s_{1} \Big\vert	\\
			&\quad \leq C\|g\|_{C^{ \alpha}(\mathbb{R})}\big(1+\|\mathfrak{R}\|_{C^{ \alpha}(\overline{B_{L_0}(0)})}\big)^{2}.
		\end{align*}
	\end{itemize}
\end{lemma}

Motivated by \cite{Chen-LeFloch-2000}, it follows from Lemmas \ref{thm6.1}--\ref{thm6.2}, \eqref{1.3},
and a direct calculation that
\begin{align}\label{6.9}
		D(\r):=&\,a_1(\r)b_1(\r)-2k(\r)^2(a_1(\r)b_2(\r)-a_2(\r)b_1(\r))\nonumber\\
		=&\,\frac{M_{\lambda_1}^2}{2(\lambda_1+1)}k(\r)^{-2\lambda_1}k'(\rho)^{-2}
		\big(k'(\rho)+(\rho k'(\rho))'\big)>0\qquad\, \text{for $\rho>0$}.
\end{align}

\begin{lemma}[{\cite[Lemmas 4.2--4.3]{Chen-LeFloch-2000}}]\label{lem7.10}
	The mollified fractional derivatives of the entropy kernel and
	the entropy flux kernel satisfy the following convergence properties{\rm :}
	\begin{itemize}
		\item[\rm (\rmnum{1})] When $0\leq \rho<\infty$,
		\begin{equation*}
			P_{2} \chi_{2}^{\tau} P_{3} \sigma_{3}^{\tau}-P_{3} \chi_{3}^{\tau} P_{2} \sigma_{2}^{\tau}
			\longrightarrow Y(\phi_{2}, \phi_{3}) Z(\rho) \sum_{\pm}(K^{\pm})^{2} \delta_{s_{1}=u \pm k(\rho)}
		\end{equation*}
		as $\tau \rightarrow 0$ weakly-star in measures in $s_{1}$ and locally uniformly in $(\rho, u)$,
		where $Y(\phi_2,\phi_3)$ satisfies \eqref{7.78-1}, $Z(\rho):=(\lambda_1+1) M_{\lambda}^{-2} k(\rho)^{2 \lambda} D(\rho)>0$
		with $D(\rho)$ defined in \eqref{6.9}, and $K^{\pm}\neq 0$ are some constants.
		
		\item[\rm (\rmnum{2})] For $j=2,3$,
		$\chi_1\,P_{j}\sigma_{j}^{\tau}-\sigma_1\,P_{j}\chi_{j}^{\tau}$ are H\"{o}lder continuous
		in $(\rho,u,s_1)$, uniformly in $\tau$,
		and there exists a H\"{o}lder continuous function $X=X(\rho,u,s_1)$, independent of the mollifying sequence $\phi_j$, such that
		\begin{equation*}
			\chi(s_1) P_{j} \sigma_{j}^{\tau}-P_{j} \chi_{j}^{\tau} \sigma(s_1) \longrightarrow X(\rho, u, s_{1})
			\qquad \text{as $\tau\rightarrow 0$}
		\end{equation*}
		uniformly in $\left(\rho, u, s_{1}\right)$ on the sets on which $\rho$ is bounded.
	\end{itemize}
\end{lemma}

\begin{lemma}\label{lem7.11}
	For any test function $\psi\in\mathcal{D}(\R)$,
	\begin{align}\label{7.85}
			&\lim\limits_{\tau\to 0}\int_{\R}\overline{\chi(s_1)}\;\overline{P_2\chi_{2}^{\tau}P_3\sigma_3^{\tau}
				-P_3\chi_3^{\tau}P_2\sigma_2^{\tau}}(s_1)\psi(s_1)\,\mathrm{d}s_1
			\nonumber\\
			&=Y(\phi_2,\phi_3)\int_{\mathcal{H}}Z(\rho)
			\sum\limits_{\pm}(K^{\pm})^2\,\overline{\chi(u\pm k(\rho))}\,\psi(u\pm k(\rho))\,\mathrm{d}\nu_{(t,r)}(\rho,u),
	\end{align}
	where $Y(\phi_2,\phi_3)$ is defined by \eqref{7.78-1} and $Z(\rho)$ is given in {\rm Lemma \ref{lem7.10}}.
\end{lemma}

\smallskip
\noindent\textbf{Proof.} $\,$ It follows from Lemma \ref{lem7.10} that,
when $\rho$ is bounded,
\begin{equation*}
	P_{2} \chi_{2}^{\tau} P_{3} \sigma_{3}^{\tau}-P_{3} \chi_{3}^{\tau} P_{2} \sigma_{2}^{\tau}
	\rightarrow Y(\phi_{2}, \phi_{3}) Z(\rho) \sum_{\pm}\left(K^{\pm}\right)^{2} \delta_{s_{1}=u \pm k(\rho)}\qquad \text{as }\tau\to 0
\end{equation*}
locally uniform in $(\rho,u)$ and hence pointwise for all $(\rho,u)$.
Therefore, we have
\begin{align}
	&\lim _{\tau \rightarrow 0} \int_{-\infty}^{\infty}
	\overline{\chi(s_{1})}\langle\,\nu_{(t,r)},\,
	(P_{2} \chi_{2}^{\tau} P_{3} \sigma_{3}^{\tau}-P_{3} \chi_{3}^{\tau} P_{2} \sigma_{2}^{\tau})
	{\mathbf{I}}_{\{\rho \leq \rho^{*}\}}\rangle \psi(s_{1})\,\mathrm{d} s_{1}\nonumber\\
	&=\lim _{\tau \rightarrow 0}\langle\,\nu_{(t,r)},\, \int_{-\infty}^{\infty} \overline{\chi(s_{1})}(P_{2}
	\chi_{2}^{\tau} P_{3} \sigma_{3}^{\tau}-P_{3} \chi_{3}^{\tau} P_{2} \sigma_{2}^{\tau}) \psi(s_{1}) \,\mathrm{d}
	s_{1}{\mathbf{I}}_{\{\rho \leq \rho^{*}\}}\rangle \nonumber\\
	&=\langle\,\nu_{(t,r)}, \, Y(\phi_{2}, \phi_{3}) Z(\rho) \sum_{\pm}(K^{\pm})^{2} \overline{\chi(u \pm k(\rho))} \psi(u \pm k(\rho))
	{\mathbf{I}}_{\{\rho \leq \rho^{*}\}}\rangle \nonumber\\
	&=Y(\phi_{2}, \phi_{3}) \sum_{\pm}(K^{\pm})^{2}\langle\,\nu_{(t,r)},
	\, Z(\rho) \overline{\chi(u \pm k(\rho))} \psi(u \pm k(\rho)){\mathbf{I}}_{\{\rho \leq \rho^{*}\}} \rangle.\label{7.87}
\end{align}

For $\rho\geq \rho^{*}$,
we notice that
\begin{equation}\label{7.89}
	P_2\chi_{2}^{\tau}\,P_3\sigma_3^{\tau}-P_3\chi_{3}^{\tau}\, P_2\sigma_2^{\tau}
	=P_2\chi_{2}^{\tau}\,P_3(\sigma_3^{\tau}-u\chi_{3}^{\tau})-P_3\chi_{3}^{\tau}\,P_2(\sigma_{2}^{\tau}-u\chi_{2}^{\tau}).	
\end{equation}
Using Lemma \ref{lem7.3}, we see that \eqref{7.89} consists of a sum of terms of the form:
$$
A_{i,\pm}(\rho)B_{j,\pm}(\rho)(s_{3}-s_{2})T_2(s_2-u\pm k(\rho))T_3(s_3-u\pm k(\rho))
$$
with $T_{2}, T_{3} \in\{\delta, \mathrm{PV},H, \mathrm{Ci}\}$, the terms of the form:
$$
A_{i,\pm}(\rho)B_{j,\pm}(\rho)
\big(T_2(s_2-u\pm k(\rho))T_3(s_3-u\pm k(\rho))-T_2(s_3-u\pm k(\rho))T_3(s_2-u\pm k(\rho))\big),
$$
and the terms of the form:
$$
A_{i,\pm}(\rho)B_{j,\pm}(\rho)(s_{3}-u)\big(T_2(s_2-u\pm k(\rho))T_3(s_3-u\pm k(\rho))-T_2(s_3-u\pm k(\rho))T_3(s_2-u\pm k(\rho))\big)
$$
with $T_{2} \in \{\delta, H, \mathrm{PV}, \mathrm{Ci}, r_{\chi}\}$ and $T_{3} \in\{H, \mathrm{Ci}, r_{\sigma}\}$.
We emphasize that, in the last two cases, when $T_2,T_3\in \{r_{\chi},r_{\sigma}\}$,
$A_{i,\pm}(\rho)B_{j,\pm}(\rho)$ or $A_{i,\pm}(\rho)B_{j,\pm}(\rho)(s_{k}-u)$ should be
replaced by $1$.

Before applying Lemma \ref{lem7.9}, we now show that $\overline{\chi(s)}$ is H\"{o}lder continuous.
In fact, it follows from Corollary \ref{cor7.1} and Lemma \ref{lem7.6} that,
for any $s,s'\in\R$ and $\alpha\in (0,\min\{\lambda_1,1\}]$,
\begin{equation}\label{7.78-2}
	\sup_{s,s'\in \R}\frac{|\overline{\chi(s)}-\overline{\chi(s')}|}{|s-s'|^{\alpha}}
	=\int_{\mathcal{H}}\frac{|\chi(s)-\chi(s')|}{|s-s'|^{\alpha}}\,\mathrm{d}\nu_{(t,r)}
 \leq \int_{\mathcal{H}}C\big(1+\rho|\ln\rho|\big)\,\mathrm{d}\nu_{(t,r)}<\infty,
\end{equation}
which implies that $\overline{\chi(s)}$ is H\"{o}lder continuous.
Hence, using Lemma \ref{lem7.9} and the fact that $|s_{j}-u|\leq k(\rho)$ for $j=2,3$, we obtain
\begin{equation*}
	\begin{aligned}
		&\Big\vert\int_{-\infty}^{\infty}\overline{\chi(s_{1})}(P_{2} \chi_{2}^{\tau} P_{3} \sigma_{3}^{\tau}
		-P_{3} \chi_{3}^{\tau} P_{2} \sigma_{2}^{\tau}) \psi(s_{1}) \,\mathrm{d}s_{1}{\mathbf{I}}_{\{\rho \geq \rho^{*}\}}\Big\vert \\
		&\leq C\max_{j, k, \pm}\Big\{|A_{j, \pm} k(\rho)|\big(|B_{k, \pm}|+ \|r_{\sigma}(\r,\cdot)\|_{C^{\alpha_1}(\overline{B_{L_0}})}\big),\nonumber\\
		&\qquad \qquad\quad  \,\|r_{\chi}(\r,\cdot)\|_{C^{ \alpha_1}(\overline{B_{L_0}})}
		\big(|B_{j, \pm}k(\r)|+\|r_{\sigma}(\r,\cdot)\|_{C^{ \alpha_1}(\overline{B_{L_0}})}\big)\Big\}
		\\&\leq  C\big(1+\rho^{2+\t_2}\big)
		= C\big(\rho^{\beta(\g_2+1)}+1\big)\qquad\, \text{for $\rho\geq \rho^{*}$},
	\end{aligned}
\end{equation*}
with $L_{0}:=|\operatorname{supp}\psi|+2$ and $\beta=\frac{\t_2+2}{\g_2+1}\in(0,1)$.
Thus, using Lemmas \ref{lem7.6} and \ref{lem7.10},
and Lebesgue's dominated convergence theorem, we obtain
\begin{equation*}
	\begin{aligned}
		&\lim _{\tau \rightarrow 0} \int_{-\infty}^{\infty}  \overline{\chi(s_{1})}
		\langle\,\nu_{(t,r)},\,(P_{2} \chi_{2}^{\tau}\, P_{3} \sigma_{3}^{\tau}-P_{3} \chi_{3}^{\tau}\, P_{2} \sigma_{2}^{\tau})
		{\mathbf{I}}_{\{\rho \geq \rho^{*}\}}\rangle \,\psi(s_{1})\,\mathrm{d}s_{1} \\
		&=Y(\phi_{2}, \phi_{3}) \sum_{\pm}(K^{\pm})^{2}\langle\nu_{(t,r)},\, Z(\rho)\, \overline{\chi(u \pm k(\rho))}
		\,\psi(u \pm k(\rho)){\mathbf{I}}_{\{\rho \geq \rho^{*}\}}\rangle,
	\end{aligned}
\end{equation*}
which, with \eqref{7.87}, yields \eqref{7.85}. This completes the proof.
$\hfill\square$

\smallskip
\begin{lemma}\label{lem7.12}
	For any test function $\psi\in \mathcal{D}(\R)$,
	\begin{align}\label{7.93}
		\lim\limits_{\tau\to 0}\int_{\R}
		&\Big(\overline{P_3\chi_3^{\tau}}\;\,\overline{P_2\chi_2^{\tau}\,\sigma(s_1)-\chi(s_1)\,P_2\sigma_2^{\tau}}\nonumber\\
		&\quad -\overline{P_2\chi_2^{\tau}}\;\,\overline{P_3\chi_3^{\tau}\,\sigma(s_1)
			+\chi(s_1)\,P_3\sigma_3^{\tau}}\Big)\,\psi(s_1)\,\mathrm{d}s_1=0.
	\end{align}
\end{lemma}

\noindent\textbf{Proof.} $\,$Fix $(\rho,u)\in \mathbb{H}$. It follows from Lemma \ref{lem7.10}  that
\begin{equation*}
	\big(\chi(s_1)\,P_3\sigma_3^{\tau}-P_3\chi_3^{\tau}\,\sigma(s_1)\big)(\rho,u,s_1)
	\longrightarrow X(\rho,u,s_1)\qquad \text{uniformly in }s_1\text{ as }\tau \to 0.
\end{equation*}
It is clear that
\begin{equation*}
	\begin{aligned}
		&\int_{\mathbb{R}} \overline{P_{2} \chi_{2}^{\tau}}\,\big(\chi(s_1)\, P_{3} \sigma_{3}^{\tau}
		-P_{3} \chi_{3}^{\tau} \, \sigma(s_1)\big)\, \psi(s_{1}) \,\mathrm{d}s_{1} \\
		&=\int_{\mathcal{H}} \int_{\mathbb{R}} (P_{2} \chi_{2}^{\tau})(\tilde{\rho}, \tilde{u}, s_{1})\,
		\big(\chi(s_1)\, P_{3} \sigma_{3}^{\tau}-P_{3} \chi_{3}^{\tau} \,\sigma(s_1)\big)(\rho, u, s_{1})\,\psi(s_{1})
		\,\mathrm{d}s_{1} \mathrm{d} \nu_{(t,r)}(\tilde{\rho}, \tilde{u}).
	\end{aligned}
\end{equation*}
It follows from Lemma \ref{lem7.3} that $P_{j} \chi_{j}^{\tau}, j=2,3$, are measures in $s_{1}$ such that
$$
\|P_{j} \chi_{j}^{\tau}(\tilde{\rho}, \tilde{u}, \cdot)\|_{\mathfrak{M}, \alpha} \leq C_{\alpha}\tilde{\r}
\qquad\,\,\mbox{for large $\tilde{\rho}$},
$$
where
$\|\mu\|_{\mathfrak{M}, \alpha}
=\sup \left\{|\langle\mu, f\rangle|\,:\, f\in C_{0}^{ \alpha}(\mathbb{R})\right.$
and $\left.\|f\|_{C^{ \alpha}(\mathbb{R})} \leq 1\right\}$ with $\alpha \in(0,1)$.
Then we  use Lemma \ref{lem7.6} and Lebesgue's dominated convergence theorem
to pass the limit inside the Young measure to obtain
\begin{align}
	&\int_{\mathbb{R}} \overline{P_{2} \chi_{2}^{\tau}}\,\big(\chi(s_1)\, P_{3} \sigma_{3}^{\tau}
	-P_{3} \chi_{3}^{\tau}\, \sigma(s_1)\big) \,\psi(s_{1}) \,\mathrm{d}s_{1}\nonumber\\
	& \longrightarrow  \int_{\mathcal{H}}\int_{\R}(P_{1} \chi)(\tilde{\rho}, \tilde{u}, s_{1}) X(\rho, u, s_{1})
	\psi(s_{1})\,\mathrm{d}s_{1}\mathrm{d}\nu_{(t,r)}(\tilde{\rho}, \tilde{u})\nonumber
\end{align}
pointwise in $(\rho,u)$ as $\tau \to 0$.
Now we are going to prove
\begin{equation}\label{7.97}
	\Big\vert\int_{\R}\overline{P_2\chi_2^{\tau}}\,\big(\chi(s_1)\,P_3\sigma_3^{\tau}
	-P_3\chi_3^{\tau}\,\sigma(s_1)\big)\,\psi(s_1)\,\mathrm{d}s_1\Big\vert\leq C\big(\rho^{\beta(\gamma_2+1)}+1\big)
\end{equation}
for some constants $C>0$ and $\beta\in (0,1]$, which are both independent of $\tau$.
Once \eqref{7.97} is proved, it follows from Lebesgue's dominated convergence theorem that
\begin{equation*}
	\begin{aligned}
		&\lim\limits_{\tau\to 0}\int_{\R}\overline{P_2\chi_2^{\tau}}(s_1)\,\overline{\chi(s_1)\,P_3\sigma_3^{\tau}
			-P_3\chi_3^{\tau}\,\sigma(s_1))}\,\psi(s_1)\,\mathrm{d}s_1
		\\&=\lim\limits_{\tau\to 0}\int_{\mathcal{H}}\int_{\R}\overline{P_2\chi_2^{\tau}}(s_1)\,\big(\chi(s_1)\,P_3\sigma_3^{\tau}
		-P_3\chi_{3}^{\tau}\,\sigma(s_1)\big)(\rho,u,s_1)\,\psi(s_1)\,\mathrm{d}s_1\mathrm{d}\nu_{(t,r)}(\rho,u)
		\\&=\int_{\mathcal{H}}\int_{\mathcal{H}}\int_{\R}(P_{1} \chi)(\tilde{\rho}, \tilde{u}, s_{1}) X(\rho, u, s_{1}) \psi(s_{1})\,\mathrm{d}s_{1}\mathrm{d}\nu_{(t,r)}(\tilde{\rho},\tilde{u})\,\mathrm{d}\nu_{(t,r)}(\rho,u).
	\end{aligned}
\end{equation*}
Since $X(\r,u,s_{1})$ is independent of the choice of the mollifying functions $\phi_{2}^{\tau}$ and $\phi_{3}^{\tau}$ from Lemma \ref{lem7.10}, we may
interchange the roles of $s_{2}$ and $s_{3}$ to conclude the proof of \eqref{7.93}.

To see the validity of \eqref{7.97}, we begin by observing that, for $j=2,3$,
$\overline{P_j\chi_{j}^{\tau}}(s_1)$ and $\psi(s_1)$ are independent of $(\rho,u)$.
Then it suffices to estimate the function:
\begin{equation}\label{7.99}
	\chi(s_1)\,P_j\sigma_j^{\tau}-P_j\chi_{j}^{\tau}\,\sigma(s_1)=\chi(s_1)\,P_{j}(\sigma_{j}^{\tau}-u\chi_{j}^{\tau})-(\sigma(s_1)-u\chi(s_1))\,P_j\chi_{j}^{\tau}.
\end{equation}
It follows from Lemmas \ref{thm6.1}--\ref{thm6.2} and \ref{lem7.3} (also see \cite[Proof of Lemma 4.2]{Chen-LeFloch-2000}) that
$$
\mbox{RHS of } \eqref{7.99} = E^{1,\tau}+E^{2,\tau}+E^{3,\tau}+E^{4,\tau},
$$
with
\begin{align*}	
	E^{1,\tau}=&\sum\limits_{\pm}\big(A_{1,\pm}(\r)b_{1}(\r)G_{\lambda_1}(s_1)+A_{1,\pm}(\r)b_2(\r)G_{\lambda_1+1}(s_1)\big)\nonumber\\[-1mm]
	&\qquad\times
	\big((s_j-s_1)\delta(s_j-u\pm k(\rho))\big)*\phi_{j}^{\tau}\\[1mm]
	&+\sum\limits_{\pm}\big(A_{3,\pm}(\r)b_1(\r)G_{\lambda_1}(s_1)+A_{3,\pm}(\r)b_2(\r)G_{\lambda_1+1}(s_1)\big)\nonumber\\[-1mm]
    &\qquad\quad \times
	\big((s_j-s_1)PV(s_j-u\pm k(\rho))\big)*\phi_{j}^{\tau}\\[1mm]
	:=&E_{1}^{1,\tau}+E_{2}^{1,\tau},\\[1mm]
	E^{2,\tau}=&-\sum_{\pm}A_{1,\pm}(\r)g_2(s_1)\big(\delta(s_{j}-u\pm k(\r))*\phi_{j}^{\tau}\big)\nonumber\\
	&-\sum_{\pm}A_{3,\pm}(\r)g_2(s_1)\big(PV(s_{j}-u\pm k(\r))*\phi_{j}^{\tau}\big)\\[1mm]	
	:=&E_{1}^{2,\tau}+E_{2}^{2,\tau},\\[2mm]
	E^{3, \tau}=& \sum_{\pm}\big(B_{1,\pm}(\r)a_1(\r)-A_{1,\pm}(\r)b_{1}(\r)\big)G_{\lambda_1}(s_1)\big((s_j-u)\delta(s_j-u\pm k(\rho))\big)*\phi_{j}^{\tau}\nonumber\\
	&+\sum_{\pm}\big(B_{1,\pm}(\r)a_2(\r)-A_{1,\pm}(\r)b_2(\r)\big)G_{\lambda_1+1}(s_1)\big((s_j-u)\delta(s_j-u\pm k(\rho))\big)*\phi_{j}^{\tau}\nonumber\\
	&+\sum_{\pm}B_{1,\pm}(\r)g_1(s_1)\big((s_j-u)\delta(s_j-u\pm k(\rho))\big)*\phi_{j}^{\tau} \nonumber\\
	&+\sum_{\pm} \big(B_{3,\pm}(\r)a_1(\r)-A_{3,\pm}(\r)b_{1}(\r)\big)G_{\lambda_1}(s_1)\big((s_j-u)PV(s_j-u\pm k(\rho))\big)*\phi_{j}^{\tau}\nonumber\\
	&+\sum_{\pm} \big(B_{3,\pm}(\r)a_2(\r)-A_{3,\pm}(\r)b_{2}(\r)\big)G_{\lambda_1+1}(s_1)\big((s_j-u)PV(s_j-u\pm k(\rho))\big)*\phi_{j}^{\tau}\nonumber\\
	&+\sum_{\pm} B_{3,\pm}(\r)g_1(s_1)\big((s_j-u)PV(s_j-u\pm k(\rho))\big)*\phi_{j}^{\tau}\nonumber\\
	:=&E_{1}^{3,\tau
	}+E_{2}^{3,\tau
	}+E_{3}^{3,\tau
	}+E_{4}^{3,\tau
	}+E_{5}^{3,\tau
	}+E_{6}^{3,\tau
	},\nonumber
\end{align*}
and $E^{4,\tau}$ is the remainder term which consists of the mollification of continuous functions,
where we have used the notation: $G_{\lambda_1}(s_1)=[k(\rho)^2-(u-s_1)^2]^{\lambda_1}$,
and $g_{i}(s_1)=g_{i}(\rho,u-s_1)$ for $i=1,2$.

We first demonstrate the uniform bound on the term involving the delta measures. By direct calculation, we have
\begin{align}\label{7.100-4}
		&\delta(s_{j}-u+k(\r))*\phi_{j}^{\tau}=\frac{1}{\tau}\phi_{j}(\frac{s_{1}-u+k(\rho)}{\tau}),\nonumber\\
		&\big((s_{j}-s_{1})\delta(s_{j}-u+k(\r))\big)*\phi_{j}^{\tau}=-\frac{s_1-u+k(\r)}{\tau}\phi_{j}(\frac{s_{1}-u+k(\rho)}{\tau}),\nonumber\\
		&\big((s_{j}-u)\delta(s_{j}-u+k(\r))\big)*\phi_{j}^{\tau}\nonumber\\
		&\,\,\,=\big((s_{j}-s_{1})\delta(s_{j}-u+k(\r))\big)*\phi_{j}^{\tau}
		+\big((s_{1}-u)\delta(s_{j}-u+k(\r))\big)*\phi_{j}^{\tau}\nonumber\\
		&\,\,\,=(s_1-u)\tau^{-1}\phi_{j}(\frac{s_{1}-u+k(\rho)}{\tau})
		-\frac{s_1-u+k(\r)}{\tau}\phi_{j}(\frac{s_{1}-u+k(\rho)}{\tau}).
\end{align}
Noting that
\begin{equation}\label{7.100-2}
	\begin{split}
		&G_{\lambda_1+1}(s_1)=G_{\lambda_1}(s_1)\,(k(\rho)-u+s_1)\,(k(\rho)+u-s_1),\\[1mm]
		& |g_{i}(s_1)|\leq \|\partial_{u}g_{i}(\rho,\cdot)\|_{L^{\infty}(\R)}|s_{1}-u\pm k(\rho)|\qquad
		\mbox{ for $i=1,2$},
	\end{split}
\end{equation}
using \eqref{7.100-4}--\eqref{7.100-2}
and Lemmas \ref{thm6.1}--\ref{thm6.2} and \ref{lem7.1}--\ref{lem7.3}, we obtain
\begin{equation}\label{7.100-3}
	E_{1}^{3,\tau}=0,\qquad |E_{1}^{1,\tau}|+|E_{1}^{2,\tau}|+|E_{2}^{3,\tau}|+|E_{3}^{3,\tau}|\leq C_{\phi}\big(1+\rho^{\frac{3}{2}+\frac{\t_2}{2}}\big).
\end{equation}

For the term involving with the principal value distribution, a direct calculation shows that
\begin{align*}
	|PV*\phi_{j}^{\tau}(x)|&=\Big\vert\int_{0}^{\infty}\frac{\phi_{j}^{\tau}(x-y)-\phi_{j}^{\tau}(x+y)}{y}\,\mathrm{d}y\Big\vert\nonumber\\
	&=\frac{1}{\tau}\Big\vert\int_{0}^{\infty}\frac{1}{y}\big(\phi_{j}(\frac{x+y}{\tau})
	-\phi_{j}(\frac{x-y}{\tau})\big)\,\mathrm{d}y\Big\vert.
\end{align*}
If $|x|\leq 2\tau$, we have
\begin{align}\label{7.103-1}
	|PV*\phi_j^{\tau}(x)|&\leq \frac{1}{\tau}\int_{0}^{4\tau}\frac{1}{|y|}\big\vert\phi_j(\frac{x-y}{\tau})
	-\phi_j(\frac{y+x}{\tau})\big\vert \,\mathrm{d}y
	\leq C\frac{1}{\tau}\|\phi_{j}'\|_{L^{\infty}}=\frac{C_{\phi}}{|x|}.
\end{align}
On the other hand, if $|x|\geq 2\tau$, we can obtain
\begin{equation}\label{7.103-2}
	|PV*\phi_{j}^{\tau}(x)|
	\leq \frac{2}{\tau}\int_{|x|-\tau}^{|x|+\tau}\frac{\|\phi\|_{L^{\infty}}}{|x|-\tau}\,\mathrm{d}y\leq \frac{C_{\phi}}{|x|}.
\end{equation}
Notice that
\begin{equation}
\begin{aligned}\label{7.103-3}
		&\big((s_{j}-s_1)PV(s_{j}-u-k(\r))\big)*\phi_{j}^{\tau}\\
		&\,\,\,=\big((s_{j}-u+k(\r))PV(s_{j}-u+k(\r))\big)*\phi_{j}^{\tau}\\
		&\,\,\,\quad
		+\big((u-k(\r)-s_1)PV(s_{j}-u+k(\r))\big)*\phi_{j}^{\tau},\\
		&\big((s_{j}-u)PV(s_{j}-u-k(\r))\big)*\phi_{j}^{\tau}\\
		&\,\,\,=\big((s_{j}-u+k(\r))PV(s_{j}-u+k(\r))\big)*\phi_{j}^{\tau}
		-k(\r)PV(s_{j}-u+k(\r))*\phi_{j}^{\tau},
\end{aligned}
\end{equation}
which, with \eqref{7.100-2}, \eqref{7.103-1}--\eqref{7.103-3}, and
Lemmas \ref{thm6.1}--\ref{thm6.2} and \ref{lem7.1}--\ref{lem7.3}, yields
\begin{align}\label{7.103-4}
	E_{4}^{3,\tau}=0,\qquad
	|E_{2}^{1,\tau}|+|E_{2}^{2,\tau}|+|E_{5}^{3,\tau}|+|E_{6}^{3,\tau}|\leq C_{\phi}\big(1+\rho^{\frac{3}{2}+\frac{\t_2}{2}}\big).
\end{align}

Combining \eqref{7.100-3} with \eqref{7.103-4} yields that there exists $\beta_1=\frac{3+\t_2}{2(\g_2+1)}\in (0,1)$ such that
\begin{equation}\label{7.104}
	|E^{1,\tau}+E^{2,\tau}+E^{3,\tau}|\leq C_{\phi}(1+\rho^{\frac{3}{2}+\frac{\t_2}{2}})\leq C_{\phi}(1+\rho^{\beta_1(\g_2+1)}).
\end{equation}

For $E^{4,\tau}$ consisting of the mollification of continuous functions,
direct calculations show that
\begin{equation}\label{7.109}
	|E^{4,\tau}|\leq C_{\phi}\, \big(1+\rho^{2+\t_2} |\ln \rho|\big)\leq C_{\phi}\,\big(1+\rho^{\beta_2(\g_2+1)}\big),
\end{equation}
with $\beta_2=\frac{4+3\t_2}{2(\g_2+1)}\in (0,1)$.
Combining \eqref{7.109} with \eqref{7.104}, we conclude the proof of \eqref{7.97}.
$\hfill\square$

\section{\,Existence of Global Finite-Energy Solutions of CEPEs}
In this section, we complete the proof of Theorem \ref{thm1.1}.
Since the proof is similar to \cite{Chen-He-Wang-Yuan-2021},
we sketch the proof for the self-containedness of this paper. We divide the proof into four steps.

\smallskip
1. Since $(\rho^{\v},m^{\v})(t,r)$ obtained in Theorem \ref{thm1.2} satisfies all
the assumptions of Theorem \ref{thm1.3}, then it follows from Theorem \ref{thm1.3}
that there exists a vector function $(\rho,m)(t,r)$ such that, up to a subsequence as $\v\to 0$,
\begin{align}
	&(\rho^{\v},m^{\v})\longrightarrow (\rho,m)
	\quad \text{{\it a.e.} $(t,r)\in \R_{+}^2$},
	\label{8.1}\\
	&(\rho^{\varepsilon}, m^{\varepsilon}) \longrightarrow(\rho, m) \quad
	\text {in $L_{\mathrm{loc}}^{p_{1}}(\mathbb{R}_{+}^{2}) \times L_{\mathrm{loc}}^{p_{2}}(\mathbb{R}_{+}^{2})$} \label{8.4}
\end{align}
for $p_{1} \in[1, \gamma_2+1)$ and $p_{2} \in [1, \frac{3(\gamma_2+1)}{\gamma_2+3})$, where $L_{\mathrm{loc}}^{p_{j}}(\mathbb{R}_{+}^{2})$
represents $L^{p_{j}}([0, T] \times K)$
for any $T>0$ and
$K \Subset (0, \infty), j=1,2$.

Noting \eqref{8.1} and $\rho^{\v}\geq 0$ {\it a.e.} from Lemma \ref{lem5.1},
it is direct to show that $\rho(t,r)\geq 0$ {\it a.e.} on $\R_{+}^2$.
Moreover, it follows from \eqref{1.29} that $\sqrt{\rho^{\v}}u^{\v}r=\frac{m^{\v}}{\sqrt{\rho^{\v}}}r$
is uniformly bounded in $L^{\infty}(0,T;L^2(\R))$.
Then using Fatou's lemma yields
\begin{equation*}
	\int_{0}^{T}\int_{0}^{\infty}\frac{|m(t,r)|^2}{\rho(t,r)}\,r^{2}\mathrm{d}r\mathrm{d}t\leq \liminf_{\v\to 0}
	\int_{0}^{T}\int_{0}^{\infty}\frac{|m^{\v}(t,r)|^2}{\rho^{\v}(t,r)}\,r^{2}\mathrm{d}r\mathrm{d}t<\infty.
\end{equation*}
Thus, $m(t,r)=0$ {\it a.e.} on $\{(t,r)\,:\,\r(t,r)=0\}$,
and we can define the limit velocity $u(t,r)$ as
\begin{equation*}
	\begin{aligned}
		&u(t,r)=\frac{m(t,r)}{\rho(t,r)}\qquad \text{{\it a.e.} on }\{(t,r)\,:\,\rho(t,r)\neq 0\},\\
		&u(t,r)=0 \qquad \text{{\it a.e.} on }\{(t,r)\,:\,\rho(t,r)=0\text{ or }r=0\}.
	\end{aligned}
\end{equation*}
Therefore, $m(t,r)=\rho(t,r)u(t,r)$ {\it a.e.} on $\R_{+}^2$.
Also, we can  define $\big(\frac{m}{\sqrt{\rho}}\big)(t,r):=\sqrt{\rho(t,r)}u(t,r)$,
which is zero {\it a.e.} on $\{(t,r)\,:\,\rho(t,r)=0\}$.
Moreover, using \eqref{1.31} and Fatou's lemma, we obtain
\begin{equation*}
	\int_{0}^{T}\int_{d}^{D}\rho|u|^3\,\mathrm{d}r\mathrm{d}t
	\leq \liminf_{\v\to 0}\int_{0}^{T}\int_{d}^{D}\rho^{\v}|u^{\v}|^3\,\mathrm{d}r\mathrm{d}t\leq C(d,D,M,E_0,T)<\infty
\end{equation*}
for any
$[d,D]\Subset (0,\infty)$.

By similar calculations as in \cite[\S 5]{Chen-He-Wang-Yuan-2021}, we obtain that, as $\v\to 0$,
\begin{equation}\label{8.3}
	\frac{m^{\varepsilon}}{\sqrt{\rho^{\varepsilon}}} \equiv \sqrt{\rho^{\varepsilon}} u^{\varepsilon}
	\longrightarrow \frac{m}{\sqrt{\rho}} \equiv \sqrt{\rho} u \qquad
	\text{strongly in }L^{2}([0, T] \times[d, D]; r^{n-1} \,\mathrm{d} r \mathrm{d} t)
\end{equation}
for any $T>0$ and
$[d, D] \Subset (0, \infty)$.

From \eqref{8.4}--\eqref{8.3}, we also obtain the convergence of the mechanical energy as $\varepsilon \rightarrow 0$:
\begin{equation}\label{8.5}
	\eta^{*}(\rho^{\varepsilon}, m^{\varepsilon}) \longrightarrow \eta^{*}(\rho, m)
	\qquad \text { in } L_{\mathrm{loc}}^{1}(\mathbb{R}_{+}^{2}).
\end{equation}
Using \eqref{8.4}, \eqref{8.5}, and Fatou's lemma,
and taking limit $\v \to 0$ in \eqref{1.28}--\eqref{1.29}, we have
\begin{equation}\label{8.6}
	\int_{t_{1}}^{t_{2}} \int_{0}^{\infty}
	\big(\eta^{*}(\rho, m)(t, r)+\rho^{\gamma_2}(t,r)+\rho(t, r)\big)\, r^{2}\mathrm{d} r \mathrm{d} t
	\leq C(M, E_{0})(t_{2}-t_{1}),
\end{equation}
which implies
\begin{equation}\label{8.7}
	\sup_{0\leq t\leq T}\int_{0}^{\infty}
	\big(\eta^{*}(\rho, m)(t, r)+\rho^{\gamma_2}(t,r)+\rho(t, r)\big)\,r^{2}\mathrm{d}r \leq C(M, E_{0}).
\end{equation}
This indicates that $\rho(t,r)\in L^{\infty}([0,T];L^{\gamma_2}(\R;r^{2}\mathrm{d}r))$,
which implies that  $\rho(t,\mathbf{x})$ is a function in $L^{\infty}([0,T];L^{\g_2}(\R^3))$
with $\g_2>1$ (rather than a measure in $(t,\mathbf{x})$).
Therefore, no delta measure ({\it i.e.}, concentration) is formed in the density in the time
interval $[0,T]$, especially at the origin: $r=0$.

\smallskip
2. For the convergence of the gravitational potential functions $\Phi^{\v}(t,r)$,
by similar calculation in \cite[\S 5]{Chen-He-Wang-Yuan-2021}, we obtain that,
as $\v\to 0$ (up to a subsequence),
\begin{equation}\label{8.7-2}
	\Phi_{r}^{\varepsilon}(t, r) r^{2}
	= \int_{0}^{r} \rho^{\varepsilon}(t, y)\,y^{2}\mathrm{d}y
	\longrightarrow  \int_{0}^{r} \rho(t, y)\,y^{2}\mathrm{d}y \qquad \text { {\it a.e.} $(t,r)\in \R_{+}^2$}.
\end{equation}
Thus, using \eqref{5.34-2}, \eqref{8.1}, \eqref{8.7-2}, Fatou's lemma,
and similar arguments as in \eqref{8.6}--\eqref{8.7}, we have
\begin{equation*}
	\int_{0}^{\infty}\Big(\int_{0}^{r} \rho(t, y)\,y^{2}\mathrm{d}y\Big) \rho(t, r)\,r \mathrm{d} r
	\leq C(M, E_{0}) \qquad \text{ for {\it a.e.} $\,t \geq 0$}.
\end{equation*}
On the other hand, it follows from \eqref{5.34-3} that there exists a function $\Phi(t, \mathbf{x})=\Phi(t, r)$
such that, as $\varepsilon \rightarrow 0$ (up to a subsequence),
\begin{equation*}
	\begin{aligned}
		&\Phi^{\varepsilon} \rightharpoonup \Phi \qquad
		\text{ weak-star in $L^{\infty}(0, T ; H_{\operatorname{loc}}^{1}(\mathbb{R}^{3}))$
			and weakly in $L^{2}(0, T ; H_{\mathrm{loc}}^{1}(\R^3)$},\\
		&\|\Phi(t)\|_{L^6(\mathbb{R}^{3})}+\|\nabla \Phi(t)\|_{L^{2}\left(\mathbb{R}^{3}\right)} \leq C(M, E_{0})
		\qquad \text { {\it a.e.} $t \geq 0$}.
	\end{aligned}
\end{equation*}
Thus, by \eqref{8.7-2} and the uniqueness of limit, we obtain
that $\Phi_{r}(t, r) r^{2}=\int_{0}^{r} \rho(t, z) z^{2} \,\mathrm{d} z$
{\it a.e.} $(t, r) \in \mathbb{R}_{+}^{2}$.
By
similar arguments in \cite[\S 5]{Chen-He-Wang-Yuan-2021}, we also have the strong convergence of the potential functions:
\begin{equation}\label{8.7-8}
	\lim _{\varepsilon \rightarrow 0} \int_{0}^{T}\int_{0}^{\infty}\big|(\Phi_{r}^{\varepsilon}-\Phi_{r})(t, r)\big|^{2} r^{2}
	\,\mathrm{d} r \mathrm{d} t=0 \qquad \text { for $\gamma_2>\frac{6}{5}$}.
\end{equation}

\smallskip
3. Now we define
\begin{equation*}
	(\rho,\mathcal{M},\Phi)(t,\mathbf{x} ):=(\rho(t,r),m(t,r)\frac{\mathbf{x}}{r}, \Phi(t,r)).
\end{equation*}
Then it follows from \eqref{1.33}, \eqref{8.7-8}, and Fatou's lemma that
\begin{align*}
	&\int_{t_{1}}^{t_{2}} \int_{\mathbb{R}^{3}}\Big(\frac{1}{2}\Big|\frac{\mathcal{M}}{\sqrt{\rho}}\Big|^{2}
	+\rho e(\rho)-\frac{1}{2}|\nabla\Phi|^2\Big)(t, \mathbf{x}) \,\mathrm{d} \mathbf{x} \mathrm{d} t\nonumber\\
	&\quad \leq(t_{2}-t_{1}) \int_{\mathbb{R}^{3}}\Big(\frac{1}{2}\Big|\frac{\mathcal{M}_{0}}{\sqrt{\rho_{0}}}\Big|^{2}
	+\rho_{0} e(\rho_{0})-\frac{1}{2}|\nabla\Phi_{0}|^2\Big)(\mathbf{x}) \,\mathrm{d}\mathbf{x},
\end{align*}
which implies that, for {\it a.e.} $t \geq 0$,
\begin{align}\label{8.10}
		&\int_{\mathbb{R}^{3}}\Big(\frac{1}{2}\Big|\frac{\mathcal{M}}{\sqrt{\rho}}\Big|^{2}
		+\rho e(\rho)-\frac{1}{2}|\nabla\Phi|^2\Big)(t, \mathbf{x}) \,\mathrm{d}\mathbf{x}\nonumber\\
		&\leq \int_{\mathbb{R}^{3}}\Big(\frac{1}{2}\Big|\frac{\mathcal{M}_{0}}{\sqrt{\rho_{0}}}\Big|^{2}+\rho_{0} e(\rho_{0})
		-\frac{1}{2}|\nabla\Phi_{0}|^2\Big)(\mathbf{x})\,\mathrm{d}\mathbf{x}.
\end{align}
On the other hand, using \eqref{1.29}, \eqref{8.7}, and \eqref{8.7-8}, we obtain
\begin{equation}\label{8.10-1}
	\int_{\mathbb{R}^{3}}\Big(\frac{1}{2}\Big|\frac{\mathcal{M}}{\sqrt{\rho}}\Big|^{2}+\rho e(\rho)
	+\frac{1}{2}|\nabla\Phi|^2\Big)(t, \mathbf{x}) \,\mathrm{d}\mathbf{x}\leq C(M,E_0).
\end{equation}
Combining \eqref{8.10} with \eqref{8.10-1}, we complete the proof of \eqref{1.16-1}.

\smallskip
4. Using \eqref{5.35}, \eqref{5.53}--\eqref{5.80}, and
similar arguments as in \cite[\S 5]{Chen-He-Wang-Yuan-2021},
we conclude the proof of \eqref{1.17}--\eqref{1.18-1}
which, along with Steps 1--3, shows that $(\rho, \mathcal{M}, \Phi)(t, \mathbf{x})$ is indeed a global weak solution
of problem \eqref{1.1} and \eqref{1.11}--\eqref{1.11-1} in sense of Definition \ref{def1.1}.
This completes the proof.
$\hfill\square$

\bmhead{Acknowledgments}

Gui-Qiang G. Chen's research is partially supported
by the UK Engineering and Physical Sciences Research Council Awards
EP/L015811/1, EP/V008854/1, and EP/V051121/1.
Feimin Huang's research is partially supported by the National Natural Science Foundation of China No. 12288201
and the National Key R\&D Program of China No. 2021YFA1000800.
Tianhong Li's research is partially supported by the National Natural Science Foundation of China No. 10931007.
Yong Wang's research is partially supported by the National Natural Science Foundation of China No. 12022114 and No. 12288201,
and CAS Project for Young Scientists in Basic Research, Grant No. YSBR-031.

\section*{Declarations}

The authors have no competing interests to declare that are relevant to the content of this article. In addition, our manuscript has no associated data.

\medskip
\begin{appendices}

\section{\,Some Inequalities}\label{AppendixB}
\subsection{\,A sharp Sobolev inequality}
In this subsection, we recall a sharp Sobolev inequality, which is used in \S\ref{BE}.
The proof can be found in  \cite[\S 8.3]{Lieb-Loss-2001}.

\smallskip
\begin{lemma}[Sobolev inequality]\label{lemC.0}
	Let $n\geq 3$ and $\nabla f\in L^2(\R^n)$ with $\lim_{|\mathbf{x}|\to \infty}f(\mathbf{x})=0$.
	Then
	\begin{equation*}
		\|f\|_{L^{\frac{2n}{n-2}}}^2\leq A_{n}\|\nabla f\|_{L^2}^2,
	\end{equation*}
	where $A_{n}=\frac{4}{n(n-2)}\omega_{n+1}^{-\frac{2}{n}}$ is the best constant
	and $\omega_{k}=\frac{2\pi^{\frac{k}{2}}}{\Gamma(\frac{k}{2})}$ is the surface area of the unit sphere
	in $\R^k$.	
\end{lemma}

\subsection{\,Some variants of the Gr\"onwall inequality}
In this subsection, we introduce some variants of the Gr\"onwall inequality, which plays an essential role in
identifying the singularities of the entropy kernel and entropy flux kernel;
see also \cite{Schrecker-Schulz-2020}.

\smallskip
\begin{lemma}[A variant of Gr{\"o}nwall inequality {\cite[Theorem 1.2.4]{Qin-2017}}]\label{lemB.1}
	Let $x(t), y(t), z(t)$, and $w(t)$ be non-negative continuous functions on $J=[t_0,t_1]$ with $t_0\geq 0$. If
	$$
	x(t)\leq y(t)+z(t)\int_{t_0}^tw(s) x(s)\,\mathrm{d} s\qquad \text{for $t\in J$},
	$$
	then
	$$
	x(t)\leq y(t)+z(t)\int_{t_0}^t w(s)y(s)\exp\Big(\int_{s}^tw(r)z(r)\,\mathrm{d} r\Big)\,\mathrm{d} s\qquad
	\text{for $t\in J$}.
	$$
\end{lemma}

\begin{lemma}\label{lemB.2}
	Let $\t\geq 0$, and let $d(s)$ be defined in \eqref{A.16}.
	Assume that $x(t)\geq 0$ is measurable and locally integrable, and satisfies
	\begin{equation}\label{B.4-1}
		x(t)\leq Ct^{\theta}+\frac{1}{t}\int_{0}^td(s)x(s)\,\mathrm{d}s\qquad \text{for $t\geq \rho^{*}$}
	\end{equation}
	for some constant $C>0$. Then there exists a possibly larger constant $\tilde{C}>0$ independent of $t$
	such that, for $t\geq \rho^{*}$,
	$$
	x(t)\leq
	\begin{cases}
		\tilde{C}t^{\theta_2} \quad &\text{if }0\leq \t <\t_2,\\
		\tilde{C}t^{\theta_2}\ln t\,\,\,&\text{if }\t=\t_2,\\
		\tilde{C}t^{\t} \quad&\text{if }\t>\t_2.
	\end{cases}
	$$
\end{lemma}

\smallskip
\noindent\textbf{Proof.}
Since $x(t)$ is positive and locally integrable, then,
using Lemma \ref{lemA.3}, there exists a constant $C>0$ that may depend on $\rho^{*}$,
but independent of $t$, such that
$$
\frac{1}{t}\int_{0}^{\rho^{*}}d(s)x(s)\,\mathrm{d} s\leq Ct^{-1}\qquad\text{for $t\geq \rho^{*}$}.
$$
This, with \eqref{B.4-1}, yields that
$x(t)\leq Ct^{\theta}+\frac{1}{t}\int_{\rho^{*}}^{t}d(s)x(s)\,\mathrm{d} s$
for $t\geq \rho^{*}$.
Applying Lemma \ref{lemB.1}, we obtain
\begin{equation}\label{B.6}
	x(t)\leq Ct^{\t}+\frac{1}{t}\int_{\r^{*}}^{t}Cs^{\t}d(s)\mathrm{exp}\Big(\int_{s}^{t}\frac{d(r)}{r}
	\,\mathrm{d} r\Big)\,\mathrm{d} s\qquad \text{for $t\geq \rho^{*}$}.
\end{equation}
It is clear that
\begin{equation}\label{B.7}
	\Big\vert\int_{s}^{t}\frac{d(r)}{r}\,\mathrm{d} r\Big\vert
	\leq \int_{s}^{t}\frac{1+\t_2}{r}\,\mathrm{d} r+\int_{s}^{t}\frac{|d(r)-(1+\t_2)|}{r}\,\mathrm{d} r.
\end{equation}
It follows from Lemma \ref{lemA.3} that
$
\frac{|d(r)-(1+\t_2)|}{r}\leq Cr^{-\epsilon-1}
$
for $r\geq \rho^{*}$
which, with \eqref{B.7}, yields
$$
\mathrm{exp}\Big(\int_{s}^{t}\frac{d(r)}{r}\,\mathrm{d} r\Big)\leq C\Big(\frac{t}{s}\Big)^{1+\t_2}
\qquad \text{for $t\geq s\geq \rho^{*}$}.
$$
Combining \eqref{B.6} with $|d(s)|\leq 3$, we obtain that, for $t\geq \rho^{*}$,
\begin{equation}\label{B.8}
	x(t)\leq Ct^{\t}+\frac{1}{t}\int_{\r^{*}}^{t}Cs^{\t}\Big(\frac{t}{s}\Big)^{1+\t_2}\,\mathrm{d} s
	\leq Ct^{\t}+Ct^{\t_2}\int_{\r^{*}}^{\t}s^{-1-\t_2+\t}\,\mathrm{d} s.
\end{equation}

{\it{Case 1}}. If $0\leq \t<\t_2$, it follows from \eqref{B.8} that
$$
x(t)\leq Ct^{\t}+Ct^{\t_2}\Big(\int_{\rho^{*}}^{\infty}s^{-1-\t_2+\t}\,\mathrm{d} s\Big)
\leq \tilde{C}t^{\t_2}\qquad \text{for $t\geq \rho^{*}$}.
$$

{\it{Case 2}}. If $ \t=\t_2$, it follows from \eqref{B.8} that
$$
x(t)\leq  Ct^{\t_2}+Ct^{\t_2}\Big(\int_{\rho^{*}}^{t}s^{-1}\,\mathrm{d} s\Big)\leq \tilde{C}t^{\t_2}\ln t
\qquad \mbox{for $t\geq \rho^{*}$}.
$$

{\it{Case 3}}. If $ \t>\t_2$, then it follows from \eqref{B.8} that
$$
x(t)\leq  Ct^{\t}+Ct^{\t_2}\Big(\int_{\rho^{*}}^{t}s^{-1-\t_{2}+\t}\,\mathrm{d} s\Big)
\leq \tilde{C}t^{\t}\qquad \text{for $t\geq \rho^{*}$}.
$$
This completes the proof.
$\hfill\square$

\smallskip
\begin{corollary}\label{corB.1}
	If $x(t)$ satisfies
	$$
	x(t)\leq Ct^{\theta}\ln t+\frac{1}{t}\int_{0}^td(s)x(s)\,\mathrm{d} s\qquad \text{for $t\geq \rho^{*}$},
	$$
	with $\t>\t_2$, then
	$\,x(t)\leq Ct^{\t} \ln t\,$
	for $t\geq \rho^{*}$.
\end{corollary}

\end{appendices}

\end{document}